\documentclass[12pt]{article} 

\usepackage[utf8]{inputenc} 

\usepackage{geometry} 
\usepackage{CJKutf8}

\usepackage{graphicx} 
\usepackage[parfill]{parskip} 
\usepackage[all]{xy}
\usepackage{amsfonts}

\usepackage{amsmath}

\usepackage{hyperref}

\usepackage{tikz}
\usetikzlibrary{shapes}

\newcommand{\bi}{\begin{itemize}}
\newcommand{\ei}{\end{itemize}}
\newcommand{\be}{\begin{enumerate}}
\newcommand{\ee}{\end{enumerate}}
\newcommand{\beq}{\begin{eqnarray*}}
\newcommand{\eeq}{\end{eqnarray*}}
\newcommand{\bc}{\begin{center}}
\newcommand{\ec}{\end{center}}
\newcommand{\bL}{\begin{LARGE}}
\newcommand{\eL}{\end{LARGE}}
\newcommand{\bH}{\begin{Huge}}
\newcommand{\eH}{\end{Huge}}

 \voffset = 0.0 mm \oddsidemargin = -5 mm
\begin{document}
\newcounter{homework} 
\newcounter{tothomework} 
\newcommand{\zth}[1]{\begin{CJK}{UTF8}{bsmi}#1\end{CJK}}
\newcommand{\at}{\makeatletter @\makeatother}
\begin{center}

\begin{huge}

{\bf Introduction to Queueing Theory and Stochastic Teletraffic
Models}

\end{huge}

\end{center}

\vspace{10mm}
\begin{Large}{\bf  Moshe Zukerman

\vspace{6 mm}

  EE Department

  \vspace{1 mm}

  City University of Hong Kong }

  \vspace{1 mm}

  email: moshezu\at gmail.com \end{Large}

\vspace{10mm}

Copyright M. Zukerman \copyright ~2000--2024. \\This book can be used for educational and research purposes under the condition that it (including this first page) is not modified in any way.


\section*{Preface}

The aim of this textbook is to provide students with basic knowledge
of stochastic models with a special focus on queueing models that may apply to telecommunications topics, such as traffic, performance evaluation, resource provisioning, and
traffic management. These topics are included in a field called
{\em teletraffic}. This book assumes
prior knowledge of a programming language and mathematics normally
taught in an electrical engineering bachelor program.

The book aims to enhance intuitive and physical understanding of the
theoretical concepts it introduces. The famous mathematician
Pierre-Simon Laplace is quoted to say that ``Probability is common sense reduced to
calculation'' \cite{BT02}; as the content of this book falls under the field of
applied probability, Laplace's quote very much applies. Accordingly, the book aims to
link intuition and common sense to the mathematical models and techniques it uses. It mainly focuses on steady-state analyses and avoids discussions of time-dependent analyzes.

 A unique feature of this
book is the considerable attention given to guided homework assignments involving
computer simulations and analyses. By successfully completing these
assignments, students learn to simulate and analyze stochastic models,
such as queueing systems and networks, and by interpreting the results, they
gain insight into the queueing performance effects and
principles of telecommunications systems modeling. Although the
book, at times, provides intuitive explanations, it still presents the
important concepts and ideas required for the understanding of teletraffic, queueing theory
fundamentals and related queueing behavior of telecommunications networks and systems.
These concepts and ideas form a strong base for the more mathematically inclined students who can follow up
with the extensive literature on probability models and queueing theory.
A small sample of it is listed at the end of this book.

The first two chapters provide background on
probability and stochastic processes topics relevant to the queueing
and teletraffic models of this book. These two chapters
provide a summary of the key topics with relevant
homework assignments that are specially tailored for understanding
the queueing and teletraffic models discussed in later chapters.
The content of these chapters is mainly based on \cite{BT02,fell68,ros2003,ross93,ross70,ross76,ross96}.
Students are encouraged to also study the original textbooks
for more explanations, illustrations, discussions,
examples, and homework assignments.

Chapter \ref{general} discusses general queueing notation and
concepts. Chapter \ref{simulations}
aims to assist the student in performing simulations of queueing
systems. Simulations are useful and important in many cases
where exact analytical results are not available.
An important learning objective of this book is to
train students to perform queueing simulations.
Chapter \ref{det} provides analyses of deterministic
queues. Many queueing theory books tend to exclude deterministic
queues; however, the study of such queues is helpful for beginners in
that it helps them better understand non-deterministic queueing
models. Chapters \ref{mm1queue} -- \ref{multiservice} provide
analyses of a wide range of queueing and teletraffic models most of which
fall under the category of continuous-time Markov-chain processes.
Chapter \ref{discrete} provides an example of a discrete-time queue
that is modeled as a discrete-time Markov chain. In Chapter
\ref{mg1}, various aspects of a single
server queue with Poisson arrivals and general service times are studied, mainly focussing on
mean value results as in \cite{BG92}. Then, in Chapter \ref{gg1q},
some selected results of a single server queue with a general
arrival process and general service times are provided. Chapter \ref{multiaccess} focuses on multi-access applications, and
in Chapter  \ref{networks}, we extend our discussion to
queueing networks. Finally,
in Chapter \ref{trafmod}, stochastic processes that have been used
as traffic models are discussed with a special focus on their
characteristics that affect queueing performance.

Throughout the book, there is an emphasis on linking the theory with telecommunications applications as
demonstrated by the following examples. Section \ref{link_dim}  describes how properties of Gaussian
distribution can be applied to link dimensioning. Section  \ref{muxmm1} shows, in the context of an
M/M/1 queueing model, how optimally to set a link service rate such that delay requirements are met
and how the level of multiplexing affects the spare capacity required to meet such delay requirements.
An application of M/M/$\infty$ queueing model to a multiple access performance problem \cite{BG92}
is discussed in Section \ref{mminfmaccess}. Then later in Chapter \ref{multiaccess} more multi-access models are presented. In Sections \ref{dim2Erlang} and \ref{dim_mmk},
discussions on dimensioning and related utilization issues of a multi-channel system are
presented. Especially important is the emphasis on the insensitivity property of models such as M/M/$\infty$, M/M/$k$/$k$, processor sharing, and multi-service that lead to practical and robust approximations as described in Chapters \ref{secmminf}, \ref{secerlang}, \ref{processorsharing}, and \ref{multiservice}.
Section \ref{mobcellsim} guides the reader to simulate a mobile cellular network.
Section \ref{ppbpmodel} describes a traffic model applicable to the Internet.

Last but not least, the author wishes to thank all the students and
colleagues who provided comments and questions that helped
develop and edit the manuscript over the years. 



 \newpage

 \tableofcontents

 \newpage
\section{Background on Relevant Probability Topics}
\label{probability}

\setcounter{homework}{1} 

\setcounter{tothomework}{1} 

Probability theory provides the foundation for queueing theory and
stochastic teletraffic models; it is, therefore, important that the
student masters the probability concepts required for the material
that follows.
Although the
coverage here is comprehensive in the sense that it discusses all the
probability concepts and techniques used in later chapters, and it includes many examples and exercises; still, readers may be aided by additional probability texts, such as
\cite{BT02} and \cite{ross76} to complement their study for improved learning.

\subsection{Events, Sample Space, and Random Variables} Consider an
experiment with an uncertain outcome. The term ``experiment'' refers
to any uncertain scenario, such as tomorrow's weather, tomorrow's
share price of a certain company, or the result of flipping a coin.
The {\em sample space} is a
 set of all possible outcomes of an experiment. An {\em event} is a subset
of the sample space. Consider, for example, an experiment that
consists of tossing a die. The sample space is
$\{1,~2,~3,~4,~5,~6\}$, and an event could be the set
 $\{2,~3\}$, or $\{6\}$, or the empty set $\{\}$ or even the entire sample space
$\{1,~2,~3,~4,~5,~6\}$. Events are called {\em mutually exclusive}
if their intersection is the empty set. A set of events is said to
be {\em exhaustive} if its union is equal to the sample space.

A {\em random variable} is a real-valued function defined on the
sample
 space. This definition appears somewhat contradictory to the wording ``random variable'' as a random variable
is not at all random because it is actually a deterministic
function, which assigns a real-valued number to each possible outcome
of an experiment. It is the outcome of the experiment that is random
and therefore the name: random variable. If we consider the flipping
a coin experiment, the possible outcomes are Head (H) and Tail (T),
hence the sample space is $\{{\rm H,~T}\}$, and a random variable
$X$ could assign $X=1$ for the outcome H, and $X=0$ for the outcome
T.

{\bf Remark} (Moran): While not every subset is an event and not every function is a random variable, for all practical purposes, we can and do assume that they are.

If $X$ is a random variable, then $Y=g(X)$ for some function $g(\cdot)$ is
also a random variable. In particular, some functions of interest are $Y=cX$
for some constant $c$ and $Y=X^n$ for some integer $n$.

If $X_1,~X_2,~X_3,~ \ldots,~X_n$ is a sequence of random variables, then
$Y=\sum_{i=1}^n X_i$ is also a
random variable.


\subsubsection*{Homework \ref{probability}.\arabic{homework}}

\addtocounter{homework}{1} \addtocounter{tothomework}{1}

Consider the experiment to be tossing a coin. What is the Sample
Space?  What are the events associated with this Sample Space?

\subsubsection*{Guide} Notice that although the sample space
includes only the outcome of the experiments, which are Head (H) and
Tail (T), the events associated with this sample space include
also the empty set, which in this case is the event $H \cap T $
and the entire sample space, which in this case is the event $H
\cup T $. $~~~\Box$

\subsection{Probability, Conditional Probability and Independence}

Consider a sample space $S$. Let $A$ be a subset of $S$, the probability of
$A$ is the function on $S$ and all its subsets, denoted P$(A)$, that
satisfies the following three axioms:
\begin{enumerate}
\item $0 \leq P(A) \leq 1$
\item $P(S)=1$
\item The probability of the union of mutually exclusive events is equal to
the sum of the probabilities of these events.
\end{enumerate}

Normally, a higher probability signifies a higher likelihood of occurrence. In
particular, if we conduct a very large number of experiments, and we generate the {\em histogram} by measuring how many times each of the possible occurrences
actually occurred, then we normalize the histogram by dividing all its values by the total number of experiments to obtain the relative frequencies. These
measurable relative frequencies can provide an intuitive interpretation of the theoretical concept of probability. Accordingly, the term {\it limiting relative frequency} is often used to interpret probability.

We use the notation $P(A \mid B) $ for the {\em conditional
probability} of $A$ given $B$, which is the probability of the
event $A$ given that we know that event $B$ has occurred. If we
know that $B$ has occurred, then it is our new sample space, and for
$A$ to occur, the relevant experiment outcomes must be in $ A
\cap B $, hence the new probability of $A$, namely the probability
$P(A \mid B) $, is the ratio between the probability of $ A \cap B
$ and the probability of $B$. Accordingly,

\begin{equation} \label{cond_prob} P(A \mid B)= \frac{P(A \cap B)} {
P( B)}. \end{equation}

Since the event $A \cap B $ is equal to the event $B \cap A$, we have that

$$ P(A \cap B) = P(B \cap A)= P(B \mid A) P( A),$$ so by
(\ref{cond_prob}) we obtain

\begin{equation} \label{cond_prob2} P(A \mid B)= \frac{P(B \mid A)
P( A)} { P( B)}. \end{equation} Eq. (\ref{cond_prob2}) is useful to
obtain the conditional probability of one event ($A$) given another
($B$) when $P(B \mid A)$ is known or easier to obtain, than   $P(A
\mid B)$.

{\bf Remark: } The intersection of $A$ and $B$ is also denoted by
$A,B$ or $AB$ in addition to $A\cap B$.

If events $A$ and $B$ are {\em independent}, which means that if one of
them
occurs, then the probability of the other to occur is not affected, then
\begin{equation}
P(A \mid B)= P(A)
\end{equation}

and hence, by Eq.\@ (\ref{cond_prob}), if $A$ and $B$ are
independent, then
\begin{equation}
\label{ind}
P(A \cap B)= P(A)   P(B).
\end{equation}

Let $B_1,~ B_2,~ B_3,~ \ldots,~B_n$ be a sequence of mutually
exclusive and exhaustive events in $S$, and let $A$ be another
event in $S$. Then,
\begin{equation}
\label{compevent}
A=\bigcup_{i=1}^{n} (A\cap B_i )
\end{equation}
and since the $B_i$s are mutually exclusive, the events $ A\cap B_i$s are
also mutually exclusive. Hence,
\begin{equation}
\label{compprob}
P(A)=\sum_{i=1}^{n} P(A\cap B_i ).
\end{equation}
Thus, by Eq.\@ (\ref{cond_prob}),
\begin{equation}
\label{compprob2}
P(A)=\sum_{i=1}^{n} P(A\mid B_i )\times P(B_i).
\end{equation}
The latter is a very useful formula for deriving the probability of a given
event by conditioning and unconditioning on a set of mutually exclusive and
exhaustive events. It is called {\em the Law of Total Probability}.
Therefore, by Eqs. (\ref{compprob2}) and (\ref{cond_prob}) (again), we
obtain the following formula for conditional probability between two events:
\begin{equation}
\label{bayes}
P(B_1 \mid A)= \frac{P(A \mid B_1) P(B_1)}{\sum_{i=1}^{n} P(A\mid B_i
)\times P(B_i)}.
\end{equation}
The latter is known as Bayes' theorem (or Bayes' law or Bayes' rule).

Note that Eq. (\ref{cond_prob2}) is also referred to as Bayes' theorem.

Set $A$ is said to be a {\em subset} of set $B$, namely, $A \subset B$ if $A\neq B$ and every element of $A$ is an element of $B$. If we do not require $A\neq B$, but still require that every element of $A$ is an element of $B$, then we say that $A$ is a subset of or equal to set $B$, which is denoted by $A\subseteq B$. If event $A$ occurs implies that event $B$ occurs, then we have that $A\subseteq B$. In such a case (and of course, under the stronger condition $A \subset B$), $A = A\cap B$, so

\begin{equation} \label{psubset} P(A) = P(A\cap B)=P(A\mid B) P(B).\end{equation}

\subsection{Probability and Distribution Functions}

Random variables are related to events. When we say that random variable $X$ takes value $x$,
this means that $x$ represents a certain outcome of an experiment, which is an event,
so $\{X=x\}$ is an event. Therefore,
 we may assign probabilities to all possible values of the
random
variable. This function denoted $P_X(x)=P(X=x)$ will henceforth be called {\em probability
function}. Other names used in the literature for a probability
function include {\em probability distribution function}, {\em probability distribution}, or simply {\em distribution}.
Because probability theory is used in many applications, in many cases, there are many alternative terms to describe the same thing.
It is important that the student is familiar with these alternative terms, so we will use these terms alternately in this book.

The {\em cumulative distribution function} of random variable $X$ is defined for all $x \in R$ ($R$ being the
set of all real numbers), is defined as
\begin{equation}
F_X(x)=P(X\leq x).
\end{equation}
Accordingly, the {\em complementary distribution function} $\bar{F}_X(x)$ is
defined by
\begin{equation}
\bar{F}_X(x)=P(X > x).
\end{equation}
Consequently, for any random variable, for every $x \in R$, $F(x) +
\bar{F}(x)=1$.
As the complementary and the cumulative distribution functions, as well as the probability function, can be obtained from each other, we will use the terms {\em distribution function} when we refer to any of these functions without being specific.

We have already mentioned that if $X$ is a random variable, then $Y=g(X)$ is also a random variable. In this case, if $P_X(x)$ is the probability function of $X$, then the probability function of $Y$ is
\begin{equation}
P_Y(y) = \sum_{x:g(x)=y} P_X(x).
\end{equation}

\subsection{Joint Distribution Functions}

In some cases, we are interested in the probability that two or
more random variables are within a certain range. For this
purpose, we define, {\em the joint distribution function} for $n$
random variables $X_1,~X_2,~ \ldots, ~X_n$, as follows:
\begin{equation}
{F}_{X_1,~X_2,~..., X_n} (x_1,~x_2,~..., x_n)=P(X_1 \leq x_1,~ X_2 \leq x_2,~ \ldots,~X_n \leq x_n).
\end{equation}
Having the joint distribution function, we can obtain the distribution
function of a single random variable, say, $X_1$, as
\begin{equation}
{F}_{X_1}(x_1)= {F}_{X_1,~X_2,~ \ldots,~ X_n}(x_1,~\infty,~ \ldots,~ \infty).
\end{equation}

A random variable is called {\em discrete} if it takes at most a
countable number of possible values. On the other hand, a {\em
continuous} random variable takes an uncountable number of possible
values.

When the random variables $X_1,~X_2,~..., X_n$ are discrete, we can use
their {\em joint probability function} which is defined by
\begin{equation}
\label{jointdis}
P_{X_1,~X_2,~ \ldots, ~X_n}(x_1,~x_2,~ \ldots, ~x_n)=P(X_1 = x_1,~ X_2 = x_2,~ \ldots,~X_n = x_n).
\end{equation}

The probability function of a single random variable can then be obtained by
\begin{equation}
\label{singldisc}
P_{X_1}(x_1)=\sum_{x_2} \cdots \sum_{x_n}  P_{X_1,~X_2,~ \ldots,~ X_n}(x_1,~x_2,~ \ldots, ~x_n).
\end{equation}

In this section and in sections \ref{cond_pr_rv},
\ref{indep_rv}, \ref{convol}, when we mention random variables or
their distribution functions, we consider them all to
be discrete. Then, in Section \ref{continRV}, we will introduce the
analogous definitions and notation relevant to their continuous
counterparts.

We have already mentioned the terms probability
function, distribution, probability distribution function, and probability distribution.
These terms apply to discrete as well as to continuous random variables.
There are, however, additional terms that are used to describe probability function only for discrete random variables; they are: {\em probability mass function}, and {\em probability mass}, and there are equivalent terms used only for continuous random variables -- they are  {\em probability density function}, {\em density function} and simply {\em density}.

\subsection{Conditional Probability for Random Variables}
\label{cond_pr_rv}
The conditional probability concept, which we defined for events, can also
apply to random variables.
Because $\{X = x \}$, namely, the random variable $X$ takes value $x$,
is an event, by the definition of conditional probability (\ref{cond_prob}), we have
\begin{equation}
\label{cond_rv}
P(X=x \mid Y=y) = \frac{P(X=x, Y=y)}{P(Y=y)}.
\end{equation}
Let  $P_{X\mid Y}(x \mid y) = P(X=x \mid Y=y)$ be the conditional
probability of a discrete random variable $X$ given $Y$, we obtain by
(\ref{cond_rv})
\begin{equation}
\label{cond_rv1}
P_{X\mid Y}(x \mid y) = \frac{P_{X,Y}(x,y)} {P_Y(y)}.
\end{equation}
Noticing that
\begin{equation}
\label{jointosing}
P_Y(y) = \sum_x P_{X,Y}(x,y),
\end{equation}
we obtain by (\ref{cond_rv1})
\begin{equation}
\label{sumofx}
\sum_x P_{X\mid Y}(x \mid y) = 1.
\end{equation}
This means that if we condition on the event $\{Y=y \}$ for a
specific $y$, then the probability function of $X$ given $\{Y=y \}$ is a
legitimate probability function. This is consistent with our
discussion above. The event $\{Y=y \}$ is the new sample space and
$X$ has a legitimate distribution function there.
By (\ref{cond_rv1})
\begin{equation}
\label{cond_rv2}
P_{X,Y}(x,y) = P_{X\mid Y}(x \mid y)  {P_Y(y)}
\end{equation}
and by symmetry
\begin{equation}
\label{cond_rv2a}
P_{X,Y}(x,y) = P_{Y\mid X}(y \mid x)  {P_X(x)}
\end{equation}
so the latter and (\ref{jointosing}) gives
\begin{equation}
\label{jointosing1}
P_Y(y) = \sum_x P_{X,Y}(x,y) = \sum_x P_{Y\mid X}(y \mid x)  {P_X(x)}
\end{equation}
which is another version of the Law of Total Probability (\ref{compprob2}).

\subsection{Independence between Random Variables}
\label{indep_rv}
The definition of independence between random variables is very much related to the definition of independence
between events because when we say that random variables $U$ and $V$ are independent, it is equivalent
to say that the events $\{U=u\}$ and $\{V=v\}$ are independent for every $u$ and $v$.
Accordingly, random variables $U$ and $V$ are said to be independent if
\begin{equation}
\label{ind_rv}
P_{U,V}(u,v) = P_U(u) P_V(v)~~~ {\rm for ~all}~ u,v.
\end{equation}
Notice that by (\ref{cond_rv2}) and (\ref{ind_rv}), we obtain an equivalent
definition of independent random variables $U$ and $V$ which is
\begin{equation}
\label{ind_rv2}
P_{U \mid V}(u \mid v) = P_U(u)
\end{equation}
which is equivalent to  $P(A\mid B) = P(A)$ which we used to define independent events $A$ and $B$.

\subsection{Convolution}
\label{convol}
Consider independent random variables $V_1$ and $V_2$ that have probability functions $P_{V_1}(v_1)$ and $P_{V_2}(v_2)$, respectively,
and their sum which is another random variable $V=V_1+V_2$. Let us now derive the probability function $P_V(v)$ of $V$.
\begin{eqnarray*}
P_V(v) & = & P(V_1+V_2=v) \\
                 & = & \sum_{v_1} P(V_1=v_1,V_2=v-v_1) \\
                 & = &  \sum_{v_1} P_{V_1}(v_1) P_{V_2}(v-v_1).
            \end{eqnarray*}
The latter is called the {\em convolution} of the probability functions $P_{V_1}(v_1)$ and $P_{V_2}(v_2)$.

Let us now extend the result from two to $k$ random variables.
Consider $k$ independent random variables  $X_i$, $i=1,2,3, ~\ldots,~
k$. Let $P_{X_i}(x_i)$ be the probability function of $X_i$, for
$i=1,2,3, ~\ldots,~ k$, and let $Y=\sum_{i=1}^k X_i$. If $k=3$, then we
first compute the convolution of $X_1$ and $X_2$ to obtain the
probability function of $V=X_1+X_2$ using the above convolution
formula and then we use the formula again to obtain the probability
function of $Y=V+X_3=X_1+X_2+X_3$. Therefore, for an arbitrary $k$,
we obtain
\begin{equation}
\label{convolution}
P_Y(y) = \sum_{x_2,~ x_3, ~\ldots, ~x_k:~x_2+ x_3+ ~\ldots ~+ x_k
\leq y}  \left( P_{X_1} ( y-\Sigma_{i=2}^k x_i  ) \prod_{i=2}^k
P_{X_i}(x_i) \right).
            \end{equation}

            If all the  random variable $X_i$, $i=1,2,3,~ \ldots,~ k$, are independent and identically distributed
            (IID) random variables, with probability function  $P_{X_1}(x)$, then the probability function  $P_Y(y)$ is
            called the k-fold convolution of $P_{X_1}(x)$.

\subsubsection*{Homework \ref{probability}.\arabic{homework}}

\addtocounter{homework}{1} \addtocounter{tothomework}{1}

Consider again the experiment of tossing a coin. Let $P(H)$ be the probability that the toss is Head and $P(T)$ the probability that the toss is Tail. Assume that
$$P(H)=P(T)=0.5.$$   Illustrate each of the Probability Axioms for this
case. $\Box$

\subsubsection*{Homework \ref{probability}.\arabic{homework}}
\addtocounter{homework}{1} \addtocounter{tothomework}{1}

Now consider an experiment involving three coin tosses. The outcome
of the experiment is now a 3-long string of Heads and Tails. Assume
that all the coin tosses have a probability of 0.5 for each of the events Head and Tail and that the coin tosses are independent events. \begin{enumerate} \item Write the sample
space where each outcome of the experiment is an ordered 3-long string of Heads and Tails. \item What is the probability of each
outcome? \item Consider the event $$A=\{Exactly~ one ~head
~occurs\}.$$ Find $P(A)$ using the additivity axiom. \end{enumerate}

{\bf Partial Answer:} $P(A)=1/8 + 1/8 + 1/8 = 3/8.$
 $\Box$

\subsubsection*{Homework \ref{probability}.\arabic{homework}}
\addtocounter{homework}{1} \addtocounter{tothomework}{1}

Now consider again three coin tosses. Find the probability
$P(A\mid B)$ where $A$ and $B$ are the events:\\
$A=$ more than one head came up\\
$B=$ 1st toss is a head.

\subsubsection*{Guide:} $P(B)=4/8$; $P(A\cap B)=3/8$; $P(A\mid
B)=(3/8)/(4/8)=3/4.$  $\Box$

\subsubsection*{Homework \ref{probability}.\arabic{homework}  (K. T. Ko)}
\addtocounter{homework}{1} \addtocounter{tothomework}{1}

Consider a medical test for a certain disease. The medical test
detects the disease with a probability of 0.99 and fails to detect the
disease with a probability of 0.01. If the disease is not present, then the
test indicates that it is present with a probability of 0.02 and that it
is not present with a probability of 0.98. Consider two cases:\\ Case a:
The test is done on a randomly chosen person from the population
where the occurrence of the disease is 1/10000.\\ Case b: The test
is done on patients who are referred by a doctor and have a prior
probability (before they do the test) of 0.3 to have the disease.

Find the conditional (posterior) probability of a person having the disease if the test
shows positive outcomes in each of these cases.

\subsubsection*{Guide:}
$A=$ person has the disease. \\
$B=$ test is positive.\\
$\bar{A}=$ person does not have the disease. \\$\bar{B}=$ test is negative.\\
We need to find $P(A\mid B)$. \\

{\bf Case a:} \\
We know: $P(A)= 0.0001$.\\
$P(\bar{A})= 0.9999$. \\
$P(B\mid A)=
0.99$.\\
$P(B\mid \bar{A})=0.02.$

By the Law of Total Probability:

$$P(B)=P(B\mid A) P(A) + P(B\mid \bar{A}) P(\bar{A}).$$

$P(B) = 0.99 \times 0.0001 + 0.02 \times 0.9999 = 0.020097.$

Now put it all together and use Eq. (\ref{cond_prob2}) to obtain:

$P(A \mid B)= 0.004926108.$

{\bf Case b:} \\$P(A)= 0.3.$

Repeat the previous derivations to show that for this case $P(A \mid
B)= 0.954983923.$
 $~~~\Box$

\subsubsection*{Homework \ref{probability}.\arabic{homework}}
\addtocounter{homework}{1} \addtocounter{tothomework}{1}
In a
multiple choice exam, there are 4 answers to a question. A student
knows the right answer with a probability of 0.8 (Case 1), with
a probability of 0.2 (Case 2), and with a probability of 0.5 (Case 3). If the
student does not know the answer, then s/he always guesses, with
a probability of success being 0.25. Given that the student marked the
right answer, what is the probability that he/she knows the answer?\\

\subsubsection*{Guide:}
$A=$ Student knows the answer. \\$B=$ Student marks correctly.\\
$\bar{A}=$ Student does not know the answer. \\$\bar{B}=$ Student marks incorrectly.\\
We need to find $P(A\mid B)$. \\

{\bf Case 1:}

We know: $P(A)= 0.8$.\\ $P(\bar{A})= 0.2$. \\$P(B\mid A)= 1$.\\
$P(B\mid \bar{A})=0.25.$

By the Law of Total Probability:

$$P(B)=P(B\mid A) P(A) + P(B\mid \bar{A}) P(\bar{A}).$$

$P(B) = 1 \times 0.8 + 0.25 \times 0.2 = 0.85.$

Now put it all together and use Eq. (\ref{cond_prob2}) to obtain:

$$P(A \mid B)= 0.941176471. $$

{\bf Case 2:}

Repeat the previous derivations to obtain:

$P(A)= 0.2$

$P(B) = 0.4$

$P(A \mid B)= 0.5. $

{\bf Case 3:}

Repeat the previous derivations to obtain:

$P(A)= 0.5$

$P(B) = 0.625$

$P(A \mid B)= 0.8. $ $~~~\Box$

\subsubsection*{Homework \ref{probability}.\arabic{homework}  \cite{Dreyfus77}}
\addtocounter{homework}{1} \addtocounter{tothomework}{1}

There are two coins. For the first coin $P(H)=0.3$ and $P(T)=0.7$, and for the second coin $P(H)=0.7$ and $P(T)=0.3$.
A coin is selected at
random with a 0.5 probability to select each coin, and then the selected coin is flipped. The result is Head. What is the probability that the coin that
was selected and flipped is the first coin?

{\bf Solution}

$A=$ The first coin was selected and flipped.\\
$B=$ The second coin was selected and flipped.\\
$H=$ The result is Head\\ $T=$ The result is Tail.

$P(A)=P(B)=0.5$\\
$P(H\mid A)=0.3$\\
$P(H\mid B)=0.7$.

$$P(A\mid H) = \frac{P(H \mid A) P(A)}{ P(H\mid A) P(A) + P(H\mid B) P(B)} = \frac{0.3 \times 0.5}{ 0.3 \times 0.5 + 0.7 \times 0.5 } =  0.3.$$

From the next question, it is clear that this result is expected.

$~~~\Box$

\subsubsection*{Homework \ref{probability}.\arabic{homework} }
\addtocounter{homework}{1} \addtocounter{tothomework}{1}

There are two coins. For the first coin $P(H)=\alpha$ and $P(T)=1-\alpha$, and for the second coin $P(H)=1-\alpha$ and $P(T)=\alpha$.
A coin is selected at
random with a 0.5 probability to select each coin, and then the selected coin is flipped. The result is Head. What is the probability that the coin that
was selected and flipped is the first coin?

{\bf Solution}

$P(A)=P(B)=0.5$\\
$P(H\mid A)=\alpha$\\
$P(H\mid B)=1-\alpha$.

$$P(A\mid H) = \frac{P(H \mid A) P(A)}{ P(H\mid A) P(A) + P(H\mid B) P(B)} = \frac{\alpha \times 0.5}{ \alpha \times 0.5 + (1-\alpha) \times 0.5 } =  \alpha.$$

$~~~\Box$

\subsubsection*{Homework \ref{probability}.\arabic{homework}  \cite{Dreyfus77}}
\addtocounter{homework}{1} \addtocounter{tothomework}{1}

Now consider the case:

$P(A)=0.2$\\
$P(B)=0.8$\\
$P(H\mid A)=0.3$\\
$P(H\mid B)=0.7$.

What is the probability that the coin that was selected is flipped and
came up Head is the first coin?

{\bf Solution}

$$P(A\mid H) = \frac{P(H \mid A) P(A)}{ P(H\mid A) P(A) + P(H\mid B) P(B)} = \frac{0.3 \times 0.2}{ 0.3 \times 0.2 + 0.7 \times 0.8 } =  0.096774194.$$

Clearly, we observe that reducing the prior probability of the event of selecting the first coin reduces its conditional (posterior) probability $P(A\mid H)$.

$~~~\Box$

\subsubsection*{Homework \ref{probability}.\arabic{homework}}
\addtocounter{homework}{1} \addtocounter{tothomework}{1}
Watch the following YouTube link on the problem known as the Monty Hall problem:\\
\href{https://www.youtube.com/watch?v=mhlc7peGlGg}{\url{https://www.youtube.com/watch?v=mhlc7peGlGg}}\\
Then study the Monty Hall problem in: \\
\href{http://en.wikipedia.org/wiki/Monty_Hall_problem}{\url{http://en.wikipedia.org/wiki/Monty_Hall_problem}}\\
and understand the different solution approaches and write a computer simulation to estimate the probability of winning the car by always switching. $\Box$


\subsection{Selected Discrete Random Variables} We present here
several discrete random variables and their corresponding
 distribution functions. Although our coverage here is
non-exhaustive, we do consider all the discrete random variables mentioned later in this book.

\subsubsection{Non-parametric}
Discrete non-parametric random variable $X$ is characterized by a set of $n$ values that $X$ can take with nonzero probability: $a_1, a_2, \ldots, a_n$, and a set of $n$ probability values $p_1, p_2, \ldots p_n$,  where  $p_i=P(X=a_i), ~~ i=1,2, \ldots n$.  The distribution of $X$, in this case, is called a {\it non-parametric distribution} because it does not depend on a mathematical function that its shape and range are determined by certain parameters of the distribution.

\subsubsection{Bernoulli}

We begin with the Bernoulli random variable. It represents
an outcome of an experiment that has only two possible outcomes.
Let us call them ``success'' and ``failure''. These two outcomes
are mutually exclusive and exhaustive events. The Bernoulli random
variable assigns the value $X=1$ to the ``success'' outcome and
the value $X=0$ to the ``failure'' outcome. Let $p$ be the
probability of the ``success'' outcome, and because ``success''
and ``failure'' are mutually exclusive and exhaustive, the
probability of the ``failure'' outcome is $1-p$. The probability
function in terms of the Bernoulli random variable is:
\begin{eqnarray}
P(X=1)& = & p \\
P(X=0)& = & 1-p. \nonumber
\end{eqnarray}

\subsubsection{Geometric}

There are two types of geometric distributions. The first one is a distribution of a random variable $X$ that represents the number of independent Bernoulli trials, each of which with $p$ being the probability of success, required until the first success. For $X$ to
be equal to $i$, we must have $i-1$ consecutive failures and then one
success in $i$ independent Bernoulli trials. Therefore, we obtain
\begin{equation} \label{geom} P(X = i) = (1 - p)^{i-1}p ~~~~ {\rm
for} ~~i = 1, 2, 3,~\ldots ~. \end{equation} The complementary
distribution function of the geometric random variable is $$ P(X >
i) = (1 - p)^i ~~~~ {\rm for} ~~i = 0, 1, 2, 3,~\ldots ~. $$

The second type is a distribution of a random variable $Y$ that represents the number of failures that occur {\it before} the first success and accordingly, it is given by $Y=X-1$ where $X$ is a random variable of the geometric distribution of the first type.

In this book, we use the name {\it geometric distribution} for the geometric distribution of the first type, and the geometric distribution of the second type we call {\it geometric distribution of failures}. Then, the corresponding random variables of the first and second type distributions will be called the {\it geometric random variable} and the {\it geometric random variable of failures}, respectively.

A geometric random variable (of either type) possesses an important property
called memorylessness. In particular, the discrete random variable $X$
is {\it memoryless} if \begin{equation}
    P(X>m+n \mid X > m)=P(X>n), ~~~m = 0, 1, 2, \ldots, {\rm and}~~
    n = 0, 1, 2, \ldots
    \end{equation}
    The Geometric random variable is memoryless because it is based
    on independent Bernoulli trials, and therefore, the fact that so
    far we had $m$ failures does not affect the probability that the
    next $n$ trials will be failures.

    The two types of geometric random variables are the only discrete random variables that are memoryless.

\subsubsection{Binomial}

Assume that $n$ independent Bernoulli trials are performed. Let $X$ be a
random
variable representing the number of successes in these $n$ trials. Such a random variable is called a  binomial random variable with parameters $n$ and
$p$.
Its probability function is:
\[
P(X=i) = {n \choose i} p^i(1-p)^{n-i} \qquad i=0,~1,~2,~\ldots, ~n.
\]
Notice that a Binomial random variable with parameters 1 and $p$
is a Bernoulli random variable. The Bernoulli and binomial random
variables have many applications. In particular, it is used as a
model for voice and data sources. Such sources alternate
between two states ``on'' and ``off''. During the ``on'' state, the
source is active and transmits at a rate equal to the

``off'' state, the source is idle. If $p$ is the proportion of
time that the source is active, and if we consider a superposition
of $n$ independent identical sources, then the binomial
distribution gives us the probability of the number of sources
that are simultaneously active, which is important for resource
provisioning.

\subsubsection*{Homework \ref{probability}.\arabic{homework}}

\addtocounter{homework}{1} \addtocounter{tothomework}{1} Consider a
state with a voter population $N$. There are two candidates in the
state election for governor and the winner is chosen based on a
simple majority. Let $N_1$ and $N_2$ be the total number of votes
obtained by Candidates 1 and 2, respectively, from voters other than
Johnny. Johnny just voted for Candidate 1, and he would like to know
the probability that his vote affects the election results. Johnny realizes that
the only way he can affect the result of the election is if the votes
are equal (without his vote) or the one he voted for (Candidate 1) had (without
his vote) one call less than Candidate 2. That is, he tries to find the probability of the event
$$0 \geq N_1-N_2 \geq -1.$$
Assume that any other voter (excluding
Johnny) votes independently for Candidates 1 and 2 with
probabilities $p_1$ and $p_2$, respectively, and also that $p_1 +
p_2 <1$ to allow for the case that a voter chooses not to vote for
either candidate. Derive a formula for the probability that Johnny's
vote affects the election results and provide an algorithm and a
computer program to compute it for the case $N=2,000,000$ and
$p_1=p_2=0.4$. \subsubsection*{Guide} By the definition of
conditional probability, $$P(N_1=n_1,N_2=n_2)=P(N_1=n_1)P(N_2=n_2
\mid N_1=n_1)$$ so $$P(N_1=n_1,N_2=n_2)= {{N-1} \choose
n_1}p_1^{n_1}(1-p_1)^{N-n_1-1}{{N-n_1-1} \choose
n_2}p_2^{n_2}(1-p_2)^{N-n_1-1-n_2}.$$ Then, as the probability of the union of mutually exclusive events is the sum of their probabilities, the required probability
is $$\sum_{k=0}^{\lfloor(N-1)/2\rfloor} P(N_1=k,N_2=k) +
\sum_{k=0}^{\lceil(N-1)/2\rceil - 1} P(N_1=k,N_2=k+1),$$ where
$\lfloor x \rfloor$ is the largest integer smaller or equal to x.
and $\lceil x \rceil$ is the smallest integer greater or equal to x.
$~~~\Box$

Next, let us derive the probability distribution of the random
variable $Y=X_1 + X_2$ where $X_1$ and $X_2$ are two independent
Binomial random variables with parameters $(N_1,p)$ and $(N_2,p)$,
respectively.  This will require the derivation of the convolution
of $X_1$ and $X_2$ as follows.

\begin{eqnarray*}
P_Y(k) & = & P(X_1+X_2=k) \\
                 & = & \sum_{i=0}^k P(\{X_1=i \}\cap \{X_2=k-i\}) \\
                 & = &  \sum_{i=0}^k P_{X_1}(i) P_{X_2}(k-i)\\
                 & = & \sum_{i=0}^k {N_1 \choose i} p^i (1-p)^{N_1
                 -i} {N_2 \choose {k-i}} p^{k-i} (1-p)^{N_2 -
                 (k-i)}\\
                   & = & p^k (1-p)^{N_1+N_2-k} \sum_{i=0}^k {N_1
                   \choose i}{N_2 \choose {k-i}}\\
   & = & p^k (1-p)^{N_1+N_2-k} {{N_1 + N_2}
                   \choose k}.\\
            \end{eqnarray*}

            We can conclude that the convolution of two independent binomial
            random variables with parameters $(N_1,p)$ and
            $(N_2,p)$ has a binomial distribution with parameters
            $N_1+N_2,p$. This is not surprising. As we recall the
            binomial random variable represents the number of
            successes of a given number of independent Bernoulli trials. Thus,
            the event of having $k$ Bernoulli successes out of $N_1
            + N_2$ trials is equivalent to the event of having some
            (or none)
            of the successes  out of the $N_1$ trials and the remaining
            out of the $N_2$ trials.

\subsubsection*{Homework \ref{probability}.\arabic{homework}}

\addtocounter{homework}{1} \addtocounter{tothomework}{1}

            In the last step of the above proof, we have used the
            equality

            $$ {{N_1 + N_2}
                   \choose k} = \sum_{i=0}^k {N_1
                   \choose i}{N_2 \choose {k-i}}. $$
Prove this equality. \subsubsection*{Guide} Consider the equality
$$(1+\alpha)^{N_1+N_2}=(1+\alpha)^{N_1}(1+\alpha)^{N_2}.$$ Then,
equate the binomial coefficients of $\alpha^k$ on both sides.
$~~~\Box$

\subsubsection{Poisson}
\label{poissonrvs}
A Poisson random variable with parameter $\lambda$ has the following
probability function: \begin{equation} \label{poissonrv}
P(X=i)=e^{-\lambda}\frac{\lambda^i} {i!} \qquad i=0,~1,~2,~3,~
\ldots~ . \end{equation}

To compute the values of $P(X=i)$, it may be convenient to use the
recursion \begin{equation} \label{poissonrec}
P(X=i+1)=\frac{\lambda}{i+1} P(X=i)\end{equation}  with
$$P(X=0)=e^{-\lambda}.$$

However, if the parameter $\lambda$ is large, other settings may be required. For example,
set $\hat{P}(X=\lfloor \lambda \rfloor) = 1$, where $\lfloor x
\rfloor$ is the largest integer smaller or equal to x, and to
compute recursively,  using (\ref{poissonrec}) (applied to the
$\hat{P}(X=i)$'s), a sufficient number of $\hat{P}(X=i)$ values for
$i > \lambda$ and $i < \lambda$ such that $\hat{P}(X=i)>\epsilon$,
where $\epsilon$ is chosen to meet a given accuracy requirement.
Clearly, the $\hat{P}(X=i)$'s do not sum up to one, and therefore
they do not form a proper probability distribution, so they will need to be normalized as follows.

Let
 ${\alpha}$ and ${\beta}$ be the lower  and upper bounds,
 respectively, of the $i$ values for which $\hat{P}(X=i)>\epsilon$.
Then, the probabilities $P(X=i)$ are approximated using the
normalization:

\begin{equation} P(X=i)=\left\{\begin{array}{ll}
\frac{\hat{P}(X=i)}{\sum_{i=\alpha}^{\beta}
\hat{P}(X=i)} & \mbox{ $\alpha \leq i \leq \beta$}\\
~~~~0 & \mbox{otherwise.} \end{array} \right. \end{equation}

The importance of the Poisson random variable lies in its property
to approximate the binomial random variable in cases when $n$ is very
large and $p$ is very small so that $np$ is not too large and not
too small. In particular, consider a sequence of binomial random
variables $X_n, ~~n=1,2, \ldots$ with parameters $(n,p)$ where
$\lambda = np$, or $p=\lambda/n$. Then, the probability function

$$\lim_{n\to\infty} P(X_n=k)$$

is a Poisson probability function with parameter $\lambda$.

To prove this we write:

$$\lim_{n\to\infty} P(X_n=k)=\lim_{n\to\infty}{n \choose k} p^k
(1-p)^{n-k}.$$

Substituting $p=\lambda/n$, we obtain

$$\lim_{n\to\infty} P(X_n=k)=\lim_{n\to\infty}{n! \over (n-k)!k!}
\left({\lambda \over n}\right)^k \left(1-{\lambda\over
n}\right)^{n-k},$$

or

$$\lim_{n\to\infty} P(X_n=k)=\lim_{n\to\infty}{n! \over (n-k)!n^k}
\left({\lambda^k \over k!}\right) \left(1-{\lambda\over
n}\right)^{n} \left(1-{\lambda\over n}\right)^{-k}.$$

Now notice that

$$\lim_{n\to\infty} \left(1-{\lambda\over n}\right)^{n} =
e^{-\lambda},$$

$$\lim_{n\to\infty}  \left(1-{\lambda\over n}\right)^{-k} = 1,$$

and

$$ \lim_{n\to\infty}{n! \over (n-k)!n^k} =1.$$

Therefore,

$$\lim_{n\to\infty} P(X_n=k)= \frac{\lambda^k e^{-\lambda}} {k!}.$$

In Subsection \ref{ztransforms}, this important limit will be shown
using Z-transform. The Poisson random variable accurately models the
number of calls arriving at a telephone exchange or Internet service
provider in a short period of time, a few seconds or a minute, say.
In this case, the population of customers (or flows) $n$ is large.
The probability $p$ of a customer making a call within a
given short period of time is small, and the calls are typically
independent because they are normally generated by independent
individual people from a large population. Therefore,
models based on Poisson random variables
have been used successfully for the design and dimensioning of
telecommunications networks and systems for many years. When we
refer to items in a queueing system in this book, they will be
called customers, jobs, messages, or packets, interchangeably.

Next, let us derive the probability distribution of the random
variable $Y=X_1 + X_2$, where $X_1$ and $X_2$ are two independent
Poisson random variables with parameters $\lambda_1$ and
$\lambda_2$, respectively.  This will require the derivation of the
convolution of $X_1$ and $X_2$ as follows.

\begin{eqnarray*}
P_Y(k) & = & P(X_1+X_2=k) \\
                 & = & \sum_{i=0}^k P(\{X_1=i \}\cap \{X_2=k-i\}) \\
                 & = &  \sum_{i=0}^k P_{X_1}(i) P_{X_2}(k-i)\\
                 & = & \sum_{i=0}^k \frac{\lambda_1^i}{i!}  e^{-\lambda_1} \frac{\lambda_2^{k-i} }{(k-i)!}e^{-\lambda_2}\\
                  & = & e^{-(\lambda_1+\lambda_2)}\sum_{i=0}^k \frac{\lambda_1^i \lambda_2^{k-i}}{i!(k-i)!}\\
                   & = & \frac{e^{-(\lambda_1+\lambda_2)}}{k!}\sum_{i=0}^k \frac{k!\lambda_1^i
                   \lambda_2^{k-i}}{i!(k-i)!}\\
                    & = & \frac{e^{-(\lambda_1+\lambda_2)}}{k!}\sum_{i=0}^k{k \choose i}  \lambda_1^i
                   \lambda_2^{k-i}\\
                     & = & \frac{e^{-(\lambda_1+\lambda_2)}(\lambda_1+\lambda_2)^k}{k!}.
            \end{eqnarray*}

            We have just seen that the random variable $Y$ has a Poisson
            distribution with parameter $\lambda_1+\lambda_2$.

\subsubsection*{Homework \ref{probability}.\arabic{homework}}

\addtocounter{homework}{1} \addtocounter{tothomework}{1} Consider a Poisson random variable $X$
with parameter $\lambda=500$. Write a program that computes the
probabilities $P(X=i)$ for $0 \leq i \leq 800$ and plot the function
$P_X(x)$. $~~~\Box$

\subsubsection*{Homework \ref{probability}.\arabic{homework}}

\addtocounter{homework}{1} \addtocounter{tothomework}{1} Let $X_1$
and $X_2$ be two independent Poisson distributed random variables
with parameters $\lambda_1$ and $\lambda_2$, respectively. Let
$Y=X_1+X_2$. Find the distribution of $(X_1\mid Y)$. In particular,
find for any given $k$, the conditional probability $P(X_1=j \mid Y=k)~{\rm for}~j=0,1,2, \ldots,
k$.

 \subsubsection*{Guide}
 $$P(X_1=j \mid Y=k) = \frac {P(X_1=j \cap Y=k)}{P( Y=k)}.$$

Now notice that $$\{X_1=j \cap Y=k\}=\{X_1=j \cap X_2=k-j\}.$$ When
you claim that the two events $A$ and $B$ are equal, you must be able to
show that $A$ implies $B$ and that $B$ implies $A$. Show both for
the present case.

Because $X_1$ and $X_2$ are independent, we have that $$P(X_1=j \cap
X_2=k-j)= P(X_1=j ) P (X_2=k-j).$$

Now recall that $X_1$ and $X_2$ are Poisson distributed random
variables with parameters $\lambda_1$ and $\lambda_2$, respectively,
and that $Y$ is the convolution of $X_1$ and $X_2$  and therefore
 has a Poisson distribution with parameter $\lambda_1 +
\lambda_2$.  Put it all together and show that $P(X_1=j \mid Y=k)$ is
a binomial probability function with parameters
$\lambda_1/(\lambda_1 + \lambda_2)$ and $k$. $~~~\Box$

\subsubsection{Pascal}

The Pascal random variable $X$ with parameters $k$ (integer and
$\geq 1$) and $p$ (real within (0,1]), represents a sum of $k$
geometric random variables each with parameter $p$. This is equivalent to having
the $k$th successful Bernoulli trial at the $i$th trial.
In other words, the Pascal random variable represents the number of independent Bernoulli trials until the $k$th success. Accordingly, for a Pascal random variable, $X$, to be equal to $i$, the $i$th trial must be the $k$th successful trial associated with the $k$th geometric random variable. Then, there must also be exactly $k-1$ successes among the
first $i-1$ trials.

The probability to have a success at the $i$th trial, as in any trial, is equal to $p$,
and the probability of having $k-1$ successes
among the first $i-1$ is equal to the probability of having a
binomial random variable with parameters $p$ and $i-1$ equal to $k-1$,
for $i\geq k \geq 1$, which is equal to
\[
{i-1 \choose k-1} p^{k-1}(1-p)^{i-k} \qquad k=1,~2,~\ldots, ~i
\]
and since the two random variables here, namely, the Bernoulli and
the Binomial, are independent (because the underlying Bernoulli
trials are independent), we obtained
\begin{equation}
P(X=i)={i-1 \choose k-1} p^{k}(1-p)^{i-k} \qquad i=k,~k+1,~k+2,~
\ldots ~.
\end{equation}

An alternative formulation of the Pascal random variable $Y$ is defined as the number of failures required to have $k$ successes. In this case, the relationship between $Y$ and $X$ is given by $Y=X-k$ for which the probability mass function is given by

\begin{equation}
P(Y=j)={k+j-1 \choose k-1} p^{k}(1-p)^{j} \qquad j=0,~1,~2,~
\ldots ~.
\end{equation}

One important property of the latter is that the cumulative distribution function can be
expressed as

\begin{equation}
\label{pascalcdf}
F_Y(j) = P(Y \leq j) = 1-I_p(j+1,k) \qquad j=0,~1,~2,~
\ldots ~,
\end{equation}
where $I_p(j+1,k)$ is the regularized incomplete beta function with parameters $p,j+1$ and $k$.

\subsubsection{Discrete Uniform}

The discrete uniform probability function with integer parameters $a$ and $b$ takes equal non-zero values for $x=a, a+1, a+2, \dots, b$.
Its probability function is given by

    $$
P_X(x)=\left\{\begin{array}{ll}
\frac{1}{b-a+1} & \mbox{if $x=a, a+1, a+2, \dots, b$}\\
0 & \mbox{otherwise.} \end{array} \right. .$$

\subsubsection*{Homework \ref{probability}.\arabic{homework}}

\addtocounter{homework}{1} \addtocounter{tothomework}{1}

Consider a group of people. Assume that their birthdays are independent and there are 365 days in a year. That is, their birthdays can be viewed as independent discrete uniform random variables with parameters 1 and 365. Let $A$ be the event that there are at least two people in a group with the same birthday. Let $P_N(A)$ be the probability of the event $A$ occurring given that the number of people in the group is $N$. Write a recursive formula to derive $P_N(A)$.  Obtain the values of $P_2(A)$, $P_5(A)$ and $P_{10}(A)$. Find the smallest values of $N$, such that  $P_N(A)\geq 0.5$ and $P_N(A)\geq 0.999$, respectively.

 \subsubsection*{Guide and Answers}
Notice the following recursive relationship for the probability of not having two people with the same birthday among $N$ people:
$$1- P_N(A)= \left(1-\frac{N-1}{365}\right)\left(1-P_{N-1}(A)\right),~~~ {\rm for} ~~~ N=2, 3, 4, ~\ldots, ~365,$$
with $P_1(A)=0$, and $P_N(A) =1$ for $N \geq 365$.

$P_2(A)=0.00274$, $P_5(A)=0.027$ and $P_{10}(A)=0.117$

For $N=23$, $P_{23}(A)= 0.5073$.

For $N=70$, $P_{70}(A)= 0.99916$.

Are you surprised? $~~~\Box$

\subsection{Continuous Random Variables and their Distributions}
\label{continRV}
Continuous random variables are related to
cases whereby the set of possible outcomes is uncountable. A continuous random variable $X$ is a function that assigns a real number to the outcome of an experiment and is characterized by the existence of a function $f(\cdot)$ defined for all $x\in R$, which
has the property that for any set $A\subset R$,
\begin{equation}
\label{contA}
P(X\in A)= \int_A f(x) dx.
\end{equation}
Such a function is the {\em probability density function} (or simply the {\em density}) of $X$.
Since the continuous random variable $X$ must take a value in $R$ with probability 1, $f$
must satisfy,
\begin{equation}
\int_{-\infty}^{+\infty} f(x) dx=1.
\end{equation}
If we consider Eq.\@ (\ref{contA}), letting $A=[a, b]$, we obtain,
\begin{equation}
\label{conta(1)a(2)} P(a \leq X \leq b) = \int_a^b f(x) dx.
\end{equation}
An interesting point to notice is that the probability of a
continuous random variable taking a particular value is equal to
zero. If we set $a=b$ in Eq.\@ (\ref{conta(1)a(2)}), we obtain
\begin{equation}
P(X=a) = \int_{a}^{a} f(x) dx=0.
\end{equation}
As a result, for a continuous random variable $X$, the cumulative distribution function $F_X(x)$ is equal to both $P(X\leq x)$
and to $P(X < x)$. Similarly, the complementary distribution function is
equal to both $P(X\geq x)$ and to $P(X > x)$.

By Eq.\@ (\ref{conta(1)a(2)}), we obtain
\begin{equation}
F_X(x)=P(X\leq x)= \int_{-\infty}^{x} f(s) ds.
\end{equation}

Hence, the probability density function is the derivative of the distribution function.

An important concept that gives rise to a continuous version of the Law of Total Probability is the continuous equivalence of Eq.\@
(\ref{jointdis}), namely, the joint distribution of continuous
random variables. Let $X$ and $Y$ be two continuous random
variables. The joint density of $X$ and $Y$ denoted $f_{X,Y}(x,y)$
is a nonnegative function that satisfies
\begin{equation}
P(\{X,Y\}\in A) = \int\!\!\!\int_{\{X,Y\}\in A} f_{X,Y}(x,y) dxdy
\end{equation}
for any set $A\subset R^2$.

The continuous equivalence of the first equality in (\ref{jointosing1}) is:
\begin{equation}
f_{Y}(y)= \int_{-\infty}^{\infty} f_{X,Y}(x,y) dx.
\end{equation}
Another important concept is the conditional density of one continuous random variable on another.
Let $X$ and $Y$ be two continuous random variables with joint density $f_{X,Y}(x,y)$. For any $y$, such
that the density of $Y$ takes a positive value at $Y=y$ (i.e. such that $f_{Y}(y)>0$),
the conditional density of $X$ given $Y$ is defined as
\begin{equation}
\label{cond_dens}
f_{X\mid Y}( x\mid y) = \frac{f_{X,Y}(x,y)}{f_{Y}(y)}.
\end{equation}
For every given fixed $y$, it is a legitimate density because
\begin{equation}
\int_{-\infty}^{\infty} f_{X\mid Y}( x\mid y)dx = \int_{-\infty}^{\infty} \frac{f_{X,Y}(x,y)dx}{f_{Y}(y)}= \frac{f_{Y}(y)}{f_{Y}(y)}=1.
\end{equation}
Notice the equivalence between the conditional probability (\ref{cond_prob}) and
the conditional density (\ref{cond_dens}).
By (\ref{cond_dens})
\begin{equation}
f_{X,Y}(x,y) = f_{Y}(y)  f_{X\mid Y}(x\mid y)
\end{equation}
so
\begin{equation}
f_{X}(x) = \int_{-\infty}^{\infty} f_{X,Y}(x,y) dy = \int_{-\infty}^{\infty} f_{Y}(y)  f_{X\mid Y}(x\mid y) dy.
\end{equation}
Recall again that $f_{X,Y}(x,y)$ is defined only for $y$ values such that $f_{Y}(y)>0$.

Let define the event $A$ as the event \{$X\in A$\} for $A \subset R$. Thus,
\begin{equation}
P(A) = P(X\in A) = \int_A f_{X}(x) dx = \int_A  \int_{-\infty}^{\infty} f_{Y}(y)  f_{X\mid Y}(x\mid y) dy dx.
\end{equation}
Hence,
\begin{equation}
P(A) = \int_{-\infty}^{\infty} f_{Y}(y) \int_A  f_{X\mid Y}(x\mid y) dx dy
\end{equation}
and therefore
\begin{equation}
\label{tot_cont}
P(A) = \int_{-\infty}^{\infty} f_{Y}(y) P(A \mid Y=y) dy
\end{equation}
which is the continuous equivalence of the Law of Total Probability
(\ref{compprob2}).

\subsubsection*{Homework \ref{probability}.\arabic{homework}}

\addtocounter{homework}{1} \addtocounter{tothomework}{1}

I always use two trains to go to work. After traveling on my first train and walking in the interchange station, I arrive at the platform (in the interchange station) at a time that is uniformly distributed between 9:01 and 9:04 AM. My second train arrives at the platform (in the interchange station) exactly at times 9:02 and 9:04 AM. Derive the density of my waiting time at the interchange station.
Ignore any queueing effects, and assume that I never incur any queueing delay.

\subsubsection*{Guide:}

Let $D$ be a random variable representing the delay in minutes. Plot $D$ as a function of the arrival time.

Then, derive the complementary distribution of $D$ to obtain:

$$P(D>t) = \left\{\begin{array}{ll}
1-\frac{2t}{3} & \mbox{~~~$0\leq t \leq 1,$} \\
\frac{2-t}{3}  & \mbox{~~ $1< t<2, $} \\
0 & \mbox{~~~otherwise.}
\end{array} \right. $$

Finally, obtain the density from this complementary distribution function.  $~~~\Box$

We will now discuss the concept of convolution as applied to continuous random variables.
Consider independent random variables $U$ and $V$ that have densities $f_U(u)$ and $f_V(v)$, respectively, and their sum which is another random variable $X=U+V$. Let us now derive the density $f_X(x)$ of $X$.
\begin{eqnarray}
\label{conv2}
f_X(x) & = & P(U+V=x) \nonumber \\
                 & = & \int_u f(U=u,V=x-u)du \nonumber \\
                 & = &  \int_u f_U(u) f_V(x-u)du.
            \end{eqnarray}
The latter is the {\em convolution} of the densities $f_U(u)$ and $f_V(v)$.

As in the discrete case the convolution $f_{Y}(y)$, of $k$ densities $f_{X_i}(x_i)$,  $i=1,2,3, ~\ldots,~ k$,
of random variables  $X_i$, $i=1,2,3, ~\ldots,~ k$, respectively, is given by

\begin{equation}
f_Y(y) = \int\!\!\!\int_{x_2, ~\ldots, ~x_k:~x_2+
~\ldots, +x_k \leq y} \left( f_{X_1} ( y-\Sigma_{i=2}^k x_i  )
\prod_{i=2}^k f_{X_i}(x_i) \right).
            \end{equation}

And again, in the special case where all the random variables $X_i$,
$i=1,2,3, ~\ldots,~ k$, are IID, the density  $f_Y$ is the k-fold
convolution of $f_{X_1}$.

\subsubsection*{Homework \ref{probability}.\arabic{homework}}

\addtocounter{homework}{1} \addtocounter{tothomework}{1}
Consider the following joint density:
\begin{equation} f_{X,Y}(x,y)= \left\{\begin{array}{ll}
2 & \mbox{ $0\leq x+y \leq 1 $},~~ x\geq 0, ~~y \geq 0,\\
0 & \mbox{~~otherwise.} \end{array} \right. \end{equation}
\begin{enumerate}
\item Show that this is a legitimate joint density by showing first that all relevant probabilities are nonnegative and that the 2-dimensional integral of this joint density over the entire state space is equal to 1.
\item Derive the marginal density $f_Y(y)$.
\item Derive the conditional density $f_{X\mid Y}(x \mid y)$.
\end{enumerate}
\subsubsection*{Guide}
To show that this is a legitimate density observe that the joint density is nonnegative and
also $$\int_0^1 \int_0^{1-x} 2 dydx=1.$$

\begin{equation} f_{Y}(y)= \left\{\begin{array}{ll}
\int_0^{1-y} f_{X,Y}(x,y)dx=2-2y  & \mbox{ $0\leq y \leq 1$}\\
~~~~~~~~~~~~0 & \mbox{~~otherwise.} \end{array} \right. \end{equation}

\begin{equation} f_{X\mid Y}(x \mid y)= \left\{\begin{array}{ll}
\frac{2}{2-2y} = \frac{1}{1-y}  & \mbox{ $0\leq x \leq 1-y $}\\
~~~~~~~0 & \mbox{~~otherwise.} \end{array} \right. \end{equation}
~~~~~~~~~~~~~~~~~~~~~~~~~~~~~~~~~~~~~~~~~~~~~~~~~~~~~~~~~~~~~~~~~~~~~~~~~~~~~~~~~~~~~~~~~~~~~~~~~~~~~~~~~~~~~~~~~~~~~~~~~~~~~~~~~~~~~~~~~~~~~$~~~\Box$

\subsection{Selected Continuous Random Variables}
\label{selectedcont}
We will now discuss several continuous random variables and their
corresponding probability distributions: uniform, exponential, hyper-exponential, Erlang, hypo-exponential
Gaussian, multivariate  Gaussian, and Pareto. These are selected because of their applicability in
teletraffic and related queueing models and consequently their relevance to the material in this book.

\subsubsection{Uniform}
\label{susecuniform}

The probability density function of the uniform random variable takes
nonnegative values over the interval $[a,b]$ and is given by
\begin{equation}
f(x)=\left\{\begin{array}{ll}
\frac{1}{b-a} & \mbox{if $a \leq x \leq b$}\\
0 & \mbox{otherwise.}
\end{array}
\right.
\end{equation}

Of particular interest is the special case - the uniform (0,1)
random variable. Its probability density function is given by
\begin{equation}
f(x)=\left\{\begin{array}{ll}
1 & \mbox{if $0 \leq x \leq 1$}\\
0 & \mbox{otherwise.}
\end{array}
\right.
\end{equation}

The uniform (0,1) random variable is very important in simulations.
Almost all computer languages have a function by which we can
generate uniform (0,1) random variates. By a simple transformation
such uniform (0,1) random variates can be translated to a sequence of
random variates of any distribution as follows. Let $U_1(0,1)$ be
the first uniform (0,1) random variate, and let $F(x)$ be a
distribution function of an arbitrary random variable. Set,
\begin{equation} \label{dev} U_1(0,1)=F(x_1) \end{equation} so $x_1=F^{-1}(
U_1(0,1))$ is the first random variate from the distribution
$F(\cdot)$. Then, generating the second uniform (0,1) random variate,
the second $F(\cdot)$ random number is obtained in the same way,
etc.

This method of generating random variates from any distribution is
known by the following names: inverse transform sampling, inverse
transformation method, inverse probability integral transform, and
Smirnov transform.

To see why this method works, let $U$ be a uniform (0,1) random
variable. Let $F(x)$ be an arbitrary cumulative distribution
function. Let the random variable $Y$ be defined by: $Y=F^{-1}(U)$.
That is, $U=F(Y)$. We will now show that the distribution of $Y$,
namely $P(Y\leq x)$, is equal to $F(x)$. Notice that $P(Y\leq x) =
P[F^{-1}(U) \leq x] = P[U \leq F(x)]$. Because $U$ is a uniform
(0,1) random variable, then $P[U \leq F(x)] = F(x)$. Thus,
$P(Y\leq x) = F(x)$.  $~~~\Box$

Please notice that if $U$ is a uniform (0,1) random variable, then $1-U$ is also a  uniform (0,1) random variable.  Therefore, instead of (\ref{dev}), we can write $$U_1(0,1)=1- F(x_1)=\bar{F}(x_1).$$ In various cases, it is more convenient to use the complementary $\bar{F}(x_1)$ instead of the cumulative distribution function. One of these cases is the exponential distribution as illustrated in the next section.

\subsubsection*{Homework \ref{probability}.\arabic{homework}}

\addtocounter{homework}{1} \addtocounter{tothomework}{1} Let $X_1, X_2, X_3, \ldots X_k$ be a sequence of $k$ independent random variables having a uniform $(0,s)$ distribution. Let $X = \min\{X_1, X_2, X_3, \ldots, X_k\}$. Prove that
\begin{equation} \label{minuniform} P(X>t) = \left\{\begin{array}{ll}
1 & \mbox{for $t \leq 0$}\\
(1-\frac{t}{s})^k & \mbox{for $0 < t < s$}\\
0 & \mbox{otherwise.}
\end{array}
\right. \end{equation}
{\bf Hint:} $P(X>t) = P(X_1>t)P(X_2>t)P(X_3>t) \cdots P(X_k>t).$
    $~~~\Box$

\subsubsection*{Homework \ref{probability}.\arabic{homework}}
\addtocounter{homework}{1} \addtocounter{tothomework}{1}
    Derive the convolution of two independent uniform (0,1) random variables.

\subsubsection*{Guide}

Since $U$ and $V$ are  uniform (0,1) random variables, for the product $f_U(u) f_V(x-u)$ in Eq. (\ref{conv2}) to be non-zero,
$u$ and $x$ must satisfy:\\
$0\leq u \leq 1$ and $0 \leq x-u \leq 1$, \\
or \\$\max(0,x-1) \leq u \leq \min(1,x)$\\ and  $0\leq x\leq 2.$

Therefore,

$$ f_X(x) = \left\{\begin{array}{ll}
u \bigg|_{\max(0,x-1)}^{\min(1,x)}  & \mbox{ $0\leq x \leq 2$} \vspace{4 mm}\\
~~~~~~~0 & \mbox{~~otherwise.} \end{array} \right. $$

or

$$ f_X(x) = \left\{\begin{array}{ll}
{\min(1,x)}-{\max(0,x-1)}   & \mbox{ $0\leq x \leq 2$} \vspace{4 mm}\\
~~~~~~~0 & \mbox{~~otherwise.} \end{array} \right. ~~~\Box $$

\subsubsection{Exponential}
\label{exponential}

The exponential
random variable has one parameter $\mu$ and its probability density function is given by,
\begin{equation}
f(x)=\left\{\begin{array}{ll}
\mu e^{-\mu x} & \mbox{if $x \geq 0$}\\
0 & \mbox{otherwise.} \end{array} \right. \end{equation}
Its
distribution function is given by \begin{equation} F(x)=
\int_{0}^{x} \mu e^{-\mu s}  ds = 1- e^{-\mu x} \qquad x\geq 0.
\end{equation}
A convenient and useful way to describe the exponential random variable is by its complementary distribution
function. It is given by, \begin{equation} \bar{F}(x)= e^{-\mu x}
\qquad x\geq 0. \end{equation}
An important application of the exponential random variable is the time until the next call (or connection request) arrives at a switch.

Interestingly, such time does not depend on how long ago the last call arrived. In other words, the exponential random variable is memoryless. In particular, a continuous random variable is called memoryless if for any $t\geq 0$ and $s\geq 0$,
\begin{equation} \label{memoryless} P(X>s+t\mid X>t)= P(X>s).
\end{equation}
If our lifetime were memoryless, then the probability
we survive at least 80 years given that we have survived 70 years is
equal to the probability that a newborn baby lives to be 10 years old.
Of course, the human lifetime is not memoryless, but, as mentioned above, inter-arrivals of phone calls at a telephone exchange are approximately memoryless. To show that the exponential random variable is memoryless, we show that Eq. ({\ref{memoryless}}) holds using the conditional probability
definition together with the complementary distribution function of
an exponential random variable as follows. \begin{eqnarray*}
P(X>s+t\mid X>t) & = & \frac{ P(X>s+t \cap X>t) }{ P(X>t)} \\
                 & = & \frac{ P(X>s+t)}{ P(X>t)} \\
                 & = & \frac{ e^{-\mu (s+t)}} { e^{-\mu t}}\\
                 & = &  e^{-\mu s} = P(X>s).
\end{eqnarray*}
Not only is the exponential random variable memoryless, but it is also the only memoryless continuous random variable.

\subsubsection*{Homework \ref{probability}.\arabic{homework}}
\addtocounter{homework}{1} \addtocounter{tothomework}{1}

Show how to apply the Inverse transform sampling to generate exponential variates.
\subsubsection*{Guide}

As discussed, we can use the complementary distribution function instead of the cumulative distribution function. Accordingly,
to generate random variates from an exponential distribution with parameter $\lambda$, follow the following steps.
\begin{enumerate}
\item Obtain a new uniform (0,1) variate. Most computer programs provide a function that enables you to do it. Call this variate $U(0,1)$.
\item Write: $U(0,1) = e^{-\lambda x} $.  Then, isolating $x$ in the latter gives: $$ x = \frac{-\ln U(0,1)}{\lambda}. $$ Using this equation, you can obtain an exponential variate $x$ from the uniform (0,1) variate $U(0,1)$.
    \item Repeat steps 1 and 2 for each required exponential variate. $~~~\Box$
    \end{enumerate}

\subsubsection*{Homework \ref{probability}.\arabic{homework}}

\addtocounter{homework}{1} \addtocounter{tothomework}{1}
Assume that the computer language that you use provides a sequence of independent uniform (0,1) variates.
Write computer programs that generate a sequence of 10,000 independent random variates from: \begin{enumerate}
\item  an exponential distribution with $\mu=1$;
\item  a discrete uniform  distribution with $a=1$ and $b=10$;
\item  a geometric distribution with $p=0.4$.
\end{enumerate}
Plot histograms to illustrate the correctness of these three cases.
$~~~\Box$



Let $X_1$ and $X_2$ be independent and exponentially distributed random variables
with parameters $\lambda_1$ and $\lambda_2$, respectively.
 We are interested to know the distribution of $X=\min[X_1,X_2]$.
In other words, we are interested in the distribution of the time
that passes until the first one of the two random variables $X_1$
and $X_2$ occurs. This is as if we have a competition between the
two, and we are interested in the time of the winner, whichever it
is.
 Then, for $t\geq 0$
\begin{equation} \label{Pmin}
P(X>t)=P(\min[X_1 X_2]>t)=P(X_1>t \cap X_2>t)= e^{-\lambda_1 t}
e^{-\lambda_2 t}= e^{-(\lambda_1+\lambda_2) t}. \end{equation}
Thus, the distribution of $X$ is exponential with parameter
$\lambda_1+\lambda_2$.

Another interesting question related to the competition between
two exponential random variables is what is the probability that
one of them, say $X_1$, wins. That is, we are interested in the
probability of $X_1 < X_2$. This is obtained using the continuous
version of the Law of Total Probability (\ref{tot_cont}) as follows:
\begin{equation} \label{Pwin} P(X_1 < X_2) = \int_{0}^{\infty}(1-e^{-\lambda_1 t})\lambda_2 e^{-\lambda_2
t}dt = \frac{\lambda_1}{\lambda_1+\lambda_2}. \end{equation}

To understand the latter, note that $X_2$ can take many
values: $t_1, t_2, t_3, \ldots$, infinitely many values ...

All these values, that $X_2$ may take, lead to the events\\ $X_2=t_1, X_2=t_2, X_2=t_3, \ldots $
that are mutually exclusive and exhaustive.

Then, using the continuous version of the Law of Total Probability, namely,
integration of the product\\ $P(X_1 < t)$ times the density of $X_2$ at $t$, will give us
the probability of $X_1<X_2$. \\By integrating over all $t$ we ``add up'' the probabilities of
infinitely many mutually exclusive and exhaustive events that make up the
event $X_1<X_2$. And this is precisely what the Law of Total Probability does!

In the following table, we point out the equivalence between the corresponding terms in the two equations (\ref{tot_cont}) and (\ref{Pwin}).

\begin{center}
\renewcommand{\arraystretch}{1.4}
\begin{tabular}{|c|c|}  \hline
{\bf term in (\ref{tot_cont})} & {\bf equivalent term in (\ref{Pwin})}  \\  \hline
event $A$ & event $\{X_1 < X_2\}$ \\  \hline
random variable $Y$  & random variable $X_2$ \\ \hline
 event $\{Y=y\}$ &  event $\{X_2=t\}$ \\  \hline
event $\{A \mid Y=y\}$ & event $\{X_1 < t\}$  \\  \hline
$P(A \mid Y=y)$  & $P(X_1<t)=1-e^{-\lambda_1 t}$  \\  \hline
density  $f_{Y}(y)$ & density $f_{X_2}(t) = \lambda_2 e^{-\lambda_2t}$ \\  \hline
\end{tabular}
\end{center}

In a similar way, \begin{equation} P(X_1
> X_2)=\frac{\lambda_2}{\lambda_1+\lambda_2}. \end{equation}
As expected, $P(X_1 < X_2)+P(X_1 > X_2)=1.$ Notice that as $X_1$
and $X_2$ are continuous-time random variables, the probability
that they are equal to each other is equal to zero.

\subsubsection{Relationship between Exponential and Geometric Random
Variables}

We have learned that the geometric random variable is the only
discrete random variable that is memoryless. We also know that the
only memoryless continuous random variable is the exponential random
variable. These facts indicate an interesting relationship between
the two. Let $X_{exp}$ be an exponential random variable with
parameter $\lambda$ and let $X_{geo}$ be a geometric  random
variable with parameter $p$.

Let $\delta$ be an ``interval'' size used to discretize  the
continuous values that $X_{exp}$ takes, and we are interested to
find $\delta$ such that $$F_{X_{exp}}(n\delta) =
F_{X_{geo}}(n),~~~n=1,2,3, \dots . $$ To find such a $\delta$, it is
more convenient to consider the complementary distributions. That
is, we aim to find $\delta$ that satisfies $$P(X_{exp}> n\delta) =
P(X_{geo} > n),~~~ n=1,2,3, \ldots, $$ or $$e^{-\lambda n\delta} =
(1-p)^n,~~~ n=1,2,3, \ldots,  $$ or $$e^{-\lambda \delta} = 1-p.$$
Thus, $$\delta = \frac{-\ln (1-p)}{\lambda} ~~~{\rm and}~~~ p=
1-e^{-\lambda \delta}.$$

We can observe that as the interval size $\delta$ approaches zero,
the probability of success $p$ also approaches zero, and under these
conditions, the two distributions approach each other.

\subsubsection{Hyper-Exponential}

Let $X_i$ for $i=1,2,3, ~\ldots,~ k$ be $k$ independent exponential
random variables with parameters $\lambda_i$, $i=1,2,3,~ \ldots,~ k$,
respectively. Let $p_i$ for $i=1,2,3, ~\ldots,~ k$ be $k$ nonnegative
real numbers such that $\sum_{i=1}^k p_i = 1$. A random variable $X$
that is equal to $X_i$ with probability $p_i$ is called
Hyper-exponential. By the Law of total probability, its density is
\begin{equation}
f_X(x)=\sum_{i=1}^k p_i f_{X_i}(x).
\end{equation}

\subsubsection{Erlang}

A random variable $X$ has an Erlang distribution with parameters
$\lambda$ (positive real) and $k$ (positive integer) if its density
is given by
\begin{equation}
\label{erlang} f_X(x)=\frac{\lambda^k x^{k-1} e^{-\lambda
x}}{(k-1)!}.
\end{equation}

This distribution is named after the Danish mathematician Agner Krarup Erlang (1878 -- 1929) who was the originator
of queueing theory and teletraffic. The term {\it Erlang} is used for other concepts in the field as we will see in later chapters.

\subsubsection*{Homework \ref{probability}.\arabic{homework}}

\addtocounter{homework}{1} \addtocounter{tothomework}{1} Let $X_i$, $i=1,~2,~ \ldots, ~k$, be $k$
independent exponentially distributed random variables each with
parameter $\lambda$, prove by induction that the random variable $X$
defined by the sum $X=\sum_{i=1}^k X_i$  has Erlang distribution
with parameters $k$ and $\lambda$. In other words, $f_X(x)$ of
(\ref{erlang}) is a k-fold convolution of $\lambda e^{-\lambda x}$.
$~~~\Box$

\subsubsection*{Homework \ref{probability}.\arabic{homework}}

\addtocounter{homework}{1} \addtocounter{tothomework}{1}

Let $X_1$ and $X_2$ be independent and Erlang distributed random variables
with parameters $(k,\lambda_1)$ and $(k,\lambda_2)$, respectively.
Find the probability of $P(X_1<X_2)$.

 \subsubsection*{Guide}

 Define the probability of success $p$ as $$p=\frac{\lambda_1}{\lambda_1+\lambda_2}$$ the probability of failure by $1-p$.
 Then, observe that the event $\{X_1<X_2\}$ is equivalent to the event: having $k$ successes before having $k$ failures, so the required probability is the probability of a Pascal random variable $Y$ with parameters $p$ and $k$ to be less or equal to $k-1$. This observation is explained as follows.
 Consider individual points that occur randomly on the time axis starting from $t=0$. The points are of two types 1 and 2. The first type 1 point occurs at time $t_1(1)$ where $t_1(1)$ is exponentially distributed with parameter $\lambda_1$. The second type 1 point  occurs at time $t_2(1)$ where $t_2(1)- t_1(1)$ is also exponentially distributed with parameter $\lambda_1$, and in general, the $n$th type 1 point occurs at time $t_n(1)$ where $t_n(1)- t_{n-1}(1)$ is exponentially distributed with parameter $\lambda_1$. Observe that $t_k(1)$ is a sum of $k$ exponentially distributed random variables with parameter $\lambda_1$, and therefore, it follows an Erlang distribution with parameters $(k,\lambda_1)$ exactly as $X_1$. Equivalently,
  we can construct the time process of type 2 points  where $t_1(2)$ and all the inter-point times $t_n(2)- t_{n-1}(2)$ are exponentially distributed with parameter $\lambda_2$. Then, $t_k(2)$ follows an Erlang distribution with parameters $(k,\lambda_2)$ exactly as $X_2$.

  Accordingly, the event $\{X_1<X_2\}$ is equivalent to the event $\{t_k(1) < t_k(2)\}$. Now consider a traveler that travels on the time axis starting from time 0. This traveler considers type 1 points as successes and type 2 points as failures, where $p$ is the probability that the next point is of type 1 (a success) and $1-p$ is
  the probability that the next point is of type 2 (a failure). The event $\{t_k(1) < t_k(2)\}$  is equivalent to having $k$ successes before having $k$ failures, which lead to the observation that $P(X_1<X_2)$ is the probability of a Pascal random variable $Y$ with parameters $p$ and $k$ to be less or equal to $k-1$.

 Based on this observation the probability $P(X_1<X_2)$ is obtained by (\ref{pascalcdf})
 as
 \begin{equation}
P(X_1<X_2) = F_Y(j) = P(Y \leq k-1) = 1-I_p(k,k).
\end{equation}
where $I_p(\cdot,\cdot)$ is the regularized incomplete beta function.
$~~~\Box$

\subsubsection*{Homework \ref{probability}.\arabic{homework}}

\addtocounter{homework}{1} \addtocounter{tothomework}{1}

Again, consider the two independent and Erlang distributed random variables $X_1$ and $X_2$ with parameters $(k,\lambda_1)$ and $(k,\lambda_2)$, respectively. Assume $\lambda_1 < \lambda_2$. Investigate the probability $P(X_1<X_2)$ as $k$ approaches infinity.
Use numerical, intuitive, and rigorous approaches.
$~~~\Box$

\subsubsection{Hypo-Exponential}

Let $X_i$, $i=1,~2,~ \ldots, ~k$ be $k$ independent exponentially
distributed random variables each with parameters $\lambda_i$,
respectively. The random variable $X$, defined by the sum $X=\sum_{i=1}^k X_i$,  is called hypo-exponential. In other words,
the density of $X$ is a  convolution of the $k$ densities $\lambda_i
e^{-\lambda_i x}$, $i=1,~2,~ \ldots, ~k$. The Erlang distribution is
a special case of hypo-exponential when all the $k$ random variables
are identically distributed.

\subsubsection{Gaussian}
\label{gaussrv}

A continuous random variable, which is commonly used in many applications, is the Gaussian (also called Normal) random variable.
We say that the random variable $X$ has a Gaussian distribution with
parameters $m$ and $\sigma^2$ if its density is given by
\begin{equation}
f_X(x)=\frac{1}{\sqrt{2\pi}\sigma} e^{-(x-m)^2/2\sigma^2} \qquad
-\infty < x <\infty.
\end{equation}
This density is symmetric and bell-shaped.

The wide use of the Gaussian random variable is rooted in the so-called
{\bf The central limit theorem}. This theorem is the most important
result in probability
theory. Loosely speaking, it says that the sum of a large number
of independent random variables (under certain conditions)
 has Gaussian (normal) distribution.
This is also true if
the distribution of these random variables is very different from Gaussian.
This theorem explains why so
many populations in nature and society have bell-shaped Gaussian
histograms and justifies the use of the Gaussian distribution as their model.
In Section \ref{limittheorems} we will further discuss the central limit theorem and demonstrate its applicability to the telecommunication
link dimensioning problem in Section \ref{link_dim}.

\subsubsection{Pareto}

Another continuous random variable often used in telecommunication
modeling is the {\bf Pareto} random variable. This random
variable, for a certain parameter range, can be useful in
modeling lengths of data bursts in data and multimedia networks
\cite{addie09,addie02,zukerman03}. We choose to define the Pareto random variable
with parameters $\gamma$ and $\delta$ by its complementary
distribution function which is given by
\[
 P \left(X>x\right) =  \left\{ \begin{array}{ll}
\left(\frac{x}{\delta}\right)^{-\gamma}, & x\ge \delta \\
1, & {\rm otherwise.} \end{array}  \right.
\]

Here $\delta > 0$  is the scale parameter representing a minimum
value for the random variable, and $\gamma > 0$ is the shape
parameter of the Pareto distribution.

\subsubsection*{Homework \ref{probability}.\arabic{homework}}

\addtocounter{homework}{1} \addtocounter{tothomework}{1} Write a computer program that generates a
sequence of 10,000 random variates from a Pareto distribution with $\gamma=1.2$ and $\delta=4$. Plot a histogram to illustrate the correctness of your results. $~~~\Box$

\subsection{Moments and Variance}
\label{moments}
\subsubsection{Mean (or Expectation)}
The {\em mean} (or the expectation) of a discrete random variable is
defined by
\begin{equation} \label{disexpec}
E[X]=\sum_{\{n:P(n)>0\}} nP_X(n).
\end{equation}

This definition can be intuitively explained by the concept of {\it limiting relative frequencies} illustrated by the following example.
Let the outcome of an experiment be a person's height.  Let $p_i$ be the probability that a person's height is $i$ cm.
Consider a sample of $N$ people.  Each one of them reports
his/her height, rounded to the nearest cm.  Let $n_i$ be the number of people reported
a height of $i$ cm.
These $n_i$ values can be graphically presented in a histogram. Then, the relative frequency $n_i/N$ approximates $p_i$.
This approximation becomes more and more accurate as $N$ increases.
This approximation is consistent with the requirement $$\sum_i p_i = 1.$$
If we set $p_i=n_i/N$, then since $\sum_i n_i = N$, we obtain
$$\sum_i p_i = 1.$$

Then, the average height is given by $$\frac{\sum_i in_i}{N}. $$

For large $N$, we set $$p_i= \frac{n_i}{N}$$.
Substituting in the above, we obtain
$${\rm The~ average~ height} = {\sum_i ip_i}$$
which is the definition of mean.

This is related to the Weak Law of Large Numbers discussed below in Section \ref{limittheorems}.

\subsubsection*{Homework \ref{probability}.\arabic{homework}}
\addtocounter{homework}{1} \addtocounter{tothomework}{1}
Consider a discrete random variable $X$ that takes the values 1, 2 and 3  with the following probabilities $P(X=1)=0.5$, $P(X=2)=0.3$, $P(X=3)=0.2$, and it takes all other values with probability zero. Find the mean of $X$.

\subsubsection*{Answer}
$E[X]=1.7$.
$~~~\Box$

\subsubsection*{Homework \ref{probability}.\arabic{homework}}
\addtocounter{homework}{1} \addtocounter{tothomework}{1}
In an American Roulette game, players place their bets, and  a ball lands on one of 38 numbered pockets in a turning wheel with the same probability of 1/38. In all the questions about Roulette, we assume that the number that the ball lands on is independent of the numbers it landed on in any other previous or future games (wheel turnings). A player has various betting options to choose from. Any of these options can be represented by a set of $n$ numbered pockets. If the player bets $B$ dollars and if the ball lands on any of the numbers in the chosen set, the payout given to the player is given by
$$\frac{1}{n}(36-n)B = \frac{36B}{n} -B ~~~ {\rm dollars}.$$
If the ball does not land in any of the chosen numbers the player loses the $B$ dollars bet.

For example, if $B=1$ and $n=18$, the player wins one dollar payout. Roulette has 18 red pockets and 18 black pockets, and there are two additional green pockets for zero and double-zero. One way to bet on a chosen set of 18 numbers is to bet on the red (on all 18 red-colored pockets).

Find the mean return for the player who bets $B=1$ in one game for a relevant range of $n$ values.

For more information on Roulette
and on the mean return, see:


{\url{http://en.wikipedia.org/wiki/Roulette}}



\subsubsection*{Guide}
For the case $B=1$ and $n=18$, with probability 18/38, the player wins one dollar, and with probability 20/38, the player loses one dollar. Let the random variable $X$ represent the return in one game. Therefore, the mean return is given by

\begin{equation}
\label{mean_win}
E[X] = (1)\frac{18}{38} + (-1) \frac{20}{38} = -0.052631579. \end{equation}

Repeat this calculation for various $n$ values.
$~~~\Box$

Equivalently, the mean of a continuous random variable is defined as
\begin{equation}
E[X]=\int_{-\infty}^{\infty} xf_X(x)dx.
\end{equation}
A useful expression for the mean of a continuous nonnegative
random variable $Z$ (i.e. a random variable $Z$ with the property
that its density $f(z)=0$ for $z<0$) is:
\begin{equation}
\label{ezc} E[Z]=\int_{0}^{\infty} P(Z>z)dz=\int_{0}^{\infty}
[1-F_Z(z)]dz.
\end{equation}
The discrete equivalence of the latter is:
\begin{equation}
\label{ezd} E[Z]=\sum_{n=0}^{\infty} P(Z>n)=\sum_{n=0}^{\infty}
[1-F_Z(n)].
\end{equation}

\subsubsection*{Homework \ref{probability}.\arabic{homework}}

\addtocounter{homework}{1} \addtocounter{tothomework}{1} Show
(\ref{ezc}) and (\ref{ezd}).

\subsubsection*{Guide}
For the case of discrete nonnegative random variables, by (\ref{disexpec}),
$$E[Z]=\sum_{n=1}^\infty nP_Z(n). $$

This can be written as

$$
E[Z]=\left\{\begin{array}{lllll}
 P_Z(1) & +  P_Z(2) & +  P_Z(3) & +  P_Z(4) & +  \dots \\
  & +  P_Z(2) & +  P_Z(3) & +  P_Z(4) & +  \dots \\
   & & +  P_Z(3) & +  P_Z(4) & +  \dots \\
     &  &  & +  P_Z(4) & +  \dots \\
      &  &  &  & +  \dots
 \end{array} \right. $$

Therefore,

$$ E[Z] = \sum_{n=1}^\infty P_Z(n) + \sum_{n=2}^\infty P_Z(n) +  \sum_{n=3}^\infty P_Z(n) +  \sum_{n=4}^\infty P_Z(n) + ~~\dots~. $$
and (\ref{ezd}) follows.

Another way to write this proof is \cite{Whitt12}

\begin{eqnarray*}
E[Z] & = & \sum_{n=1}^\infty nP_Z(n) \\
  &   = & \sum_{n=1}^\infty \left( \sum_{j=1}^n 1 \right)P_Z(n) \\
     &  = & \sum_{j=1}^\infty  \sum_{n=j}^\infty 1 P_Z(n) \\
      & = & \sum_{n=0}^{\infty} P(Z>n).
 \end{eqnarray*}

To show (\ref{ezc}) use  an analogous approach (defining $y=f_Z(z)$) as follows \cite{Whitt12}

\begin{eqnarray*}
E[Z] & = & \int_{0}^{\infty} zf_Z(z)dz\\
& = & \int_{0}^{\infty} \left( \int_0^z 1 dy \right) f_Z(z)dz\\
& = & \int_{0}^{\infty} 1 \left( \int_y^\infty  f_Z(z) dz \right) dy\\
& = & \int_{0}^{\infty} P(Z>y)dy\\
& = & \int_{0}^{\infty} P(Z>z)dz.
 \end{eqnarray*}
Check carefully and understand all the steps and operations that led to the above proofs.
 $~~~\Box$

\subsubsection*{Homework \ref{probability}.\arabic{homework}}

\addtocounter{homework}{1} \addtocounter{tothomework}{1} Let $X_1, X_2, X_3, \ldots X_k$ be a
sequence of $k$ independent random variables having a
uniform (0,s) distribution. Let $X = \min\{X_1, X_2, X_3, \ldots, X_k\}$. Prove that $$E[X]=\frac{s}{k+1}.$$ {\bf Hint:}
Use (\ref{minuniform}) and (\ref{ezc}). $~~~\Box$

As mentioned above, a function of a random variable is also a random
variable. The mean of a function of random variables denoted $g(\cdot)$ by
\begin{equation}
\label{gxdis} E[g(X)]=\sum_{\{k:P_X(k)>0\}} g(k)P_X(k)
\end{equation}
for a discrete random variable
and
\begin{equation}
\label{gxcont} E[g(X)]=\int_{-\infty}^{\infty} g(x)f_X(x)dx
\end{equation}
for a continuous random variable.
If $a$ and $b$ are constants, then for a random variable $X$ (either
discrete or continuous) we have:
\begin{equation}
\label{Eax}
E[aX]= aE[X],
\end{equation}
\begin{equation}
E[X-b]= E[X]-b,
\end{equation}
and
\begin{equation}
E[aX-b]= aE[X]-b.
\end{equation}

\subsubsection{Moments}

The {\bf $n$th moment} of the random variable $X$ is defined by
$E[X^n]$. Substituting $g(X)=X^n$ in (\ref{gxdis}) and in
(\ref{gxcont}), the $n$th moment of $X$ is given by:
\begin{equation}
E[X^n]=\sum_{\{k:P_X(k)>0\}} k^n P_X(k)
\end{equation}
for a discrete random variable
and
\begin{equation}
E[X^n]=\int_{-\infty}^{\infty} x^n f_X(x)dx
\end{equation}
for a continuous random variable. Similarly, the $n$th central
moment of the random variable $X$ is defined by
 $E[(X-E[X])^n]$. Substituting $g(X)=(X-E[X])^n$ in (\ref{gxdis})
and in (\ref{gxcont}), the {\bf $n$th central moment} of $X$ is
given by:
\begin{equation}
E[(X-E[X])^n]=\sum_{\{k:P(k)>0\}} (k-E[X])^n P_X(k)
\end{equation}
for a discrete random variable
and
\begin{equation}
E[(X-E[X])^n]=\int_{-\infty}^{\infty} (x-E[X])^n f_X(x)dx
\end{equation}
for a continuous random variable.
By definition, the first moment is the mean.
\subsubsection{Variance}

The second central moment is
called
the {\bf variance}. It is defined as
\begin{equation}
\label{vardef}
Var[X]=E[(X-E[X])^2].
\end{equation}
The variance of a random variable $X$ is
given by
\begin{equation}
Var[X]=\sum_{\{k:P(k)>0\}} (k-E[X])^2 P_X(k)
\end{equation}
if $X$ is discrete, and by
\begin{equation}
Var[X]=\int_{-\infty}^{\infty} (x-E[X])^2 f_X(x)dx
\end{equation}
if it is continuous.

By (\ref{vardef}) we obtain
\begin{equation}
\label{vardef2}
Var[X]=E[(X-E[X])^2] = E[X^2 - 2XE[X] + (E[X])^2] = E[X^2] - (E[X])^2.
\end{equation}

To understand the latter, notice that
$$E[X^2 - 2XE[X] + (E[X])^2] = E[X^2] - E[2XE[X]] + E[(E[X])^2] = E[X^2] - 2E[X]E[X] + (E[X])^2. $$

 \subsubsection*{Homework \ref{probability}.\arabic{homework}}
\addtocounter{homework}{1} \addtocounter{tothomework}{1}
Consider two independent and identically distributed random variables $X$ and $Y$ that obey the following probability function:

\begin{equation} P(X=i)=\left\{\begin{array}{ll}
1/5 & \mbox{ ~~for~~$i=-2, -1, 0, 1, {\rm and} ~2$ }\\
~~~~0 & \mbox{otherwise.} \end{array} \right. \end{equation}

Let $Z= |X|$,    $U=\max(X,Y)$ and $V=\min(X,Y)$.
Find the probability function, the mean and variance of $Z, U$, and $V$.
$~~~\Box$

 \subsubsection*{Homework \ref{probability}.\arabic{homework}}
\addtocounter{homework}{1} \addtocounter{tothomework}{1}
For random variable $X$ and constant $c$, show the following:  $$Var[X+c] = Var[X]$$ and $$Var[cX] = c^2 Var[X].$$

 \subsubsection*{Guide}
 $
 Var[X+c]= E[(X+c)^2] - (E[X+c])^2=E[X^2+2Xc+c^2]-(E[X])^2-2cE[X] - c^2$\\$= E[X^2] - (E[X])^2=Var[X]$

 and

 $Var[cX] =E[(cX)^2] - (E[cX])^2=  c^2E[X^2] - c^2(E[X])^2= c^2 Var[X].$
$~~~\Box$

While the mean provides the average according to the limiting relative frequencies concept, the variance is a measure of the
level of variation of the possible values of the random variable.
Another measure of such variation is the {\bf standard deviation}
denoted $\sigma_X$, or simply $\sigma$, and defined by
\begin{equation}
\sigma_X = \sqrt{Var[X]}.
\end{equation}
Hence the variance is often denoted by $\sigma^2$.

Notice that the first central moment $E[x-E[X]]$ is not very useful
because it is always equal to zero. Therefore, the second central moment
$E[(x-E[X])^2]$, which is the variance, and its square root, the
standard deviation, are used for measuring the level of variation of
a random variable.

The standard deviation of the random variable $X$, namely $\sigma_X$, has always the same units as the random variable $X$. For example, if the random variable $X$ represents the number of bytes in an IP packet, then $\sigma_X$ is also in bytes. If $X$ represents the service time in seconds of an IP packet, then $\sigma_X$ is also in seconds. Since $Var[X] = \sigma_X^2$, the units of the variance are in the units of $\sigma_X^2$ or $X^2$. Thus, for the above two examples, of being in bytes or in seconds, the units of variance are ${\rm bytes}^2$ or in ${\rm seconds}^2$, respectively.

The mean of the sum of random variables is always the sum of their means, namely,
\begin{equation}
E\left[\sum_{i=1}^n X_i\right] = \sum_{i=1}^n E[X_i]
\end{equation}
but the variance of a sum of random variables is not always equal to the sum of
their variances. This is true for independent random variables. That is, if the
random variables $X_1,~ X_2, ~X_3, ~ \ldots, ~X_n$ are independent, then
\begin{equation}
Var\left[\sum_{i=1}^n X_i \right] = \sum_{i=1}^n Var[X_i].
\end{equation}
Also, if $X_1,~ X_2, ~X_3, ~ \ldots, ~X_n$ are independent, then
\begin{equation}
\label{etimes}
E[\Pi_{i=1}^n X_i] = \Pi_{i=1}^n E[X_i].
\end{equation}

\subsubsection*{Homework \ref{probability}.\arabic{homework}}

\addtocounter{homework}{1} \addtocounter{tothomework}{1}
\label{HWnumber}
Consider an experiment of tossing a 6-sided die. Assume that
the die is fair, i.e., each side has the same probability (1/6) to
occur.  Consider a random variable $X$ that takes the value $i$ if
the outcome of the toss is $i$, for $i=1,2,3, \cdots, 6.$ Find
$E[X]$, $Var[X]$ and $\sigma_X$.

 \subsubsection*{Answers}

$E[X] = 3.5$; $E[X^2] = 15.16666667 $; $Var[X] =  2.916666667 $;
$\sigma_X = 1.707825128.$
(See Eq. (\ref{vardu}) below.) $~~~\Box$


\subsubsection{Conditional Mean and the Law of Iterated Expectation}

In many applications, it is helpful to use the concept of {\bf
Conditional Expectation (or Mean)} to derive moments of unknown
distributions.

Let $E[X \mid Y]$ be the conditional expectation of
random variable $X$ given the event $\{Y=y\}$ (random variable $Y$ is equal to $y$) for each relevant value of $y$.

The conditional expectation of two discrete random variables is defined by
\begin{equation}
\label{condindis} E[X \mid Y=j] = \sum_{i} i P (X=i \mid Y=j ).
\end{equation}
If $X$ and $Y$ are continuous, their conditional expectation is
defined as
\begin{equation}
\label{condincont} E[X \mid Y=y] = \int_{x=-\infty}^{\infty}
x f_{X \mid Y} (x \mid y) dx.
\end{equation}
It is important to realize that $E[X \mid Y]$ is a random variable
which is a function of the random variable $Y$. Therefore, if we
consider its mean (in the case that $X$ and $Y$ are discrete) we obtain
\begin{eqnarray}
E_Y [ E [X\mid Y]] & = &   \sum_{j} E [X \mid Y=j ] P(Y=j)   \nonumber \\
& = & \sum_{j} \sum_{i} i P (X=i \mid Y=j ) P(Y=j)   \nonumber \\
& = & \sum_{i} i \sum_{j} P (X=i \mid Y=j ) P(Y=j)   \nonumber \\
& = & \sum_{i} i P (X=i ) = E[X].
\end{eqnarray}
Thus, we have obtained the following formula for the mean $E[X]$
\begin{equation}
\label{meancondind} E[X] =  E_Y [ E [X\mid Y]].
\end{equation}
The latter is called {\it the law of iterated expectations}. It is also known by many other names, including the law of iterated expectations,
the law of total expectation, and Adam's law.

The law of iterated expectations also applies to continuous random variables. In this case, we have:
\begin{eqnarray*}
E_Y [ E [X\mid Y]] & = &  \int_{y=-\infty}^{\infty} E [X \mid Y=y ] f_Y(y)dy  \\
& = &  \int_{y=-\infty}^{\infty} \int_{x=-\infty}^{\infty}  x f_{X \mid Y} (x \mid y) dx f_Y(y)dy  \\
& = &  \int_{x=-\infty}^{\infty} x \int_{y=-\infty}^{\infty}  f_{X \mid Y} (x \mid y) f_Y(y)dy dx  \\
& = &  \int_{x=-\infty}^{\infty} x f_X (x) dx = E[X].
\end{eqnarray*}

\subsubsection*{Homework \ref{probability}.\arabic{homework}}

\addtocounter{homework}{1} \addtocounter{tothomework}{1} Show that $E[X] =  E_Y [ E [X\mid Y]]$
holds also for the case where $X$ is discrete and $Y$ is continuous
and vice versa. $~~~\Box$

{\bf Remark:} $P (X=x \mid Y=y)$ is itself a random variable that is a function of the values $y$ taken by the random variable $Y$.  Therefore, by definition
$ E_Y[P (X=x \mid Y=y)] = \sum_y P (X=x \mid Y=y) P(Y=y)$ which leads to another way to express the Law of Total Probability:
\begin{equation}
\label{lotpmean} P_X(x) = E_Y[P (X=x \mid Y=y)] .
\end{equation}

\subsubsection{Conditional Variance and the Law of Total Variance}

Define the {\bf Conditional Variance} as
\begin{equation}
\label{varcondind} Var[X \mid Y] = E [(X - E[X\mid Y])^2 \mid Y ].
\end{equation}
This gives rise to the following useful formula for the variance of
a random variable known as {\it The Law of Total Variance}:
\begin{equation}
\label{evve} Var[X] = E[Var[X\mid Y]] + Var[E [ X\mid Y]].
\end{equation}
Due to its form, it is also called {\it EVVE}, or {\it Eve's law}. In this book, we will mainly refer to it as EVVE.

To show EVVE, we recall (\ref{vardef2}): $Var[X] = E[X^2] -
(E[X])^2$, and (\ref{meancondind}): $E[X] =  E_Y [ E [X\mid Y]]$, we
obtain
\begin{equation}
\label{evvestep1} Var[X] = E[E[X^2\mid Y]] - (E[E[X\mid Y]])^2.
\end{equation}
Then, using $ E[X^2]  = Var[X] + (E[X])^2$ gives
\begin{equation}
\label{evvestep2} Var[X] = E[Var[X \mid Y] + (E[X \mid Y])^2 ] - (E[E[X\mid Y]])^2
\end{equation}
or
\begin{equation}
\label{evvestep3} Var[X] = E[Var[X \mid Y]] + E[E[X \mid Y]]^2  - (E[E[X\mid Y]])^2.
\end{equation}
Now considering again the formula $Var[X] = E[X^2] - (E[X])^2$, but instead of the random variable
$X$ we put the random variable $ E[X\mid Y]$, we obtain
\begin{equation}
\label{evvestep4}
Var[E[X\mid Y]] = E[E[X \mid Y]^2]  - (E[E[X\mid Y]])^2,
\end{equation}
observing that the right-hand side of (\ref{evvestep4}) equals to
the last two terms in the right-hand side of (\ref{evvestep3}), we
obtain EVVE.

To illustrate the use of conditional mean and variance, consider
the following example. Every second, the number of Internet flows
that arrive at a router, denoted $\phi$, has mean $\phi_e$ and variance $\phi_v$. The
number of packets in each flow, denoted $\varsigma$, has mean $\varsigma_e$ and variance $\varsigma_v$.
Assume that the number of packets in each flow and the number of
flows arriving per second are independent. Let $W$ be the total number of packets
arriving at the router per second, which has mean $W_e$ and variance
$W_v$. Assume $W=\varsigma \phi.$
The network designer,
aiming to meet certain quality of service (QoS) requirements,
makes sure that the router serves the arriving packets at the rate
of $s_r$ per second, such that $s_r=W_e + 4 \sqrt{W_v}$.
 To compute $s_r$
one needs to have the values of $W_e$ and $W_v$.
Because $\phi$ and $\varsigma$ are independent
$E[W|\phi]=\phi \varsigma_e$
and by (\ref{meancondind})
$$W_e=E[W]=E[E[W|\phi]]=E[\phi]E[\varsigma]=\phi_e\varsigma_e.$$
Note that the relationship
\begin{equation}
\label{pe} W_e=\phi_e\varsigma_e
\end{equation}
is also obtained
directly by (\ref{etimes}). In fact, the above proves (\ref{etimes}) for the case
of two random variables.

Also, by EVVE,
$$Var[W]=E[Var[W|\phi]] + Var[E[W|\phi]]
=\varsigma_v E[\phi^2]+(\varsigma_e)^2Var[\phi].$$
Therefore,
\begin{equation}
\label{pv} W_v=\phi_v \varsigma_v + \varsigma_v \phi_e^2  + \phi_v \varsigma_e^2.
\end{equation}

\subsubsection*{Homework \ref{probability}.\arabic{homework}}
\addtocounter{homework}{1} \addtocounter{tothomework}{1}
\begin{enumerate}
\item Provide detailed derivations of Equations
(\ref{pe}) and (\ref{pv}) using (\ref{meancondind}) and
(\ref{evve}).

\subsubsection*{Guide} Observe that because $W=\phi \varsigma$, derivation of $Var[W|\phi]$ is obtained by $Var[\phi \varsigma]$ when $\phi$ is a given constant, which is simplified by  $$Var[\phi \varsigma]=\phi^2 Var[\varsigma]=\phi^2 \varsigma_v.$$
Therefore, $$Var[W|\phi] = \phi^2 \varsigma_v.$$
Now it is clear that $Var[W|\phi]$ is a random variable which is a function of the random variable $\phi$ (where $\varsigma_v$ is clearly a constant). Then, we can obtain the mean of $Var[W|\phi]$, by
$$E[Var[W|\phi]]=E[\phi^2 \varsigma_v]=\varsigma_v E[\phi^2]. $$
The last equality is based on Equation (\ref{Eax}).

To obtain $E[\phi^2]$, notice that
$$ \phi_v = Var[\phi] = E[\phi^2] - (E[\phi])^2 = E[\phi^2] - \phi_e^2.$$
Therefore, $$E[\phi^2]= \phi_v + \phi_e^2. $$

Similarly, the derivation of $E[W|\phi]$ is obtained by $E[\phi \varsigma]$ when $\phi$ is a given constant, which is simplified by  $$E[\phi \varsigma]=\phi E[\varsigma]=\phi \varsigma_e.$$

Accordingly, $$E[W|\phi] = \phi \varsigma_e.$$
 In a similar way to the previous elaboration, it is clear that $E[W|\phi]$ is a random variable which is a function of the random variable $\phi$ (where $\varsigma_e$ is  a constant). Then, we can obtain the variance of $E[W|\phi]$, by
$$Var[E[W|\phi]]=Var[\phi \varsigma_e]=\varsigma_e^2 Var[\phi]=\varsigma_e^2 \phi_v. $$

Now put it all together, and obtain Equation (\ref{pv}).

\item Derive Equations (\ref{pe}) and (\ref{pv}) in
a different way, considering the independence of the number of
packets in each flow and the number of flows arriving per second.
$~~~\Box$
\end{enumerate}

\subsubsection*{Homework \ref{probability}.\arabic{homework}}
\addtocounter{homework}{1} \addtocounter{tothomework}{1}
Let $Y$ be a random variable representing the outcome of rolling a 6-sided die. Namely, $Y$ is a discrete uniform random variable with parameters $a=1$ and $b=6$, so $P(Y=i)=1/6$ for $i=1, 2, 3,~ \ldots,~ n$. If $Y=i$, then the random variable $X$ is defined as the sum of $i$ independent Bernoulli random variables, each of which with parameter $p=0.5$. That is, $X$ is a binomial random variable with parameters $n=i$ and $p=0.5$. Find the mean and variance of $X$.
\subsubsection*{Guide}

$$E[Y]=3.5.$$
$$Var(Y) = 35/12=2.916666667.$$
(See Homework 1.29, and Eq. (\ref{vardu}) below. )
$$E[X \mid Y]=0.5Y$$ and $$Var[X \mid Y]=0.25Y.$$
Recall that the variance of a Bernoulli random variable with parameter $0<p<1$ is equal to $p(1-p)$. In our case, $Var[X \mid Y]$ is a variance of a sum of $Y$ independent Bernoulli random variables, each of which with parameter $p=0.5$.

To obtain $E[X]$, we use the Law of Iterated Expectations. Namely, $$E[X]= E[E[X \mid Y]]=E[0.5Y]=0.5E[Y].$$

Notice that $$Var[E[X \mid Y]]=Var[0.5Y]=(0.5)^2Var[Y],$$ and $$E[Var[X \mid Y]]=E[0.25Y]=0.25E[Y].$$

Then, by EVVE,

$$Var[X] = E[Var[X\mid Y]] + Var[E [ X\mid Y]].$$

$~~~\Box$

\subsubsection*{Homework \ref{probability}.\arabic{homework}}
\addtocounter{homework}{1} \addtocounter{tothomework}{1}

A traveler considering buying travel insurance would like to assess the statistics of the cost that may be associated with unforeseen insurable incidents or incidents that may be associated with his/her planned travel. The probability that the traveler experiences incident(s) is $p_{inc}$. That is, the probability that nothing happens during travel is $1-p_{inc}$, and the probability that some incidents (e.g., accident, theft, illness, etc.) happen is $p_{inc}$. If incidents happen, the total unforeseen cost of the
incidents has a mean of $M_{inc}$ and a standard deviation of $\sigma_{inc}$. In this document, the term ``incidents'' refers to one or more incidents.
What is the mean loss and the standard deviation of the loss incurred by the traveler?

\subsubsection*{Guide}

Let $Y$ be a random variable that takes the value 1 if accidents happen and the value 0 if no incidents happen. Accordingly, $Y$ is Bernoulli distributed with parameter $p_{inc}$.
That is,  $P(Y=1)=p_{inc}$, and $P(Y=0)=1-p_{inc}$.
Let $X$ be a random variable representing the total cost of the incidents. We know that $X$ depends on $Y$. We need to obtain $E[X]$ and $Var[X]$, which are obtainable by the law of iterated expectations and EVVE, respectively.

Therefore, according to these definitions and the given parameters, we obtain,

$$E_X[X\mid Y=1]=M_{inc},$$

and

$$E_X[X\mid Y=0]=0,$$

which imply

$$E_X[X\mid Y]=M_{inc}Y.$$

The latter makes it very clear that $E_X[X\mid Y]$ is a random variable which is a function of the random variable $Y$.

In addition, we obtain,

$$Var_X[X\mid Y=1]=\sigma_{inc}^2,$$
and
$$Var_X[X\mid Y=0]=0.$$

Also, since $Y$ is a Bernoulli random variable with parameter $p_{inc}$, its variance is equal to $(1-p_{inc})p_{inc}$.

Then, by the law of iterated expectations,

\begin{eqnarray*}
E[X] & = &  E_Y [ E [X\mid Y]] \\
& = &  P(Y=1) E_X[X\mid Y=1] + P(Y=0) E_X[X\mid Y=0]  \\
 & = &  p_{inc}M_{inc} + 0 \\
 & = &  p_{inc}M_{inc},
\end{eqnarray*}

and by EVVE,

\begin{eqnarray*}
Var[X] & = &  E[Var[X\mid Y]] + Var[E [ X\mid Y]]  \\
& = & P(Y=1)Var[X\mid Y=1] + P(Y=0) Var[X\mid Y=0] + Var[M_{inc}Y] \\
& = & p_{inc}\sigma_{inc}^2 + 0 + M_{inc}^2 p_{inc}(1-p_{inc}) \\
& = & p_{inc}\sigma_{inc}^2 + M_{inc}^2 p_{inc}(1-p_{inc}). ~~~\Box
\end{eqnarray*}

\subsubsection*{Homework \ref{probability}.\arabic{homework}}
\addtocounter{homework}{1} \addtocounter{tothomework}{1}

Now consider the following special case for the previous problem: $p_{inc}=0.04$, $\sigma_{inc}=4000$ dollars, and $M_{inc}=3000$ dollars. Find the mean, variance, and standard deviation of the cost due to unforeseen insurable incidents.

\subsubsection*{Answers}
mean: 120 dollars; variance: 985600 [dollar$^2$]; standard deviation: 992.77 dollars (rounded).
$~~~\Box$


\subsubsection*{Homework \ref{probability}.\arabic{homework}}
\addtocounter{homework}{1} \addtocounter{tothomework}{1}

Consider the above two problems and assume that if incidents occur and the traveler does not have insurance, then the traveler must pay the cost incurred as a result of the incidents. Otherwise, if the traveler has insurance, the insurance company pays it fully. Also, assume that the event that incidents happen to a traveler is independent of whether or not the traveler buys insurance. In this way, if a traveler buys insurance, the  cost to the traveler has a mean that is equal to the insurance premium and a zero variance.
John is a traveler that considering buying travel insurance.
The insurance premium costs 200 dollars, but considering his time associated with claiming back the insurance, he is willing to buy travel insurance only if the mean plus half the standard deviation of the total costs of the insurable incidents is higher than 300 dollars. Will he buy or not buy the travel insurance?

\subsubsection*{Guide}

To answer this question, notice that the mean cost is 120 dollars and half the standard deviation is $992.77/2= 496.39$ dollars, so mean + 0.5 $\times$ standard deviation is equal to $616.39$ dollars which is more than $300$ dollars, so he will buy the travel insurance.
$~~~\Box$

\subsubsection*{Homework \ref{probability}.\arabic{homework}}
\addtocounter{homework}{1} \addtocounter{tothomework}{1}

Now consider an insurance company that has sold $N_p$  insurance policies as described in the above three problems. Assume $N_p = 1,000$.  Also, assume that the incidents that occur for different customers are independent.  What is the mean and standard deviation of the total claims payouts of the insurance company?
Compare them to $200N_p$, which is the income generated by the premium paid by the customers.
What is the mean and standard deviation of the difference between the income generated by the premium paid by the customers, namely $200N_p$, and the total claims payouts of the insurance company? What do you think is the distribution of this difference? What is the probability that this difference is negative? What is the probability is less than half of its mean?

\subsubsection*{Answers}

The mean of the total claims payouts is $p_{inc} M_{inc} N_p =0.04 \times 3000 \times 1000 = 120,000 $

The variance of the total claims payouts is $N_p Var[X] = 1,000 \times 985,600 = 985,600,000 $.

The standard deviation of the total claims payouts is $\sqrt{985,600,000}=31,394.267$ dollars.

Notice that the mean + $0.5 \times$ standard deviation of the payouts is $120,000 + 31,394.267/2 = 135,697.13$ dollars is far smaller than the income of $200N_p=200,000$ dollars.

The mean of this difference is $N_p(200 - E[X]) = 1,000(200 -120) = 80,000$ dollars.

The variance of the difference is $N_p Var[X] = 1,000 \times 985,600 = 985,600,000 $ dollars$^2$.

The standard deviation of the difference is $\sqrt{985,600,000}=31,394.267$ dollars.

Notice that the results for the variance and standard deviation of the difference are the same as the results for those of the total claims. This is because
$$Var[X+c]=Var[X].$$

Based on the central limit theorem, the distribution of the difference between income and payouts is approximately normal with the above mean and standard deviation.

The probability that this difference is negative is approximately equal to 0.005413 using the Excel function: =NORM.DIST(0,80000,31394.267,TRUE).

The probability that this difference is less than 40,000  is approximately equal to 0.1013 using the Excel function: =NORM.DIST(40000,80000,31394.267,TRUE).
$~~~\Box$

\subsubsection*{Homework \ref{probability}.\arabic{homework}}
\addtocounter{homework}{1} \addtocounter{tothomework}{1}

Now consider the case, $N_p = 10,000$, and answer all the questions of the previous Problem. Then, compare the two cases and discuss the economy of scale implications.

\subsection*{Answers}

The mean of this difference is $N_p(200 - E[X]) = 10,000(200 -120) = 800,000$ dollars.

The variance of the difference is $N_p Var[X] = 10,000 \times 985,600 = 9,856,000,000 $ dollars$^2$.

The standard deviation of the difference is $\sqrt{9,856,000,000}=99,277.39$ dollars.

Based on the central limit theorem, again the distribution of the difference between income and payouts is approximately normal with the above mean and standard deviation.

The probability that this difference is negative is approximately equal to $3.9 \times 10^{-16}$ using the Excel function: =NORM.DIST(0,800000,99277.39,TRUE).

The probability that this difference is less than 40,000  is approximately equal to $2.8 \times 10^{-5}$ using the Excel function: =NORM.DIST(400000,800000,99277.39,TRUE).

We observe that for the larger insurance company with 10,000 customers, it is easier to guarantee that the difference between income and payouts will not be negative and will be over half its mean value - in this case, over 40 dollars per policy. Notice also that the payouts are mean $\pm$ a certain constant $\times$ the standard deviation. Therefore the standard deviation to mean ratio can be viewed as a measure of the risk or variability of income. Notice that in the case of  $N_p = 1,000$, this ratio is equal to 0.39, while in the case of $N_p = 10,000$, it is equal to approx 0.011, which further demonstrates a far lower risk in the case of the larger company.
$~~~\Box$

\subsubsection*{Homework \ref{probability}.\arabic{homework}}
\addtocounter{homework}{1} \addtocounter{tothomework}{1}

In insurance, the term {\it loss ratio} is the ratio of the total payouts in claims to the total income from premiums. Find the mean and standard deviations of the loss ratio in the two cases of  $N_p = 1,000$ and $N_p = 10,000$.


\subsubsection*{Answers}
In the case $N_p = 1,000$, the mean of the loss ratio is $$\frac{1,000\times 120}{1,000 \times 200} = 0.6;$$ the variance of the loss ratio is $$\frac{985,600,000}{(1,000\times 200)^2}=0.02464;$$
and the standard deviation is $$\sqrt{0.02464}=0.157. $$

In the case $N_p = 10,000$, the mean of the loss ratio is $$\frac{10,000\times 120}{10,000 \times 200} = 0.6;$$ the variance of the loss ratio is $$\frac{9,856,000,000}{(10,000\times 200)^2}=0.002464;$$
and the standard deviation is $$\sqrt{0.002464}=0.05. ~~~\Box$$

\subsubsection*{Homework \ref{probability}.\arabic{homework}}

A question that is often asked is the following. Can insurance companies use my generic DNA tests in their decisions about my health insurance premium, coverage, or eligibility? In the US, a law called the Genetic Information Nondiscrimination Act (GINA) limits such use by insurance companies. The insurance companies may argue that using DNA information reduces their risk and, therefore, they will be able to provide lower insurance premiums {\bf on average}. Analyze such a claim based on what you have learned.

\subsubsection*{Guide}

Consider the following basic arguments to develop your answer. In particular, try to provide a bridge between theory and practice. Notice that for very large insurance companies based on economy of scale, the overall risk may already be low, so having DNA information may not reduce such risk significantly. (This last statement relies on the simple examples illustrated in this Chapter and may not be true in practice.) However, for small companies, accurate information based on DNA results may be more useful in reducing their risk and may improve competition, which in turn may lead to overall cost reduction.

\subsection{Mean and Variance of Specific Random Variables}

If $X$ is a Bernoulli random variable with parameter $p$, its mean and variance are obtained by
$$E[X]=0 \times P(X=0) + 1 \times P(X=1) = 0(1-p)+1(p) = p$$
and
$$Var[X] = E[X^2] - (E[X])^2 = p - p^2 = p(1-p),$$
or
\begin{eqnarray*}
Var[X] & = & E[(X-E[X])^2] \\
 & = & (0-p)^2  P(X=0) + (1-p)^2 P(X=1)\\
  & = &  p^2(1-p) + (1-p)^2p \\
   & = & p(1-p).
   \end{eqnarray*}
Since the binomial random variable is a sum of $n$ IID
 Bernoulli random variables, its mean and its variance are $n$ times the mean
and variance of these Bernoulli random variables, respectively. Notice also
that
by letting $p\rightarrow 0$, and $np \rightarrow \lambda$, both the mean
and the variance of the binomial random variable approach $\lambda$, which is the
value of both the mean and variance of the Poisson random variable.

\subsubsection*{Homework \ref{probability}.\arabic{homework}}

\addtocounter{homework}{1} \addtocounter{tothomework}{1}

Consider the previous problem and plot the probability function,
distribution function, and complementary distribution function of
$X$. $\Box$

\subsubsection*{Homework \ref{probability}.\arabic{homework}}

\addtocounter{homework}{1} \addtocounter{tothomework}{1}

Let $X$ be a geometric random variable with parameter $p$. Derive its mean and variance using the law of iterated expectations and EVVE.

\subsubsection*{Guide}

Consider a Bernoulli random variable $Y$ that takes the value 1 (success) with probability $p$ and the value 0 (failure) with probability $1-p$. Let $Y$ be the first Bernoulli trial associated with the geometric random variable $X$. Therefore, $E[X|Y=1]=1$ and $E[X|Y=0]=1+E[X]$. The former represents the case of having a success in the first trial, so $X=1$, and the second represents a failure in the first trial, so due to the memoryless property of the geometric random variable, the expected number of required trials will be one (the first failure) plus the mean of another independent geometric random variable that has mean of $E[X]$.

Then, $E[X|Y]$ can be written as $$E[X|Y]=1+(1-Y)E[X].$$
Accordingly,
$$E[X]=E_YE[X|Y]=E[1+(1-E[Y])E[X]] = 1+(1-p)E[X].$$
Thus,
$$E[X]=\frac{1}{p}.$$
Also,
$$
Var[X|Y]=\left\{\begin{array}{ll}
0 & \mbox{if $Y=1$}\\
Var(X) & \mbox{if $Y=0.$} \end{array} \right. $$

Thus,
$$Var[X|Y]=(1-Y)Var[X],$$
so
$$E_Y[Var[X|Y]]=Var[X](1-p).$$
Also,
$$Var \left[1+(1-Y)\frac{1}{p}\right]=\frac{1}{p^2}p(1-p).$$
Then, by EVVE,
$$Var[X] = Var[X](1-p) + \frac{p(1-p)}{p^2}.$$
Hence,
$$Var[X]=\frac{1-p}{p^2}. ~~~\Box $$

\subsubsection*{Homework \ref{probability}.\arabic{homework}}
\addtocounter{homework}{1} \addtocounter{tothomework}{1}

Let $Y$ be a geometric random variable of failures with parameter $p$. Obtain its mean and variance.

\subsubsection*{Guide}

We know that $Y=X-1$, where

$$E[X]=\frac{1}{p},~~~~~~{\rm and}~~~~~~Var[X]=\frac{1-p}{p^2}.$$

Therefore,

$$E[Y]=E[X] - 1 = \frac{1}{p}-1= \frac{1-p}{p}.$$
and
$$Var[Y]=Var[X]=\frac{1-p}{p^2}. ~~~\Box $$

\subsubsection*{Homework \ref{probability}.\arabic{homework}}

\addtocounter{homework}{1} \addtocounter{tothomework}{1}

Consider a discrete uniform probability function with parameters $a$ and $b$ that takes equal non-zero values for $x=a, a+1, a+2, \dots, b$. Derive its mean $E[X]$ and variance $Var[X]$.

\subsubsection*{Guide}

First, show that $$E[X] = \frac{a+b}{2}.$$

To obtain $Var[X]$, notice first that the variance of this discrete uniform distribution over the $x$ values of $a, a+1, a+2, \ldots b$ has the same value as the variance of a discrete uniform distribution over the $x$ values of $1, 2, \ldots b-a+1$. Then, for convenience, set $n=b-a+1$ and derive the variance for the probability function,
$$
P_X(x)=\left\{\begin{array}{ll}
\frac{1}{n} & \mbox{if $x=1, 2, \dots, n$}\\
0 & \mbox{otherwise.} \end{array} \right. $$

Next, by induction on $n$ show that

$$E[X^2] = \frac{(n+1)(2n+1)}{6}.$$

Finally, use the equation $Var[X]= E[X^2] - (E[X])^2$ to show that

\begin{equation}
\label{vardu}
Var[X] = \frac{n^2-1}{12}=\frac{(b-a+1)^2-1}{12}.
\end{equation}

$~~~\Box$

\subsubsection*{Homework \ref{probability}.\arabic{homework}}

\addtocounter{homework}{1} \addtocounter{tothomework}{1}

Consider an exponential random variable with parameter $\lambda$. Derive its mean and Variance.

 \subsubsection*{Guide}

Find the mean by $$E[X] = \int_0^\infty x \lambda e^{-\lambda x}
dx.$$

Use integration by parts to show that:

$$E[X] = \left[ -xe^{-\lambda x }  \right ]_0^\infty + \int_0^\infty
 e^{-\lambda x} dx = \frac{1}{\lambda}.$$

 Now notice how  (\ref{ezc}) simplifies the derivation of the mean:

$$E[X] = \int_0^\infty e^{-\lambda x}dx = \left[ \frac{1}{-\lambda} e^{-\lambda x} \right]_0^\infty =\frac{1}{\lambda}.$$

Then, use integration by parts to derive the second moment.
Understand and verify the following derivations:

\begin{eqnarray*}
E[X^2] & = & \int_0^\infty x^2 \lambda e^{-\lambda x} dx \\
                 & = & \left[ -x^2e^{-\lambda x }  \right ]_0^\infty  + 2\int_0^\infty
 xe^{-\lambda x} dx \\
                 & = & \left[ \left(-x^2e^{-\lambda x }  -\frac{2}{\lambda}xe^{-\lambda x }
                 -\frac{2}{\lambda^2}e^{-\lambda x } \right) \right ]_0^\infty\\
                 & = &  \frac{2}{\lambda^2}.
\end{eqnarray*}

$$Var[X] = E[X^2] - (E[X])^2 = \frac{2}{\lambda^2} -
\frac{1}{\lambda^2} = \frac{1}{\lambda^2}~~~\Box $$

{\bf The mean of the Pareto random variable} is given by

\begin{equation}
\label{meanpareto} 
E[X]=\left\{\begin{array}{ll}
\infty & \mbox{if $0 < \gamma \le 1$}\\
{\frac{\delta\gamma}{\gamma-1}} & \mbox{$\gamma > 1$.}
\end{array}
\right.
\end{equation}
 For $0< \gamma
\leq 2$, the variance $Var[X]=\infty$.

The following table provides the mean and the variance of some of the above-mentioned random variables.
\begin{center}
\renewcommand{\arraystretch}{1.4}
\begin{tabular}{|c|c|c|c| }  \hline
random variable & parameters & mean & variance  \\  \hline
Bernoulli & $0 \leq p \leq 1 $ & $p$ & $p(1-p)$ \\  \hline
geometric  & $0 \leq p \leq 1 $ & $1/p$ & $(1-p)/p^2$ \\  \hline
geometric (of failures) & $0 \leq p \leq 1 $ & $(1-p)/p$ & $(1-p)/p^2$ \\  \hline
binomial & $n$ and $0 \leq p \leq 1$ & $np$ &  $np(1-p)$ \\ \hline
Poisson & $\lambda>0$ & $\lambda$ & $\lambda$ \\  \hline
discrete uniform & $a$ and $b$ & $(a+b)/2$ & $[(b-a+1)^2-1]/12$ \\  \hline
uniform & $a$ and $b$ & $(a+b)/2$ & $(b-a)^2/12$ \\  \hline
exponential & $\mu>0$ & $1/\mu$ & $1/\mu^2$ \\  \hline Gaussian &
$m$ and $\sigma$ & $m$ & $\sigma^2$ \\  \hline Pareto & $\delta>0$
and $1<\gamma\leq 2$ & ${{\delta\gamma}/(\gamma-1)}$
 & $\infty$ \\  \hline
\end{tabular}
\end{center}

\subsection{Sample Mean and Sample Variance}

If we are given a sample of $n$ realizations of a random variable $X$, denoted $X(1), X(2), ~\ldots, X(n)$ we will use the {\bf Sample Mean}
defined by \begin{equation}
\label{samplemean} S_m = \frac{\sum_{i=1}^n X(i)}{n}
\end{equation}
as an estimator for the mean of $X$. For example, if we run a simulation of a queueing system and observe $n$ values of customer delays for $n$ different customers, the Sample Mean will be used to estimate a customer delay.

If we are given a sample of $n$ realizations of a random variable $X$, denoted $X(1), X(2), ~\ldots, X(n)$
we will use the {\it Sample Variance} defined by \begin{equation}
\label{samplemean1} S_v = \frac{\sum_{i=1}^n [X(i)-S_m]^2}{n-1}~~~~~~~~~~n>1
\end{equation}
as an estimator for the variance of $X$.
The sample standard deviation is then $\sqrt{S_v}$.

\subsubsection*{Homework \ref{probability}.\arabic{homework}}

\addtocounter{homework}{1} \addtocounter{tothomework}{1} Generate 10 variates from an exponential distribution of a given mean and compute the Sample Mean and Sample Variance. Compare them with the real mean and variance. Then, increase the sample to 100, 1000, $~\ldots,$ 1,000,000. Observe the difference between the real mean and variance and the sample mean and variance. Repeat the experiment for a Pareto variates that have the same mean. Discuss differences. $~~~\Box$

\subsection{Covariance and Correlation}

When random variables are positively dependent, namely, if when
one of them obtains high values, the others are likely to obtain
high values also, then the variance of their sum may be much higher than the sum of the individual variances. This is very significant
for bursty multimedia traffic modeling and resource provisioning.
For example, let time be divided into consecutive small-time
intervals, if $X_i$ is the amount of traffic that arrives during the
$i$th interval, and assume that we use a buffer that can store
traffic that arrives in many intervals, the probability of buffer
overflow will be significantly affected by the variance of the
amount of traffic arrives in a time period of many intervals,
which in turn is strongly affected by the dependence between the
$X_i$s. Therefore, there is a need to define a quantitative
measure for dependence between random variables. Such measure is
called the {\bf covariance}. The covariance of two random
variables $X_1$ and $X_2$, denoted by $Cov[X_1,X_2]$, is defined
by
\begin{equation}
\label{cov}
Cov[X_1,X_2] = E[(X_1 - E[X_1])(X_2 - E[X_2])].
\end{equation}
Intuitively, by Eq.\@ (\ref{cov}), if high value of $X_1$ implies
high value of $X_2$, and low value of $X_1$ implies low value of
$X_2$, the covariance is high. By Eq.\@ (\ref{cov}),
\begin{equation}
\label{cov2} Cov[X_1,X_2] = E[X_1 X_2] - E[X_1] E[X_2].
\end{equation}
Hence, by (\ref{etimes}), if $X_1$ and $X_2$ are independent, then
$Cov[X_1,X_2]=0$. The
variance of the sum of two random variables $X_1$ and $X_2$ is given by
\begin{equation}
\label{var2}
Var[X_1+X_2] = Var[X_1] + Var[X_2] + 2 Cov[X_1,X_2].
\end{equation}
This is consistent with our comments above. The higher the dependence
between the two random variables, as measured by their covariance, the
higher the variance of their sum, and if they are independent, hence
$Cov[X_1,X_2]=0$, the variance of their sum is equal to the sum of their
variances. Notice that the reverse is not always true: $Cov[X_1,X_2]=0$
does not necessarily imply that $X_1$ and $X_2$ are independent.

Notice also that negative covariance results in a lower value for the
variance of their sum than the sum of the individual variances.

\subsubsection*{Homework \ref{probability}.\arabic{homework}}
\addtocounter{homework}{1} \addtocounter{tothomework}{1} Prove that
$Cov[X_1,X_2]=0$ does not necessarily imply that $X_1$ and $X_2$ are
independent.
\subsubsection*{Guide}
The proof is by a counter-example. Consider two random variables $X$
and $Y$ and assume that both have Bernoulli distribution with
parameter $p$. Consider random variable $X_1$ defined by $X_1=X+Y$
and another random variable $X_2$ defined by $X_2=X-Y.$ Show that
$Cov[X_1,X_2]=0$ and that $X_1$ and $X_2$ are not independent.
$~~~\Box$

Let the sum of the random variables $X_1, X_2, X_3, \ldots, X_k$ be
denoted by $$S_k=X_1+ X_2+ X_3+ \ldots + X_k.$$ Then,
\begin{equation}
\label{varsum} Var[S_k] = \sum_{i=1}^k Var[X_i] + 2\sum_{i<j}
Cov[X_i,X_j]
\end{equation}
where $\sum_{i<j} Cov[X_i,X_j]$ is a sum over all $Cov[X_i,X_j]$
such that $i$ and $j$ is a pair selected without repetitions out of
$1,2,3, \ldots k$ so that $i<j$.

\subsubsection*{Homework \ref{probability}.\arabic{homework}}
\addtocounter{homework}{1} \addtocounter{tothomework}{1} Prove Eq.\@
(\ref{varsum}).
\subsubsection*{Guide}
First show that $S_k-E[S_k]=\sum_{i=1}^k (X_i-E[X_i])$ and that
$$(S_k-E[S_k])^2=\sum_{i=1}^k (X_i-E[X_i])^2 +2\sum_{i<j} (X_i-E[X_i])(X_j-E[X_j]). $$ Then, take the expectations of both sides of the
latter. $~~~\Box$

If we consider $k$ independent random variables denoted $X_1, X_2,
X_3, \ldots, X_k$, then by substituting $Cov[X_i,X_j]=0$ for all
relevant $i$ and $j$ in (\ref{varsum}), we obtain
\begin{equation}
\label{varsumind} Var[S_k] = \sum_{i=1}^k Var[X_i].
\end{equation}
\subsubsection*{Homework \ref{probability}.\arabic{homework}}
\addtocounter{homework}{1} \addtocounter{tothomework}{1} Use Eq.\@
(\ref{varsum}) to explain the relationship between the variance of
a Bernoulli random variable and a binomial random variable.
\subsubsection*{Guide}
Notice that a binomial random variable with parameters $k$ and $p$ is
a sum of $k$ independent Bernoulli random variables  with parameter
$p$. $~~~\Box$

The covariance can take any value between $-\infty$ and $+\infty$,
and in some cases, it is convenient to have a normalized dependence
measure - a measure that takes values between -1 and 1. Such a measure
is the {\bf correlation}. Noticing that the covariance is bounded by
\begin{equation}
\label{covbound}
Cov[X_1,X_2] \leq \sqrt{Var[X_1]Var[X_2]},
\end{equation}
the correlation of two random variables $X$ and $Y$ denoted by
$Corr[X,Y]$ is defined by
\begin{equation}
\label{corr} Corr[X,Y] = \frac{Cov[X,Y]}{\sigma_{X}\sigma_{Y}},
\end{equation}
assuming $Var[X]\neq 0$ and $Var[Y]\neq 0$.

\subsubsection*{Homework \ref{probability}.\arabic{homework}}
\addtocounter{homework}{1} \addtocounter{tothomework}{1} Prove that $\mid Corr[X,Y] \mid \leq 1$.
\subsubsection*{Guide}
Let $C=Cov[X,Y]$, and show
$C^2-\sigma_{X}^2\sigma_{Y}^2\leq 0$, by noticing that
$C^2-\sigma_{X}^2\sigma_{Y}^2$ is a discriminant of the quadratic
$a^2\sigma_{X}^2 + 2aC + \sigma_Y^2$ which must be nonnegative
because $E[a(X-E[X])+(Y-E[Y])]^2$ is nonnegative. $~~~\Box$

\subsection{Transforms}
\label{transf}
Transforms are useful in analyses of probability models and queueing systems. We will first consider the following general
definition \cite{BT02}
for a transform function $\Gamma$ of a random variable $X$:
\begin{equation}
\label{transforms}
\Gamma_X(\omega) = E[e^{\omega X}]
\end{equation}
where $\omega$ is a complex scalar.
Transforms have two important properties:
\begin{enumerate}
\item There is a one-to-one correspondence between transforms and probability distributions.
This is why they are sometimes called {\em characteristics} functions.
This means that for any distribution function, there is a unique transform function that
characterizes it and for each transform function there is a unique probability distribution
it characterizes. Unfortunately, it is not always easy to convert a transform to its probability
distribution, and therefore, in some cases where we are able to obtain the transform but not its
probability function, we use it as a means to characterize the random variable statistics instead of the probability distribution.
\item Having a transform function of a random variable we can generate its moments. This is why
transforms are sometimes called {\em moment generating} functions. In many cases, it is easier to
obtain the moments having the transform than having the actual probability distribution.
\end{enumerate}
We will now show how to obtain the moments of a continuous random
variable $X$ with density
function $f_X(x)$ from its transform function $\Gamma_X(\omega)$.
By definition,
\begin{equation}
\label{transtomom1}
\Gamma_X(\omega) = \int_{-\infty}^{\infty} e^{\omega x} f_X(x) dx.
\end{equation}
Taking the derivative with respect to $\omega $ leads to
\begin{equation}
\label{transtomom2}
\Gamma_X^{'}(\omega) = \int_{-\infty}^{\infty} x e^{ \omega x} f_X(x) dx.
\end{equation}
Letting   $\omega \rightarrow 0$, we obtain
\begin{equation}
\label{transtomom3}
\lim_{\omega \rightarrow 0} \Gamma_X^{'}(\omega) = E[X],
\end{equation}
and in general, taking the $n$th derivative and  letting   $\omega \rightarrow 0$, we obtain
\begin{equation}
\label{transtomom4}
\lim_{\omega \rightarrow 0} \Gamma_X^{(n)}(\omega ) = E[X^n].
\end{equation}
\subsubsection*{Homework \ref{probability}.\arabic{homework}}
\addtocounter{homework}{1} \addtocounter{tothomework}{1} Derive
Eq.\@ (\ref{transtomom4}) using (\ref{transtomom1}) --
(\ref{transtomom3}) completing all the missing steps. $~~~\Box$

Consider for example the exponential random variable $X$ with
parameter $\lambda$ having density function $f_X(x) = \lambda
e^{-\lambda x} $ and derive its transform function. By definition,
\begin{equation}
\label{transform_exp}
\Gamma_X(\omega) = E[e^{\omega X} ] = \lambda
\int_{x=0}^{\infty}e^{ \omega x}e^{-\lambda x} dx,
\end{equation}
which gives after some derivations
\begin{equation}
\label{transf_exp}
\Gamma_X(\omega) = \frac{\lambda}{\lambda - \omega}.
\end{equation}
\subsubsection*{Homework \ref{probability}.\arabic{homework}}
\addtocounter{homework}{1} \addtocounter{tothomework}{1} Derive
Eq.\@ (\ref{transf_exp}) from (\ref{transform_exp}) $~~~\Box$

Let $X$ and $Y$ be random variables and assume that $Y = a X + b$.
The transform of $Y$ is given by
\begin{equation}
\label{transform_axplusb}
\Gamma_Y(\omega) = E[e^{\omega Y}] = E[e^{ \omega (a X + b) }] = e^{\omega b} E[e^{\omega a X } ]= e^{\omega b} \Gamma_X(\omega a).
\end{equation}

Let random variable $Y$ be the sum of independent random variables $X_1$ and $X_2$, i.e., $Y=X_1+X_2$.
 The transform of $Y$  is given by
\begin{equation}
\label{transform_x1plusx2}
\Gamma_Y(\omega) = E[e^{\omega Y}] = E[e^{ \omega (X_1 + X_2) }] =
E[e^{\omega  X_1 } ]E[e^{\omega  X_2 } ]=
\Gamma_{X_1}(\omega)\Gamma_{X_2}(\omega).
\end{equation}
This result applies to a sum of $n$ independent random variables, so the transform of
a sum of independent random variables is equal to the product of their transforms.
If $Y=\sum_{i=1}^n X_i$ and all the $X_i$s are $n$ independent and
identically distributed (IID) random variables, then
\begin{equation}
\label{transform_iids}
\Gamma_Y(\omega) = E[e^{\omega Y}] = [\Gamma_{X_1}(\omega)]^n.
\end{equation}

Let us now consider a Gaussian random variable $X$ with parameters $m$ and
$\sigma$ and density
\begin{equation}
f_X(x)=\frac{1}{\sqrt{2\pi}\sigma} e^{-(x-m)^2/2\sigma^2} \qquad -\infty <
x <\infty.
\end{equation}
Its transform is derived as follows
\begin{eqnarray*}
\Gamma_X(\omega) & = &   E[e^{\omega X}]\\
& = & \int_{-\infty}^\infty
\frac{1}{\sqrt{2\pi}\sigma} e^{-(x-m)^2/2\sigma^2} e^{\omega x}dx  \\
& = & e^{(\sigma^2 \omega^2/2)+m \omega} \int_{-\infty}^\infty
\frac{1}{\sqrt{2\pi}\sigma} e^{-(x-m)^2/2\sigma^2} e^{\omega x}
e^{-(\sigma^2 \omega^2/2)- m \omega} dx\\
& = & e^{(\sigma^2 \omega^2/2)+m \omega} \int_{-\infty}^\infty
\frac{1}{\sqrt{2\pi}\sigma} e^{-(x-m-\sigma^2\omega)^2/2\sigma^2} dx \\
& = & e^{(\sigma^2 \omega^2/2)+m \omega}.
\end{eqnarray*}
Let us use the transform just derived to obtain the mean and variance
of a Gaussian random variable with parameters $m$ and $\sigma$.
Taking the first derivative and putting $\omega=0$, we obtain
\begin{equation}
E[X] = \Gamma_X^{'}(0) = m.
\end{equation}
Taking the second derivative and setting $\omega=0$, we obtain
\begin{equation}
E[X^2] = \Gamma_X^{(2)}(0) = \sigma^2+m^2.
\end{equation}
Thus,
\begin{equation}
Var[X] = E[X^2] - (E[X])^2 = \sigma^2 + m^2 - m^2 = \sigma^2.
\end{equation}

A Gaussian random variable with a mean equal to zero
and variance equal to one is
called {\em standard Gaussian}. It is well known that if
$Y$ is Gaussian with mean $m$ and standard deviation
$\sigma$, then the random variable $X$ defined as
\begin{equation}
\label{standardGauss}
X=\frac{Y-m}{\sigma}
\end{equation}
is standard Gaussian.

Substituting $\sigma=1$ and $m=0$ in the above transform of a Gaussian
random variable, we obtain that
\begin{equation}
\label{trans_G_standard}
\Gamma_X(\omega)= e^{(\omega^2/2)}
\end{equation}
is the transform of a standard Gaussian random variable.

\subsubsection*{Homework \ref{probability}.\arabic{homework}}
\addtocounter{homework}{1} \addtocounter{tothomework}{1} Show the consistency between the results
obtained for the transform of a Gaussian random variable, (\ref{transform_axplusb}),
(\ref{standardGauss}) and (\ref{trans_G_standard}). $~~~\Box$

Let $X_i, ~~i=1, 2, 3,~ \ldots,~ n$ be $n$ independent
random variables and let
$Y$ be a random variable that equals $X_i$ with probability $p_i$
for $i=1, 2, 3,~ \ldots,~ N$.  Therefore, by the Law of Total Probability,
\begin{equation}
P(Y=y) = \sum_{i=1}^{N} p_i P(X_i =y)
\end{equation}
or for continuous densities
\begin{equation}
f_Y(y) = \sum_{i=1}^{n} p_i f_{X_i}(y).
\end{equation}
Its transform is given by
\begin{eqnarray*}
\Gamma_Y(\omega) & = &   E[e^{\omega Y}]\\
& = & \int_{-\infty}^\infty f_Y(y) e^{\omega y}dx \\
& = & \int_{-\infty}^\infty [\sum_{i=1}^{n} p_i f_{X_i}(y)] e^{\omega y}dx \\
& = & \int_{-\infty}^\infty \sum_{i=1}^{n} p_i f_{X_i}(y) e^{\omega y}dx \\
& = & \sum_{i=1}^{n} p_i \Gamma_{X_i}(\omega).
\end{eqnarray*}
Notice that if the $X_i$ are exponential random variables, then by definition, $Y$ is hyper-exponential.

Particular transforms include the Z, the Laplace, and the Fourier transforms.

The {\bf Z-transform} $\Pi_X(z)$ applies to integer valued random variable $X$ and is defined by $$\Pi_X(z)=E[z^X].$$ This is a special case
of (\ref{transforms}) by setting $z=e^\omega$.

The {\bf Laplace transform} applies to nonnegative valued random variable $X$ and is
defined by $$L_X(s)=E[e^{-sX}]~~ {\rm for} s\geq 0.$$ This is a special case
of (\ref{transforms}) by setting $\omega=-s$.

The {\bf Fourier transform}
applies to both nonnegative and negative valued random variable $X$ and is
defined by $$\Upsilon_X(s)=E[e^{i\theta X}],$$ where $i = \sqrt{-1}$ and $\theta$ is real. This is a special case
of (\ref{transforms}) by setting $\omega=i\theta$.

We will only use the Z and Laplace transforms in this book.

\subsubsection{Z-transform}
\label{ztransforms}

Consider a discrete and nonnegative random variable $X$, and let
$p_i=P(X=i)$, $i=0,1,2,~ \ldots$ with  $\sum_{i=0}^{\infty} p_i=1$. The
Z-transform of $X$ is defined by
\begin{equation}
\label{Ztrans}  \Pi_X(z)=E [z^X] = \sum_{i=0}^{\infty} p_i z^i,
\end{equation}
where $z$ is a real number that satisfies $0 \leq z \leq 1$. Note that in many applications the Z-transform is defined for complex $z$.
However, for the purpose of this book, we will only consider real $z$ within $0 \leq z \leq 1$.

\subsubsection*{Homework \ref{probability}.\arabic{homework}}

\addtocounter{homework}{1} \addtocounter{tothomework}{1} Prove the following properties of the
Z-transform $\Pi_X(z)$:
\begin{enumerate}
\item $\lim_{z\rightarrow 1^{-}} \Pi_X(z)=1$ ($z\rightarrow 1^{-}$
is defined as $z$ approaches 1 from below). \item
$p_i=\Pi_X^{(i)}(0)/i!$ where $\Pi_X^{(i)}(z)$ is the $i$th
derivative of $\Pi_X(z)$. \item $E[X]=\lim_{z\rightarrow
1^{-}}\Pi_X^{(1)}(z)$. $~~~\Box$
\end{enumerate}

For simplification of notation, in the following, we will use
$\Pi_X^{(i)}(1)= \lim_{z\rightarrow 1^{-}}\Pi_X^{(i)}(z)$, but the reader must keep in mind that a straightforward substitution of
$z=1$ in $\Pi_X^{(i)}(z)$ is not always possible and the limit needs
to be derived. An elegant way to show the 3rd property is to consider  $
\Pi_X(z)=E [z^X]$, and exchanging the operation of derivative and
expectation, we obtain $ \Pi_X^{(1)}(z)=E [Xz^{X-1}]$, so
$\Pi_X^{(1)}(1)=E [X]$. Similarly,
\begin{equation}
\label{Ztransderivat1}  \Pi_X^{(i)}(1)=E[X(X - 1) \ldots (X - i + 1)].
\end{equation}

\subsubsection*{Homework \ref{probability}.\arabic{homework}}
\addtocounter{homework}{1} \addtocounter{tothomework}{1} Show that the variance $Var[X]$ is given
by
\begin{equation}
\label{Ztransvar}  Var[X] = \Pi_X^{(2)}(1) + \Pi_X^{(1)}(1) -
(\Pi_X^{(1)}(1) )^2.
\end{equation}
~~~~~~~~~~~~~~~~~~~~~~~~~~~~~~~~~~~~~~~~~~~~~~~~~~~~~~~~~~$~~~\Box$

\subsubsection*{Homework \ref{probability}.\arabic{homework}}

\addtocounter{homework}{1} \addtocounter{tothomework}{1} Derive a formula for $E[X^i]$ using the
Z-transform. $~~~\Box$

As a Z-transform is a special case of the transform $\Gamma_Y(\omega) = E[e^{\omega Y}]$, the following results hold.

If random variables $X$ and $Y$ are related by $Y=aX+b$ for real numbers $a$ and $b$, then
\begin{equation}
\label{Ztransaxplusb}  \Pi_Y(z) = z^b\Pi_X(z a).
\end{equation}
Let random variable $Y$ be the sum of independent random variables
$X_1$,  $X_2$, ~\ldots,~ $X_n$ ($Y=\sum_{i=1}^n X_i$),
 The Z-transform of $Y$ is given by
\begin{equation}
\label{transform_sumrvs}
\Pi_Y(z) = \Pi_{X_1}(z) \Pi_{X_2}(z)\Pi_{X_3}(z) ~\ldots ~\Pi_{X_n}(z).
\end{equation}
If $X_1, X_2, ~ \ldots, ~ X_n $ are also identically distributed, then
\begin{equation}
\label{transform_niids}
\Pi_Y(z) = [\Pi_{X_1}(z)]^n.
\end{equation}
Let us now consider several examples of Z-transforms of nonnegative
discrete random variables. If $X$ is a Bernoulli random variable
with parameter $p$, then its Z-transform is given by
\begin{equation}
\label{Ztransbernoulli}  \Pi_X(z) = (1-p)z^0 + p z^1 = 1-p +pz.
\end{equation}
Its mean is $E[X] = \Pi_X^{(1)}(1) = p$ and by (\ref{Ztransvar}) its variance is $p(1-p)$.

If $X$ is a Geometric random variable with parameter $p$, then its Z-transform is given by
\begin{equation}
\label{Ztransgeom}  \Pi_X(z) = p\sum_{i=1}^\infty    (1 - p)^{i-1}z^i=\frac{pz}{1 - (1 - p)z}.
\end{equation}
Its mean is $E[X] = \Pi_X^{(1)}(1) = 1/p$ and by (\ref{Ztransvar}) its variance is $(1-p)/p^2$.

If $X$ is a Binomial random variable with parameter $p$, then we can obtain its Z-transform
either by definition or by realizing that a Binomial random variable is a sum of $n$ IID Bernoulli random variables.
Therefore its Z-transform is given by
\begin{equation}
\label{Ztransbinom}  \Pi_X(z) = (1-p +pz)^n = [1+(z-1)p]^n.
\end{equation}

\subsubsection*{Homework \ref{probability}.\arabic{homework}}

\addtocounter{homework}{1} \addtocounter{tothomework}{1} Verify that the latter is consistent with
the Z-transform obtained using $\Pi_X(z) = \sum_{i=0}^{\infty} p_i
z^i$. $~~~\Box$

The mean of $X$ is $E[X] = \Pi_X^{(1)}(1) = np$ and by
(\ref{Ztransvar}) its variance is $np(1-p)$.

If $X$ is a Poisson random variable with parameter $\lambda$, then its Z-transform is given by
\begin{equation}
\label{Ztranspoisson}  \Pi_X(z) = \sum_{i=0}^{\infty} p_i z^i = e^{-\lambda} \sum_{i=0}^{\infty} \frac{\lambda^iz^i}{i!}=e^{(z-1)\lambda}.
\end{equation}
Its mean is $E[X] = \Pi_X^{(1)}(1) = \lambda$ and by (\ref{Ztransvar}) its variance is also equal to $\lambda$.

We can now see the relationship between the Binomial and the Poisson random variables.
If we consider the Z-transform of the Binomial random variable $\Pi_X(z) = (1-p +pz)^n$, and set
$\lambda = np$ as a constant so that $\Pi_X(z) = (1+ (z-1)\lambda/n)^n$ and let $n \rightarrow \infty$, we obtain
\begin{equation}
\label{binompoisson}  \lim_{n \rightarrow \infty} (1-p +pz)^n = \lim_{n \rightarrow \infty} [1+ (z-1)\lambda/n]^n = e^{(z-1)\lambda}
\end{equation}
which is exactly the Z-transform of the Poisson random variable. This proves the convergence of the binomial
 to the Poisson random variable if we keep $np$ constant and let $n$ go to infinity.

\subsubsection{Laplace Transform}
\label{Lap_transforms}

The Laplace transform of a non-negative random variable $X$ with   density
 $f_X(x)$ is defined as
\begin{equation}
\label{laplacedef} {\mathcal{L}}_X(s)= E[e^{-sX}] = \int_{0}^{\infty} e^{-sx}f_X(x) dx.
\end{equation}
As it is related to the transform $\Gamma_X(\omega) = E[e^{\omega X}]$ by setting $\omega=-s$, similar
derivations to those made for $\Gamma_X(\omega)$ above give the following.

If $X_1$, $X_2$, ~ \ldots,~ $X_n$  are $n$ independent random variables, then
\begin{equation}
\label{Laplacesumofnrvs}
{\mathcal{L}}_{X_1 + X_2 +~ \ldots ~+X_n}(s) = {\mathcal{L}}_{X_1}(s){\mathcal{L}}_{X_2}(s)~ \ldots ~ {\mathcal{L}}_{X_n}(s).
\end{equation}
Let $X$ and $Y$ be random variables and $Y = a X + b$. The Laplace
transform of $Y$  is given by
\begin{equation}
\label{Ltransform_axplusb}
{\mathcal{L}}_Y(s) = e^{-sb} {\mathcal{L}}_X(s a).
\end{equation}

The $n$th moment of random variable $X$ is given by
\begin{equation}
\label{Ltransform_nthmom}
E[X^n] = (-1)^n {\mathcal{L}}_X^{(n)}(0)
\end{equation}
where ${\mathcal{L}}_Y^{(n)}(0)$ is the $n$th derivative of ${\mathcal{L}}_X(s)$ at $s=0$ (or at the limit $s \rightarrow 0$). Therefore,
\begin{equation}
\label{Ltransform_var}
Var[X] = E[X^2]-(E[X])^2 = (-1)^2 {\mathcal{L}}_X^{(2)}(0) - ((-1) {\mathcal{L}}_X^{(1)}(0))^2 = {\mathcal{L}}_X^{(2)}(0) - ({\mathcal{L}}_X^{(1)}(0))^2.
\end{equation}
Let $X$ be an exponential random variable with parameter $\lambda$.
Its Laplace transform is given by
\begin{equation}
\label{Lap_transf_exp}
{\mathcal{L}}_X (s) = \frac{\lambda}{\lambda + s}.
\end{equation}
\subsubsection*{Homework \ref{probability}.\arabic{homework}}
\addtocounter{homework}{1} \addtocounter{tothomework}{1} Derive
(\ref{Laplacesumofnrvs})--(\ref{Lap_transf_exp}) using the
derivations made for $\Gamma_X(\omega)$ as a guide. $~~~\Box$

Now consider $N$ to be a nonnegative discrete (integer) random variable of a probability
distribution that has the Z-transform $\Pi_N(z)$, and let $Y=X_1+X_2+ ~\ldots~+ X_N$, where
$ X_1, ~X_2,  ~\ldots,~ X_N$ are nonnegative IID random variables with a common distribution that has the Laplace transform ${\mathcal{L}}_X(s)$ (i.e., they are exponentially distributed).
Let us derive the Laplace transform of $Y$.
Conditioning and unconditioning on $N$, we obtain
\begin{equation}
\label{Ltransform_sumN}
{\mathcal{L}}_Y(s) = E[e^{-sY}]= E_N [E[e^{-s(X_1+X_2+~\ldots~ +X_N \mid N)}]].
\end{equation}
Therefore, by the independence of the $X_i$,
\begin{equation}
\label{Ltransform_sumN2}
{\mathcal{L}}_Y(s) =  E_N [ E[e^{-sX_1} + E[e^{-sX_2}  + ~\ldots~ + E[e^{-sX_N} ]] = E_N [({\mathcal{L}}_X(s))^N ].
\end{equation}
Therefore
\begin{equation}
\label{Ltransform_sumN3}
{\mathcal{L}}_Y(s) =  \Pi_N  [({\mathcal{L}}_X(s)) ].
\end{equation}
An interesting example of (\ref{Ltransform_sumN3}) is the case where
the $X_i$ are IID exponentially distributed each with parameter
$\lambda$, and $N$ is geometrically distributed with parameter $p$.
In this case, we already know that since $X$ is an exponential
random variable, we have ${\mathcal{L}}_X(s) = \lambda/(\lambda +
s)$, so
\begin{equation}
\label{Ltransform_sumN4}
{\mathcal{L}}_Y(s) =  \Pi_N  \left( \frac{ \lambda}{\lambda + s}  \right).
\end{equation}
We also know that $N$ is geometrically distributed,
so $\Pi_N (z) = pz/[1-(1-p)z]$. Therefore,
from (\ref{Ltransform_sumN4}), we obtain,
\begin{equation}
\label{Ltransform_sumN5}
{\mathcal{L}}_Y(s) =  \frac{\frac{p\lambda}{\lambda +s}}
{ 1-\frac{(1-p)\lambda}{\lambda + s}}
\end{equation}
and after some algebra, we obtain
\begin{equation}
\label{Ltransform_sumN6}
{\mathcal{L}}_Y(s) =  \frac{p\lambda}{s + p\lambda}.
\end{equation}
This result is interesting. We have shown that $Y$ is
exponentially distributed with parameter $p\lambda$.

Note that if $N$ is a geometric random variable of failures, then this result will not apply because $Y$ will have an atom at zero, and therefore it cannot be exponentially distributed.

\subsubsection*{Intuitive interpretation}

As above, let $Y$ be a positive continuous random variable which is equal to the sum $ X_1 + X_2 + \ldots + X_N$ where the $X_i, ~~~ i=1, 2, \ldots, N$ are exponential IID random variables (also called phases), and $N$ is a geometric random variable with parameter $p$. Assume that you are part of a simulation that starts at time $t=0$ and at time $t=Y$, you are about to win a very significant award. Assume that someone in the background is generating the exponential variates and the Bernoulli experiments that make up the geometric variate. That is, after every period $X_i$, a Bernoulli experiment is taking place to decide if this was the last exponential phase. At any point in time $t$, for $0 < t < Y$, you can see the length of the time interval $(0,t)$, has nothing to do with the remaining time until you receive the award, i.e., the time interval $(t,Y)$, because the remaining time of your phase at time $t$ is memoryless (because it is exponential) and the remaining number of phases is memoryless (because it is geometric). Accordingly, $Y$ must be memoryless, and because the exponential distribution is the only continuous random variable that is memoryless,  $Y$ must be exponential.

\subsubsection*{Homework \ref{probability}.\arabic{homework}}
\addtocounter{homework}{1} \addtocounter{tothomework}{1}
Let $X_1$ and $X_2$ be exponential random variables with parameters $\lambda_1$ and $\lambda_2$, respectively.
Consider a random variable $Y$ defined by the following algorithm.
\begin{enumerate}
\item Initialization: $Y=0$.
\item Conduct an experiment to obtain the values of $X_1$ and $X_2$. If $X_1<X_2$, then $Y=Y+X_1$ and Stop. Else, $Y=Y+X_2$ and repeat 2.
\end{enumerate}
Show that $Y$ is exponentially distributed with parameter $\lambda_1$.
\subsubsection*{Hint}
Notice that $Y$ is a geometric sum of exponential random variables.
 $~~~\Box$

\subsubsection*{Homework \ref{probability}.\arabic{homework}}
\addtocounter{homework}{1} \addtocounter{tothomework}{1}

Derive the density, the Laplace transform,
the mean and the variance of $Y$ in the following three cases.

\begin{enumerate}
\item Let $X_1$ and $X_2$ be exponential random variables with parameters
$\mu_1$ and $\mu_2$, respectively. In this case, $Y$ is
a hyperexponential random variable with density $f_Y(y)=pf_{X_1}(y) + (1-p)f_{X_2}(y)$.
\item Let $X_1$ and $X_2$ be exponential random variables with parameters
$\mu_1$ and $\mu_2$, respectively. The hypoexponential random variable $Y$ is defined by $Y=X_1+X_2$.
\item  Let $Y$ be an Erlang random variable, namely,
$Y=\sum_{i=1}^k X_i$ where the $X_i$s are IID
exponentially distributed  random variables with
parameter $\mu$.
\end{enumerate}
Now plot the standard deviation to mean ratio for the cases of hyperexponential and Erlang
random variables over a wide range of parameter values and discuss
implications. For example, show that for Erlang(k) the standard
deviation to mean ratio approaches zero as k approaches infinity.
$~~~\Box$

\subsection{Multivariate Random Variables and Transform}

A {\em multivariate random variable} is a vector $X = (X_1, X_2,~
\ldots,~ X_k)$ where each of the $k$ components is a random variable.
A multivariate random variable is also known as {\em random vector}.
These $k$ components of a random vector are related to events
(outcomes of experiments) on the same sample space and they can be
continuous or discrete.
They also have a legitimate well-defined
joint distribution (or density) function. The distribution of each
individual component $X_i$ of the random vector is its marginal
distribution. A transform of a random vector $X = (X_1, X_2, ~\ldots,~
X_k)$ is called  {\em multivariate transform} and is defined by
\begin{equation}
\Gamma_{X}(\omega_1, \omega_2, ~\ldots,~ \omega_k)  = E[s^{\omega_1
X_1, \omega_2 X_2, ~\ldots,~ \omega_k X_k}].
\end{equation}

\subsection{Probability Inequalities and Their Dimensioning Applications}
In the course of the design of telecommunications networks, a fundamentally important problem is how much capacity
a link should have. If we consider the demand as a non-negative random variable $X$ and the link capacity as a fixed scalar
$C>0$, we will be interested in the probability that the demand exceeds the capacity $P(X>C)$.
The more we know about the distribution the more accurate our estimation of $P(X>C)$ will be.

If we know only the mean, we use the so-called {\bf Markov
inequality}:
\begin{equation}
\label{markovineq}
P(X>C) \leq \frac{E[X]}{C}.
\end{equation}

\subsubsection*{Homework \ref{probability}.\arabic{homework}}
\addtocounter{homework}{1} \addtocounter{tothomework}{1} Prove Eq.\@
(\ref{markovineq}).
\subsubsection*{Guide}
Define a new random variable $U(C)$ a function of $X$ and $C$
defined by: $U(C)=0$ if $X<C$, and $U(C)=C$ if $X\geq C$. Notice
$U(C)\leq X$, so $E[U(C)]\leq E[X]$. Also, $E[U(C)] = CP(U(C)=C) =
CP(X\geq C)$, and Eq.\@ (\ref{markovineq}) follows.
 $~~~\Box$

  If we know the mean and the variance of $X$, then we can use the so-called {\bf Chebyshev inequality}:
\begin{equation}
\label{Chebyshevineq}
P(\mid X-E[X] \mid >C) \leq \frac{Var[X]}{C^2}.
\end{equation}
\subsubsection*{Homework \ref{probability}.\arabic{homework}}
\addtocounter{homework}{1} \addtocounter{tothomework}{1} Prove Eq.\@
(\ref{Chebyshevineq}).
\subsubsection*{Guide}
Define a new random variable $(X-E[X])^2$ and apply the Markov inequality
putting $C^2$ instead of C obtaining: $$P((X-E[X])^2\geq C^2) \leq \frac{E[(X-E[X])^2]}{C^2}=\frac{Var[X]}{C^2}.$$ Notice that the two events $(X-E[X])^2\geq C^2$ and $\mid X-E[X] \mid \geq C$ are identical.
 $~~~\Box$

Another version of Chebyshev inequality is
\begin{equation}
\label{Chebyshevineq2}
P(\mid X-E[X] \mid >C^*\sigma) \leq \frac{1}{(C^*)^2}
\end{equation}
for $C^*>0$.
  \subsubsection*{Homework \ref{probability}.\arabic{homework}}
\addtocounter{homework}{1} \addtocounter{tothomework}{1} Prove and
provide interpretation to Eq.\@ (\ref{Chebyshevineq2}).
\subsubsection*{Guide}
Observe that the right-hand side of (\ref{Chebyshevineq2}) is equal to $\frac{Var[X]}{Var[X](C^*)^2}$.
   $~~~\Box$
     \subsubsection*{Homework \ref{probability}.\arabic{homework}}
\addtocounter{homework}{1} \addtocounter{tothomework}{1} For a wide
range of parameter values, study numerically how tight the bounds
provided by Markov versus Chebyshev inequalities are. Discuss the
differences and provide interpretations.
   $~~~\Box$

   A further refinement of the Chebyshev inequality is the following {\bf
   Kolmogorov inequality}.  Let $X_1, X_2, X_3, ~\ldots,~ X_k$ be a sequence of mutually independent random
   variables (not necessarily identically distributed)
   and let $S_k=X_1+ X_2+ X_3+~ \ldots~+ X_k$ and $\sigma(S_k)$ be the standard deviation of $S_k$. Then, for
   every $\epsilon >0$,
\begin{equation}
\label{Kolmogorov} P(\mid S_k-E[S_k] \mid < \theta \sigma(S_k) {\rm
~~for~~all~~} k=1,2, \ldots, n) \geq 1-\frac{1}{\theta^2}.
\end{equation}
 The interested reader may consult Feller \cite{fell68} for the
 proof of the Kolmogorov inequality. We are however more interested in
 its teletraffic implications. If we let time be divided into
 consecutive intervals and we assume that $X_i$ is the number of
 packets arrive during the $i$th interval, and if the number of packets that arrive during the different intervals are mutually independent, then it is rare that we will have within
 a period of $n$ consecutive intervals any period of $k$ consecutive
 intervals ($k\leq n$) during which the number of packets arriving
 is significantly more than the average.

   \subsection{Limit Theorems}
   \label{limittheorems}
   Let $X_1, X_2, X_3, ~\ldots,~ X_k$ be a sequence of IID random
   variables with mean $\lambda$ and variance $\sigma^2$.
   Let $\bar{S}_k$ be the {\it sample mean} of these $k$ random variables defined by
   $$ \bar{S}_k =\frac{X_1+ X_2+ X_3+~ \ldots~+ X_k}{k}.$$
   This gives
 $$ E[\bar{S}_k] =\frac{E[X_1]+ E[X_2]+ E[X_3]+~ \ldots ~+ E[X_k]}{k}=\frac{k\lambda}{k}=\lambda.$$
 Recalling that the  $X_i$s are independent, we obtain
\begin{equation} \label{varsamplemean} Var[\bar{S}_k] = \frac{\sigma^2}{k}.\end{equation}
  \subsubsection*{Homework \ref{probability}.\arabic{homework}}
\addtocounter{homework}{1} \addtocounter{tothomework}{1} Prove Eq.\@
(\ref{varsamplemean}).
    $~~~\Box$

 Applying Chebyshev's inequality, we obtain
    \begin{equation}
    \label{PSbarminuslambda} P(\mid \bar{S}_k - \lambda \mid \geq \varepsilon)\leq \frac{\sigma^2}{k \varepsilon^2}~~{\rm for~~all}~~
    \varepsilon>0.
    \end{equation}
    Noticing that as $k$ approaches infinity, the right-hand side of  (\ref{PSbarminuslambda})
    approaches zero which implies that the left-hand side approaches zero as well. This leads to the
    so-called {\bf the weak law of large numbers} that states the following. Let $X_1,X_2,X_3,~ \ldots, ~X_k$ be $k$ IID random variables with common mean $\lambda$. Then,       \begin{equation} \label{lawlargenumbers} P\left(\left| \frac{X_1+ X_2+ X_3 +~\ldots ~+ X_k}{k} - \lambda \right|  \geq \varepsilon)\right) \rightarrow 0 ~~{\rm as}~~ k\rightarrow \infty ~~{\rm for~~ all}~~ \varepsilon>0. \end{equation}
    What the weak law or large number essentially says is that the sample mean approaches the mean as the sample size increases.

    Next we state the central limit theorem that we have mentioned in Section \ref{gaussrv}.
    Let $X_1,X_2,X_3, ~\ldots,~ X_k$ be $k$ IID random variables with common mean $\lambda$ and variance $\sigma^2$. Let random variable $Y_k$ be defined as
       \begin{equation} \label{clt1} Y_k= \frac{X_1+ X_2+ X_3 +~\ldots~ + X_k-k\lambda}{\sigma\sqrt{k}}. \end{equation}
       Then,
         \begin{equation} \label{clt2} \lim_{k\rightarrow \infty} P(Y_k \leq y) = \Phi(y) \end{equation}
     where $\Phi(\cdot)$ is the distribution function of a standard Gaussian random variable given by
$$\Phi(y)=\frac{1}{\sqrt{2\pi}}\int_{-\infty}^y e^{-t^2/2}dt.$$
 \subsubsection*{Homework \ref{probability}.\arabic{homework}}
\addtocounter{homework}{1} \addtocounter{tothomework}{1} Prove that $E[Y_k]=0$ and that $Var[Y_k]=1$ from first principles without using the central limit theorem.
    $~~~\Box$

As we mentioned in Section \ref{gaussrv}, the central limit theorem is considered the most important result in probability. Notice that it implies that the sum of $k$ IID random variable with common mean $\lambda$ and variance $\sigma^2$ is approximately Gaussian with mean $k\lambda$ and variance $k\sigma^2$ \emph{regardless} of the distribution of these variables.

Moreover, under certain conditions, the central limit theorem also applies in the case of sequences that are not identically distributed, provided one of a number of conditions apply. One of the
cases where the central limit theorem also applies in the case of non-IID random variables is due to Lyapunov described as follows.
Consider $X_1,X_2,X_3,~ \ldots, ~X_k$ to be a sequence of independent random variables. Let $\lambda_n=E[X_n], ~n=1,2,~ \ldots, ~k$ and $\sigma^2_n=Var[X_n],~n=1,2,~ \ldots, ~k$, and assume that all $\lambda_n$ and  $\sigma^2_n$ are finite.
Let
$$\hat{S}^2_n = \sum_{i=1}^n \sigma_i^2,$$
$$\hat{R}^3_n = \sum_{i=1}^n E[|X_i - \lambda_i|^3],$$
and assume that $\hat{S}^2_n$ and $\hat{R}^3_n$ are finite for all $n=1,2,~ \ldots, ~k$. Further, assume that
$$\lim_{k \rightarrow \infty} \frac{\hat{R}}{\hat{S}} = 0.$$
The latter is called the ``Lyapunov condition''.

If these conditions hold, then the random variable $\sum_{i=1}^k X_i $ has Gaussian distribution with mean
$\sum_{i=1}^k \lambda_i$ and variance $\sum_{i=1}^k  \sigma^2_i $.
This generalization of the central limit theorem to non-IID random variables, based on the Lyapunov condition,
is called ``Lyapunov's central limit theorem''.

\subsection{Link Dimensioning}
\label{link_dim} Before we end this chapter on probability, let us
demonstrate how the probability concepts discussed so far can be
used to provide simple means for link dimensioning. We will consider several scenarios of sources (individuals
or families) sharing a communication link. Each of the sources has certain requirements for capacity and the common link must be
dimensioned in such a way that minimizes the cost for the telecommunications provider but still meets the individual QoS requirements.
The link dimensioning procedures that we described below apply to user requirements for capacity. These requirements apply to
transmissions from the sources to the network as well as to downloads from the networks to the user
or to a combination of downloads and transmissions. We are not concerned with specific directions of transmission. We assume that the capacity of the common
link can be used in either direction. When we say a source ``transmits'' it should always be read as ``transmits and/or downloads''.

\subsubsection{Case 1: Homogeneous Individual Sources}

Consider $N$
independent sources (end-terminals), sharing a transmission link of
capacity $C$ [Mb/s]. Any one of the sources transmits data in accordance with an on-off process. That is, a source alternates between two states: 1) the on state during which the source transmits at a rate $R$ [Mb/s], and 2) the off state during which the source is
idle. Assume that the proportion of time the source is in the on-state
is $p$, so it is in the off-state $1-p$ of the time. The question is
how much capacity should the link have so it can serve all $N$
sources such that the probability that the demand exceeds the total
link capacity is not higher than $\alpha$.

We first derive the distribution of the total traffic demanded by the $N$
sources. Without loss of generality, let us normalize the traffic generated
by a source during an on-period by setting $R=1$. Realizing that the demand generated by a single source
is Bernoulli distributed with parameter $p$, we obtain that
the demand generated by all $N$ sources has Binomial distribution with parameters $p$ and $N$.
Accordingly, finding the desired capacity
is reduced to finding the smallest $C$ such that
\begin{equation}
\label{binomialdim}
\sum_{i=C+1}^N {N \choose i} p^i(1-p)^{N-i} \leq \alpha.
\end{equation}
Since the left-hand side of (\ref{binomialdim}) increases as $C$ decreases,
and since its value is zero if $C=N$, all we need to do to find the
optimal $C$ is to
compute the value of the left-hand side of (\ref{binomialdim}) for $C$  values of $N-1$, $N-2$,
$\dots$ until we find the first $C$ value for which
the inequality (\ref{binomialdim}) is violated.
Increasing that $C$ value by
one will give us the desired optimal $C$ value.

If $N$ is large we can use the central limit theorem and approximate the
Binomial distribution by a Gaussian distribution.
Accordingly, the demand can be approximated by a Gaussian random variable with mean
$Np$ and variance $Np(1-p)$ and simply find $C_G$ such that
the probability of our Gaussian random variable
to exceed $C_G$ is $\alpha$.

It is well known that Gaussian random variables obey the so-called
68-95-99.7\% Rule which means that the following apply to a
random variable $X$ with mean $m$ and standard deviation $\sigma$.
\begin{eqnarray*}
P(m-\sigma \leq X \leq m+\sigma) & = &   0.68 \\
P(m-2\sigma \leq X \leq m+2\sigma) & = & 0.95  \\
P(m-3\sigma \leq X \leq m+3\sigma) & = & 0.997.
\end{eqnarray*}
Therefore, if $\alpha=0.0015$, then $C_G$
should be three standard deviations above the mean, namely,
\begin{equation}
\label{binomialdimG}
C_G= Np + 3 \sqrt{Np(1-p)}.
\end{equation}
Note that $\alpha$ is a preassigned QoS  measure representing the proportion of time that the
demand exceeds the supply, and under the zero buffer approximation during such time, some
traffic is lost. If  it is required that $\alpha$
is lower than 0.0015, then more than three standard deviations above the mean are required.
Recall that for our original problem, before we introduced the Gaussian approximation, $C=N$ guarantees that there is sufficient capacity to serve
all arriving traffic without losses. Therefore, we set our dimensioning rule for the optimal $C$ value as follows:
\begin{equation}
\label{optcbinomialdimG}
C_{opt}= \min{ \left[ N, Np + 3 \sqrt{Np(1-p)} \right] }.
\end{equation}

\subsubsection{Case 2: Non-homogeneous Individual Sources}
Here we generalize the above scenario to the case where the traffic and the
peak rates of different sources can be different. Consider $N$ sources where the $i$th source transmits  at rate
$R_i$ with probability $p_i$, and at rate $0$ with probability $1-p_i$.
In this case, where the sources are non-homogeneous, we must invoke a generalization of
the central limit theorem that allows for non-IID random variables
(i.e., the ``Lyapunov's central limit theorem'').
Let $R_X(i)$ be a random variable representing the rate transmitted by source $i$.
We obtain:  $$E[R_X(i)]=p_iR_i.$$ and $$Var[R_X(i)] = R_i^2 p_i-(R_ip_i)^2 =R_i^2p_i(1-p_i).$$
The latter is consistent with the fact that $R_X(i)$ is equal to $R_i$ times a Bernoulli random variable.
We now assume that the random variable $$\Sigma_R=\sum_{i=1}^N R_X(i) $$
has a Gaussian distribution with mean $$E[\Sigma_R] = \sum_{i=1}^N E[R_X(i)] =\sum_{i=1}^N p_iR_i$$ and variance
$$Var[\Sigma_R] = \sum_{i=1}^N Var[R_X(i)] =\sum_{i=1}^N  R_i^2p_i(1-p_i).$$
Notice that the allocated capacity should not be more than the total sum
of the peak rates of the individual sources. Therefore, in this more general case, the QoS requirement
$\alpha=0.0015$, our optimal
$C$ value  is set to:
\begin{equation}
\label{optcgen}
C_{opt}= \min{ \left[ \sum_{i=1}^N R_i , E[\Sigma_R] + 3 \sqrt{Var[\Sigma_R]} \right] }.
\end{equation}

 \subsubsection*{Homework \ref{probability}.\arabic{homework}}
\addtocounter{homework}{1} \addtocounter{tothomework} {1}
There are 20 sources, each transmits at a peak rate of 10 Mb/s with a probability of 0.1 and is idle with a probability of 0.9, and there are other 80 sources,  each transmits at a peak rate of 1 Mb/s with a probability of 0.05 and is idle with probability 0.95.
A service provider aims to allocate the minimal capacity $C_{opt}$ such that no more than 0.0015 of the time, the demand of all these 100 sources
exceeds the available capacity.
Set $C_{opt}$ using the above-described approach.

{\bf Answer:} $C_{opt} =$ 64.67186 Mb/s.

Notice the difference in contributions to the total variance of sources from the first group versus such contributions
of sources from the second group.

Consider a range of examples where the variance is the dominant part of $C_{opt}$ versus examples where the variance is not
the dominant part of $C_{opt}$.
    $~~~\Box$

    \subsubsection{Case 3: Capacity Dimensioning for a Community}

    In many cases, the sources are actually a collection of sub-sources. A source could be a family of several members, and at any given point in time,
    one or more of the family members are accessing the link. In such a case, we assume that source $i$, $i=1,2,3, \ldots, N$,
    transmits at rate $R_j(i)$ with probability $p_{ij}$ for $j=0, 1,2,3, \ldots, J(i)$.  For all $i$,
    $R_0(i)\equiv 0$ and $R_{J(i)}(i)$ is defined to be the peak rate of source $i$.
For each source (family) $i$, $R_j(i)$ and $p_{ij}$ for $j= 1,2,3, \ldots, J(i)-1$, are set based on measurements
for the various rates reflecting the total rates transmitted by active family members and their respective proportion of time used. For example, for a certain family $i$,  $R_1(i)$
could be the rate associated with one individual
 family member browsing the web,  $R_2(i)$ the rate associated with one individual
 family member using Voice over IP, $R_3(i)$ the rate associated with one individual
 family member watching video, $R_4(i)$ the rate associated with one individual
 family member watching the video and another browsing the web, etc.  The $p_{ij}$ is the proportion
 of time during the busiest period of the day that $R_i(j)$ is used.

    Again, defining $R_X(i)$ as
    a random variable representing the rate transmitted by source $i$, we have
$$E[R_X(i)]=\sum_{j=0}^{J(i)} p_{ij}R_j(i)~~{\rm for}~~ i=1,2,3, \ldots, N.$$
and
$$Var[R_X(i)] =\sum_{j=0}^{J(i)} \{R_j(i)\}^2 p_{ij}-\{E[R_X(i)]\}^2~~{\rm for}~~ i=1,2,3, \ldots, N.$$
Again, assume that the random variable
$$\Sigma_R=\sum_{i=1}^N R_X(i) $$
has a Gaussian distribution with mean
$$E[\Sigma_R] = \sum_{i=1}^N E[R_X(i)]$$
and variance
$$Var[\Sigma_R] = \sum_{i=1}^N Var[R_X(i)].$$
Therefore, in this general case, for the QoS requirement
$\alpha=0.0015$,
our optimal
$C$ value  is again set by
\begin{equation}
\label{optcgen1}
C_{opt}= \min{ \left[ \sum_{i=1}^N R_{J(i)}(i) , E[\Sigma_R] + 3 \sqrt{Var[\Sigma_R]} \right] }.
\end{equation}

\newpage
\section{Relevant Background on Stochastic Processes}
\label{stochastic}

\setcounter{homework}{1} 

Aiming to understand behaviors of various natural and artificial
processes, researchers often model them as collections of random
variables where the mathematically defined statistical
characteristics and dependencies of such random variables are fitted
to those of the real processes. The research in the field of
stochastic processes has therefore three facets:
\begin{description}
\item [Theory:] mathematical explorations of stochastic processes
models that aim to better understand their properties.
\item [Measurements:] taken on the real process in order to identify its statistical characteristics.
\item [modeling:] fitting the measured statistical characteristics of the real process with those of a model and development of new models of stochastic processes
that well match the real process.
\end{description}
This chapter provides background on the basic theoretical aspects of
stochastic processes, which form a basis for queueing theory and
teletraffic models discussed in the later chapters.

\subsection{General Concepts}
\label{gencons}

For a given {\em index set} $T$, a {\em stochastic process}
$\{X_t,~t\in T\}$ is an indexed collection of random variables. They may or may
not be identically distributed. In many applications, the index $t$
is used to model time. Accordingly, the random variable $X_t$ for
a given $t$ can represent, for example, the number of telephone
calls that have arrived at an exchange by time $t$.

If the index set $T$ is countable, the stochastic process is called
a {\em discrete-time} process, or a {\em time series} \cite{ande76,
Box70, nets73}. Otherwise, the stochastic process is called a {\em
continuous-time} process. Considering our previous example, where
the phone calls arriving at an exchange by time $t$ are modeled as a continuous-time process $\{X_t,~t\in T\}$, we can
alternatively, use a discrete-time process to model, essentially,
the same thing. This can be done by defining the discrete-time
process $\{X_n,~n=1,~2,~3,~ \ldots~\}$, where $X_n$ is a random variable
representing, for example, the number of calls arriving within the $n$th minute.

A stochastic process $\{X_t,~t\in T\}$ is called a {\em discrete space} stochastic process if the random variables $X_t$ are discrete, and it is
called {\em continuous space} stochastic process if it is continuous. We
therefore have four types of stochastic processes:
\begin{enumerate}
\item Discrete Time Discrete Space
\item Discrete Time Continuous Space
\item Continuous Time Discrete Space
\item Continuous Time Continuous Space.
\end{enumerate}

A discrete-time stochastic process $\{X_n,~n=1,~2,~3,~ \ldots~\}$ is
 {\em strictly stationary}
 if for any subset of $\{X_n\}$, say, $\{
X_{n(1)},~X_{n(2)},~ X_{n(3)},~\ldots,~X_{n(k)}\}$, for any integer $m$,
the joint probability function
$P(X_{n(1)},~X_{n(2)},~X_{n(3)},~ \ldots,~X_{n(k)})$, is equal to the
joint probability function $P(X_{n(1)+ m},~X_{n(2)+ m},~X_{n(3)+
m},~ \ldots,~X_{n(k)+m})$. In other words,\\ $P(X_{n(1)+ m},~X_{n(2)+
m},~X_{n(3)+ m},~ \ldots,~X_{n(k)+m})$ is independent of $m$. In this
case, the probability structure of the process does not change with
time. An equivalent definition for strict stationarity is applied
also to a continuous-time process, accept that in that case, $m$ is
non-integer. Notice that for the process to be strictly stationary,
the value of $k$ is unlimited as the joint probability should be
independent of $m$ for any subset of $\{X_n,~n=1,~2,~3,~ \ldots\}$. If
$k$ is limited to some value $k^*$, we say that the process is {\em
stationary of order} $k^*$.

An equivalent definition applies to a continuous-time stochastic
process. A continuous-time stochastic process $X_t$ is said to be
strictly stationary if its statistical properties do not change with
a shift of the origin. In other words, the process $X_t$ is
statistically the same as the process $X_{t-d}$ for any value of
$d$.

An important stochastic process is the {\em Gaussian Process}
defined as a process that has the property that the joint
probability function (density) associated with any set of times is
multivariate Gaussian. The importance of the Gaussian process lies
in its property to be an accurate model for the superposition of many
independent processes. This makes the Gaussian process a useful
model for heavily multiplexed traffic that arrive at switches or
routers deep in a major telecommunications network. Fortunately, the Gaussian process is not only useful, but it is also simple and
amenable to analysis. Notice that for a multivariate Gaussian
distribution, all the joint moments of the Gaussian random variables
are fully determined by the joint first- and second-order moments of
the variables.  Therefore, if the first- and second-order moments do
not change with time, the Gaussian random variables themselves are
stationary. This implies that for a Gaussian process, stationarity
of order two (also called {\em weak stationarity}) implies strict
stationarity.

For a time series $\{X_n,~n=1,~2,~3,~ \ldots~\}$, weak stationarity
implies that, for all $n$, $E[X_n]$ is constant, denoted $ E[X]$,
independent of $n$. Namely, for all $n$,
\begin{equation}
\label{Eind} E[X]=E[X_n].
\end{equation}
Weak stationarity (because it is stationarity of order two) also
implies that the covariance between $X_n$ and $X_{n+k}$, for any
$k$, is independent of $n$, and is only a function of $k$, denoted
$U(k)$. Namely, for all $n$,
\begin{equation}
\label{autocov} U(k) = Cov[X_n, X_{n+k}].
\end{equation}
Notice that, the case of $k=0$ in Eq.\@ (\ref{autocov}), namely,
\begin{equation}
\label{autocov0} U(0) = Cov[X_n, X_{n}]= Var[X_n]
\end{equation}
implies that the variance of $X_n$ is also independent of $n$. Also
for all integer $k$,
\begin{equation}
\label{negautocov}
U(-k) = U(k)
\end{equation}
because $ Cov[X_n,X_{n+k}]=Cov[X_{n+k},X_{n}]=Cov[X_n,X_{n-k}]$. The
function $U(k)$, $k = 0, 1, 2,~ \ldots$,   is called the {\em
autocovariance function}. The value of the autocovariance function
at $k$, $U(k)$, is also called the autocovariance of lag $k$.

Important parameters are the so-called Autocovariance Sum, denoted
$S$, and Asymptotic Variance Rate (AVR) denoted $v$
\cite{azatr94,azitc94}. They are defined by:
\begin{equation}
\label{autosum} S= \sum_{i=1}^{\infty}U(i)
\end{equation}
and
\begin{equation}
\label{AVR} v= \sum_{i=-\infty}^{\infty}U(i).
\end{equation}
Notice that
\begin{equation}
\label{AVRS} v= 2S+Var[X_n].
\end{equation}
Another important definition of the AVR that justifies its name is
\begin{equation}
\label{AVR2} v= \lim_{n\rightarrow \infty } \frac{Var[S_n]}{n}.
\end{equation}
We will further discuss these concepts in Section \ref{parfit}.

\subsubsection*{Homework \ref{stochastic}.\arabic{homework}}
\addtocounter{homework}{1} \addtocounter{tothomework}{1} Prove that
the above two definitions are equivalent; namely, prove that
\begin{equation}
\label{AVR3} \lim_{n\rightarrow \infty } \frac{Var[S_n]}{n}= 2S+Var[X_n]
\end{equation}
where $$S_n=\sum_{i=1}^n X_i.$$
\subsubsection*{Guide}
Define $$S(k^*)=\sum_{i=1}^{k^*} U(i)$$
and notice that $$\lim_{j\rightarrow \infty} S(j) = S.$$
Let $$S_{k^*}=\sum_{i=1}^{k^*} X_i$$
and notice that

$$\sum_{i<j}Cov[X_i,X_j]= \sum_{n=1}^{k^*-1} \sum_{k=1}^{k^* - n} Cov[X_n,X_{n+k}]= \sum_{n=1}^{k^*-1} \sum_{k=1}^{k^* - n} U(k).$$
Noticing that by the weak stationarity property, we have that $Var[X_i] = Var[X_j]$ and
$Cov[X_i,X_{i+k}] = Cov[X_j,X_{j+k}]$ for all  pairs $i,j$ ,
and letting $k^* \rightarrow \infty $, by  (\ref{varsum}), we obtain
$$Var[S_k^*] = k^* Var[X_n] + 2k^* S $$ which leads to (\ref{AVR3}).
$~~~\Box$

The {\em autocorrelation function} at lag $k$, denoted $C(k)$, is
the normalized version of the autocovariance function, and since by
weak stationarity, for all $i$ and $j$, $Var[X_j] = Var[X_i]$, it is
given by:
\begin{equation}
\label{autocovck} C(k) = \frac{U(k)}{Var[X_n]}.
\end{equation}
A stochastic process is called {\em ergodic} if every realization
contains sufficient information on the probabilistic structure of
the process. For example, let us consider a process that can be in
either one of two realizations: either $X_n=1$ for all $n$, or
$X_n=0$ for all $n$. Assume that each one of these two realizations
occurs with a probability of 0.5. If we observe any one of these
realizations, regardless of the duration of the observations, we
shall never conclude that $E[X_n]=0.5$. We shall only have the
estimations of either $E[X_n]=0$ or $E[X_n]=1$, depending on which
realization we happen to observe. Such a process is not ergodic.

Assuming $\{X_n,~n=1,~2,~3,~ \ldots~\}$ is ergodic and stationary, and we
observe $m$ observations of this $\{X_n\}$ process, denoted by
$\{\hat{A}_n,~n=1,~2,~3,~ \ldots,~m\}$, then the mean of the process
$E[A]$ can be estimated by
\begin{equation}
\label{Eest} \hat{E}[A]=\frac{1}{m} \sum_{n=1}^m \hat{A}_n,
\end{equation}
and the autocovariance function $U(k)$ of the process can be estimated by
\begin{equation}
\label{Uest} \hat{U}(k)=\frac{1}{m-k} \sum_{n=k+1}^m
(\hat{A}_{n-k} - E[A])(\hat{A}_n - E[A]).
\end{equation}
\subsection{Two Orderly and Memoryless Point Processes}
In this section, we consider a very special class of stochastic
processes called {\em point} processes that also possess two
properties: {\em orderliness} and {\em memorylessness}. After
providing, somewhat intuitive, definitions of these concepts, we
will discuss two processes that belong to this special class: one is
discrete-time - called the {\em Bernoulli process} and the other is
continuous-time - called the {\em Poisson process}.

We consider here a physical interpretation, where a {\em point
process} is a sequence of events which we call {\em arrivals}
occurring at random in points of time $t_i$, $i= 1, 2,~ \ldots$ ,
$t_{i+1}>t_i$, or $i=\ldots,~ -2, -1, 0, 1, 2, ~\ldots$,
$t_{i+1}>t_i$. The index set, namely, the time, or the set where the
$t_i$ get their values from, can be continuous or discrete, although
in most books, the index set is considered to be the real line or its non-negative part. We call our events arrivals to relate is to
the context of queueing theory, where a point process typically
corresponds to points of arrivals, i.e.,  $t_i$ is the time of the
$i$th arrival that joints a queue. A point process can be defined by
its {\em counting process}  $\{ N(t),~t\geq 0\}$, where $N(t)$ is
the number of arrivals occurred within  $[0, t)$. A counting process
$\{N(t)\}$ has the following properties:
\begin{enumerate}
\item
$N(t)\geq 0$,
 \item $N(t)$  is integer, \item if $s > t$, then $N(s) \geq N(t)$ and
$N(s)-N(t)$ is the number of occurrences within $(t,s]$.
\end{enumerate}
Note that $N(t)$ is not an independent process because, for example,
if $t_2 > t_1$, then $N(t_2)$ is dependent on the number of arrivals
in $[0,t_1)$, namely, $N(t_1)$.

Another way to define a point process is by the stochastic process
of the inter-arrival times ${\Delta_i}$ where $\Delta_i=t_{i+1}-t_i$.

One important property of a counting process is the so-called {\em
Orderliness} which means that the probability that two or more
arrivals happen at once is negligible. Mathematically, for a
continuous-time counting process to be {\em orderly}, it should
satisfy: \begin{equation} \lim_{\Delta t \rightarrow 0} P(N(t+
\Delta t) - N(t) > 1 \mid N(t+ \Delta t) - N(t) \geq 1) = 0.
\end{equation} Another very important property is the {\em
memorylessness}. A stochastic process is {\em memoryless} if at any
point in time, the future evolution of the process is statistically
independent of its past.

\subsubsection{Bernoulli Process}
The Bernoulli process is a discrete-time stochastic process made up
of a sequence of IID Bernoulli distributed random variables $\{X_i,
~i=0, 1, 2, 3,~\ldots\}$ where for all $i$, $P(X_i=1)=p$ and
$P(X_i=0)=1-p$. In other words, we divide time into consecutive
equal time slots. Then, for each time slot $i$, we conduct an independent Bernoulli
experiment. If $X_i=1$, we say that there was an {\em arrival} at
time slot $i$. Otherwise, if $X_i=0$, we say that there was no
arrival at time slot $i$.

The Bernoulli process is both orderly and memoryless. It is orderly
because, by definition, no more than one arrival can occur at any
time slot as the Bernoulli random variable takes values of more than
one with probability zero. It is also memoryless because the
Bernoulli trials are independent, so at any discrete point in time
$n$, the future evolution of the process is independent of its past.

The counting process for the Bernoulli process is another
discrete-time stochastic process $\{ N(n), n\geq 0\}$ which is a
sequence of Binomial random variables $N(n)$ representing the total
number of arrivals occurring within the first $n$ time slots. Notice
that since we start from slot 0, $N(n)$ does not include slot $n$ in
the counting. That is, we have
\begin{equation}
P[N(n)=i] =  {n \choose i} p^i(1-p)^{n-i} \qquad i=0,~1,~2,~\ldots
,~n.
\end{equation}
The concept of an inter-arrival time for the Bernoulli process can be
explained as follows. Let us assume without loss of generality that there was an arrival at time slot $k$, the inter-arrival time
will be the number of slots between $k$ and the first  time slot to
have an arrival following $k$. We do not count the time slot $k$, but we
do count the time slot of the next arrival. Because the Bernoulli process is memoryless, the inter-arrival times are IID, so we can
drop the index $i$ of ${\Delta_i}$, designating the $i$ inter-arrival
time, and consider the probability function of the random variable
${\Delta}$ representing any inter-arrival time. Because ${\Delta}$
represents a number of independent Bernoulli trials until a success, it is
geometrically distributed, and its probability function is given by
\begin{equation}
P(\Delta = i) =   p(1-p)^{i-1} \qquad i=1,~2,~\ldots ~.
\end{equation}
Another important statistical measure is the time it takes $n$ until
the $i$th arrival. This time is a sum of $i$ inter-arrival times
which is a sum of $i$ geometric random variables which we already know has a Pascal distribution with parameters $p$ and $i$, so we
have
\begin{equation}
P[{\rm the~} i {\rm th~arrival ~occurs ~in ~time ~slot~} n] = {n-1
\choose i-1} p^{i}(1-p)^{n-i} \qquad i=i,~i+1,~i+2,~\ldots~.
\end{equation}
The reader may notice that the on-off sources discussed in Section
\ref{link_dim} could be modeled as Bernoulli processes where the on-periods
are represented by consecutive successes of Bernoulli trials
and the off-periods by failures. In this case, for each on-off
process, the length of the on- and the off-periods are both
geometrically distributed. Accordingly, the {\em superposition} of
$N$ Bernoulli processes with parameter $p$ is another discrete-time stochastic process where the number of arrivals during the different
slots are IID and binomial distributed with parameters $N$ and $p$.

\subsubsection*{Homework \ref{stochastic}.\arabic{homework}}
\addtocounter{homework}{1} \addtocounter{tothomework}{1} Prove the last statement. $~~~\Box$

Another important concept is {\em merging} of processes which is
different from superposition. Let us use a sensor network example to
illustrate it. Consider $N$ sensors that are spread around a country
to detect certain events. Time is divided into consecutive fixed-length
time slots and a sensor is silent if it does not detect an event in a given time slot and active (transmitting an alarm signal) if it does. Assume that the
probability that sensor $i$ detects an event in a given time slot
is equal to $p_i$, and that the probability of such detection is
independent from time slot to time slot. Accordingly, we have a Bernoulli process associated with sensor $i$. We also assume that the $N$
Bernoulli processes associated with the $N$ sensors are independent. Assume that an alarm is sounded during a time slot when at least one of the sensors is active. We
are interested in the discrete-time process representing alarm
sounds. The probability that an alarm is sounded in a given time slot
is the probability that at least one of the sensors is active which
is one minus the probability that they are all silent. Therefore the
probability that the alarm is sounded is given by
\begin{equation}
P_a=1-\prod_{i=1}^N (1-p_i).
\end{equation}
Now, considering the independence of the processes, we can realize
that the alarms follow a Bernoulli process with parameter $P_a$.

In general, the arrival process is the process that results from the merging of
$N$ Bernoulli processes, which is the process of time slots during which at least one of the $N$ processes records an arrival. Unlike
superposition in which we are interested in the total number of
arrivals, in merging we are only interested to know if there were at least one arrival within a time slot without any interest in how
many arrivals were there in total.

Let us now consider {\em splitting}. Consider a Bernoulli process
with parameter $p$ and then color each arrival, independently of all
other arrivals, in red with probability $q$ and in blue with
probability $1-q$. Then, in each time slot, we have a red arrival with
probability $pq$ and a blue one with probability $p(1-q)$. Therefore,
the red arrivals follow a Bernoulli process with parameter $pq$ and
the blue arrivals follow a Bernoulli process with parameter
$p(1-q)$.

\subsubsection{Poisson Process} The Poisson process is a
continuous-time point process that is also memoryless and orderly.
It applies to many cases where a certain event occurs at different
points in time. Such occurrences of the events could be, for
example, arrivals of phone call requests at a telephone exchange. As
mentioned above such a process can be described by its {\em counting
process} $\{N(t), ~t\geq 0\}$ representing the total number of
occurrences by time $t$.

From now on, let us call the occurrences of the events {\it arrivals}. Using the term ``arrivals'' makes it relevant to modeling arrivals of customers to a queue, phone calls to a telephone exchange, or packets to a router which are common applications of the queueing models that are discussed in this book. Accordingly, $\{N(t), ~t\geq 0\}$ is the number of arrivals in the interval $[0, t]$.

When we consider a point process in general or a Poisson process in particular, we stumble on a problem because we cannot consider the probability of arrival at a certain point in time because such probability is always equal to zero. Therefore, we consider the probability of having a certain number of arrivals in a certain time interval or consider the probability of having an arrival during an infinitesimal time interval.

A counting process $\{N(t)\}$ is defined as {\em a Poisson process}
with rate $\lambda>0$ if it satisfies the following three
conditions.
\begin{enumerate}
\item $N(0)=0$.
\item The numbers of arrivals in two non-overlapping intervals are
independent. That is, for any $s>t>u>v>0$, the random variable
$N(s)-N(t)$ (representing the number of arrivals between time $s$ and time $t$), and the random variable $N(u)-N(v)$ (representing the number of arrivals between time $v$ and time $u$) are independent.
This means that the Poisson process has what is called {\it
independent increments}.
\item The number of arrivals in an
interval of length $t$, for $t>0$, has a Poisson distribution with mean $\lambda
t$.
\end{enumerate}
These three conditions will be henceforth called the {\em
three Poisson process conditions}.

By definition, the Poisson process $N(t)$ has what is called {\em
stationary increments} \cite{papo84,ross93}, that is, for any
$t_2>t_1$, the random variable $N(t_2)-N(t_1)$, and the random
variable $N(t_2+u)-N(t_1+u)$ have the same distribution for any
$u>0$. In both cases, the distribution is Poisson with parameter
$\lambda(t_2-t_1)$. Intuitively, if we choose the time interval
$\Delta=t_2-t_1$ to be arbitrarily small (almost a ``point'' in
time), then the probability of having an arrival there is the
same regardless of where the ``point'' is. Loosely speaking, every
point in time has the same chance of having an arrival. Therefore,
arrivals are equally likely to happen at all times. This property
is also called {\em time-homogeneity} \cite{BT02}. The Poisson process is called also a ``pure chance process'' because it represents arrivals in time that are happening by pure chance independently of each other.

Another important property of the Poisson process is that the
inter-arrival times is exponentially distributed with
parameter $\lambda$. This is shown by considering $s$ to be an
occurrence and $T$ the time until the next occurrence, noticing that
$P(T>t)=P(N(t)=0)=e^{-\lambda t}$, and recalling the properties of
independent and stationary increments. As a result, the mean
inter-arrival time is given by
\begin{equation}
\label{ET}
E[T]=\frac{1}{\lambda}.
\end{equation}
This property can serve as another definition of the Poisson process. A Poisson process with rate $\lambda$ is a point process where the inter-arrival times are independent and exponentially distributed random variables with mean $1/\lambda$.

By the memoryless property of the exponential distribution, the time
until the next occurrence is always exponentially distributed and
therefore, at any point in time, not necessarily at points of
occurrences, the future evolution of the Poisson process is
independent of the past, and is always probabilistically the same.
The Poisson process is therefore memoryless. Actually, the
independence of the past can be explained also by the Poisson
process property of {\em independent increments} \cite{ross93}, and
the fact that the future evolution is probabilistically the same can
also be explained by the stationary increments property.

An interesting paradox emerges when one considers the Poisson
process. If we consider a random point in time, independent of a
given Poisson process, the time until the next occurrence event has
exponential distribution with parameter $\lambda$. The
Poisson process is known to be {\it reversible}, namely, a Poisson process in reverse is also a Poisson process with the same rate.
In later sections, we will further discuss the reversibility concept of stochastic processes.
Then, at any
point in time, the time passed from the last Poisson occurrence
event also has an exponential distribution with parameter $\lambda$.
Therefore, if we pick a random point in time, the mean length of the
interval between two consecutive Poisson occurrences must be
$1/\lambda + 1/\lambda  = 2/\lambda$. How can we explain this
phenomenon, if we know that the time between consecutive Poisson
occurrences must be exponentially distributed with mean $1/\lambda$?
The explanation is that if we pick a point of time at random we are
likely to pick an interval that is longer than the average.

\subsubsection*{Homework \ref{stochastic}.\arabic{homework}}
\addtocounter{homework}{1} \addtocounter{tothomework}{1} Demonstrate the above paradox as follows. Generate
a Poisson process with rate $\lambda=1$ for a period of time of length $T\geq 10,000$. Pick
a point in time from a uniform distribution within the interval [1,10000].
Record the length of the interval (between two consecutive Poisson occurrences)
that includes the chosen point in time. Repeat the experiment 1000 times.
Compute the average length of the intervals you recorded. $~~~\Box$

A {\em superposition} of a number of Poisson processes is another
point process that comprises all the points of the different
processes. Another important property of the Poisson process is that
superposition of two Poisson processes with parameters $\lambda_1$
and $\lambda_2$ is a Poisson process with parameter $\lambda_1 +
\lambda_2$. Notice that in such a case, at any point in time, the
time until the next occurrence is a competition between two
exponential random variables one with parameter $\lambda_1$ and the
other with parameter $\lambda_2$. Let $T$ be the time until the
winner of the two occurs, and let $T_1$ and $T_2$ be the time until
the next occurrence of the first process and the second process,
respectively. Then, by (\ref{Pmin})\\
\begin{equation} P(T>t)= e^{-(\lambda_1+\lambda_2) t}. \end{equation} Thus, the inter-arrival
time of the superposition is exponentially distributed with
parameter $\lambda_1+\lambda_2$. This is consistent with the fact
that the superposition of the two processes is a Poisson process
with parameter $\lambda_1 + \lambda_2$.

\subsubsection*{Homework \ref{stochastic}.\arabic{homework}}
\addtocounter{homework}{1} \addtocounter{tothomework}{1} Prove that a superposition of $N$ Poisson
processes with parameters $\lambda_1$,
$\lambda_2,~\dots,~\lambda_N$, is a Poisson process with parameter
$\lambda_1+ \lambda_2 +,~\dots,+\lambda_N$.
Provide two proofs: one based on the three
Poisson process conditions, and for the other one use the
definition based on the inter-arrival
times being independent and identical
exponential random variables.
 $~~~\Box$

Another interesting question related to the superposition of Poisson
processes is the question of what is the probability that the next
event that occurs will be of a particular process. This is
equivalent to the question of having say two exponential random
variables $T_1$ and $T_2$ with parameters $\lambda_1$ and
$\lambda_2$, respectively, and we are interested in the probability
of $T_1 < T_2$. By (\ref{Pwin}), \begin{equation} P(T_1 < T_2) =
\frac{\lambda_1}{\lambda_1+\lambda_2}. \end{equation}

Before we introduce further properties of the Poisson process, we
shall introduce the following definition: a function $g(\cdot)$ is
$o(\Delta t)$ if
\begin{equation}
\label{oh}
\lim_{\Delta t \rightarrow 0} \frac{g(\Delta t)}{\Delta t} =0.
\end{equation}

Examples of functions that are $o(\Delta t)$ are $g(x)=x^v$ for $v>1$.
The sum or product of two functions which are $o(\Delta t)$ is also $o(\Delta t)$, and a constant times a function which is $o(\Delta t)$ is $o(\Delta t)$.

If a counting process $\{N(t)\}$ is {\em a Poisson process}, then
for a small (infinitesimal) interval $\Delta t$, 
we have: \begin{enumerate} \item
$P(N(\Delta t)=0) = 1- \lambda \Delta t + o(\Delta t)$ \item
$P(N(\Delta t)=1) = \lambda \Delta t + o(\Delta t)$ \item
$P(N(\Delta t)\geq 2) = o(\Delta t)$. \end{enumerate}   The above
three conditions will henceforth be called {\em small interval
conditions}. To show the first, we know that $ N(\Delta t)$ has a
Poisson distribution, therefore \begin{equation} \label{oh0}
P(N(\Delta t)=0)= e^{-\lambda \Delta t} \end{equation} and
expanding it in a Taylor series gives, \begin{equation} \label{oh1}
P(N(\Delta t)=0)=1- \lambda \Delta t + o(\Delta t). \end{equation}
The second is shown by invoking (\ref{poissonrec}) to obtain $P(N(\Delta t)=1)= \lambda
\Delta t P(N(\Delta t)=0)$ and (\ref{oh1}), or alternatively, expanding $P(N(\Delta t)=1)=e^{-\lambda \Delta t} \lambda \Delta t $ directly in a Taylor series. The third
is obtained by observing that  $ P(N(\Delta t)\geq 2)=1- P(N(\Delta t)= 1)-
P(N(\Delta t)=0)$. In fact, these three small interval conditions
plus the stationarity and independence properties together with
$N(0)=0$, can serve as an alternative definition of the Poisson
process. These properties imply that the number of occurrences per
interval has a Poisson distribution. Also, observe that the right-hand side of the above three equations (the three small interval conditions) are legitimate probabilities for infinitesimal $\Delta t$. 

\subsubsection*{Homework \ref{stochastic}.\arabic{homework}}
\addtocounter{homework}{1} \addtocounter{tothomework}{1} Prove the last statement. Namely, show
that the three small-interval conditions plus the stationarity and
independence properties together with $N(0)=0$ are equivalent to the
three Poisson process conditions.

\subsubsection*{Guide}

Define $$P_n(t) = P(N(t) = n).$$

Using the assumptions of stationary and independent increments,
show that

$$P_0(t+\Delta t) = P_0(t) P_0(\Delta t).$$ Therefore
$$\frac{P_0(t+\Delta t)-P_0(t)}{\Delta t} = P_0(t) \frac{P_0(\Delta
t)-1}{\Delta t}.$$

From the small interval conditions, we know that $P_0(\Delta t) =
1-\lambda \Delta t + o(\Delta)$, so let $\Delta t \rightarrow 0$ in
the above and obtain the differential equation: $$ {P_0}^{\prime}
(t) = - \lambda P_0(t).$$

Consider the boundary condition $P_0(0) =1$ due to the condition
$N(0) = 0$, and solve the differential equation to obtain
$$P_0(t)=e^{-\lambda t}.$$

This proves the Poisson distribution for the case $n=0$. Now
continue the proof for $n>0$. This will be done by induction, but as
a first step, consider $n=1$ to show that

$$P_1(t)=\lambda t e^{-\lambda t}.$$

Notice the $P_n(t+\Delta t)$ can be obtained by conditioning and
un-conditioning (using the Law of Total Probability) on the number
of occurrences in the interval $(t,t+\Delta t)$. The interesting
events are: \begin{enumerate} \item no occurrences with probability
$ 1 - \lambda \Delta t + o(\Delta t)$, \item one occurrence with
probability $ \lambda \Delta t + o(\Delta t)$, \item two or more
occurrences with probability $o(\Delta t)$. \end{enumerate}

Considering these events show that

$$P_n(t+\Delta t) = P_n(t)(1-\lambda \Delta t) + P_{n-1}(t)\lambda
\Delta t  + o(\Delta t)  $$

which leads to

$$\frac{P_n(t+\Delta t)-P_n(t)}{\Delta t} = - \lambda P_{n}(t) +
\lambda P_{n-1}(t) + \frac{o(\Delta t)}{\Delta t}.$$

Let $\Delta t \rightarrow 0$ in the above and obtain the
differential equation: $$ {P_n}^{\prime} (t) = - \lambda P_n(t) +
\lambda P_{n-1}(t).$$

Multiply both sides by $e^{\lambda t}$ and rearrange to obtain

\begin{equation} \label{ddt}\frac{d\{e^{\lambda
t}P_{n}(t)\}}{dt}=\lambda e^{\lambda t} P_{n-1}(t). \end{equation}

Then, use the result for $P_0(t)$ and the boundary condition of
$P_1(0)=0$ to obtain $$P_1(t)=\lambda t e^{-\lambda t}.$$

To show that the Poisson distribution holds for any $n$, assume it
holds for $n-1$, i.e., $$P_{n-1}(t)=\frac{e^{-\lambda t} (\lambda
t)^{n-1}}{(n-1)!}.$$

Now by the latter and (\ref{ddt}) obtain

$$ \frac{d\{e^{\lambda t}P_{n}(t)\}}{dt}=\frac{\lambda^n
t^{n-1}}{(n-1)!}. $$

Then, use the latter plus the condition $P_n(0) = 0$ to obtain

$$P_{n}(t)=\frac{e^{-\lambda t} (\lambda t)^{n}}{(n)!}.~~~\Box $$


A Poisson process can be viewed as a continuous-time version of the Bernoulli process \cite{Gallager11}. If we consider a Bernoulli process where time is divided into fix-length time intervals each of which is of length $\Delta t$, and the probability of having an arrival in such an interval is equal to $p= \lambda \Delta t + o(\Delta t)$, and the probability of having no arrivals is equal to $1-p$,  then as   $\Delta t \rightarrow 0$, this Bernoulli process will approach a Poisson process.

In many networking applications, it is of interest to study the
effect of {\em splitting} of packet arrival processes. In
particular, we will consider two types of splitting: {\em random
splitting} and {\em regular splitting}. To explain the difference
between the two, consider an arrival (counting) process $X(t), t \ge 0$, of packets to a certain
switch called Switch X. This packet arrival process is assumed to
follow a Poisson process with parameter $\lambda$. Some of these
packets are then forwarded to Switch A and the others to Switch B. The counting processes of packets forwarded to A and B are denoted as $X_A(t)$ and $X_B(t)$, $t\geq 0$, respectively.

Under random splitting, every packet that arrives at Switch X is
forwarded to A with probability $p$ and to B with probability $1-p$
independently of any other event associated with other packets. In
this case, the process $X_A(t), t \ge 0$  follows a Poisson process
with parameter $p\lambda$, and the process $X_B(t), t \ge 0$  follows a Poisson process
with parameter $(1-p)\lambda$. Note that $X(t)$,  $X_A(t)$, and $X_B(t)$ represent the number of occurrences during the time interval $[0, t)$ for each of these processes. We assume that the forwarding is done instantly upon arrival at Switch X, so
\begin{equation}
\label{split}
X(t) = X_A(t) + X_B(t).
\end{equation}

\subsubsection*{Homework \ref{stochastic}.\arabic{homework}}
\addtocounter{homework}{1} \addtocounter{tothomework}{1} Prove that under random splitting,
the process $X_A(t), t \geq 0$  is a Poisson
process with parameter $p\lambda$. Prove it in three ways: (1) by showing that the three Poisson process conditions hold for the $X_A(t)$ process, (2) by proving that the inter-arrival times of the $X_A$ process are IID and exponentially distributed, and (3) by showing that the small interval conditions hold for the $X_A$ Process.
\subsubsection*{Guide}
First, show that the three Poisson process conditions (with parameter $p\lambda$) hold for the $X_A(t)$ process. Let event A be defined as the event $\{ \{X_A(t)=m\} \cap \{X_B(t)=n\} \}$ for some positive integers $m$ and $n$ and let event B be the event $\{X(t)= m+n\}$. Notice that if event A occurs, then event B must occur because of (\ref{split}). Notice also that the reverse is not necessarily true.

Then, by Eq. (\ref{psubset}), we have
\begin{equation}
\label{pxa}
P(X_A(t)=m \cap X_B(t)=n) = P(X_A(t)=m \cap X_B(t)=n \mid X(t)= m+n ) P (X(t)= m+n ).
\end{equation}





Notice that the first factor in the product at the right-hand side of Equation (\ref{pxa}) is the probability of a binomial random variable with parameters $p$ and $m+n$ of having $m$ successes out of $m+n$ trials, and the second factor in the product at the right-hand side of Equation (\ref{pxa}) is the probability of a Poisson random variable with parameter $p$ to obtain value $m+n$. Therefore,

\begin{eqnarray*}
P(\{X_A(t)=m \} \cap \{ X_B(t)=n)\} \\ &
= & \left[ {m+n \choose m} p^m(1-p)^n \right] \left[ e^{-\lambda t}  \frac{(\lambda t)^{m+n}}{(m+n)!} \right] \\ &
= & \frac{(m+n)!}{m!n!}p^m(1-p)^n  e^{-\lambda t}  \frac{(\lambda t)^{m+n}}{(m+n)!}  \\ &
= & \frac{p^m(\lambda t)^m}{m!}   \frac{(\lambda t)^n (1-p)^n }{n!} e^{-\lambda t}   \\ &
= & \frac{(p \lambda t)^m}{m!}   \frac{((1-p)\lambda t)^n }{n!} e^{-p\lambda t} e^{-(1-p)\lambda t} \\ &
= & e^{-p\lambda t} \frac{(p \lambda t)^m}{m!}    e^{-(1-p)\lambda t}   \frac{((1-p)\lambda t)^n }{n!}.
\end{eqnarray*}

Then, to obtain the probability $P(X_A(t)=m), m=0.1,2, \ldots$, we use the additivity axiom that says that ``the probability of the union of mutually exclusive events is equal to the sum of the probabilities of these events'' as follows

\begin{eqnarray*}
P(\{X_A(t)=m\}) \\ &
 = & \sum_{n=0}^\infty \{ P(X_A(t)=m\} \cap\{ X_B(t)=n\} ) \\ &
 = &  \sum_{n=0}^\infty  e^{-p\lambda t} \frac{(p \lambda t)^m}{m!}    e^{-(1-p)\lambda t}   \frac{((1-p)\lambda t)^n }{n!} \\ &
  = &  e^{-p\lambda t} \frac{(p \lambda t)^m}{m!}  \sum_{n=0}^\infty  e^{-(1-p)\lambda t}   \frac{((1-p)\lambda t)^n }{n!} \\ &
  = &  e^{-p\lambda t} \frac{(p \lambda t)^m}{m!}.
 \end{eqnarray*}

 This proves that the number of occurrences in the process $X_A(t)$ during the time interval $[0, t)$ is Poisson distributed with parameter $p \lambda t$. Because of the property that the process has independent increments, this Poisson distribution with parameter $p \lambda t$ applies to any time interval of length $t$.

Notice that in the derivation above, we have used the relationship

$$ \sum_{n=0}^\infty  e^{-(1-p)\lambda t}   \frac{((1-p)\lambda t)^n }{n!} =1. $$

which holds true because a sum of probabilities of a Poisson random variable with parameter $(1-p)\lambda$ taking values from 0 to $\infty$ must be equal to 1.

 The second way to prove that the $X_A(t)$ process is a Poisson process with parameter $p\lambda$, notice that the inter-arrival times of the $X_A(t)$ process are IID random variables, each of which is a geometric (with parameter $p$) sum of exponential random variables (each with parameter $\lambda$). Such random variables are known to be exponentially distributed with a mean equal to $(1/p)(1/\lambda)$. Therefore, the $X_A(t)$ process is Poisson with rate $p\lambda$.

 The third way to prove that the $X_A(t)$ process is a Poisson process with parameter $p\lambda$, we show that the small interval conditions hold as follows.

$P (X_A(\Delta t)=0)=P(X(\Delta t)=0)+(1-p) P(X(\Delta t)=1) + o(\Delta t)=1- \lambda \Delta t + (1-p)\lambda \Delta t+ o(\Delta t) =1- p\lambda \Delta t + o(\Delta t)$,

$P(X_A(\Delta t)=1)=pP(X(\Delta t)=1)+ o(\Delta t)=p\lambda \Delta t + o(\Delta t)$,

$P(X_A(\Delta t)>1)=o(\Delta t)$,

and the stationarity and independence properties together with $X(0)=0$, follow from the same properties of the Poisson counting process $X(t)$.
$~~~\Box$

The case of splitting a Poisson process to two Poisson processes by random selection of each occurrence can be generalized to splitting a Poisson process into $N$ processes in a similar way. In general, if we consider {\bf splitting} of a Poisson process with rate $\lambda$ into $N$ point processes, where each occurrence of the original process becomes an occurrence of one of the $N$ point processes, namely, of the $i$th process with probability $p_i\geq 0$, $i=1,2,3 \ldots, N$, where $p_1 + p_2 + p_3 + \ldots + P_N = 1$. Then, each one of these $N$ processes is a Poisson process, where the rate of the $i$th process $i=1,2,3 \ldots, N$ is equal to $p_i\lambda$. To show that this is the case, we can consider one of the $N$ processes at the time, by, for example, coloring by red all occurrences with probability $p_1$ and by blue all the other occurrences (with probability $1-p_1$). Then, by the proof of the last homework problem, the process comprises the red occurrences is Poisson. We can repeat this argument to each of the processes $i=1,2,3, \ldots, N$ and show that the statement is correct.

The reader must notice that the splitting of a Poisson process into multiple Poisson processes, or splitting a Bernoulli process into Bernoulli processes, requires random selection for each occurrence. If, for example, we consider the last homework question, but the splitting is not random but regular, e.g., the first packet that arrives at Switch X
is forwarded to A, the second to B, the third to A, the fourth to B, etc. In this case, the packet stream from X to A (the X-A Process) will not follow a Poisson. It
will follow a stochastic process which is a point process where the inter-arrival times are Erlang distributed with parameters $\lambda$ and 2.

\subsubsection*{Homework \ref{stochastic}.\arabic{homework}}

\addtocounter{homework}{1} \addtocounter{tothomework}{1}
\begin{enumerate}
\item Prove the last statement.
\item Derive the mean and the variance of the inter-arrival times of
the X-A process in the two cases above: random splitting and regular
splitting.
\item Consider now 3-way splitting. Derive and compare the
mean and the variance of the inter-arrival times for the regular and
random splitting cases.
\item Repeat the above for $n$-way splitting and let $n$ increase
arbitrarily. What can you say about the burstiness/variability of
regular versus random splitting? $~~~\Box$
\end{enumerate}

The properties of the Poisson process, namely, independence and
time-homogeneity, make the Poisson process able to randomly inspect
other continuous-time stochastic processes in a way that the sample
it provides gives us enough information on what is called {\em
time-averages}. In other words, its inspections are not biased.
Examples of time averages are the proportion of time
a process $X(t)$ is in state $i$, i.e., the proportion of time
during which $X(t)=i$. Or the overall mean of the process defined by
\begin{equation}
E[X(t)] = \frac{\int_{0}^{T}X(t)dt}{T}
\end{equation}
for an arbitrarily large $T$. These properties that an occurrence
can occur at any time with equal probability, regardless of times of
past occurrences, gave rise to the expression a {\em pure chance
process} for the Poisson process.
\subsubsection*{Homework \ref{stochastic}.\arabic{homework}}
\addtocounter{homework}{1} \addtocounter{tothomework}{1} Consider a
Poisson process with parameter $\lambda$. You know that there was exactly one occurrence during the interval [0,1].
Prove that the time of the occurrence is uniformly distributed within [0,1].
\subsubsection*{Guide}
For $0\leq t \leq 1$, consider
$$P({\rm occurrence ~within~}[0,t)\mid {\rm  exactly ~one ~occurrence~within~}[0,1])$$
and use the definition of conditional probability. Notice that the latter is equal to:
$$\frac{P({\rm one~ occurrence ~within~}[0,t){\rm~and~}{\rm no ~occurrence~within~}[t,1])}{P({\rm  exactly ~one ~occurrence~within~}[0,1])}$$
or
$$\frac{P({\rm one~ occurrence ~within~}[0,t))P({\rm no ~occurrence~within~}[t,1])}{P({\rm exactly ~one ~occurrence~within~}[0,1])}.$$
Then, recall that the number of occurrences in any interval of size $T$ has
Poisson distribution with parameter $\lambda T$.
 $~~~\Box$

 \subsubsection*{Extension to multiple occurrences within an interval}
The pure chance nature of the Poisson process implies that at any point in time (or, more precisely, for any arbitrarily small interval of length $\Delta t$), the probability of an occurrence is the same, and it is independent of any other occurrences. This gives rise to the following interesting and useful property, which is an extension of the result of the last homework to multiple occurrences. If we generate $k$ uniform independent random variates within time interval [0, $T$] and we order them in increasing order of their values, the resulting point process is statistically equivalent to  a Poisson process that has exactly $k$ occurrences during the interval [0, $T$]. This implies that if one is interested in generating a realization of a Poisson process with exactly $k$ occurrences in a time interval [0, $T$], it is not necessary to repeatedly generate occurrences from a Poisson process over many such intervals until there are exactly $k$ occurrences within [0, $T$]. It is much faster to generate only $k$ independent uniform random variates. This is especially true if $k$ is large and the Poisson arrival rate is small, or vice versa, so the probability of having exactly $k$ arrivals in [0, $T$] is small.

 \subsubsection*{Mixing processes}

In addition to the Poisson process, there are other processes, the so-called {\it mixing processes}, that also have the
 property of inspections without bias. In particular, Baccelli et al. \cite{bacc06,bacc07}
promoted the use of a point process where the inter-arrival times are IID Gamma (a generalization of the Erlang random variable)
distributed for probing and measuring packet loss and delay over the Internet. Such a point-process
is a mixing process and thus can ``see time averages'' with no bias.

\subsection{Markov Modulated Poisson Process}
\label{mmpp}

The stochastic process called Markov modulated Poisson process
(MMPP) is a point process that behaves as a Poisson process with
parameter $\lambda_i$ for a period of time that is exponentially
distributed with parameter $\delta_i$. Then, it moves to mode (state)
$j$ where it behaves like a Poisson process with parameter
$\lambda_j$ for a period of time that is exponentially distributed
with parameter $\delta_j$. The $\delta$ parameters $\delta$ are called {\it mode duration parameters} \cite{ZR86a}\cite{ZR86b},\cite{ZR86c}.
In general, the MMPP can have an arbitrary
number of modes, so it requires a transition probability matrix as
an additional set of parameters to specify the probability that it
moves to mode $j$ given that it is in mode $i$. However, we are
mostly interested in the simplest case of MMPP -- the two-mode MMPP
denoted MMPP(2) and defined by only four parameters: $\lambda_0$,
$\lambda_1$, $\delta_0$, and $\delta_1$. The MMPP(2) behaves as a
Poisson process with parameter $\lambda_0$ for a period of time that
is exponentially distributed with mode duration parameter $\delta_0$. Then, it moves to
mode $1$, where it behaves like a Poisson process with parameter
$\lambda_1$ for a period of time that is exponentially distributed
with parameter $\delta_1$. Then, it switches back to mode 0, etc., alternating between the two modes 0 and 1.

A special case of MMPP(2) where $\lambda_0=0$ is called the {\em Interrupted Poisson Process} (IPP) \cite{kuczura73}. The IPP is characterized by three parameters $\lambda$, $\delta_0$, and $\delta_1$. It  behaves as a
Poisson process with parameter $\lambda>0$ during an active mode for a period of time that
is exponentially distributed with parameter $\delta_1$. Then, it is interrupted with no arrivals at all during a nonactive mode for a period of time that is exponentially distributed
with parameter $\delta_0$. Then, it switches back to an active mode, etc., alternating between the two modes active and nonactive.

 \subsection{Discrete-time Markov chains}

  \subsubsection{Definitions and Preliminaries}

 Markov chains are certain discrete space
stochastic processes which are amenable for analysis and hence are
very popular for analysis, traffic characterization, and modeling of
queueing and telecommunications networks and systems. They can be
classified into two groups:
 discrete-time Markov chains discussed here and
 continuous-time Markov chains discussed in
 the next section.

A discrete-time Markov chain is a discrete-time stochastic process
$\{X_n, n=0,~1,~2,~\ldots\}$ with the Markov property; namely, that
at any point in time $n$, the future evolution of the process is
dependent only on the state of the process at time $n$, and is independent
of the past evolution of the process. The state of the process can
be a scalar or a vector. In this section, for simplicity, we will
mainly discuss the case where the state of the process is a scalar,
but we will also demonstrate how to extend the discussion to a
multiple dimension case.

The discrete-time Markov chain $\{X_n, n=0,~1,~2,~\ldots\}$ at any
point in time may take many possible values. The set of these
possible values is finite or countable and it is called the state
space of the Markov chain, denoted by $\Theta$. A {\em
time-homogeneous Markov chain} is a process in which
$$P(X_{n+1}=i \mid X_n=j) = P(X_{n}=i \mid X_{n-1}=j)\ ~~{\rm for~
all}~n.$$ We will only consider, in this section, Markov chains
which are time-homogeneous.

A discrete-time time-homogeneous Markov chain is characterized by
the property that, for any $n$, given $X_{n}$, the distribution of
$X_{n+1}$ is fully defined regardless of states that occur before
time $n$. That is,
\begin{equation}
 P( X_{n+1} = j \mid X_n = i  ) = P( X_{n+1} = j \mid X_n =
i, X_{n-1} = i_{n-1}, X_{n-2} = i_{n-2},~\ldots).
\end{equation}

  \subsubsection{Transition Probability Matrix}

A Markov chain is characterized by the so-called {\em Transition
Probability Matrix} {\bf P} which is a matrix of one-step transition
probabilities ${P_{ij}}$ defined by
\begin{equation}
\label{onestep} P_{ij} = P( X_{n+1} = j \mid X_n = i  )~~ {\rm for~
all}~n.
\end{equation}
We can observe in the latter that the event $ \{X_{n+1} = j \}$
depends only on the state of the process at $X_n$ and the transition
probability matrix {\bf P}.

Since the $P_{ij}s$ are probabilities and since when you transit out
of state $i$, you must enter some state, all the entries in {\bf P}
are non-negatives, less or equal to 1, and the sum of entries in
each row of {\bf P} must add up to 1.

  \subsubsection{Chapman-Kolmogorov Equation}

Having defined  the one-step transition probabilities ${P_{ij}}$ in
(\ref{onestep}), let us define the $n$-step transition probability
from state $i$ to state $j$  as
\begin{equation}
\label{nstep} P_{ij}^{(n)} = P( X_{n} = j \mid X_0 = i  ).
\end{equation}
The following is known as the Chapman-Kolmogorov equation:
\begin{equation}
\label{ckeq} P_{ij}^{(n)} = \sum_{k\in \Theta} P_{ik}^{(m)}
P_{kj}^{(n-m)},
\end{equation}
 for any $m$,
such that $0<m<n$.

Let {\bf P}$^{(n)}$  be the matrix that its entries are the
$P_{ij}^{(n)}$ values.

\subsubsection*{Homework \ref{stochastic}.\arabic{homework}}

\addtocounter{homework}{1} \addtocounter{tothomework}{1}
 First, prove the Chapman-Kolmogorov equation and then
use it to prove:
\begin{enumerate}
\item {\bf P}$^{(k+n)}$ = {\bf P}$^{(k)} \times$ {\bf P}$^{(n)}$
\item {\bf P}$^{(n)}$ = {\bf P}$^n$. $~~~\Box$
\end{enumerate}

  \subsubsection{Marginal Probabilities}

Consider the marginal distribution $\pi_n(i)=P(X_n = i)$ of the
Markov chain at time $n$, over the different states $i \in \Theta$.
Assuming that the process started at time $0$, the initial
distribution of the Markov chain is $\pi_0(i)=P(X_0 = i)$, $i \in
\Theta$. Then, $\pi_n(i)$, $i \in \Theta$, can be obtained based on
the marginal probability $\pi_{n-1}(i)$ as follows
\begin{equation}
\label{onestep11} \pi_n(j) = \sum_{k \in \Theta} P_{kj}
\pi_{n-1}(k).
\end{equation}

Let the vector $\Pi_n$ be defined by $\Pi_n=
\{\pi_n(j),~j=0,~1,~2,~3,~ \ldots \}$ . The vector $\Pi_n$ can be
obtained by
\begin{equation}
\Pi_n= \Pi_{n-1} \textbf{P}=\Pi_{n-2}\textbf{P}^2=~\ldots~=\Pi_0
\textbf{P}^n.
\end{equation}

  \subsubsection{Properties and Classification of States}

One state $i$ is said to be {\em accessible} from a second state $j$
if there exists $n$, $n=0,~1,~2,~\ldots$, such that
\begin{equation}
P_{ji}^{(n)} > 0.
\end{equation}
This means that there is a positive probability for the Markov chain
to reach state $i$ at some time in the future, if it is now in state $j$.

A state $i$ is said to {\em communicate} with state $j$ if $i$ is
accessible from $j$ and $j$ is accessible from $i$.

\subsubsection*{Homework \ref{stochastic}.\arabic{homework}}

\addtocounter{homework}{1} \addtocounter{tothomework}{1}
 Prove the following:
\begin{enumerate}
\item A state communicates with itself.
\item If state $a$ communicates with $b$, then $b$ communicates with
$a$.
\item If state $a$ communicates with $b$, and $b$ communicates with
$c$, then $a$ communicates with $c$. $~~~\Box$
\end{enumerate}

A {\em communicating class} is a set of states that every pair of
states in it communicate with each other.

\subsubsection*{Homework \ref{stochastic}.\arabic{homework}}
\addtocounter{homework}{1} \addtocounter{tothomework}{1} Prove that a state cannot belong to two
different classes. In other words, two different classes must be
disjoint. $~~~\Box$

The latter implies that the state space $\Theta$ is divided into a
number (finite of infinite) of communicating classes.

\subsubsection*{Homework \ref{stochastic}.\arabic{homework}}
\addtocounter{homework}{1} \addtocounter{tothomework}{1} Provide an example of a Markov chain with
three communicating classes. $~~~\Box$

A communicating class is said to be {\em closed} if no state outside
the class is accessible from a state that belongs to the class.

A Markov chain is said to be {\em irreducible} if all the states in
its state space are accessible from each other. That is, the entire
state space is one communicating class.

A state $i$ has {\em period} $m$ if $m$ is the greatest common
divisor of the set $\{ n: P(X_n = i | X_0 = i) > 0\}$. In this case,
the Markov chain can return to state $i$ only in a number of steps
that is a multiple of $m$. A state is said to be {\em aperiodic} if
it has a period of one.

\subsubsection*{Homework \ref{stochastic}.\arabic{homework}}
\addtocounter{homework}{1} \addtocounter{tothomework}{1} Prove that in a communicating class, it
is not possible that there are two states that have different
periods. $~~~\Box$

Given that the Markov chain is in state $i$, define  {\em return
time} as the random variable representing the next time the Markov chain returns to state $i$. Notice that the return time is a random
variable $R_i$, defined by
\begin{equation}
R_i=\min \{n: X_n=i \mid X_0=i\}.
\end{equation}
A state is called {\em transient} if, given that we start in it,
there is a positive probability that we will never return back to
it. In other words, state $i$ is transient if $P(R_i < \infty) < 1$.
A state is called {\em recurrent} if it is not transient. That is, state $i$ is recurrent if
$P(R_i < \infty) = 1$.

Because in the case of a recurrent state $i$, the probability to
return to state $i$ in finite time is one, the process will visit
state $i$ infinitely many times. However, if $i$ is
transient, then the process will visit state $i$ only a
geometrically distributed number of times with parameter. (Notice
that the probability of ``success'' is $1-P(R_i < \infty)$.) In this
case the number of visits in state $i$ is finite with probability 1.

\subsubsection*{Homework \ref{stochastic}.\arabic{homework}}
\addtocounter{homework}{1} \addtocounter{tothomework}{1} Show that state $i$ is recurrent if and
only if
$$\sum_{n=0}^{\infty} p_{ii}^{(n)} = \infty.$$

\subsubsection*{Guide}
This can be shown by showing that if the condition
holds, the Markov chain will visit state $i$ an infinite number of
times, and if it does not hold, the Markov chain will visit state
$i$ a finite number of times. Let $Y_n=J_i(X_n)$, where $J_i(x)$ is
a function defined for $x = 0,1,2, ~\ldots $, taking the value 1 if
$x=i$, and 0 if $x \neq i$. Notice that $E[J_i(X_n)\mid X_0 = i] =
P(X_n=i\mid X_0 = i)$, and consider summing up both sides of the
latter. $~~~\Box$

\subsubsection*{Homework \ref{stochastic}.\arabic{homework}}
\addtocounter{homework}{1} \addtocounter{tothomework}{1} Prove that
if state $i$ is recurrent, then all the states in a class that $i$
belongs to are recurrent. In other words, you are asked to prove that recurrence is a
class property. \subsubsection*{Guide}
 Consider $m$ and $n$, such that $P_{ji}^{(m)}>0$ and
$P_{ij}^{(n)}  >0$, and argue that
$P_{ji}^{(m)}P_{ii}^{(k)}P_{ij}^{(n)}>0$ for some $m,k,n$. Then, use
the ideas and result of the previous proof. $~~~\Box$

\subsubsection*{Homework \ref{stochastic}.\arabic{homework}}
\addtocounter{homework}{1} \addtocounter{tothomework}{1} Provide an example of a Markov chain
where $P(R_i < \infty) = 1$, but $E[R_i]=\infty$. $~~~\Box$

State $i$ is called {\em positive recurrent} if $E[R_i]$ is finite.
A recurrent state that is not positive recurrent is called {\em null
recurrent}. In a finite-state Markov chain, there are no null
recurrent states, i.e.,  all recurrent states must be positive
recurrent. We say that a Markov chain is {\it stable} if all its
states are positive recurrent. This notion of stability is not
commonly used for Markov chains or stochastic processes in general
and it is different from other definitions of stability. It is
however consistent with the notion of stability of queueing systems
and this is the reason we use it here.

A Markov chain is said to be aperiodic if all its states are
aperiodic.

  \subsubsection{Steady-state Probabilities}

Consider an irreducible, aperiodic, and stable Markov chain. Then, the following limit exists.
\begin{equation}
\label{stat} \mathbf{\Pi} = \lim_{n \rightarrow \infty} \Pi_n =
\lim_{n \rightarrow \infty} \Pi_0 \textbf{P}^n
\end{equation}
and it satisfies
\begin{equation}
\label{stat1} \mathbf{\Pi} = {row ~of ~}\lim_{n \rightarrow \infty}
\textbf{P}^n
\end{equation}
where $row ~of$ ~$\lim_{n \rightarrow \infty} \textbf{P}^n$ is any
row of the matrix $\textbf{P}^n$ as $n$ approaches $\infty$. All the
rows are the same in this matrix at the limit. The latter signifies
the fact that the limit $\mathbf{\Pi}$ is independent of the initial
distribution. In other words, after the Markov chain runs for a long
time, it forgets its initial distribution and converges to
$\mathbf{\Pi}$.

We denote by $\pi_j$, $j=0, ~1,~2,~ \dots$, the components of the
vector $\mathbf{\Pi}$. That is, $\pi_j$ is the steady-state
probability of the Markov chain to be at state $j$. Namely,
\begin{equation}
\label{stat2} \pi_j = \lim_{n \rightarrow \infty} \pi_n(j)~~ {\rm
for ~all}~j.
\end{equation}
By equation (\ref{onestep11}), we obtain
\begin{equation}
\label{onestep2} \pi_n(j) = \sum_{i =0}^\infty P_{ij} \pi_{n-1}(i),
\end{equation}
then by the latter and (\ref{stat2}), we obtain
\begin{equation}
\label{stat3} \pi_j = \sum_{i=0}^\infty P_{ij} \pi_i.
\end{equation}
Therefore, recalling that $\mathbf{\Pi}$ is a proper probability
distribution, we can conclude that for an irreducible, aperiodic and
stable Markov chain, the steady-state probabilities can be obtained
by solving the following steady-state equations:
\begin{equation}
\label{sspij} \pi_j= \sum_{i=0}^{\infty} \pi_i P_{ij} ~~{\rm for~
all}~ j,
\end{equation}
\begin{equation}
\label{dsumto1} \sum_{j=0}^{\infty} \pi_j =1
\end{equation}
and
\begin{equation}
\label{nonneg}
 \pi_j \geq 0 ~~~~{\rm for~ all}~ j.
\end{equation}
In this case:
\begin{equation}
\label{piander}
 \pi_j = \frac{1}{E[R_j]}.
\end{equation}
To explain the latter, consider a large number of sequential
transitions of the Markov chain denoted $\bar{N}$,
and let $R_j(i)$ be the $i$th return time to state $j$. We assume that $\bar{N}$ is large
enough so we can neglect
edge effects. Let $N(j)$ be the number of times the process visits state $j$ during the
$\bar{N}$ sequential transitions of the Markov chain.
Then, $$
 \pi_j \approx \frac{N(j)}{\bar{N}} \approx \frac{N(j)}{\sum_{k=1}^{N(j)} R_j(k)}
 \approx \frac{N(j)}{E[R_j]N(j)}= \frac{1}{E[R_j]}.
$$

When the state space $\Theta$ is finite, one of the equations in
(\ref{sspij}) is redundant and replaced by (\ref{dsumto1}).

In matrix notation equation (\ref{sspij}) is written as:
$\mathbf{\Pi}=\mathbf{\Pi} \textbf{P}$.

Note that the steady-state vector $\mathbf{\Pi}$ is called the
stationary distribution of the Markov chain, because if we set
$\mathbf{\Pi}_0 = \mathbf{\Pi}$, $\mathbf{\Pi}_1=\mathbf{\Pi}
\textbf{P}=\mathbf{\Pi}$, $\mathbf{\Pi}_2=\mathbf{\Pi}
\textbf{P}=\mathbf{\Pi}$, $\ldots$, i.e., $\mathbf{\Pi}_n =
\mathbf{\Pi}$ for all $n$.

We know that for an irreducible, aperiodic and stable
Markov chain,$$\sum_{i=0}^\infty {P_{ji}}=1.$$ This is because we
must go from $j$ to one of the states in one step. Then, $$\pi_j\sum_{i=0}^\infty {P_{ji}}=\pi_j.$$
Then, by (\ref{sspij}), we obtain the following steady-state
equations:

\begin{equation}
\label{glob_bal} \pi_j\sum_{i=0}^\infty {P_{ji}}=
\sum_{i=0}^{\infty} \pi_i P_{ij} ~~{\rm for}~ j=0,1,2 \ldots~.
\end{equation}
These equations are called {\it global balance equations}. Equations
of this type are often used in queueing theory. Intuitively, they
can be explained as requiring that the long-term frequency of
transitions out of state $j$ should be equal to the long term
frequency of transitions into state $j$.

\subsubsection*{Homework \ref{stochastic}.\arabic{homework}}
\addtocounter{homework}{1} \addtocounter{tothomework}{1}
 \begin{enumerate}
\item Show that a discrete-time Markov chain (MC) with two states where the
rows of the transition probability matrix are identical is a
Bernoulli process.
\item Prove that in any finite MC, at least one state must be
recurrent.
\item Provide examples of MCs defined by their
transition probability matrices that their states (or some of the
states) are periodic, aperiodic, transient, null recurrent and
positive recurrent. Provide examples of irreducible and reducible
(not irreducible) and of stable and unstable MCs. You may use as
many MCs as you wish to demonstrate the different concepts.
\item For different $n$ values, choose an $n\times n$
transition probability matrix $\textbf{P}$ and an initial vector
$\mathbf{\Pi}_0$. Write a program to compute $\mathbf{\Pi}_1$,
$\mathbf{\Pi}_2$, $\mathbf{\Pi}_3$, $\ldots$ and demonstrate
convergence to a limit in some cases and demonstrate that the limit
does not exist in other cases.
\item Prove equation (\ref{piander}).
\item Consider a binary communication channel between a transmitter
and a receiver where $B_n$ is the value of the $n$th bit at the
receiver. This value can be either equal to $0$, or equal to $1$.
Assume that the event [a bit to be erroneous] is independent of the
value received and only depends on whether or not the previous bit
is erroneous or correct. Assume the following:\\
$P(B_{n+1}$ is erroneous $\mid B_n$ is correct) = 0.0001\\
$P(B_{n+1}$ is erroneous $\mid B_n$ is erroneous) = 0.01\\
$P(B_{n+1}$ is correct $\mid B_n$ is correct) = 0.9999\\
$P(B_{n+1}$ is correct $\mid B_n$ is erroneous) = 0.99\\
Compute the steady-state error probability. $~~~\Box$
 \end{enumerate}

\subsubsection{Birth and Death Process}

In many real-life applications, the state of the system either stays unchanged or sometimes
increases by one, and at other times decreases by one, and no other
transitions are possible. Such a discrete-time Markov chain
$\{X_n\}$ is called a {\em birth-and-death process}. In this case,
$P_{ij}=0$ if $|i-j|>1$ and $P_{ij}>0$ if $|i-j| \leq 1$, for $i\geq 0$ and $j\geq 0$.

Then, by the first equation of (\ref{sspij}), we obtain,
$$\pi_0 P_{01} = \pi_{1} P_{10}. $$
Then, substituting the latter in the second equation of
(\ref{sspij}), we obtain
$$\pi_1 P_{12} = \pi_{2} P_{21}. $$
Continuing in the same way, we obtain
\begin{equation}
\label{local_balance}
 \pi_i P_{i,i+1} = \pi_{i+1} P_{i+1,i}, ~~i=0,1,2, \ldots ~.
\end{equation}
These equations are called {\it detailed balance equations}. They 
together with the normalizing equation $$\sum_{i=0}^{\infty} \pi_i=1$$
constitute a set of steady-state equations for the steady-state
probabilities which are far simpler than (\ref{sspij}).

\subsubsection*{Homework \ref{stochastic}.\arabic{homework}}
\addtocounter{homework}{1} \addtocounter{tothomework}{1} Solve the
detailed balance equations together with the normalizing equations for
the $ \pi_i$, $i=0, 1, 2, \ldots~.$
\subsubsection*{Guide}
Recursively, write all $\pi_i$, $i=0, 1, 2, 3, \ldots$ in terms of
$\pi_0$. Then, use the normalizing equation and isolate $\pi_0$.
 $~~~\Box$

   \subsubsection{Reversibility}
   \label{reversibility_DT}

Consider an irreducible, aperiodic, and stable Markov chain
$\{X_n\}$. Assume that this Markov chain has been running for a long
time to achieve stationarity with the transition probability matrix
$\textbf{P}=[P_{ij}]$, and consider the process
$X_n,~X_{n-1},~X_{n-2},~\ldots$, going back in time. This reversed
process is also a Markov chain because $X_n$ has dependence
relationship only with $X_{n-1}$ and $X_{n+1}$ and conditional on $X_{n+1}$, it is independent of
$X_{n+2},~X_{n+3},~X_{n+4},~\ldots$~~. Therefore, $$P( X_{n-1} = j
\mid X_n = i  ) = P( X_{n-1} = j \mid X_n = i, X_{n+1} = i_{n+1},
X_{n+2} = i_{n+2},~\ldots).$$

In the following, we derive the transition probability matrix,
denoted $\textbf{Q}=[Q_{ij}]$ of the process $\{X_n\}$ in reverse.
Accordingly, define \begin{equation} \label{Qij} Q_{ij} = P ( X_n=j
\mid X_{n+1}=i). \end{equation} By the definition of conditional
probability, we obtain, \begin{equation} \label{Qij1} Q_{ij} = \frac
{P (X_n=j \cap X_{n+1}=i )}{P(X_{n+1}=i)} \end{equation} or
\begin{equation} \label{Qij2} Q_{ij} = \frac {P(X_n=j) P( X_{n+1}=i
\mid X_n=j ) } {P(X_{n+1}=i)} \end{equation} and if $\pi_j$ denotes
the steady-state probability of the Markov chain $\{X_n\}$ to be in
state $j$, and let $n\rightarrow \infty$, we obtain \begin{equation}
\label{Qij3} Q_{ij} = \frac {\pi_j P_{ji}  }{ \pi_i}. \end{equation}
A Markov chain is said to be {\em time reversible} if
$Q_{ij}=P_{ij}$ for all $i$ and $j$. Substituting $Q_{ij}=P_{ij}$ in
(\ref{Qij3}), we obtain, \begin{equation} \label{Qij4} \pi_i P_{ij}
= \pi_j P_{ji} ~{\rm for~all} ~i,j. \end{equation} The set of
equations (\ref{Qij4}) is also a necessary and sufficient condition
for time reversibility. This set of equations is called the {\it
detailed balance conditions}. In other words, a necessary and
sufficient condition for reversibility is that there exists a solution that sums up to unity for the detailed balance conditions.
Furthermore, if such a solution exists, it is the stationary
probability of the Markov chain, namely, it also solves the
global balance equations.

Intuitively, a Markov chain $X_n$ is time-reversible if for a large
$k$ (to ensure stationarity) the Markov chain $X_k,~ X_{k+1},
X_{k+2} ~ \ldots$ is statistically the same as the process $X_k,~
X_{k-1}, X_{k-2} ~ \ldots$. In other words, by considering the
statistical characteristics of the two
processes, you cannot tell which one is going forward and which is
going backward.

\subsubsection*{Homework \ref{stochastic}.\arabic{homework}}
\addtocounter{homework}{1} \addtocounter{tothomework}{1} Provide an example of a Markov chain that
is time reversible and another one that is not time reversible.
$~~~\Box$

  \subsubsection{Multi-dimensional Markov Chains}

So far, we discussed single-dimensional Markov chains. If the state
space is made of finite vectors instead of scalars, we can easily
convert them to scalars and proceed with the above-described
approach. For example, if the state-space is {(0,0) (0,1) (1,0)
(1,1)} we can simply change the names of the states to {0,1,2,3} by
assigning the values 0, 1, 2 and 3 to the states (0,0), (0,1), (1,0) and (1,1),
respectively. In fact, we do not even have to do it explicitly. All
we need to do is to consider a $4\times 4$ transition probability
matrix as if we have a single dimension Markov chain. Let us now
consider an example of a multi-dimensional Markov chain.

Consider a bit stream transmitted through a channel. Let $Y_n=1$ if
the $n$th bit is received correctly, and let $Y_n=0$ if the $n$th
bit is received incorrectly. Assume the following

$ P( Y_{n} = i_n \mid Y_{n-1} = i_{n-1},  Y_{n-2} = i_{n-2}  ) $\\
 $~~~~~~~~~~~~~~~~~~=P( Y_{n} = i_n \mid Y_{n-1} = i_{n-1},  Y_{n-2} = i_{n-2},
 Y_{n-3} = i_{n-3}, Y_{n-4} = i_{n-4},~\ldots  ).$

$P( Y_{n} = 0 \mid Y_{n-1} = 0,  Y_{n-2} = 0  ) = 0.9$\\
$P( Y_{n} = 0 \mid Y_{n-1} = 0,  Y_{n-2} = 1  ) = 0.7$\\
$P( Y_{n} = 0 \mid Y_{n-1} = 1,  Y_{n-2} = 0  ) = 0.6$\\
$P( Y_{n} = 0 \mid Y_{n-1} = 1,  Y_{n-2} = 1  ) = 0.001.$

By the context of the problem, we have

$P( Y_{n} = 1) = 1 -  P( Y_{n} = 0)$

so,

$P( Y_{n} = 1 \mid Y_{n-1} = 0,  Y_{n-2} = 0  ) = 0.1$\\
$P( Y_{n} = 1 \mid Y_{n-1} = 0,  Y_{n-2} = 1  ) = 0.3$\\
$P( Y_{n} = 1 \mid Y_{n-1} = 1,  Y_{n-2} = 0  ) = 0.4$\\
$P( Y_{n} =
1 \mid Y_{n-1} = 1,  Y_{n-2} = 1  ) = 0.999.$

\subsubsection*{Homework \ref{stochastic}.\arabic{homework}}
\addtocounter{homework}{1} \addtocounter{tothomework}{1}
 Explain why the process $\{Y_{n}\}$ is
not a Markov chain. $~~~\Box$

Now define the $\{X_n\}$ process as follows: \\
$X_{n} = 0$ if $ Y_{n} = 0$ and $Y_{n-1} = 0$.\\
$X_{n} = 1$ if $ Y_{n} = 0$ and $Y_{n-1} = 1$.\\
$X_{n} = 2$ if $ Y_{n} = 1$ and $Y_{n-1} = 0$.\\
$X_{n} = 3$ if $ Y_{n} = 1$ and $Y_{n-1} = 1$.

\subsubsection*{Homework \ref{stochastic}.\arabic{homework}}
\addtocounter{homework}{1} \addtocounter{tothomework}{1} Explain why the process $\{X_{n}\}$ is a
Markov chain, produce its transition probability matrix, and compute
its steady-state probabilities. $~~~\Box$

\subsection{Continuous Time Markov chains}
\label{ctmc}

  \subsubsection{Definitions and Preliminaries}

A continuous-time Markov chain is a continuous-time stochastic
process $\{X_t, t \geq 0\}$. At any point in time $t>0$, $X_t$ describes the
state of the process which is discrete. We will consider only
a continuous-time Markov chain where $X_t$ takes values that are
nonnegative integers. The time between changes in the state of the
process is exponentially distributed. In other words, the process
stays constant for an exponential time duration before changing to
another state.

In general, a continuous-time Markov chain $\{X_t, t \geq 0\}$ is defined by
the property that for all real numbers $s \geq 0$, $t \geq 0$ and
$v\geq 0$, and integers $i \geq 0$, $j \geq 0$ and $k \geq 0$,
\begin{equation}
P(X_{t+s} = j \mid X_t = i, X_v = k_v, v\leq t) = P(X_{t+s} = j \mid
X_t = i ).
\end{equation}
That is, the probability distribution of the future values of the
process $X_t$, represented by $X_{t+s}$, given the present value of
$X_t$ and the past values of $X_t$ denoted $X_v$, is independent of
the past and depends only on the present.

A general continuous-time Markov chain can also be defined as a
continuous-time discrete space stochastic process with the following
properties.
\begin{enumerate}
\item  Each time the process enters state $i$,  $i\geq 0$, it stays at that state for
an amount of time which is exponentially distributed with parameter
$\delta_i$ before making a transition into a different state.
\item  When the process leaves state $i$, it enters state $j$,  $j\geq 0$, with
probability denoted $P_{ij}$. The set of $P_{ij}$s must satisfy the
following:
\begin{eqnarray}
(1) & P_{ii}=0 ~~{\rm for~ all}~ i \nonumber \\
(2) &  \sum_j P_{ij} = 1. \nonumber
\end{eqnarray}
\end{enumerate}

An example of a continuous-time Markov chain is a Poisson process
with rate $\lambda$. The state at time $t$, $X_t$ can be the
number of occurrences by time $t$ which is the counting process
$N(t)$. In this example of the Poisson counting process
$\{X_t\}=N(t)$ increases by one after every exponential time
duration with parameter $\lambda$.

Another example is the so-called {\em pure birth process} $\{X_t, t \geq 0\}$.
It is a generalization of the counting Poisson process. Again
$\{X_t, t \geq 0\}$ increases by one every exponential amount of time but
here, instead of having a fixed parameter $\lambda$ for each of
these exponential intervals, this parameter depends on the state of
the process and it is denoted $\delta_i$. In other words, when
$X_t=i$, $i\geq 0$, the time until the next occurrence in which the process $\{X_t, t \geq 0\}$
increases from $i$ to $i+1$ is exponentially distributed with
parameter $\delta_i$. If we set $\delta_i=\lambda$ for all $i\geq 0$, we
have the Poisson counting process.

\subsubsection{Birth and Death Process}

As in the case of the discrete-time Markov chain, in many real-life
applications, such as various queueing systems, that lend themselves
to continuous-time Markov-chain modeling, the state of the system
at one point in time sometimes increases by one and at other times
decreases by one, but never increases or decreases by more than one at
one-time instance. Such a continuous-time Markov chain $\{X_t, t \geq 0\}$, as
its discrete-time counterpart, is called a {\em birth-and-death
process}. In such a process, the time between occurrences in the state
$i, ~i=0, 1, 2, \ldots,$ is exponentially distributed, with parameter $\delta_i$, and at
any point of occurrence, the process increases by one (from its
previous value $i$ to $i+1$) with probability $\upsilon_i$, for $i>0$, and with probability 1 for $i=0$, and it
decreases by one (from $i$ to $i-1$) with probability
$\vartheta_i=1-\upsilon_i$, for $i > 0$, and with probability 0, for $i=0$. The transitions from $i$ to $i+1$ are
called {\em births} and the transitions from $i$ to $i-1$ are called
{\em deaths}. Recall that the mean time between occurrences, when in
state $i$, is $1/\delta_i$. Hence, the birth rate in state $i$,
denoted $b_i$, is given by

\begin{equation}
\label{brate}
b_i=\left\{\begin{array}{ll}
\delta_i\upsilon_i & \mbox{for $i \geq 1$}\\
\delta_i & \mbox{for $i = 0$,}
\end{array}
\right.
\end{equation}

and the death rate ($d_i$) is given by

\begin{equation}
\label{drate}
d_i=\left\{\begin{array}{ll}
\delta_i\vartheta_i & \mbox{for $i \geq 1$}\\
0 & \mbox{for $i = 0$.}
\end{array}
\right.
\end{equation}


Summing up these two equations
gives the intuitive result that the total rate at state $i$ is equal
to the sum of the birth-and-death rates. Namely,
$$\delta_i=b_i+d_i, ~\mbox{for $i \geq 0$,} $$ and therefore the mean time between
occurrences is
$$\frac{1}{\delta_i}=\frac{1}{b_i+d_i},~\mbox{for $i \geq 0$}.$$

\subsubsection*{Homework \ref{stochastic}.\arabic{homework}}
\addtocounter{homework}{1} \addtocounter{tothomework}{1} Show the following:
$$\vartheta_i= \frac{d_i}{b_i+d_i}, ~\mbox{for $i \geq 0$}$$
and
$$\upsilon_i= \frac{b_i}{b_i+d_i}, ~\mbox{for $i \geq 0$}. ~~~\Box $$

Birth-and-death processes apply to queueing systems where customers
arrive one at a time and depart one at a time. Consider for example
a birth-and-death process with the death rate higher than the birth
rate. Such a process could model, for example, a stable single-server queueing system.

\subsubsection{First Passage Time}

An important problem that has applications in many fields, such as
biology, finance, and engineering, is how to derive the distribution
or moments of the time it takes for the process to transit from
state $i$ to state $j$.  In other words, given that the process is
in state $i$ find the distribution of a random variable representing
the time it takes to enter state $j$ for the first time. This random
variable is called the {\em first passage time from i to j}. Let us
derive the mean of the first passage time from $i$ to $j$ in a birth-and-death process
for the case $i<j$. To solve this problem we start
with a simpler one. Let $U_i$ be the mean passage time to go from
$i$ to $i+1$. Then, \begin{equation}
\label{mptu0} U_0=\frac{1}{b_0}.
\end{equation}
and
\begin{equation}
\label{mptui}U_i=\frac{1}{\delta_i} + \vartheta_i[U_{i-1}+U_i], ~~\mbox{for $i \geq 1$}.
\end{equation}

\subsubsection*{Homework \ref{stochastic}.\arabic{homework}}
\addtocounter{homework}{1} \addtocounter{tothomework}{1} Explain equations (\ref{mptu0}) and
(\ref{mptui}).
\subsubsection*{Guide}
Notice that $U_{i-1}$ is the mean passage time to go from
$i-1$ to $i$, so $U_{i-1}+U_i$ is the mean passage time to go from
$i-1$ to $i+1$. Equation (\ref{mptui}) essentially says that $U_i$ the mean passage time to go from
$i$ to $i+1$ is equal to the mean time the process stays in state $i$ (namely ${1}/{\delta_i}$),
plus the probability to move from $i$ to $i-1$, times
the mean passage time to go from
$i-1$ to $i+1$. Notice that the probability of moving from $i$ to $i+1$ is not considered because
if the process moves from $i$ to $i+1$ when it completes its sojourn in state $i$,
then the process reaches the target (state $i+1$), so no
further time needs to be considered.
$~~~\Box$

Therefore,
\begin{equation}
\label{mptui1}U_i=\frac{1}{b_i+d_i} +
\frac{d_i}{b_i+d_i}[U_{i-1}+U_i], ~~~i\geq 1,
\end{equation}
or
\begin{equation}
\label{mptui2} U_i = \frac{1}{b_i} + \frac{d_i}{b_i}U_{i-1}, ~~~ i\geq 1.
\end{equation}

Now we have a recursion by which we can obtain $U_0, U_1, U_2,~
\ldots$, and the mean first passage time between $i$ and $j$, for $j>i$, is
given by the sum $$\sum_{k=i}^{j-1} U_k.$$

\subsubsection*{Homework \ref{stochastic}.\arabic{homework}}
\addtocounter{homework}{1} \addtocounter{tothomework}{1} Let $b_i=\lambda$ for $i\geq 0$ and $d_i=\mu$ for
$i\geq 1$, derive a closed-form expression for $U_i$. $~~~\Box$

\subsubsection{Transition Probability Function}

For a continuous-time Markov chain $\{X_t, t \geq 0\}$, define
the {\it transition probability function} $P_{ij}(t)$ as the probability that given that the process is in state $i$ at time $t_0$, then a time $t$ later, it will be in state $j$. That is,
\begin{equation}
P_{ij}(t) = P[X(t_0+t)=j \mid X(t_0)=i].
\end{equation}
The continuous time version of the Chapman-Kolmogorov equations are
\begin{equation}
\label{chapko} P_{ij}(t+\tau)=\sum_{n=0}^\infty
P_{in}(t)P_{nj}(\tau)~~{\rm for~ all}~t\geq 0, \tau\geq 0.
\end{equation}
Using the latter to derive the limit
$$\lim_{\triangle t \rightarrow 0} \frac{P_{ij}(t+\triangle
t)-P_{ij}(t)} {\triangle t}$$ we obtain the so-called Kolmogorov's
Backward Equations:
\begin{equation}
\label{koback} P_{ij}^\prime(t)=\sum_{n\neq i} \delta_i P_{in}
P_{nj}(t) - \delta_i P_{ij}(t) ~{\rm for~all}~i,j ~{\rm and}~ t\geq
0.
\end{equation}
For a birth-and-death process, the latter becomes
\begin{equation}
\label{koback_bad} P_{0j}^\prime(t)=b_0\{ P_{1j}(t)- P_{0j}(t)\}
\end{equation}
and
\begin{equation}
\label{koback_bad2} P_{ij}^\prime(t)= b_i P_{i+1,j}(t) +d_i
P_{i-1,j} (t) - (b_i +d_i) P_{ij}(t) ~{\rm for~all~} i>0.
\end{equation}

\subsubsection{Steady-State Probabilities}
\label{prob_flux}

As in the case of the discrete-time Markov chain, define a
continuous-time Markov chain to be called {\em irreducible} if there
is a positive probability for any state to reach every state, and we
define a state in a continuous-time Markov chain to be called {\em positive
recurrent}, if the process visits and then leaves
that state, the random variable that represents the time it returns
to that state has a finite mean. As for discrete-time Markov chains, a
continuous-time Markov chain is said to be {\it stable} if all its
states are positive recurrent.

Henceforth, we only consider continuous-time Markov chains that are
irreducible, aperiodic and stable. Then, the limit of $P_{ij}(t)$ as
$t$ approaches infinity exists, and we define
\begin{equation}
\pi_j = \lim_{t\rightarrow \infty} P_{ij}(t).
\end{equation}
The $\pi_j$ values are called steady-state probabilities or
stationary probabilities of the continuous-time Markov chain. In
particular,
 $\pi_j$ is the steady-state probability of
the continuous-time Markov chain to be at state $j$. We shall now
describe how the steady-state probabilities $\pi_j$ can be obtained.

We now construct the matrix {\bf Q} which is called the {\em infinitesimal
generator} of the continuous-time Markov chain. The matrix {\bf Q} is a
matrix of one step infinitesimal rates ${Q_{ij}}$ defined by
\begin{equation}
\label{qij} Q_{ij}= {\delta_i}{P_{ij}}   {\rm ~for~} i \not= j
\end{equation}
and
\begin{equation}
\label{qii} Q_{ii}=-\sum_{j\not=i} { Q_{ij} }.
\end{equation}
{\bf Remarks:}
\begin{itemize}
\item The state-space can be finite or infinite and hence the
matrices {\bf P } and {\bf Q } can also be finite or infinite.
\item In Eq.\@ (\ref{qij}), $Q_{ij}$ is the product of the rate to
leave state $i$ and the probability of transition to state $j$
from state $i$ which is the rate of transitions from $i$ to $j$.
\end{itemize}

To obtain the steady-state probabilities $\pi_j$s, we solve the
following set of steady-state equations:
\begin{equation}
\label{ssqij}
0 = \sum_i \pi_i Q_{ij} ~~{\rm for~ all}~ j
\end{equation}
and the {\it normalization equation} that ensures that the sum of the steady-state probabilities is equal to one is
\begin{equation}
\label{sumto1}
\sum_j \pi_j =1.
\end{equation}
Denoting $ \mathbf{\Pi}=[\pi_0, \pi_1, \pi_2,~ \ldots ] $, Eq.\@
(\ref{ssqij}) can be written as
\begin{equation}
\label{matrixqij} 0=\mathbf{\Pi Q}.
\end{equation}

To explain Eqs. (\ref{ssqij}), notice that, by (\ref{qij}) and
(\ref{qii}), for a particular $j$, the steady-state equation
\begin{equation}
\label{ssqj} 0 = \sum_i \pi_i Q_{ij}
\end{equation}
is equivalent to the steady-state equations
\begin{equation}
\pi_j \sum_{i\not=j} Q_{ji} =\sum_{i\not=j} \pi_i Q_{ij}
\end{equation}
 or
\begin{equation}
\label{equilib} \pi_j\sum_{i\not=j} { \delta_j P_{ji}
}=\sum_{i\not=j} \pi_i \delta_i P_{ij}
\end{equation}

which give the following global balance equations if we consider all
$j$.

\begin{equation}
\label{gb_equilib} \pi_j\sum_{i\not=j} { \delta_j P_{ji}
}=\sum_{i\not=j} \pi_i \delta_i P_{ij} ~~{\rm for~ all}~ j,
\end{equation}
or using the $Q_{ij}$ notation,
\begin{equation}
\label{gb_equilib1} \pi_j\sum_{i\not=j} Q_{ji} =\sum_{i\not=j} \pi_i
Q_{ij} ~~{\rm for~ all}~ j.
\end{equation}

The quantity $\pi_i Q_{ij}, i \neq j$ which is the steady-state probability of
being in state $i$ times the infinitesimal rate of a transition from
state $i$ to state $j$ is called the {\em probability flux} from
state $i$ to state $j$. Eq.\@ (\ref{ssqj}) says that the total
probability flux from all states into state $j$ is equal to the
total probability flux out of state $j$ to all other states. To
explain this equality, consider a long period of time $L$. Assuming
the process returns to all states infinitely many times, during a
long time period $L$, the number of times the process moves into
state $j$ is equal (in the limit $L\rightarrow \infty$) to the
number of times the process moves out of state $j$. This leads to
Eq.\@ (\ref{equilib}) with the factor $L$ in both sides. The concept
of probability flux is equivalent to the concept of the long-term
frequency of transitions discussed above in the context of discrete-time Markov chains.

Similar to the case of discrete-time Markov chains, the set of steady-state equations
(\ref{ssqij}) and (\ref{sumto1}) is dependent and one of the equations in
(\ref{ssqij}) is redundant in the finite state space case.

For continuous-time birth-and-death processes, $Q_{ij}=0$ for
$|i-j|>1$. As in the discrete-time case, under this special
condition, the global balance equations (\ref{gb_equilib}) can be
simplified to the detailed balance equations.   We start with the first
equation of (\ref{gb_equilib}) and using the condition $Q_{ij}=0$
for $|i-j|>1$, we obtain
\begin{equation}
\label{lb_1}
\pi_0 Q_{01} = \pi_1 Q_{10}.
\end{equation}
The second equation is
\begin{equation}
\label{lb_2}
\pi_1 [Q_{10} + Q_{12}] =\pi_0 Q_{01} + \pi_2 Q_{21}.
\end{equation}
Then, Eq.\@ (\ref{lb_2}) can be simplified using (\ref{lb_1}) and we
obtain
\begin{equation}
\label{lb_2_1} \pi_1 Q_{12} = \pi_2 Q_{21}.
\end{equation}
In a similar way, by repeating the process, we obtain the following
detailed balance equations.
\begin{equation}
\label{lb} \pi_i Q_{i,i+1} = \pi_{i+1} Q_{i+1,i} ~~i=0,1,2,
~\ldots~.
\end{equation}

  \subsubsection{Multi-Dimensional Continuous-time Markov Chains}

The extension discussed  earlier, regarding multi-dimensional
discrete-time Markov chains, applies also to the case of
continuous-time Markov chains. If the state-space is made of finite
vectors instead of scalars, as discussed, there is a one-to-one
correspondence between vectors and scalars, so a multi-dimensional
continuous-time Markov chain can be converted to a single-dimension
continuous-time Markov chain and we proceed with the above-described
approach that applies to the single dimension. Note that the above-described MMPP is a Multi-Dimensional continuous-time Markov chain.
For more details on Multi-Dimensional Continuous-time Markov chains see \cite{Iver15}.

\subsubsection{Solutions by Successive Substitutions}
\label{succsub}
When the matrix {\bf Q} is large but finite, there is a need to solve a set of steady-state equations efficiently. An efficient way that normally works well for solving  such equations is to use a method of successive substitutions that applies to both discrete-time and continuous-time Markov chains, but in our description of the method we
consider for example a finite set of steady-state equations of a continuous-time Markov chain of the form

\begin{equation} \label{sseq} 0= \mathbf{\Pi}\mathbf{Q}
\end{equation}
where $\mathbf{\Pi}=[\pi_0, \pi_{1}, \pi_{2}, \pi_{3}, \dots,  \pi_{k}]$
and $\mathbf{Q}$ is the infinitesimal generator $(k+1)\times (k+1)$ matrix of the continuous-time Markov chain,
and we also have the normalization equation (\ref{sumto1}) that ensures that the sum of the steady-state probabilities is equal to one.

The method can be described as follows.
First, isolate the first element of the vector $\mathbf{\Pi}$; in this case, it is the variable $\pi_{0}$ in the first equation.
 Next, isolate  the second element of the vector $\mathbf{\Pi}$, namely, $\pi_{1}$ in the second equation, and then keep isolating all the variables of the vector $\mathbf{\Pi}$. This leads to the
following vector equation for $\mathbf{\Pi}$
\begin{equation}
\label{fixedp} \mathbf{\Pi} = \mathbf{\Pi}\mathbf{\hat{Q}},
\end{equation}
 where $\mathbf{\hat{Q}}$ is different from the original $\mathbf{Q}$ because of the algebraic operations we performed when we isolated the elements of the  $\mathbf{\Pi}$ vector.
Then, perform the successive substitution operations by first setting arbitrary
initial values to the vector $\mathbf{\Pi}$; substitute them in the right-hand side of (\ref{fixedp}) and obtain different values at the left-hand side which are then substituted back in the right-hand side, etc.
For example, the initial setting can be $\mathbf{\Pi} = 1$ without
any regard to the normalization equation (\ref{sumto1}). When the values obtained
for $\mathbf{\Pi}$ are sufficiently close to those obtained in the previous subsection, say, within a distance no
more than $10^{-6}$, stop. Finally, normalize the vector $\mathbf{\Pi}$
obtained in the last iteration using (\ref{sumto1}). This is the desired solution.
There are several versions of this method called Jacobi \cite{MAA}, Gauss–Seidel \cite{cooper81}, and  successive over-relaxation \cite{cooper81}.  In the Jacobi method, in each iteration, all the new $\pi_i$ values are obtained based only on the $\pi_i$ values of the previous iteration, while in Gauss–Seidel, every time a new $\pi_i$ value is computed it is fed back to replace the $\pi_i$ of the previous iteration to compute $\pi_{i+1}, \pi_{i+2}, \pi_{i+3}, \dots,  \pi_{k}$. The successive over-relaxation method is used to speed up the computation. See \cite{cooper81,MAA} for more details on these methods.

\subsubsection{The Curse of Dimensionality}

In many applications, the {\bf Q} matrix is too large, so it may not
be possible to solve the steady-state equations (\ref{ssqij}) in
a reasonable time. Actually, the case of a large state space (or large
{\bf Q} matrix) is common in practice.

This often occurs when the application leads to a Markov-chain model that is of high dimensionality. Consider for example a 49 cell
GSM mobile network, and assume that every cell has 23 voice
channels. Assuming Poisson arrivals, exponential call time duration, and exponential cell sojourn times. Then, this cellular mobile network can be modeled as a
continuous time Markov chain with each state representing  the
number of busy channels in each cell. In this case, the number of
states is equal to $24^{49}$, so a numerical solution of the
steady-state equations is computationally prohibitive.

\subsubsection{Simulations}

When an exact numerical solution is not attainable, we often rely on
simulations. Fortunately, due to the special structure of the
continuous-time Markov chain together with a certain property of the
Poisson process called PASTA (Poisson Arrivals See Time Averages),
simulations of continuous-time Markov-chain models can be simplified
and expedited so they lead to accurate results. To explain the PASTA
property, consider a stochastic process for which steady-state
probabilities exist. If we are interested in obtaining certain
steady-state statistical characteristics of the process (like the
$\pi_i$ in a continuous-time Markov chain), we could inspect the entire evolution of the process (in practice, for a long enough time
period), or we could use an independent Poisson inspector. (We
already discussed the property of the Poisson process to see
time-averages.) The PASTA principle means that if the arrivals
follow a Poisson process, we do not need a separate Poisson
inspector, but we could inspect the process at occurrences of points
in time just before points of arrivals.

Note that in practice, since we are limited to a finite number of
inspections, we should choose a Poisson process that will have
a sufficient number of occurrences (inspections) during the simulation
of the stochastic process that we are interested in obtaining its
steady-state statistics.

In many cases, when we are interested in steady-state statistics of
a continuous time Markov chain, we can conveniently find a Poisson
process which is part of the continuous-time Markov chain that we are
interested in and use it as a Poisson inspector. For example, if we
consider a queueing system in which the arrival process follows a
Poisson process, such process could be used for arrival times of
the inspector if it, at any inspection, does not count (include) its
own particular arrival. In other words, we consider a Poisson inspector that arrives just before its own arrival occurrences.

See Chapter \ref{simulations} for more information and examples on simulations.

\subsubsection{Reversibility}
\label{reversibility_CT}

We have discussed the {\bf time reversibility} concept in the
context of discrete-time Markov chains. In the case of a
continuous-time Markov chain, the notion of time reversibility is
similar. If you observe the process $X_t$ for a large $t$ (to ensure
stationarity) and if you cannot tell from its statistical behavior
if it is going forward or backward, it is time reversible.

Consider a stationary continuous-time Markov chain that has a unique steady-state
solution.  Its $[P_{ij}]$ matrix characterizes a discrete-time
Markov chain. This discrete-time Markov chain, called the {\em
embedded chain} of our continuous-time Markov chain, has $[P_{ij}]$
as its transition probability matrix. This embedded chain is
the sequence of states that our original continuous-time chain
visits where we ignore the time spent in each state during each
visit to that state. We already know the condition for the time reversibility of the embedded chain, so consider our continuous-time
chain and assume that it has been running for a long while, and
consider its reversed process going backward in time. In the
following, we show that also the reversed process spends an
exponentially distributed amount of time in each state. Moreover, we
will show that the reverse process spends an exponentially
distributed amount of time with parameter $\delta_i$ when in state
$i$, which is equal to the time spent in state $i$ by the original
process.
\begin{eqnarray*}
P\{X(t)=i, ~{\rm for}~ t\in [u-v,u] \mid  X(u)=i \} &= &\frac{P\{X(t)=i, ~{\rm for~} t\in [u-v,u] \cap X(u)=i \}}{P[X(u)=i]} \\
&=& \frac{P[X(u-v)=i]e^{-{\delta_i v}}}{P[X(u)=i]} ~=~ e^{-{\delta_i v}}.\\
\end{eqnarray*}
The last equality is explained by reminding the reader that the process is in steady-state, so the probability that the process is in
state $i$ at time $(u-v)$ is equal to the probability that the
process is in state $i$ at time $u$.

Since the continuous-time Markov chain is composed of two parts, its
embedded chain and the time spent in each state, and since we have
shown that the reversed process spends time in each state, which is
statistically the same as the original process, a condition for time
reversibility of a continuous-time Markov chain is that its embedded
chain is time reversible.

As we have learned when we discussed
reversibility of stationary discrete-time Markov chains, the condition for
reversibility is the existence of positive  $\hat{\pi}_i$ for all
states $i$ that sum up to unity and satisfy the detailed balance
equations:
\begin{equation}
\label{rev1} \hat{\pi}_i P_{ij} = \hat{\pi}_j P_{ji} ~{\rm
for~all~adjacent~} i,j.
\end{equation}
Recall that this condition is necessary and sufficient for reversibility and that if
such a solution exists, it is the stationary probability of the process.
The equivalent condition in the case of a stationary continuous-time
Markov chain is the existence of positive  ${\pi}_i$ for all states $i$
that sum up to unity that
satisfy the detailed balance equations of a continuous-time
Markov chain, defined as:
\begin{equation} \label{rev2} \pi_i Q_{ij} = \pi_j Q_{ji}
~{\rm for~all~~adjacent~} i,j.
\end{equation}
\subsubsection*{Homework \ref{stochastic}.\arabic{homework}}
\addtocounter{homework}{1} \addtocounter{tothomework}{1} Derive (\ref{rev2}) from (\ref{rev1}).
$~~~\Box$

It is important to notice that for a birth-and-death process, its
embedded chain is time-reversible. Consider a very long period of time $L$
during that time, the number of transitions from state $i$ to state
$i+1$, denoted $T_{i,i+1}(L)$,  is equal to the number of
transitions, denoted $T_{i+1,i}(L)$, from state $i+1$ to $i$ because
every transition from $i$ to $i+1$ must eventually follow by a
transition from $i+1$ to $i$. Actually, there may be a last
transition from $i$ to $i+1$ without the corresponding return from
$i+1$ to $i$, but since we assume that $L$ is arbitrarily long, the
number of transitions is arbitrarily large, so the effect of being off by one
transition is negligible.

Therefore, for arbitrarily long $L$,
\begin{equation}
\frac{T_{i,i+1}(L)}{L} = \frac{T_{i+1,i}(L)}{L}.
\end{equation}
Since for a birth-and-death process $Q_{ij}=0$ for $\mid i - j \mid
> 1$ and for $i=j$, and since for arbitrarily long $L$, we have
\begin{equation}
\pi_i Q_{i, i+1} = \frac{T_{i,i+1}(L)}{L} = \frac{T_{i+1,i}(L)}{L} =
\pi_{i+1} Q_{i+1,i},
\end{equation}
so our birth-and-death process is time reversible. This is an important result for the present context because many of the
queueing models discussed in this book involve birth-and-death processes.
\subsection{Renewal Process}
An informal way to describe a {\it renewal process} is a generalization of the Poisson process where the inter-arrival times are not necessarily exponentially distributed but are still positive IID random variables with finite mean. Because of this generalization, a renewal process does not have to be memoryless. There are discrete-time and continuous-time renewal processes. In a discrete-time renewal process, the inter-arrival times take only positive integer values, while in a continuous-time renewal process, the inter-arrival times are positive real-valued.

{\bf Note:}
The definition of a renewal process is sometimes extended to cases that allow a zero value for the inter-arrival times to occur with positive probability. This implies that more than one arrival can occur at the same point in time.

\subsubsection*{Homework \ref{stochastic}.\arabic{homework}}
\addtocounter{homework}{1} \addtocounter{tothomework}{1} Provide examples of renewal processes that are not memoryless.
\subsubsection*{Guide} Since the exponential random variable is the only continuous random variable with the memoryless property, take any other positive continuous IID random variables for the inter-arrival times, and you will have a renewal process that is not memoryless.
$~~~\Box$

\subsubsection*{Homework \ref{stochastic}.\arabic{homework}}
\label{ipphw}
\addtocounter{homework}{1} \addtocounter{tothomework}{1} Show that the IPP is a renewal process. Derive the mean and distribution of its inter-arrival times.
\subsubsection*{Guide} The IPP is a renewal process because the time from an arrival (called arrival A) until the next arrival (called arrival B) is independent of the evolution of the process before event A. Therefore, the inter-arrival times are independent. In fact, because of the memoryless property of the exponential distribution of the inter-arrival times during the active mode, the distribution of the inter-arrival times is equal to the distribution of the time from any moment the system is in the active mode until the next arrival. Therefore the inter-arrival times are IID. Because of the underlying structure of the IPP, the inter-arrival times are also continuous and positive.

Let $X$ be a random variable representing the inter-arrival time of the IPP. As mentioned above, because of the memoryless property of exponential distribution, $X$ is equal to the time until the next arrival, given that now the system is in active mode. Then, $X$ can be described by the following recursive equation.

\begin{equation} \label{XbyYZ} X=Y + \alpha (Z +X), \end{equation}

where $Y$ is an exponentially distributed random variable with parameter $\lambda+\delta_1$ representing the time from the last arrival until the next event that can be either an arrival or change of mode from active to nonactive, $Z$ is an exponentially distributed random variable with parameter $\delta_0$ representing the duration of a nonactive period, $X$ on the right-hand side represents the time until the next arrival from the moment the system is again in the active mode, and $\alpha$ is the probability that the next event after an arrival is a change of mode from active to nonactive, given by

\begin{equation} \label{alpha1} \alpha = \frac{\delta_1}{\lambda+\delta_1}. \end{equation}
Isolating $X$ in (\ref{XbyYZ}) yields
\begin{equation} \label{Xbyalpha} X = \frac{Y + \alpha Z}{1-\alpha}. \end{equation}
Substituting (\ref{alpha1}) in (\ref{Xbyalpha}, we obtain
\begin{equation}
X =  \frac{Y+\frac{\delta_1}{\lambda+\delta_1} Z}{1-\frac{\delta_1}{\lambda+\delta_1}} =   \frac{(\lambda+\delta_1)Y+\delta_1 Z}{ \lambda}.
\end{equation}
Thus,
\begin{equation}
\label{finalX}
X = \frac{\lambda+\delta_1}{\lambda}Y+\frac{\delta_1}{ \lambda}Z.
\end{equation}
Taking expectations on both sides of (\ref{finalX}) gives
\begin{equation} \label{EofX}  E[X] = \frac{\lambda+\delta_1}{\lambda}E[Y]+\frac{\delta_1}{ \lambda}E[Z]. \end{equation}
By the definitions of the random variables $Y$ and $Z$, we have
\begin{equation} \label{EofY} E[Y] = \frac{1}{\lambda + \delta_1}, \end{equation}
and
\begin{equation} \label{EofZ} E[Z] = \frac{1}{\delta_0}. \end{equation}

Substituting (\ref{EofY}) and (\ref{EofZ}) in (\ref{EofX}), we obtain

\begin{equation}
\label{ex1}
E[X] =  \frac{\lambda+\delta_1}{\lambda}\left(\frac{1}{\lambda + \delta_1}\right)+\frac{\delta_1}{ \lambda}\left(\frac{1}{\delta_0}\right)
=  \frac{\delta_0 + \delta_1}{ \lambda\delta_0}.
\end{equation}
Let us now obtain $E[X]$ in a different way.

Consider the two-state continuous-time Markov chain of the mode process. This process alternates between state 0 (non-active mode) and state 1 (active mode). Let $\pi_i$ be the probability of being in state $i$, $i=0,1$. The two-by-two ${\bf Q}$ matrix is $q_{00}=-\delta_0$, $q_{01}=\delta_0$, $q_{10}=\delta_1$, $q_{11}=-\delta_1$. That is,
$${\bf Q} =
\begin{bmatrix}
    -\delta_0 & \delta_0  \\
    \delta_1 & -\delta_1
    \end{bmatrix}.
$$

This gives rise to the steady state equations ${\bf \Pi Q}=0$,
which are two dependent equations, so from one of them we have

$$\pi_0 \delta_0 = \pi_1 \delta_1.$$

We also have the normalizing equation

$$\pi_0 + \pi_1 = 1.$$

Solving these equations, we obtain

$$ \pi_0 = \frac{\delta_1}{\delta_0 + \delta_1}$$

and

$$\pi_1 = \frac{\delta_0}{\delta_0 + \delta_1}.$$

The arrival rate of the IPP is given by

$$\lambda_{IPP} = \lambda \pi_1 + 0 \times \pi_0 = \frac{\lambda \delta_0}{\delta_0 + \delta_1}.$$

Therefore,  the mean inter-arrival time $E[X]$ is obtained by

$$E[X] = \frac{1}{\lambda_{IPP}} = \frac{\delta_0 + \delta_1}{\lambda \delta_0}.$$
This is consistent with the result obtained in (\ref{ex1}).

Finally, obtain the distribution of $X$ using (\ref{finalX}) by deriving the convolution of two random variables.
$~~~\Box$

\newpage
\section{General Queueing and Teletraffic Concepts}
\label{general}

\setcounter{homework}{1} 

In general, a queueing system may be characterized by a complex input
process, service time distribution, number of servers (or channels),
buffer size (or waiting room) and queue discipline. In practice,
such queueing processes and disciplines are often not amenable to
analysis. Nevertheless, insight can be often gained using simpler
queueing models. Modeling simplification is often made when the aim
is to analyze a complex queueing system or network, such as the
Internet, where
 packets on
their ways to their destinations arrive at a router where they are
stored and then forwarded according to addresses in their
headers. One of the most fundamental elements in this process is the
single-server queue (SSQ). One of the aims of telecommunications research
is to explain traffic and management processes and their effect on
queueing performance. In this section, we briefly cover basic
queueing theory concepts. We shall bypass mathematically rigorous
proofs and rely instead on simpler intuitive explanations.

\subsection{Notation}

A commonly used shorthand notation, called Kendall notation \cite{kendall53}, for
a single queue model describes the arrival process, service
distribution, the number of servers and the buffer size (waiting
room) as follows.

\begin{center}
\{arrival process\}/\{service distribution\}/\{number of
servers\}/\{buffer size\}-\{queue discipline\}
\end{center}

Commonly used characters for the first two positions in this
shorthand notation are: D (Deterministic), M (Markovian - Poisson
for the arrival process or Exponential for the service time distribution required by each customer), G (General), GI (General and independent),
and Geom (Geometric). The fourth position is used for the number of buffer places, including the buffer spaces available at the servers.
This means that if there are $k$ servers and no additional waiting room is available, then
$k$ will be written in the fourth position. The fourth position is not
used if the waiting room is unlimited. The fifth
position is used for the queue discipline.
Namely, the order in which the customers are served
in the queue.  For example: First In
First Out (FIFO), Last In First Out (LIFO), Processor Sharing (PS)
where all the customers in the queue obtain service simultaneously, and random service order
(random). The fifth position is not used for the case of the FIFO
queue discipline. Notice that the dash notation ``-'' before the fifth
position is used to designate the fifth position. This ``-'' designation avoids ambiguity
in case the fourth position is missing.

For example, D/D/1 denotes an SSQ where both the inter-arrival times and the service times are deterministic. This means that the inter-arrival times are all equal to each other, and the service times are all equal to each other. The notation M/M/1 represents an SSQ with a Poisson arrival process and independent and exponentially distributed service times with an infinite buffer and FIFO service order. Next, GI/M/1 denotes a generalization of M/M/1 where the arrival process is a renewal process and not necessarily Poisson. Then, G/M/1 is a further generalization of GI/M/1 where the arrival process is not restricted to be a renewal process, i.e., the inter-arrival times do not have to be IID, as no independence is assumed. However, in many cases in the literature and this book, G is used as GI. This is especially common in the second position, namely,  for the service times. For example, M/G/1 always represents an SSQ queue with Poisson arrivals and independent and generally distributed service times that are also independent of the arrival process. Another example is the use of the GI/G/1 notation, e.g., in \cite{kingman62,kingman66}, to represent a queueing model where the arrival process is a renewal process, namely, the inter-arrival times are generally distributed and independent, and the service times are also independent and are independent of the arrival process. Even the notation G/G/1 is often used for the case where both the service times and the inter-arrival times are IID (generally distributed) and independent of each other (e.g., \cite{cohen82,Shanthikumar07}).

Given the inconsistencies and ambiguities in using G, it is important to specify the queueing model under consideration clearly and not rely only on notation such as G/G/1. In this book, G/G/1  represents a single server queue where both the service times and the inter-arrival times are generally distributed, and no independence of any kind is assumed.

M/G/$k$/$k$ denotes a k-server queue with Poisson arrivals and generally distributed IID service times that are also independent of the arrival process without additional waiting room except at the servers. Then, the $k$-server queue M/G/$k$/$N$ (with $N \geq k$) also with Poisson arrivals and generally distributed IID service times is a generalization of M/G/$k$/$k$ where a waiting room that can accommodate up to $N-k$ customers.
M/G/1-PS denotes a single server processor sharing queue with Poisson arrivals and generally distributed IID customer service time requirements. Notice that in an M/G/$1$-PS queue, although the service time of a customer/packet starts immediately upon arrival, it may continue for a longer time than its service requirement because the server capacity is always shared among all customers/packets in the system. Accordingly, M/G/1/$N$-PS is a single server processor sharing queue with Poisson arrivals and a generally distributed IID service time requirements that apply to all the customers, and the waiting room can accommodate up to $N$ customers, including the customer in service. An important special case of the latter is the M/M/1/$N$-PS queue, where the service time distribution is exponential.

\subsection{Utilization}
\label{utiliz}

An important measure for queueing systems performance is the
utilization, denoted $\hat{U}$. It is the proportion of time that a server
is busy on average. In many systems, the server is paid for its time
regardless if it is busy or not. Normally, the time that
transmission capacity is not used is the time during which money is
spent, but no revenue is earned. It is, therefore, important to
design systems that will maintain high utilization.

If you have two identical servers and one is busy 0.4 of the time and the
other 0.6. Then, the utilization is 0.5. We always have that $0\leq \hat{U}
\leq 1$. If we consider an M/M/$\infty$ queue (Poisson arrivals,
exponentially distributed service times, and infinite servers) and
the arrival rate is finite, then the utilization is zero because the mean number of busy servers is finite, and the mean number of idle servers is infinite.

Consider a G/G/1 queue (that is, an SSQ with
an arbitrary arrival process and arbitrary service time distribution,
with an infinite buffer). Let $S$ be a random variable representing the service time and
let $E[S]=1/\mu$, i.e., $\mu$ denotes the service rate.
Further, let $\lambda$ be the arrival rate. Assume that $\mu > \lambda$ so that
the queue is {\it stable}, namely, that it will not keep growing
forever, and that whenever it is busy, eventually it will reach
the state where the system is empty. For a stable G/G/1 queue, we
have that $\hat{U}=\lambda/\mu$. To show the latter let $L$ be a
very long period of time. The average number of customers arrived within the time period $L$ is: $\lambda L$. The
average number of customers that have been served
during time period $L$ is equal to $\mu \hat{U}L$. Note that during $L$, customers are being served only during $\hat{U}$ proportion of $L$ when the system is not empty. Since $L$ can be made arbitrarily long
and the queue is stable, these two values can be considered equal. Thus, $\mu
\hat{U}L=\lambda L $. Hence,

\begin{equation}
\label{utilgg1}
\hat{U}=\frac{\lambda}{\mu}.
\end{equation}

Another Explanation for (\ref{utilgg1}) is the following.
Again, consider an arbitrarily long period of time $L$. The assumption of $L$ being arbitrarily long is again required to ensure that the edge effects are negligible - namely, effects of a customer being in service at the beginning and end of the time period $L$.
Again, there are on average $\lambda L $ arrivals during $L$. These arrivals require on average a total service time of $$ \lambda L \times \frac{1}{\mu}.$$
Therefore, the proportion of time the server is busy ($\hat{U}$) is obtained by $$\hat{U} = \frac{\lambda L \times \frac{1}{\mu}}{L} = \frac{\lambda}{\mu}.$$

Often, we are interested in the distribution of the number of
customers (or jobs or packets) in the system. Consider a G/G/1 queue
and let $\pi_n$ be the probability that there are $n$ customers in the system.
Having the utilization, we can readily obtain $\pi_0$ the probability
that the G/G/1 queue is empty. Specifically,
\begin{equation}
\label{gg1}
\pi_0=1-\hat{U}=1-\frac{\lambda}{\mu}.
\end{equation}

In the case of a multi-server queue, which in the most general case is denoted by G/G/$k$/$k+n$, the
utilization will be defined as the overall average
utilization of the individual servers. That is, each server will
have its own utilization defined by the proportion of time it is busy,
and the utilization of the entire multi-server system will be the
average of the individual server utilizations.

\subsection{Little's Formula}
\label{littles}

Another important and simple queueing theory result that applies
to a stable G/G/1 queue (and to other systems) is known as {\em Little's
formula} \cite{Little61,Stidham72,Stidham74}. It has two forms. The first form is:
\begin{equation}
\label{Litt1} E[Q]=\lambda E[D]
\end{equation}
where $E[Q]$ and $E[D]$ represent the stationary mean queue-size including
the customer in service and the mean delay (system waiting time)
of a customer from the moment it arrives until its service is
completed, respectively. In the remainder of this book, when we use terms such as
{\it mean queue-size} and {\it mean delay},
we refer to their values in steady-state,
i.e., stationary mean queue size and delay, respectively.

The second form is:
\begin{equation}
\label{Litt2} E[N_Q]=\lambda E[W_Q]
\end{equation}
where $E[N_Q]$ and $E[W_Q]$ represent the mean number of customers in the queue in steady-state
excluding the customer in service and the mean delay of a customer, in steady-state,
from the moment it arrives until its service commences (waiting
time in the queue), respectively.

An intuitive (non-rigorous) way to explain Eq. (\ref{Litt1}) is by
considering a customer that just left the system (completed
service). This customer sees behind his/her back on average $E[Q]$
customers. Who are these customers? They are the customers that had
been arriving during the time that our customer was in the system.
Their average number is $\lambda E[D]$.

Another explanation is based on the amusement park analogy of \cite{Bertsekas02}. Consider an amusement park where
customers arrive at a rate of $\lambda$ per time unit. Assume that the park is in a stationary condition which implies that it is open 24 hours a day and that the arrival rate $\lambda$ does not change in time. After arriving at the park, a customer spends time at various sites and then leaves. The mean time a customer spends in the park is represented by $E[D]$. Assume that the park charges every customer
one dollar per unit of time spent in the park. The mean queue size $E[Q]$, in this case, represents the mean number of customers in the park in steady state.

Under these assumptions, the rate at which the park earns its income in steady state is $E[Q]$  per unit of time because each of the $E[Q]$ customers pays one dollar per unit of time. Now let $L$ be an arbitrarily long period of time. The mean number of customers that arrive during $L$ is $\lambda L$. Because $L$ is arbitrarily long, the mean time that a customer spends in the park $E[D]$ is negligible relative to $L$, so we can assume that all the customers that arrive during $L$ also left during $L$. Since a customer on average pays $E[D]$ dollars for its visit in the park, and there are on average $\lambda L$ customers visiting during $L$, the total income earned by the park on average during $L$ is $\lambda L E[D]$. Therefore, the rate at which the park earns its income in steady state per unit of time  is  $$\frac{\lambda L E[D]}{L}=\lambda E[D].$$
We also know that the rate that the park earns its income is also equal to $E[Q]$. Therefore, $\lambda E[D]=E[Q]$.

For a graphical proof of Little's formula for the case of G/G/1
queue see \cite{BG92}. The arguments there may be summarized as
follows. Consider a stable G/G/1 queue that starts at time $t=0$
with an empty queue. Let $A(t)$ be the number of arrivals up to
time $t$, and let $D(t)$ be the number of departures up to time
$t$. The queue size (number in the system) at time $t$ is denoted
$Q(t)$ and is given by $Q(t)=A(t)-D(t)$, $t\geq 0$. Let $L$ be an arbitrarily long period of time. Then, the mean queue-size $E[Q]$
is given by
\begin{equation} E[Q]=\frac{1}{L}\int_0^L Q(t) dt. \end{equation}
Also, notice that
\begin{equation} \int_0^L Q(t) dt = \sum_{i=1}^{A(L)} D_i \end{equation}
where $D_i$ is the time spent in the system by the $i$th customer.
(Notice that since $L$ is arbitrarily large, there have been
arbitrarily large number of events during $[0,L]$ where our stable
G/G/1 queue became empty, so $A(L)=D(L)$.) Therefore,
\begin{equation}
\frac{1}{L}\int_0^L Q(t) dt = \frac{1}{L}\sum_{i=1}^{A(L)} D_i
\end{equation} and realizing that
\begin{equation}
\lambda=\frac{A(L)}{L},
\end{equation}
 and
\begin{equation} E[D]=\frac{1}{A(L)}\sum_{i=1}^{A(L)} D_i, \end{equation}
we obtain
\begin{equation}
E[Q]= \frac{1}{L}\int_0^L Q(t) dt =
\frac{A(L)}{L}\frac{1}{A(L)}\sum_{i=1}^{A(L)} D_i=\lambda E[D].
\end{equation}

Little's formula applies to many systems. Its applicability is not
limited to single-server queues, or single-queue systems, or
systems with infinite buffers. However, the system must be in steady state for Little's formula to apply.

Little's formula is applicable to almost any queueing system in
steady state. The system may consist of more than one queue,
more than one server, the order does not need to be FIFO, the
arrivals do not need to follow a Poisson process and the service time
does not need to be exponential.

Interestingly, the result $\hat{U}=\lambda/\mu$ for a G/G/1 queue can
also be obtained using Little's formula. Let us consider a system
to be just the server with space for at most one customer served at the server
 (excluding the waiting room in the buffer).
 The mean time a customer spends in this server
system is $1/\mu$ because this is the mean service time.
The mean arrival rate into that system must be equal to $\lambda$
because all the customers that arrive at the queue eventually
arrive at the server - nothing is lost. Let us now consider the
 number of customers at the server, denoted $N_s$. Clearly,
$N_s$ can only take the values zero or one because no more than
one customer can be at the server at any point in time. We also
know that the steady-state probability $P(N_s=0)$ is equal to
$\pi_0$. Therefore, $$E[N_s]=0\pi_0+1(1-\pi_0)=1-\pi_0=\hat{U}.$$ By Little's
formula, we have $$E[N_s]=\lambda(1/\mu),$$ so $$\hat{U}=\frac{\lambda}{\mu}.$$

Conventional notations in queueing theory for a k-server queue are $$A = \frac{\lambda}{\mu}$$ and $$\rho = \frac{A}{k}.$$
Thus, for a G/G/1 queue

\begin{equation}
\label{ENsGG1} E[N_s]=\frac{\lambda}{\mu}=\hat{U}=\rho.\end{equation}

To obtain (\ref{Litt2}) from (\ref{Litt1}), notice that
\begin{equation}
\label{EQEQW} E[Q]=E[N_Q]+E[N_s]=E[N_Q]+\frac{\lambda}{\mu}
\end{equation}
and
\begin{equation}
\label{EDEW} E[D]=E[W_Q]+1/\mu.
\end{equation}
Substituting (\ref{EQEQW}) and (\ref{EDEW}) in (\ref{Litt1}),
(\ref{Litt2}) follows.

Another interesting application of Little's formula relates the
blocking probability $P_b$ of a G/G/1/$k$ queue with its server utilization \cite{Hubner90,roberts96}. Again, consider
the server as an independent system. Since the mean number of
customers in this server system is $\hat{U}$, and the arrival rate into this server
system is $(1-P_b)\lambda$, we obtain by Little's formula:
\begin{equation} \hat{U}=(1-P_b)\lambda \mu^{-1}, \end{equation} where $\mu^{-1}$ is the mean service
time. Having $\rho=\lambda/ \mu$, we obtain
\begin{equation} P_b=1-\frac{\hat{U}}{\rho}. \end{equation}

Now observe that as $k\rightarrow \infty$ and the  G/G/1/$k$ queue  approaches a G/G/1 queue, where ${\hat{U}}={\rho}$, we obtain $P_b=1-\frac{\hat{U}}{\rho}=1-1=0$, as expected.

\subsection{Delay and Loss Systems}

We use the term {\it delay system} to describe  a queueing system with an infinite buffer, such as G/G/1 or G/G/2, where an arriving customer that finds all servers busy will wait in the queue until it is served, or system such as G/G/1-PS where the service time and overall delay is increased with increasing congestion. Increasing delay, which adversely affects users' QoS, is used by the system to save resources.

Another measure used by the system to save resources is rejecting customers at times of congestion. This is done in a system such as M/G/$k$/$k$. An arriving customer that finds all servers busy is blocked and cleared from the system. Systems, where customers are blocked and lost during congestions, are called {\it loss systems}. In loss systems, the term {\it holding time}, is often used for the time spent in the system by a call if it is served by a dedicated server. In such systems, calls that are not blocked immediately start their service upon their arrivals. Let $h$ be the mean holding time. As our holding time definition is associated with loss systems without waiting time, $h$ is also the mean service time, namely, $$h=E[S]=\frac{1}{\mu}.$$ The term holding time has often been used in relation to telephony systems, referring to the average time the phone call holds a circuit.

Systems like M/G/1/$k$ and M/G/1/$k$-PS are called {\it delay-loss systems}, such systems are characterized by both adverse effects of delay and loss taking place. If a customer arrives and all servers are busy, but the buffer (waiting room) is not full, it may suffer some additional delay, but if the buffer is full, the customer will be rejected and lost.

Finally, there are systems that are neither delay nor loss systems. These are systems where the number of servers is infinite. In such systems, an arriving customer will always receive service at the full-service rate upon arrival and will never be delayed or blocked.

\subsection{Traffic}
The concept of {amount of traffic} is very important for the telecommunications industry. It is important for the network designer that networks are properly designed and dimensioned such that they have sufficient capacity to route the traffic to the destinations and meet the required customers' quality of service. As traffic on its way from source to destination may traverse several domains (parts of the Internet), each of which is controlled and managed by different operators, the operators must negotiate and agree on the level of service they provide and the amount of traffic that is carried between them. To measure the quantity of traffic, we consider the arrival rate of customer calls as well as the amount of network resources the calls require. Clearly, a short phone call will require less resources than a long one. In particular, traffic is measured in units called {\it erlangs} and it is defined as $\lambda/\mu$, namely, if a quantity of traffic $A$ [erlangs] are offered to a system, then $A$ is given by $$A=\frac{\lambda}{\mu}.$$ In other words, the quantity of traffic in erlangs is the product of the call arrival rate ($\lambda$) and the mean time that a single server will serve a call ($1/\mu$), i.e., the mean time that calls occupy a server. Note that this is the second term mentioned in this book named after the Danish mathematician Agner Krarup Erlang. The first was the Erlang distribution.

Intuitively, the traffic is $\lambda$ times the mean service time. It represents the number of arrivals per mean service time. One erlang represents the traffic load of one arrival, on average, per mean service time. This means that the traffic load of one erlang, if admitted to the system, will require a resource of one server on average. For example, this could be provided by service from one server continuously busy, or from two servers each
of which is busy 50\% of the time. Accordingly, if the traffic is $A$, then $A$ is the mean number of servers required in steady state by this traffic load.

The total traffic offered to a system is normally generated by many users. If we have $N$ users and the $i$th user generates $A_i$ erlangs, the total traffic generated by the $N$ users is $$A=\sum_{i=1}^N A_i~~~{\rm [erlang]}.$$

\subsubsection*{Homework \ref{general}.\arabic{homework}}
\addtocounter{homework}{1} \addtocounter{tothomework}{1} Consider a wireless system that provides channels each of which can serve one phone call. There are 100 users making phone calls. Each user makes on average one phone call per hour and the average duration of a phone call is three minutes. What is the total traffic in erlangs that the 100 users generate? \subsubsection*{Guide} In such questions, it is important first to choose a consistent time unit. Here it is convenient to choose minutes. Accordingly, the arrival rate of each user is $\lambda_i = 1/60 {\rm ~~calls~ per~ minute}. $ The mean call duration (or holding time) is 3 minutes, so $A_i = 3/60 = 1/20~~{\rm [erlang]}$,
 and the total traffic is $100 \times (1/20) = 100/20 = 5$ [erlang]. $~~~\Box$

 The quantity of traffic measured in erlangs is also called {\it traffic intensity}. See Section \ref{offered_carried} for more details. Then, another related concept is {\it traffic volume}. Traffic volume is measured in units of erlang-hour (or erlang-minute, or call-hour or call-minute, etc.) and it is a measure of the traffic processed by a facility during a given period of time.  Traffic volume is the product of the traffic intensity and the given time period, namely,

 \begin{center} {\it Traffic Volume = Traffic Intensity $\times$ Time Period.} \end{center}

Let $A, \lambda, \mu, {\rm and}~ h$ be the traffic intensity, arrival rate, service rate, and holding time, respectively. Then,
 $$A=\frac{\lambda}{\mu}=\lambda h.$$

Let $a_T$ be the number of arrivals during the given period of time $T$. Then, the relevant arrival rate can be estimated by

$$\lambda=\frac{a_T}{T}.$$

Let $V_T$ be the traffic volume during $T$. Then,
$$V_T = AT = \frac{a_T h T}{T} = a_T h.  $$

 This leads to a different definition of traffic volume, namely, the product of the number of calls during $T$ and the mean holding time, and this explains the traffic volume units of call-hour or call-minute. The latter result for $V_T$ can be illustrated by the following example.

 Consider a period of time of three hours and during this period of time,  120 calls have arrived and their average holding time is three minutes, then the traffic volume is $3\times 120 =360$ call-minute, or 360 erlang-minute, or $360/60=6$ erlang-hour. Then, the traffic intensity in erlangs is obtained by dividing the traffic volume in erlang-hour by the number of hours in the given period of time. In our case, the traffic intensity is $6/3=2$ [erlangs].

\subsection{Offered and Carried Traffic}
\label{offered_carried}

There are two traffic-related concepts called: {\it offered traffic} and {\it carried traffic}.
The offered traffic is defined as the mean number of arrivals (of customers,
calls, or packets) per mean
service time. Accordingly, it is equal to the ratio $\lambda/\mu$ which is identical to the definition of traffic discussed in the previous section.
It is common to use the notation $A$ that we used for traffic, for the offered traffic as well.

The carried traffic is the average number of calls that are being served in the system. In a system such as M/G/$k$/$k$, every call is served by one server, so the carried traffic is also the mean number of servers required to serve the customers,
calls, or packets that are admitted to the system. It is also equal to the proportion of the offered traffic that is admitted to the system.
Both offered and carried traffic are measured in {\em erlangs}. The relationship between offered and carried traffic is given by

$$ [{\rm offered ~traffic}](1-P_b) = [{\rm carried ~traffic}]. $$

Therefore, in systems where blocking does not occur, such as delay systems and systems where the number of servers is infinite, we have that $P_b=0$, the offered and carried traffic are equal to each other.

As mentioned above, another term to describe traffic that is often used is the {\em traffic intensity} \cite{ITUT98}. Then, {\em offered traffic intensity} and {\em carried traffic intensity} are synonymous with offered traffic and carried traffic, respectively. In cases of infinite buffer/capacity systems, such as M/M/1 and M/M/$\infty$, the term traffic intensity is used to describe both offered traffic and carried traffic. Accordingly, $\rho$ is called traffic intensity in the M/M/1 context and $A$ is the traffic intensity in an M/M/$\infty$ system. Others use the term traffic intensity in a multiserver system for the offered load per server \cite{Whitt02}. To avoid confusion, we will only use the term traffic intensity in the context of a single server queue with an infinite buffer, in which case, traffic intensity is always equal to $\rho=\lambda/\mu$.

\subsection{Work Conservation}

Another important concept in queuing theory is the concept of {\em
work conservation}. A queuing system
is said to be {\em work conservative} if a server is never idle whenever there
is still work to be done.
For
example, M/M/1 and M/M/1/$7$ are work conservative. However, a stable M/M/$3$
is not work conservative because a server
can be idle  while there are customers served by
other servers.

\subsection{PASTA}

Many of the queueing models we consider in this book involve Poisson arrival processes. The PASTA property discussed in the previous chapter is important for analyzes
and simulations of such queueing models. Let us further explain and
prove this important property.

The PASTA property implies that arriving customers in a steady  state
will find the number of customers in the system obeying its steady-state distribution.
In other words, the statistical characteristics
(e.g., mean, variance, distribution) of the number of customers in
the system observed by arrivals is the same as those observed by
an independent Poisson inspector.

In addition to the assumption of Poisson arrivals, for PASTA to be valid
we also need the condition that arrivals after time $t$ are
independent of the queue size at time $t$, $Q(t)$. For example, if
we have a single-server queue with Poisson arrivals and the
service times have the property that the service of a customer must
always terminate before the next arrival, then the arrivals always
see an empty queue, and, of course, an independent arrival does not.
However, in all the queueing systems that we study, this condition holds because normally the server cannot predict the exact time of the next arrival because
of the pure chance nature of the Poisson process.

To prove PASTA we consider the limit
$$A_k(t) = \lim_{\Delta t \rightarrow 0} P [ Q(t)=k \mid {\rm ~an ~arrival~ occurs~
within~} (t,t + \Delta t)]. $$

Using Bayes' theorem and the condition that arrivals after time $t$
are independent of $Q(t)$, we obtain that
\begin{equation}
\label{pastaproof} A_k(t) = P[Q(t)=k].
\end{equation}
Then, by taking the limit of both sides of (\ref{pastaproof}), we
complete the proof that the queue size seen by arrivals is
statistically identical to the queue size seen by an independent
observer.$~~~\Box$

\subsubsection*{Homework \ref{general}.\arabic{homework}}
\addtocounter{homework}{1} \addtocounter{tothomework}{1} Prove Eq.\@
(\ref{pastaproof}). $~~~\Box$

\subsubsection*{Homework \ref{general}.\arabic{homework}~\cite{BG92}}
\addtocounter{homework}{1} \addtocounter{tothomework}{1}

Consider a queueing system with $k$ servers and the total waiting room for customers in the system (including a waiting room in the queue and at the server) is $N$, such that $N > k$.
This system is always full. That is, it always has $N$ customers in the system. In the beginning, there are $N$ customers in the system. When a customer leaves the system immediately another customer arrives. You are given that the mean service time is $1/\mu$.

Find the mean time from the moment a customer arrives until it leaves the system as a function of $N, k$ and $\mu$.

\subsubsection*{Solution}

Let $\lambda$ be the arrival rate into the system (which is also equal to the departure rate).
We use the term {\it System S} to denote the system composed of only the servers without the waiting room in the queue outside the servers, and we use the term {\it System E}  to denote the system composed of both the servers and the waiting room in the queue outside the servers. Therefore, $\lambda$ is the arrival rate into {\it System S} and into {\it System E}, and it is also equal to the departure rate.

Then, by applying Little's formula to {\it System S}, we have:

$$\lambda \left (\frac {1}{\mu} \right )= k$$

so

$$\lambda = k \mu.$$

Let $E[D]$ be the mean time from the moment a customer arrives until it leaves the system.

Then, by applying Little's formula to {\it System E}, we have:

$$\lambda E[D] = N.$$

Substituting $\lambda$, we obtain

 $$E[D] = \frac{N}{k\mu}.$$

 Interestingly, observe that we have $k$ servers working all the time. Each server serves customers at a rate of $\mu$ customers per unit time. Therefore, the service rate of the entire system is $k\mu$. This gives inter-departure time of 1/$k\mu$. When a new customer joins the system, there are $N$ customers in the system. Therefore, the result of $N$ times 1/$k\mu$ that we obtained is intuitively justifiable.

 $~~~\Box $

\subsubsection*{Homework \ref{general}.\arabic{homework}}
\addtocounter{homework}{1} \addtocounter{tothomework}{1}

Consider an M/G/1/2 queue with an arrival rate of one packet per millisecond  and the mean service time is one millisecond. It is also given that $E[N_Q]$ the mean number of customers in the queue (not including a customer in the service)  is equal to 0.2. What is the blocking probability? A numerical answer is required.

\subsubsection*{Solution}

Let $\pi_0$, $\pi_1$, and $\pi_2$ be the steady-state probabilities that the total number of customers in this  M/G/1/2 queueing system (including those in service and those waiting in the queue) is equal to 0, 1, and 2, respectively. Then, $N_Q=1$ when there are 2 customers in the system, and $N_Q=0$ when there are either 0 or 1 customers in the system. Therefore:

$$E[N_Q] = \pi_2 (1) +(\pi_0+\pi_1)0=\pi_2. $$

It is also given that $E[N_Q] = 0.2.$ Therefore, $\pi_2=0.2$, and since the arrival process is Poisson, by PASTA it is also the blocking probability.
$~~~\Box$

\subsection{Bit-rate Versus Service Rate}

In many telecommunications design problems, we know the bit rate of an output link and we can estimate the load on the system in terms of the arrival rate of items of interest that require service. Examples of such items include packets, messages, jobs, or calls. To apply queueing theory for evaluation of quality of service measures, it is necessary, as an intermediate step, to calculate the {\it service rate} of the system which is the number of such items that the system can serve per unit time. For example, if we know that the average message size is 20 MBytes [MB] and the bit-rate of the output link of our system (the capacity in bit-rate of the output link) is 8 Gigabits per second [Gb/s] and we are interested in the service rate of the link in terms of messages per second, we first calculate the message size in bits which is $20 \times 8 =160 $ Mega bits [Mb] or $160 \times 10^{6}$ bits. Then, the service rate is $$ \frac{8 \times 10^{9}}{160 \times 10^{6}} = 50~ {\rm messages~ per ~ second}. $$

\subsubsection*{Homework \ref{general}.\arabic{homework}}
\addtocounter{homework}{1} \addtocounter{tothomework}{1}

Consider a single server queue with a limited buffer. The arrival rate is one job per second  and the mean job size is 0.1 Gigabyte. The server serves the jobs according to a FIFO discipline at a bit rate of  one Gb/s (Gigabit per second). Jobs that arrive and find the buffer full are blocked and cleared of the system. The proportion of jobs blocked is 2\%.  Find the bit rate in [Gb/s] of the arriving jobs, the service rate in jobs per second, the queue utilization, and throughput rates in [Gb/s] and in jobs per second. Note that the service rate is the capacity of the server to render service and the throughput is the actual output in either jobs per second or Gb per second.

\subsubsection*{Solution}

Since the mean job size is 0.1 Gbytes = 0.8 Gbits, and the arrival rate is one job per second, then the bit-rate of the arriving jobs in [Gb/s] is 0.8 [Gb/s].

The service rate is 1/0.8 = 1.25 jobs per seconds.

The queue utilization is obtained by using Little's formula considering the server as the system and only the traffic that enters the buffer and reaches the server is the input to this server system. Therefore, the queue utilization is the mean queue size in this system (of only the server) and it is obtained by 1(1-0.02)/1.25 = 0.784.

The throughput is obtained by the product of the arrival rate and (1-0.02), so it is equal to 0.98 jobs per second or 0.8(1-0.02) = 0.784 [Gb/s].
$~~~\Box$

\subsection{Queueing Models}

In this book, we discuss various queueing models that are
amenable to analysis. The analysis is simplest for D/D/ type queues
where the inter-arrival and service times are deterministic (fixed
values). They will be discussed in Chapter 5. Afterward, we
will consider the so-called Markovian queues. These queues are
characterized by the Poisson arrival process, independent
exponential service times, and independence between the arrival
process and the service times. They are denoted by M in the first
two positions (i.e., ${\rm M/M/}\cdot/\cdot$). Because of the
memoryless property of Markovian queues, these queues are amenable
to analysis. In fact, they are all continuous-time Markov chains
with the state being the {\em queue-size} defined as the number in
the system $n$, and the time between state transitions is
exponential. The reason that these time periods are exponential is
that at any point in time, the remaining time until the next
arrival, or the next service completion, is a competition between
various exponential random variables.

\newpage
\section{Simulations}
\label{simulations}

\setcounter{homework}{1} 

 In many
cases, analytical solutions are not available, so simulations
are used to estimate performance measures.  Simulations are also
used to evaluate the accuracy of analytical approximations.

\subsection{Confidence Intervals}
\label{conf}

Regardless of how long we run a simulation involving random processes, we will never obtain the
the exact mathematical result of a steady-state measure we are
interested in. To assess the error of our simulation, we begin by
running a certain number, say $n$, of simulation experiments and
obtain $n$ observed values, denoted $a_1, a_2, ~ \ldots, ~a_n$, of the
measure of interest.

Let $\bar{a}$ be the observed mean and $\sigma_a^2$ the observed
variance of these $n$ observations. Their values are given by

\begin{equation}
\label{mean_observ} \bar{a} = \frac{1}{n}  \sum_{i=1}^{n} a_i
\end{equation}
and
\begin{equation}
\label{var_observ} \sigma_a^2 = \frac {1}{n-1} \sum_{i=1}^{n}( a_i
- \bar{a})^2.
\end{equation}

Then, the confidence interval of $\bar{a}$, with confidence
$\alpha$, $0 \leq \alpha \leq 1$, is given by $(\bar{a} - U_r,
\bar{a} + U_r)$, where
\begin{equation}
\label{U_rad} U_r = \{t_{(1-\alpha)/2,(n-1)}\}
\frac{\sigma_a}{\sqrt{n}}
\end{equation}
where $t_{(1-\alpha)/2,(n-1)}$ is the appropriate percentage point
for Student's t-distribution with $n-1$ degrees of freedom. The
$t_{(1-\alpha)/2,(n-1)}$ values are available in standard tables.
For example: $t_{0.025,5}=2.57$ and $t_{0.025,10}=2.23$. That is,
if we are interested in 95\% confidence and we have
$n=6$ observations, we will use $t_{0.025,5}=2.57$ to obtain the
confidence interval, and if we have $n=11$ observations, we will use $t_{0.025,10}=2.23$.

Microsoft (MS) Excel$^{\rm TM}$ provides the function $TINV$ whereby $TINV(1-\alpha,n-1)$ gives the
appropriate constant based in t-distribution for  confidence
$\alpha$ and $n-1$ degrees of freedom. Then, the confidence interval of $\bar{a}$, with confidence
$\alpha$, $0 \leq \alpha \leq 1$, is given by $$(\bar{a} - U_r,
\bar{a} + U_r),$$ where $$ U_r = TINV(1-\alpha,n-1)
\frac{\sigma_a}{\sqrt{n}}. $$

Let us now consider the above-mentioned two examples of $n=6$ and $n=11$. Using MS Excel$^{\rm TM}$,  $TINV(0.05,5)=2.57$ and $TINV(0.05,10)=2.23$.
That is, if we are interested in 95\% confidence and we have $n=6$
observations, we will use $TINV(0.05,5)=2.57$  to obtain the
confidence interval, and if we have $n=11$ observations, we will use
$TINV(0.05,10)=2.23$.

More values for $t_{(1-\alpha)/2,(n-1)}$, or equivalently, $TINV(1-\alpha,n-1)$,  for 95\% confidence ($\alpha = 0.95$) and various $n$ values, are provided in the following table.

\begin{center}
\renewcommand{\arraystretch}{1.4}
\begin{tabular}{|c|c|c|}  \hline
$\alpha$ & $n$ & $t_{(1-\alpha)/2,(n-1)}~~ (= TINV(1-\alpha,n-1))$  \\  \hline
0.95	&	5	&	2.78	\\	\hline
0.95	&	6	&	2.57	\\	\hline
0.95	&	7	&	2.45	\\	\hline
0.95	&	8	&	2.36	\\	\hline
0.95	&	9	&	2.31	\\	\hline
0.95	&	10	&	2.26	\\	\hline
0.95	&	11	&	2.23	\\	\hline
0.95	&	12	&	2.20	\\	\hline
0.95	&	13	&	2.18	\\	\hline
0.95	&	14	&	2.16	\\	\hline
0.95	&	15	&	2.14	\\	\hline
0.95	&	16	&	2.13	\\	\hline
0.95	&	17	&	2.12	\\	\hline
0.95	&	18	&	2.11	\\	\hline
0.95	&	19	&	2.10	\\	\hline
0.95	&	20	&	2.09	\\	\hline
0.95	&	21	&	2.09	\\	\hline
\end{tabular}
\end{center}

Generally, the larger the number of observations, the smaller the 95\% confidence interval. As certain simulations are very time-consuming, a decision needs to be made on the tradeoff between time and accuracy. In many cases, when the simulations are not very time-consuming, we can increase the number of observations until required
accuracy (length of confidence interval) is achieved.

\subsubsection*{Homework \ref{simulations}.\arabic{homework}}
\addtocounter{homework}{1} \addtocounter{tothomework}{1} Choose a set of 12 different real numbers to represent outcomes of multiple measurements of the same quantity, and  calculate the confidence interval based on Student's t-distribution.
$~~~\Box$

\subsubsection*{Homework \ref{simulations}.\arabic{homework}}
\addtocounter{homework}{1} \addtocounter{tothomework}{1} Generate 10 uniform (0,1)
variates compute their average. Repeat this 11 times. These repetitions will result in 11 different estimations of $E[X]$ where $X$ is a uniform (0,1) random variable. Use these 11 estimations to obtain the confidence interval for the estimation of  $E[X]$ based on the 10 variates. Then, repeat this exercise five more times, increasing the sample size (number of variates) using  100, 1,000, 10,000, 100,000 and 1,000,000 uniform (0,1) variates. Observe the length of the confidence interval as you increase the sample size.

Then, repeat this homework for the following cases: (1) $X$ is a uniform $(a,b)$ random variable for a range of values $a$ and $b$, (2) $X$ is an exponential random variable with parameter $\lambda$ (consider different values for $\lambda$), and (3) $X$ is a Pareto distributed random variable with parameters  $\gamma$ and $\delta$  (again choose a range of parameter values for $\gamma$ and $\delta$ - in particular, consider the range $0< \gamma
\leq 2$, where the variance is infinite). To consider different random variables, use the inverse transformation method (inverse transform sampling). Observe and discuss the results as they vary for the different cases.
$~~~\Box$


\subsubsection*{Homework \ref{simulations}.\arabic{homework}}
\addtocounter{homework}{1} \addtocounter{tothomework}{1}
From an assignment in Section \ref{moments} (see also Equation (\ref{mean_win})), we have learned that the mean return in every roulette game is negative. In other words, in the long run, one loses money by playing roulette. Recall that we have made the independence assumption that the results of all roulette games are independent. Now consider the following strategy. Place a bet of $B=1$ dollar on the red. If you lose, double the bet in the next game. As long as you keep losing, keep doubling and redoubling your bets until you win. The argument for this strategy is that the number of times the player needs to bet until a win is geometrically distributed, and since the latter has a finite mean, eventually, the player wins. Henceforth, we will use the following terminology: a {\em game} is one bet -- one turn of the wheel, and a {\em sequence}  represents the number of games required to achieve a win or to run out of money. The drawback is twofold. Firstly, you may run out of money, and secondly, you may reach the {\em table limit} (an upper bound on $B$) set by the Casino management. Here we assume for simplicity that there is no table limit. However, we do take into consideration the case of running out of money. In particular, assume that you have \$10,000 to bet with at the beginning of a sequence. This means that a sequence cannot be longer than 13 games (why?) and if you run out of money, you lose \$8,191 (why?).

Recall the discussion and assignments in Section \ref{selectedcont} on how to generate variates from any probability distribution and use computer simulations to show consistency with the results of the previous assignment in Section \ref{moments}.
In particular, run at least six independent runs, each of which of several million sequences where in each sequence the player has fresh starts with \$10,000 regardless of the result of the previous sequence. In each run (of millions of sequences) calculate the total net winnings (total won minus total lost), and the total number of games played. The ratio between the two will give you an estimate for the mean win/loss per game. Repeating these runs multiple times, will provide a confidence interval for this mean value. Then, compare these results with the mean value obtained by Equation (\ref{mean_win}) and demonstrate consistency between the results. $~~~\Box$

We use here Student's t-distribution (and not Gaussian) because it
is the right distribution to use when we attempt to estimate the
mean of a population which is normally distributed  when we have a
small sample size. In fact, the need to estimate such a mean based on
a small sample gave rise to the development of the Student's
t-distribution.
In the next section, we will guide the reader on how to write queueing simulations for a G/G/1 queue. $~~~\Box$

\subsection{Simulation of a G/G/1 Queue}
\label{gg1simul}
We will now present an example of how to simulate a G/G/1
queue using an approach called {\it Discrete Event Simulation} \cite{Fishman01}.
Although the example presented here is for a G/G/1 queue,
the principles can be easily extended to multi-server and/or finite buffer queues.
The first step is to generate a sequence of inter-arrival
times and service times in accordance with the given statistical characteristics of the arrival process and service times.
Note the discussion in Section \ref{susecuniform} regarding the generation of random variates, but be aware that this discussion is limited to the generation of variates associated with IID random variables.
In our example, starting at time 0, let us consider
the following inter-arrival times: 1, 2, 1, 8, 4, 5,  \ldots, and the
following sequence of service times: 4, 6, 4, 2, 5, 1,  \ldots~.

In writing a computer simulation for G/G/1, we aim to fill in the following table for several 100,000s or millions of arrivals (rows).

\begin{center}
\renewcommand{\arraystretch}{1.2}
\begin{tabular}{|c|c|c|c|c|c| }  \hline
arrival time & service duration & queue-size on arrival & service starts & service ends & delay \\
\hline
 1 & 4 & 0& 1& 5 & 4\\
\hline
 3  & 6 & 1& 5& 11 & 8\\
\hline
 4 &  4& 2& & &  \\
\hline
12 & 2 & & & &  \\
\hline
 16 & 5 & & & &  \\
\hline
 21 &  1& & &  & \\
\hline
\end{tabular}
\end{center}

The following comments explain how to fill in the table.
\begin{itemize}
\item The arrival times and the service durations values are readily
obtained from the inter-arrival and service time sequences. \item
Assuming that the previous rows are already filled in, the
``queue-size on arrival'' is obtained by comparing the arrival
time of the current arrivals and the values in the ``service
ends'' column of the previous rows. In particular, the queue size
on arrival is equal to the number of customers that arrive
before the current customer (previous rows) that their ``service
ends'' time values are greater than the arrival time value of the
current arrival.
\item The ``service starts'' value is the maximum
of the ``arrival time'' value of the current arrival and the
``service end'' value of the previous arrival. Also, notice that if the
queue size on the arrival of the current arrival is equal to zero,
the service start value is equal to the ``arrival time''
value of the current arrival and if the
queue size on the arrival of the current arrival is greater than
zero the service start value is equal to the
``service end'' value of the previous arrival.
\item The ``service
ends'' value is simply the sum of the ``service starts'' and the
``service duration'' values of the current arrival.
\item The
``delay'' value is the difference between the ``service ends''
 and the ``arrival time'' values.
 \end{itemize}

 Using the results obtained in the last column, we can estimate the delay distribution and
 moments in steady-state. However, the  ``queue-size on arrival''
 values for all the customers do not, in general, provide directly the
 steady-state queue-size  distribution and moments. To estimate accurately the
 steady-state queue-size distribution, we will need to have
 inspections performed by an independent Poisson inspector. Fortunately, due to PASTA, for  M/G/1 (including M/M/1
 and M/D/1) the ``queue-size on arrival'' values can be used
 directly to obtain the steady-state queue-size distribution and
 moments, and a separate Poisson inspector is not required. Observing the
 queue size just before the arrivals provides the right inspections
 for steady-state queue-size statistics. However, if the arrival
 process does not follow a Poisson process, a separate independent Poisson
 inspector is required.  In such a case, we generate a Poisson
 process: $t_1$, $t_2$, $t_3$, ~ \ldots, and for each $t_i$, $i=1,~2,~3, ~ \ldots$ we can
 invoke the queue-size at time $t_i$, denoted $Q_i$, in a similar way to the one we obtained
 the ``queue-size on arrival'' values. The $Q_i$ values are then
 used to evaluate the queue-size distribution and moments.

 An alternative way to evaluate the queue size distribution of a
 G/G/1 queue is to record the total time spent in each state.
 If there was an event (arrival or departure) at time $t_j$ when the G/G/1 queue entered
 state $i$ and the next event (arrival or
departure) at $t_k$ when the G/G/1 queue exited
 state $i$, then the period $t_k - t_j$ is added
to a counter recording the total time spent in the state $i$.

\subsubsection*{Homework \ref{simulations}.\arabic{homework}}
\addtocounter{homework}{1} \addtocounter{tothomework}{1} First, fill in the above table  by hand. Then, modify the values in the columns ``arrival time'' and ``service duration'' and again fill in the table.
$~~~\Box$
\subsubsection*{Homework \ref{simulations}.\arabic{homework}}
\addtocounter{homework}{1} \addtocounter{tothomework}{1} Write a computer simulation for a P/P/1
queue (a single-server queue with Pareto inter-arrival and service
time distributions) to derive estimates for the mean and distribution of the delay and of the
queue size.  Perform the simulations for a wide range of parameter
values. Compute confidence interval as described in Section
\ref{simulations}. $~~~\Box$
\subsubsection*{Homework \ref{simulations}.\arabic{homework}}
\addtocounter{homework}{1} \addtocounter{tothomework}{1} Repeat the
simulations, of the previous homework, for a wide range of parameter
values, for a U/U/1 queue, defined as a single-server queue with
Uniform inter-arrival and service time distributions, and for an
M/M/1 queue. For the M/M/1 queue, verify that your simulation
results are consistent with respective analytical results. For the
U/U/1 queue, use the Poisson inspector approach and the ``time
recording'' approach and verify that the results are consistent.
$~~~\Box$
\subsubsection*{Homework
\ref{simulations}.\arabic{homework}} \addtocounter{homework}{1} \addtocounter{tothomework}{1}
Discuss the accuracy of your estimations in the different cases.
$~~~\Box$

\subsubsection*{Homework \ref{simulations}.\arabic{homework}}
\addtocounter{homework}{1} \addtocounter{tothomework}{1} Use the principles presented here for a
G/G/1 queue simulation to write a computer simulation for a G/G/$k$/$k$
queue. In particular, focus on the cases of an M/M/$k$/$k$ queue and a U/U/$k$/$k$ queue, defined as a $k$-server system without additional waiting room where the inter-arrival and service
times are uniformly distributed, and compute results for the blocking probability for these two cases.
For a meaningful comparison use a wide range of parameter values.$~~~\Box$

\newpage
\section{Deterministic Queues}
\label{det}

\setcounter{homework}{1} 

We consider here the simple case where inter-arrival and service
times are deterministic. To avoid ambiguity,
we assume that if an arrival and a departure occur at the same time,
the departure occurs first. Such an assumption is not required for
Markovian queues where the queue size process
follows a
continuous-time Markov chain because the probability of two events occurring at the same time is zero,
but it is needed for deterministic queues. Unlike many of the Markovian queues that we study in this book,
for deterministic queues steady-state queue size distribution does not exist because the queue size deterministically fluctuates
according to a certain pattern. Therefore, for deterministic queues, we will use the notation $P(Q=n)$,
normally designating the steady-state probability of the queue size to be equal to $n$ in cases where such steady-state probability exists,
for the proportion of time that there are $n$ customers in the queue,
or equivalently, $P(Q=n)$ is the probability of having $n$ in the queue at a randomly
(uniformly) chosen point in time. Accordingly, the mean queue size $E[Q]$ will be defined by $$E[Q]=\sum_{n=0}^\infty nP(Q=n).$$
We will use the term blocking probability
$P_b$ to designate the proportion of packets that are blocked. To derive performance measures such as mean queue size, blocking probability and  utilization,
in such deterministic queues, we follow the
queue-size process, for a certain transient period, until we
discover a pattern (cycle) that repeats itself. Then, we focus on a single cycle and obtain the desired measures of that cycle.

\subsection{D/D/1}

If we consider the  case $\lambda > \mu$, the D/D/1 queue is unstable.
In this case, the queue size constantly grows and approaches infinity as $t \rightarrow \infty$,
and since there are always packets in the queue waiting for service, the
server is always busy, thus the utilization is equal to one.

Let us consider now a stable D/D/1 queue, assuming $\lambda <
\mu$. Notice that for D/D/1, given our above assumption that if an arrival and a departure occur at the same time,
the departure occurs first, the case  $\lambda =
\mu$ will also be stable. Assume that the first arrival occurs at time $t=0$. The
service time of this arrival will terminate at $t=1/\mu$. Then, another arrival will occur at time $t=1/\lambda$ which will be completely
served at time $t=1/\lambda + 1/\mu$, etc. This gives rise to a
deterministic cyclic process where the queue size takes two
values: 0 and 1 with transitions from 0 to 1 in points of time
$n(1/\lambda)$, $n = 0, 1, 2, ~ \ldots$, and transitions from 1 to 0 in
points of time $n(1/\lambda)+1/\mu$, $n = 0, 1, 2, ~ \ldots~$. Each
cycle is of time-period $1/\lambda$ during which there is a customer to
be served for a time period of $1/\mu$ and there is no customer for a time period of $1/\lambda
- 1/\mu$. Therefore, the utilization is given by
$\hat{U}=(1/\mu)/(1/\lambda)=\lambda/\mu$ which is consistent with what we know about the utilization of
G/G/1.

As all the customers that enter the system are served before the next one arrives,
the mean queue size of D/D/1 must be equal to the mean queue-size
at the server, and therefore, it is also equal to the utilization.
In other words, the queue-size alternates between the values 1 and
0, spending a time-period of $1/\mu$ at state 1, then a time-period of $1/\lambda - 1/\mu$ at
state 0, then again $1/\mu$ time at state 1, etc. If we pick a
random point in time, the probability that there is one in the
queue is given by $P(Q=1)=(1/\mu)/(1/\lambda)$, and the probability that there are no customers in the queue is given by
$P(Q=0)=1- (1/\mu)/(1/\lambda)$. Therefore, the mean queue-size is
given by $E[Q]=0P(Q=0)+1P(Q=1)=(1/\mu)/(1/\lambda)=\hat{U}$.

Moreover, we can show that out of all possible stable G/G/1 queues,
 with $\lambda$ being the arrival rate and
$\mu$ the service rate, no one will have a lower mean queue size
than its equivalent D/D/1 queue. This can be shown using Little's formula $E[Q]=\lambda
E[D]$. Notice that for each of the relevant G/G/1 queues
$E[D]=1/\mu + E[W_Q]\geq 1/\mu$, but for D/D/1 $E[W_Q]=0$. Thus,
$E[D]$ for any G/G/1 queue must be equal or greater than that of
its equivalent D/D/1 queue, and consequently by Little's formula, $E[Q]$ for any G/G/1
queue must be equal to or greater than that of its equivalent D/D/1 queue.

\subsection{D/D/$k$}

Here we consider deterministic queues with multiple servers. The
inter-arrival times are again always equal to $1/\lambda$, and the
service time of all messages is equal to $1/\mu$.
Again if we consider the  case $\lambda > k\mu$, the D/D/$k$ queue is unstable.
In this case, the queue size constantly increases and approaches infinity as $t \rightarrow \infty$,
and since there are always more than $k$ packets in the queueing system, all $k$
servers are constantly busy, thus the utilization is equal to one.

Now consider the stable case of $\lambda<k\mu$, so that the arrival
rate is below the system capacity. Notice again that given our above assumption that if an arrival and a departure occur at the same time,
the departure occurs first, the case  $\lambda = k
\mu$ will also be stable.
Extending the D/D/1 example to a general number of servers, the behavior of the D/D/$k$ queue is analyzed as follows. As $\lambda$
and $\mu$ satisfy the stability condition $\lambda\leq k\mu$, there
must exist an integer $\hat{n}$, $1 \leq \hat{n} \leq k$ such that
\begin{equation} \label{intn}
(\hat{n}-1)\mu
< \lambda \leq \hat{n}\mu,
\end{equation}
or equivalently
\begin{equation} \label{intn2}
\frac{\hat{n}-1}{\lambda}
< \frac{1}{\mu} \leq \frac{\hat{n}}{\lambda}.
\end{equation}

\subsubsection*{Homework \ref{det}.\arabic{homework}}
\addtocounter{homework}{1} \addtocounter{tothomework}{1}
Show that
\begin{equation} \label{int_n}
\hat{n}=\left\lceil \frac{\lambda}{\mu}\right\rceil
\end{equation}
satisfies  $1 \leq \hat{n} \leq k$ and (\ref{intn2}). Recall that $\lceil x \rceil$
designates the smallest integer greater or equal to $x$.
\subsubsection*{Guide}
Notice that
$$ \frac{\hat{n}}{\lambda} = \frac{\left\lceil \frac{\lambda}{\mu} \right\rceil}{\lambda}\geq  \frac{\frac{\lambda}{\mu}}{\lambda} =\frac{1}{\mu}. $$
Also,
 $$\frac{\hat{n}-1}{\lambda}= \frac{\left\lceil \frac{\lambda}{\mu} \right\rceil-1}{\lambda} <
 \frac{\frac{\lambda}{\mu}}{\lambda}= \frac{1}{\mu}.    ~~~\Box  $$

The inequalities (\ref{intn2}) imply
that if the first arrival arrives at $t=0$, there will be
additional $\hat{n}-1$ arrivals before the first customer leaves the
system. Therefore, the queue size increases incrementally taking
the value $j$  at time $t=(j-1)/\lambda$, $j=1, ~2, ~3,~ \ldots, ~\hat{n}$.
When the queue reaches $\hat{n}$ for the first time, which happens at time $(\hat{n}-1)/\lambda$, the cyclic behavior
starts. Then, at time $t=1/\mu$ the queue-size reduces to $\hat{n}-1$
when the first customer completes its service. Next, at time
$t=\hat{n}/\lambda$, the queue-size increases to $\hat{n}$ and decreases to
$\hat{n}-1$ at time $t=1/\lambda + 1/\mu$ when the second customer
completes its service. This cyclic behavior continues forever
whereby the queue-size increases from $\hat{n}-1$ to $\hat{n}$ at time points
$t=(\hat{n}+i)/\lambda$,  and decreases from $\hat{n}$ to $\hat{n}-1$ at time points
$t=i/\lambda+1/\mu$, for $i=0, ~1, ~2,~ \ldots~$. The cycle length is
$1/\lambda$ during which the queue-size process is at state $\hat{n}$,
$1/\mu - (\hat{n}-1)/\lambda$ of the cycle time, and it is at state $\hat{n}-1$,
$\hat{n}/\lambda - 1/\mu$ of the cycle time. Thus,
$$P(Q=\hat{n})=\frac{\lambda}{\mu} -
(\hat{n}-1)$$ and $$P(Q=\hat{n}-1)=\hat{n}-\frac{\lambda}{\mu}.$$
The mean queue-size $E[Q]$,
can be obtained by
$$E[Q]=(\hat{n}-1)P(Q=\hat{n}-1)+\hat{n}P(Q=\hat{n})$$
which after some
algebra gives
\begin{equation}
\label{eqmmk} E[Q]=\frac{\lambda}{\mu}.
\end{equation}
\subsubsection*{Homework \ref{det}.\arabic{homework}}
\addtocounter{homework}{1} \addtocounter{tothomework}{1}
Perform the algebraic operations that lead to (\ref{eqmmk}). $~~~\Box$.

This result is consistent with Little's formula. As customers are
served as soon as they arrive, the time each of them spends in the
system is the service time $1/\mu$ - multiplying it by $\lambda$,
gives by Little's formula the mean queue size. Since $E[Q]$ in
D/D/$k$ gives the number of busy servers, the utilization is given
by
\begin{equation}
\label{Ummk} \hat{U}=\frac{\lambda}{k\mu}.
\end{equation}

Notice that Equations (\ref{eqmmk}) and (\ref{Ummk}) apply also to
D/D/$\infty$ for finite $\lambda$ and $\mu$. Eq.\@ (\ref{eqmmk})
gives the mean queue size of D/D/$\infty$ (by Little's formula, or
by following the arguments that led to Eq.\@ (\ref{eqmmk})) and for
D/D/$\infty$, we have that $\hat{U}=0$ by (\ref{Ummk}). Also notice
that in D/D/$\infty$ there are an infinite number of servers and the
number of busy servers is finite, so the average utilization per server must be equal to zero.

\subsection{D/D/$k$/$k$}
\label{ddkk}
In D/D/$k$/$k$ there is no waiting room beyond those available at
the servers. Recall that to avoid ambiguity, we assume that if
an arrival and a departure occur at the same time,
the departure occurs first. Accordingly, if $\lambda \leq k\mu$, then we have the same queue
behavior as in D/D/$k$ as no losses will occur. The interesting case
is the one where $\lambda > k\mu$ and this is the case we focus on.
Having $\lambda>k\mu$, or $1/\mu > k/\lambda$, implies that
$$\tilde{n}=\left\lceil \frac{\lambda}{\mu} \right\rceil-k$$
satisfies $$  \frac{k+\tilde{n}-1}{\lambda} < \frac{1}{\mu} \leq \frac{k+\tilde{n}}{\lambda}.$$
\subsubsection*{Homework \ref{det}.\arabic{homework}}
\addtocounter{homework}{1} \addtocounter{tothomework}{1}
Prove the last statement.
\subsubsection*{Guide} Notice that
$$ \frac{k+\tilde{n}}{\lambda} = \frac{\left\lceil \frac{\lambda}{\mu} \right\rceil}{\lambda}\geq  \frac{\frac{\lambda}{\mu}}{\lambda} =\frac{1}{\mu}. $$
Also,
 $$\frac{k+\tilde{n}-1}{\lambda}= \frac{\left\lceil \frac{\lambda}{\mu} \right\rceil-1}{\lambda} <
 \frac{\frac{\lambda}{\mu}}{\lambda}= \frac{1}{\mu}. ~~~\Box $$

\subsubsection{The D/D/$k$/$k$ Process and Its Cycles}

Again, consider an empty system with the first arrival occurring
at time $t=0$. There will be additional $k-1$ arrivals before all
the servers are busy. Notice that because $1/\mu > k/\lambda$, no
service completion occurs before the system is completely
full. Then, $\tilde{n}$ additional arrivals will be blocked before the
first customer completes its service at time $t=1/\mu$ at which time the
queue-size decreases from $k$ to $k-1$. Next, at time
$t=(k+\tilde{n})/\lambda$, the queue-size increases to $k$ and reduces to
$k-1$ at time $t=1/\lambda + 1/\mu$ when the second customer
completes its service. This behavior of the queue size alternating
between the states $k$ and $k-1$ continues until all the first
$k$ customers complete their service which happens at time
$t=(k-1)/\lambda + 1/\mu$ when the $k$th customer completes its
service, reducing the queue-size from $k$ to $k-1$. Next, an
arrival at time $t=(2k+\tilde{n}-1)/\lambda$ increased the queue-size from
$k-1$ to $k$. Notice that the point in time $t=(2k+\tilde{n}-1)/\lambda$
is an end-point of a cycle that started at $t=(k-1)/\lambda$. This
cycle comprises two parts: the first is a period of time where
the queue size stays constant at $k$ and all the arrivals are blocked,
and the second is a period of time during which no losses occur and the queue-size
alternates between $k$ and $k-1$.
Then, a new cycle of duration $(k+\tilde{n})/\lambda$ starts and this new cycle
ends at $t=(3k+2\tilde{n}-1)/\lambda$. In general, for each $j=1,~2,~3,~ \ldots$,
a cycle of duration $(k+\tilde{n})/\lambda$ starts at
$t=(jk+(j-1)\tilde{n}-1)/\lambda$ and ends at $t=((j+1)k+j\tilde{n}-1)/\lambda$.

\subsubsection{Blocking Probability, Mean Queue-size and Utilization}

In every cycle, there are $k+\tilde{n}$ arrivals out of which $\tilde{n}$ are blocked.
The blocking probability is therefore $$P_b=\frac{\tilde{n}}{k+\tilde{n}}.$$
Since
$$k+\tilde{n}=\left\lceil \frac{\lambda}{\mu} \right\rceil,$$
 the blocking probability is given by
\begin{equation}
\label{bpddkk}
P_b=\frac  {\left \lceil \frac{\lambda}{\mu}
\right \rceil-k}{\left \lceil \frac{\lambda}{\mu}  \right \rceil}.
\end{equation}
Let $A=\lambda/\mu$, the mean-queue size is obtained using Little's formula to be given by
\begin{equation}
\label{mqsddkk}
E[Q]=\frac{\lambda}{\mu}  (1-P_b) = \frac{kA}{\lceil A \rceil}.
\end{equation}
As in D/D/$k$, since every customer that enters a D/D/$k$/$k$ system does not
wait in a queue, but immediately enters service,  the utilization is given by
\begin{equation}
\label{utilddkk}
\hat{U}=\frac{E[Q]}{k}=\frac{A}{\lceil A \rceil}.
\end{equation}

\subsubsection{Proportion of Time Spent in Each State}

Let us now consider a single cycle and derive the proportion of
time spent in the states $k-1$ and $k$, denoted $P(Q=k-1)$ and $P(Q=k)$, respectively.
In particular, we consider the first cycle of duration
$$\frac{k+\tilde{n}}{\lambda} = \frac{\lceil A \rceil}{\lambda}$$
that starts at time $$t_s=\frac{k-1}{\lambda}$$
and ends at time $$t_e=\frac{2k+\tilde{n}-1}{\lambda}.$$
We define the first part of this cycle (the part during which arrivals are blocked)
to begin at $t_s$ and to end at the point in time when the $\tilde{n}$th arrival of this cycle is blocked which is
$$t_{\tilde{n}}=t_s + \frac{\tilde{n}}{\lambda} = \frac{k-1+\tilde{n}}{\lambda}=\frac{\lceil A \rceil-1}{\lambda}.$$
The second part of the cycle starts at $t_{\tilde{n}}$ and ends at $t_e$. The queue size is equal to $k$ for the entire duration
of the first part of the cycle. However, during the second part of the cycle, the queue-size alternates between the values $k$ and $k-1$ creating
a series of $k$ mini-cycles each of duration $1/\lambda$. Each of these mini-cycles is again composed of two parts.
During the first part of each mini-cycle, $Q=k$, and during  the second part of each mini-cycle, $Q=k-1$.
The first mini-cycle starts at time $t_{\tilde{n}}$ and ends at
$$t_{1e}=t_{\tilde{n}} + \frac{1}{\lambda}=\frac{\lceil A \rceil}{\lambda}.$$
The first part of the first mini-cycle starts at time $t_{\tilde{n}}$ and ends at time $1/\mu$, and the second part starts at $1/\mu$ and ends at
time $t_{1e}$. Thus, the time spent in each mini-cycle at state $Q=k-1$ is equal to
$$t_{1e}-\frac{1}{\mu} =  \frac{\lceil A \rceil}{\lambda} - \frac{1}{\mu} =  \frac{\lceil A \rceil}{\lambda} - \frac{\frac{\lambda}{\mu}}{\lambda} = \frac{\lceil A \rceil - A}{\lambda}.$$
Because there are $k$ mini-cycles in a cycle, we have that the total time spent in state $Q=k-1$ during a cycle is
$$\frac{k(\lceil A \rceil - A)}{\lambda}.$$
Because $P(Q=k-1)$ is the ratio of the latter to the total cycle duration, we obtain,
\begin{equation}
\label{pQequalk-1}
P(Q=k-1)=\frac{\frac{k(\lceil A \rceil - A)}{\lambda}}{\frac{\lceil A \rceil}{\lambda}}.
\end{equation}
The time spent in state $Q=k$ during each cycle is the total cycle duration minus the time spent in state $Q=k-1$. Therefore, we obtain
\begin{equation}
\label{pQequalk}
P(Q=k)=\frac{\frac{\lceil A \rceil}{\lambda}-\frac{k(\lceil A \rceil - A)}{\lambda}}{\frac{\lceil A \rceil}{\lambda}}.
\end{equation}

\subsubsection*{Homework \ref{det}.\arabic{homework}}
\addtocounter{homework}{1} \addtocounter{tothomework}{1}

\begin{enumerate}
\item
Show that the results for the queue-size probabilities $P(Q=k-1)$
and  $P(Q=k)$ in (\ref{pQequalk-1}) and (\ref{pQequalk}) are consistent with the result for the mean queue-size in
(\ref{mqsddkk}). In other words, show that

$$ (k-1)P(Q=k-1) + k P(Q=k) = E[Q] $$
or equivalently
$$(k-1)\left\{ \frac{\frac{k(\lceil A \rceil - A)}{\lambda}}{\frac{\lceil A \rceil}{\lambda}}\right\} + k \left\{   \frac{\frac{\lceil A \rceil}{\lambda}-\frac{k(\lceil A \rceil - A)}{\lambda}}{\frac{\lceil A \rceil}{\lambda}}    \right\} = \frac{kA}{\lceil A \rceil}.$$
\item Consider a D/D/3/3 queue with $1/\mu =
5.9$ and $1/\lambda = 1.1$. Start with the first arrival at $t=0$
and produce a two-column table showing the time of every arrival and
departure until $t=20$, and the corresponding queue-size values
immediately following each one of these events.
\item Write a general simulation program for a D/D/$k$/$k$ queue and use it to validate (\ref{mqsddkk}) and the results for $P(Q=k-1)$ and  $P(Q=k)$
in (\ref{pQequalk-1}) and (\ref{pQequalk}). Use it also to confirm the results you obtained for the D/D/3/3 queue.
\item Consider a D/D/1/$n$ queue for $n>1$. Describe the evolution of its queue-size
process and derive formulae for its mean queue size, mean delay,
utilization,
 and blocking probability. Confirm your results by simulation $~~~\Box$.
\end{enumerate}
\subsection{Summary of Results}
The following table summarizes the results on D/D/1, D/D/$k$ and D/D/$k$/$k$. Note that we do not
consider the case $\lambda=k\mu$ for which the results of the case
$\lambda<k\mu$ are applicable assuming that if a departure and an arrival occur
at the same time, the departure occurs before the arrival.
\begin{center}
\renewcommand{\arraystretch}{1.4}
\begin{tabular}{|c|c|c|c| }  \hline
Model & Condition & $E[Q]$ & $\hat{U}$  \\  \hline
D/D/1 & $\lambda<\mu$ & $\lambda/\mu$ & $\lambda/\mu$ \\  \hline
D/D/1 & $\lambda>\mu$ & $\infty$ & 1 \\  \hline
D/D/$k$ & $\lambda<k\mu$ & $A=\lambda/\mu$ &  $A/k$ \\ \hline
D/D/$k$ & $\lambda>k\mu$ & $\infty$ &  1\\ \hline
D/D/$k$/$k$ & $\lambda<k\mu$  & $A$  & $A/k$ \\  \hline
D/D/$k$/$k$ & $\lambda>k\mu$  & ${kA}/{\lceil A \rceil}$ & ${A}/{\lceil A \rceil}$ \\  \hline
\end{tabular}
\end{center}

\subsubsection*{Homework \ref{det}.\arabic{homework}}
\addtocounter{homework}{1} \addtocounter{tothomework}{1}
Justify the following statements.
\begin{enumerate}
\item D/D/1 is work conservative.
\item D/D/$k$ is work conservative (following a certain finite initial period) if $\lambda>k\mu$.
\item D/D/$k$ is not work conservative if $\lambda<k\mu$.
\item D/D/$k$/$k$ is not work conservative for all possible values of the parameters $\lambda$ and $\mu$ if we assume that if arrival and departure occur at the same time, then the arrival occurs before the departure.

\end{enumerate}
\subsubsection*{Guide}
Notice that D/D/$k$ is work conservative only if there are always more than $k$ customers in the system.
Notice that for  D/D/$k$/$k$ (under the above assumption) there are always periods
of time during which less than $k$ servers are busy.
$~~~\Box$.

\newpage
\section{M/M/1}
\label{mm1queue}

\setcounter{homework}{1} 

Having considered the straightforward cases of deterministic queues,
we will now discuss queues where the inter-arrival and
service times are non-deterministic. We will begin with cases where
the inter-arrival and service times are independent and
exponentially distributed (memoryless). Here we consider the M/M/1
queue where the arrival process follows a Poisson process with
parameter $\lambda$ and service times are assumed to be IID and
exponentially distributed with parameter $\mu$, and are independent
of the arrival process. As M/M/1 is a special case of G/G/1, all the
results that are applicable to G/G/1 are also applicable to M/M/1.
For example, $\hat{U}=\lambda/\mu$, $\pi_0=1-\lambda/\mu$ and Little's
formula. It is the simplest Markovian queue; it has only a single server
and an infinite buffer. It is equivalent to a continuous-time Markov chain on the
states: 0, 1, 2, 3, \ldots. Assuming that the M/M/1 queue-size process starts at state 0, it will stay in state 0 for a period of
time that is exponentially distributed with parameter $\lambda$ then
it moves to state 1. The time the process stays in state $n$, for
$n\geq 1$, is also exponentially distributed, but this time, it is a
competition between two exponential random variables, one of which is
the time until the next arrival - exponentially distributed with
parameter $\lambda$, and the other is the time until the next
departure - exponentially distributed with parameter $\mu$.
As discussed in Section \ref{exponential}, the
minimum of the two is therefore also exponential with parameter
$\lambda+\mu$ and this minimum is the time the process stays in
state $n$, for $n\geq 1$. We also know from the discussion in  Section \ref{exponential} that after spending
an exponential amount of time with parameter $\lambda+\mu$, the process will move to state $n+1$ with probability
$\lambda/(\lambda+\mu)$ and to state $n-1$ with probability
$\mu/(\lambda+\mu)$.

\subsection{Steady-State Queue Size Probabilities}
\label{mm1steadystate}

As the M/M/1 queue-size process increases by only one,  decreases
by only one, and stays an exponential amount of time at each state,
it is equivalent to a birth-and-death process. Therefore,
by Eqs. (\ref{qij}) and (\ref{qii}), the infinitesimal
generator (the ${\bf Q}$ matrix) for the M/M/1 queue-size process is given by

$Q_{i,i+1}=\lambda $ for $i=0, 1, 2, 3, ~ \ldots$\\
$Q_{i,i-1}=\mu $ for $i= 1, 2, 3, 4, ~ \ldots$\\
$Q_{0,0}=-\lambda $ \\
$Q_{i,i}=-(\lambda +\mu) $ for $i=1, 2, 3,~\ldots$~.

That is,

$${\bf Q} =
\begin{bmatrix}
    -\lambda & \lambda  & 0 & 0 & 0 & \ldots \\
    \mu & -(\lambda + \mu) & \lambda  & 0 & 0 & \ddots \\
    0 & \mu & -(\lambda + \mu) & \lambda  & 0 &  \ddots \\
    0 & 0 & \mu & -(\lambda + \mu) & \lambda  &  \ddots \\
    \vdots & \hspace{-8 mm} \ddots & \hspace{-16 mm}  \ddots & \hspace{-16 mm}  \ddots & \ddots &  \ddots
    \end{bmatrix}.
$$

Notice that each column in the infinitesimal generator is associated with one global balance equation. In particular, the $i$th column provides the infinitesimal rates of transitions into and out of state
$i$. The total rates out of state $i$ are all in the main diagonal of the ${\bf Q}$ matrix with a negative sign, and the rates into state $i$ are in positive matrix entries that are not on the main diagonal. For example, in the case of the first column of the ${\bf Q}$ matrix that is associated with the state $i=0$, we have the value of entry on the main diagonal is given by $Q_{0,0}=-\lambda$ which indicates that the total rates out of state $i=0$ is equal to $\lambda$, and the value of an entry in the first column which is not on the main diagonal is $Q_{1,0}=\mu $, and it indicates that the rate into state $i=0$ is equal to $\mu$ (from state $i=1$). In the case of the second column of the ${\bf Q}$ matrix associated with the rate $i=1$, we have the value of entry on the main diagonal is given by $Q_{1,1}=-(\lambda+\mu)$ which indicates that the total rates of $\lambda+\mu$ out of state $i=1$, and the value of  entries in the second column which is not on the main diagonal are given by $Q_{2,1}=\mu $ and $Q_{0,1}=\lambda $  that indicate that the rates into state $i=1$ are  $\mu$ (from state $i=2$) and $\lambda$ (from state $i=0$). In general, for the $i+1$ column of the ${\bf Q}$ matrix, for $i> 0$, associated with the rate $i$, we have the value of entry on the main diagonal  $Q_{i,i}=-(\lambda+\mu)$ implying total rates out of state $i$ to be equal to $\lambda+\mu$, and the value of  entries in the $i+1$ column not on the main diagonal are given by $Q_{i+1,i}=\mu $ and $Q_{i-1,i}=\lambda $. They indicate that the rates into state $i$ are  $\mu$ (from state $i+1$) and $\lambda$ (from state $i-1$).

Substituting this infinitesimal generator in Eq.\@ (\ref{ssqij})
we readily obtain the following global balance
steady-state
equations for the M/M/1 queue: $${\bf \Pi Q}=0.$$ That is,

$\pi_0 \lambda = \pi_1 \mu$\\
$\pi_1 (\lambda+\mu) = \pi_2 \mu + \pi_0 \lambda$\\
and in general for $i \geq 1$:
\begin{equation}
\label{sseqmm1global}
\pi_{i} (\lambda+\mu) = \pi_{i+1} \mu + \pi_{i-1} \lambda.
\end{equation}
To explain (\ref{sseqmm1global}) intuitively,
Let $L$ be a very long time.
During $L$, the total time that
the process stays in state $i$ is equal to $\pi_iL$.
For the case $i \geq 1$, since the arrival process is a Poisson process, the mean number of
transitions out of state $i$ is equal to $(\lambda + \mu)\pi_iL$. This
can be explained as follows.
For the case $i \geq 1$, the mean number of events that
occur during $L$ in state $i$ is $(\lambda+\mu)\pi_iL$ because
as soon as the process enters state $i$ it stays there on
average an amount of time equal to $1/(\mu + \lambda)$ and
then it moves out of
state $i$ to either
state $i+1$, or to state $i-1$. Since during time $\pi_iL$, there are, on average,
$(\lambda+\mu)\pi_iL$ interval times of size $1/(\mu + \lambda)$, then
$(\lambda+\mu)\pi_iL$ is also the mean number of events (arrivals and
departures) that occur in state $i$ during $L$. In a similar way we
can explain that for the case $i=0$, the mean number of transition
out of state $i=0$ is equal
to $(\lambda)\pi_iL$ for $i=0$. It is the mean number of events that occur during $L$
(because there are no departures at state $0$).

Recalling the notion of probability flux, introduced in Section
\ref{prob_flux}, we notice that the global balance equations (\ref{sseqmm1global})
equate for each state the total probability flux out of the state
and the total probability flux into that state.

A solution of the global balance equations (\ref{sseqmm1global})
together with the following normalizing equation that will guarantee
that the sum of the steady-state probabilities must be equal to one:
\begin{equation}
\label{sumto1mm1}
\sum_{j=0}^{\infty} \pi_j =1
\end{equation}
will give the steady-state probabilities of M/M/1.

However, the global balance equations (\ref{sseqmm1global}) can be
simplified recursively as follows.
We first write the first equation:
$\pi_0 \lambda = \pi_1 \mu$\\
Then, we write the second equation
$\pi_1 (\lambda+\mu) = \pi_2 \mu + \pi_0 \lambda$\\
Then, we observe that these two equations yield
$\pi_1 \lambda = \pi_2 \mu$\\
Then, recursively using all the equations (\ref{sseqmm1global}), we
obtain:
\begin{equation}
\label{sseqmm1} \pi_i \lambda = \pi_{i+1} \mu, ~{\rm for}
~i=0,~1,~2,~\ldots~.
\end{equation}

Notice that the steady-state equations (\ref{sseqmm1}) are the detailed balance
equations of the continuous-time Markov chain that describes the
stationary behavior of the
queue-size process of M/M/1. What we have noticed here is that the
global balance equations, in the case of the M/M/1 queue, are equivalent to the detailed balance
equations. In this case, a solution of the detailed balance
equations (\ref{sseqmm1}) and (\ref{sumto1mm1}) that sum up to unity will give the steady-state
probability distribution of the queue size. Recalling the discussion
we had in Section \ref{reversibility_CT}, this implies that the M/M/1
queue is reversible, which in turn implies that the output process
of M/M/1 is also a Poisson process. This is an important result that
will be discussed later in Section \ref{output}.

Another way to realize that the queue size process of M/M/1 is reversible
is to recall that this process is a birth-and-death process.
And we already know from Section \ref{reversibility_CT} that
birth-and-death processes are reversible.

Let $\rho=\lambda/\mu$, by (\ref{sseqmm1})we obtain,

$\pi_1 = \rho \pi_0 $\\
$\pi_2 = \rho \pi_1 = \rho^2 \pi_0$\\
$\pi_3 = \rho \pi_2 = \rho^3 \pi_0$

and in general:

\begin{equation}
\label{mm1piis}
\pi_i = \rho^i \pi_0  ~{\rm for} ~i=0,~1,~2,~\ldots~.
\end{equation}
As M/M/1 is a special case of G/G/1, we can use Eq.\@ (\ref{gg1}) to
obtain $\pi_0=1-\rho$, so
\begin{equation}
\label{mm1}
\pi_i = \rho^i (1-\rho)  ~{\rm for} ~i=0,~1,~2,~\ldots~.
\end{equation}

\subsubsection*{Homework \ref{mm1queue}.\arabic{homework}}
\addtocounter{homework}{1} \addtocounter{tothomework}{1}

Show that $\pi_0=1-\rho$ by summing up the $\pi_i$s in (\ref{mm1piis}) and equating the sum to 1. $~~~\Box$

\subsection{State Transition Diagram of M/M/1}

In general, state transition diagrams are used to represent a system as a collection of states and activities associated with various relationships
among the states. Such diagrams show how the system moves from one state to another, and the rates of movements between states.
State transition diagrams have many applications related to design and analysis  of real-time and object-oriented systems. Queueing systems that are modeled by continuous-time Markov chains are often described by their state transition diagram that provides the complete information of their  detailed balance equations. In particular, the state transition diagram of M/M/1 is\footnote{The author would like to thank Yin Chi Chan for his help in producing the various state transition diagrams in this book.}:

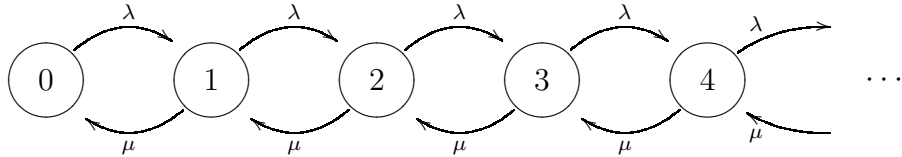
$$\xymatrix@C=32pt{
*++++[o][F]{0}\ar@/^{20pt}/[r]^{\lambda} &
*++++[o][F]{1}\ar@/^{20pt}/[l]^{\mu}\ar@/^{20pt}/[r]^{\lambda} &
*++++[o][F]{2}\ar@/^{20pt}/[l]^{\mu}\ar@/^{20pt}/[r]^{\lambda} &
*++++[o][F]{3}\ar@/^{20pt}/[l]^{\mu}\ar@/^{20pt}/[r]^{\lambda} &
*++++[o][F]{4}\ar@/^{20pt}/[l]^{\mu}\ar@/^{6pt}/[r]!<-20pt,20pt>^{\lambda} &
*++++[]{\cdots}\ar@/^{6pt}/[]!<-20pt,-20pt>;[l]^{\mu}
}
$$

The states are the numbers in the circles: 0, 1, 2, 3, $\ldots$, and the rates downwards and upwards are $\mu$ and $\lambda$, respectively.
We observe that the rates of
transitions between the states in the state transition diagram of M/M/1 are consistent with the rates in the
detailed balance equations of M/M/1 (\ref{sseqmm1}).

\subsection{Queue Performance Measures: $E[Q]$, $E[N_Q]$, $E[D]$, and $E[W_Q]$}
\label{perf_measure_mean}

Let $Q$ be a random number representing the
queue size in steady-state. Its mean is obtained by $E[Q]=\sum_{i=0}^\infty i\pi_i$.
This leads to:
\begin{equation}
\label{meanmm1}
E[Q]=\frac{\rho}{1-\rho}.
\end{equation}
\subsubsection*{Homework \ref{mm1queue}.\arabic{homework}}
\addtocounter{homework}{1} \addtocounter{tothomework}{1}
Perform the algebraic operations that lead to (\ref{meanmm1}). $~~~\Box$

Now, by (\ref{EQEQW}),
\begin{equation}
\label{EQEQNMM1} E[Q]=E[N_Q]+\rho,
\end{equation}
so, the mean number of customers in the queue (excluding the one in service) is obtained to be given by
\begin{equation}
\label{EQNMM1}
E[N_Q]=\frac{\rho}{1-\rho} -\rho = \frac{\rho^2}{1-\rho}
\end{equation}
and by Little's formula, the mean delay is obtained by
\begin{equation}
\label{EDMM1} E[D]=\frac{E[Q]}{\lambda}=\frac{\rho}{(1-\rho)\lambda}=\frac{1}{\mu-\lambda}.
\end{equation}
Then, the mean time of waiting in the queue (excluding the time in service) is obtained by (\ref{EDEW}) to be given by
\begin{equation}
\label{EWMM1} E[W_Q]=E[D] - 1/\mu = \frac{1}{\mu-\lambda} - \frac{1}{\mu} = \frac{\rho}{\mu - \lambda}.
\end{equation}
We already know that $\rho$ is equal to the ratio $\lambda/\mu$, and that it is equal to the server utilization (the proportion of time the server is busy). Now we also observe by (\ref{EDMM1}) and (\ref{EWMM1}) that $\rho$ is also equal to the ratio $E[W_Q]/E[D]$ in the case of M/M/1 (Credit: Ivanovich).
\subsubsection*{Homework \ref{mm1queue}.\arabic{homework}}
\addtocounter{homework}{1} \addtocounter{tothomework}{1}
Show that (\ref{EWMM1}) can be also obtained from (\ref{EQNMM1}) by using Little's formula. $~~~\Box$

\subsection{Using Z-Transform}
The Z-transform defined in Section \ref{transf}, also known as Probability Generating Function, is a powerful tool to derive
statistics of queueing behavior.

As an example, we will now demonstrate how the Z-transform is used to
derive the mean queue size of M/M/1.

Let us multiply the $n$th equation of (\ref{sseqmm1}) by $z^n$.
Summing up both sides  will give
\begin{equation}
\label{Ztrmm1}  \frac{\Psi(z)-\pi_0}{z}=\rho \Psi(z)
\end{equation}
where $\Psi(z)=\sum_{i=0}^{\infty} \pi_i z^i$.
Letting $z$
approach 1 (from below) gives
\begin{equation}
\pi_0=1-\rho
\end{equation}
 which is consistent with what we know already.
Substituting it back in (\ref{Ztrmm1}) gives after simple
algebraic manipulation:
\begin{equation}
\label{Ztramm1}  \Psi(z)=\frac{1-\rho}{1-\rho z}.
\end{equation}
Taking the derivative and substituting $z=1$, after some algebra we
obtain
\begin{equation}
\label{Ztreq}  E[Q]=\Psi^{(1)}(1)=\frac{\rho}{1-\rho}
\end{equation}
which is again consistent with what we know about M/M/1 queue.

\subsubsection*{Homework \ref{mm1queue}.\arabic{homework}}
\addtocounter{homework}{1} \addtocounter{tothomework}{1}
\begin{enumerate}
\item Derive equations (\ref{Ztrmm1}) -- (\ref{Ztreq}).
\item Derive the variance of the M/M/1 queue size using Z-transform.
$~~~\Box$
\end{enumerate}

\subsection{Mean Delay of Delayed Customers}
\label{meandelayed}
So far we were interested in the delay statistics of all customers. Now suppose that we are interested in the mean delay of only those customers that found the server busy upon
their arrivals and had to wait in the queue before they commence service. We assume that
the arrival rate $\lambda$ and the service rate $\mu$ are given,
then the  mean number of customers in the queue $E[N_Q]$ is given by
$$E[N_Q] =  E[Q] - E[N_s] = \frac{\rho}{1-\rho} - \rho = \frac{\rho^2}{1-\rho}.$$

Denote:\\
 $\hat{D}$ = The delay of a delayed customer including the service time\\
 $\hat{W_Q}$ = The delay of a delayed customer in the queue excluding the service time.

 To obtain $E[\hat{W_Q}]$, we use Little's formula where we consider the queue (without the server) as the system
and the arrival rate of the delayed customers which is $\lambda \rho$.
Thus $$E[\hat{W_Q}]= \frac{E[N_Q]}{\lambda \rho}=\frac{1}{\mu-\lambda},$$
and $$E[\hat{D}] = E[\hat{W_Q}] + \frac{1}{\mu}= \frac{1}{\mu-\lambda} + \frac{1}{\mu}.$$

Now, let us check the latter using the Law of Iterated Expectation as follows:
\begin{eqnarray*}
E[D] & = &  (1-\rho) [\mbox{Mean  delay of a non-delayed customer}] \\
& & + ~~~~~\rho [\mbox{Mean  delay of a delayed customer}]\\
& = & (1-\rho)\frac{1}{\mu} + \rho \left( \frac{1}{\mu-\lambda} + \frac{1}{\mu}\right) = \frac{1}{\mu-\lambda}.
\end{eqnarray*}
and we observe that consistency is achieved.   Notice that this consistency check is an alternative way to obtain $E[\hat{D}]$.

\subsubsection*{Homework \ref{mm1queue}.\arabic{homework}}
\addtocounter{homework}{1} \addtocounter{tothomework}{1}
Derive  $E[\hat{D}]$ using the Law of Iterated Expectation. $~~~\Box$

\subsection{Delay Distribution}

By (\ref{mm1}), and by the PASTA principle, an arriving customer from the moment it arrives until it leaves the
system will have to stay in the system a geometric number of IID phases (representing service times of the customers ahead of it in the queue, as well as its own service time) each of which is exponentially distributed with parameter $\mu$.
We have already shown that a geometrically distributed sum
of an IID exponentially distributed random variables is
exponentially distributed (see Eq. (\ref{Ltransform_sumN6}) in
Section \ref{Lap_transforms}). Therefore, the total delay of any
arriving customer in an M/M/1 system must be exponentially
distributed.

Accordingly, to derive the density of the delay, all that is left to do is
to obtain its mean which can be derived by
(\ref{meanmm1}) invoking Little's formula. Another way to obtain
the mean delay is by noticing from (\ref{mm1}) that the number of phases is
geometrically distributed with mean $1/(1-\rho)$. Observe that the mean number of phases
must equal $E[Q] + 1$ which is the mean queue size observed by an arriving
customer plus one more phase which is the service time of the
arriving customer. Thus, the mean number of phases is $$E[Q] + 1 = \frac{\rho}{1-\rho} +1 = \frac{1-\rho+\rho}{1-\rho}=\frac{1}{1-\rho}.$$

\subsubsection*{Homework \ref{mm1queue}.\arabic{homework}}
\addtocounter{homework}{1} \addtocounter{tothomework}{1}
Prove that the number of phases is geometrically distributed with mean $1/(1-\rho)$.
\subsubsection*{Guide}
Let $P_h$ be the number of phases. We know that in steady-state an arriving customer will find $Q$ customers in the system, where
$$P(Q=i)=\pi_i=\rho^i(1-\rho).$$
Since $P_h=Q+1$, we have $$P(P_h=n)=P(Q+1=n)=P(Q=n-1)=\rho^{n-1}(1-\rho). ~~~\Box $$

The
mean delay equals the mean number of phases times the mean service time $1/\mu$. Thus,
\begin{equation}
\label{meanDmm1} E[D]=\frac{1}{(1-\rho)\mu}= \frac{1}{\mu-\lambda}.
\end{equation}
This result is consistent with the result obtained in (\ref{EDMM1}) using Little's formula.

Substituting $1/E[D]=\mu-\lambda$ as the parameter of
exponential density, the density of the delay distribution
is obtained to be given by
\begin{equation}
\label{delaymm1} \delta_D(x)=\left\{\begin{array}{ll}
(\mu-\lambda) e^{(\lambda-\mu) x} & \mbox{if $x \geq 0$}\\
0 & \mbox{otherwise.}
\end{array}
\right.
\end{equation}

Having derived the distribution of the total delay (in the queue and in service),
let us now derive the distribution of the queueing delay (excluding the service time). That is, we are interested in deriving
$P(W_Q>t),~~~ t\geq 0$. By the Law of Total Probability, we obtain:

\begin{eqnarray}
\label{wq}
P(W_Q>t)& = & \rho P(W_Q>t|{\rm server~busy}) + (1-\rho)P(W_Q>t|{\rm server~ not~ busy}) \nonumber \\
        & = &\rho P(W_Q>t|{\rm server~ busy})~~ t\geq 0.
  \end{eqnarray}
To find $P(W_Q>t|{\rm server~ busy})$, let us find $P(N_q=n|{\rm server~ busy})$.

\begin{eqnarray*}
P(N_q=n|{\rm server~ busy}) & = & \frac{P(N_q=n \cap {\rm server~busy})}{P({\rm server ~ busy})} \\
                 & = & \frac{P(Q=n+1 \cap n+1\geq 1)}{\rho} \\
                 & = &  \frac{\rho^{n+1}(1-\rho)}{\rho} ~~ n=0,1,2,3 \ldots\\
                 & = & \rho^{n}(1-\rho)~~ n=0,1,2,3 \ldots .
            \end{eqnarray*}
Note that this is the same geometric distribution as that of $P(Q=n)$. Therefore, the random variable
$\{W_Q>t|{\rm server~ busy}\}$ is a geometric sum of exponential random variables and therefore has
exponential distribution. As a result,
$$  P(W_Q>t|{\rm server~ busy}) = e^{-(\mu-\lambda)t} $$
and by (\ref{wq}) we obtain
$$ P(W_Q>t) = \rho  e^{-(\mu-\lambda)t}   ~~ t\geq 0. $$

By the fact that the delay distribution in the cased of the M/M/1 queue is exponentially distributed with parameter $1/E[D]=\mu-\lambda$, expressed explicitly by Equation (\ref{delaymm1}), we can also obtain the variance of the delay in the M/M/1 queue to be given by

\begin{equation}
\label{vardelaymm1}
Var[D] = \frac{1}{(\mu-\lambda)^2}= \frac{1}{\mu^2(1-\rho)^2}.
\end{equation}

\subsection{The Departure Process}
\label{output}

We have already mentioned in Section \ref{mm1steadystate} the fact that the output process of an
M/M/1 is Poisson. This  is one of the results of the so-called Burke's theorem \cite{Burke56}.
In steady-state, the departure
process of a stable M/M/1, where $\rho<1$, is a
Poisson process with parameter $\lambda$ and is independent of the
number in the queue after the departures occur.
If we already know that the output process is Poisson, given that
the arrival rate is $\lambda$ and given that there are no losses,
all the traffic that enters must depart. Therefore the rate of the
output process must also be equal to $\lambda$.

We have shown in Section \ref{mm1steadystate} the reversibility of
M/M/1 queue-size process showing that the detailed balance equations
and the normalizing equation yield
the steady-state distribution of the queue-size process.

By reversibility, in steady-state, the arrival process of the reversed process must also follow a
Poisson process with parameter $\lambda$ and this process is the
departure process of the forward process.  Therefore the
departures follow a Poisson process and the
inter-departure times are  independent of the
number in the queue after the departures occur in the same way that inter-arrival
times are independent of the queue size before the arrivals.

Now that we know that in steady-state the departure process of a stable M/M/1 queue
is Poisson with parameter $\lambda$, we also know that, in steady-state, the
inter-departure times are also exponentially distributed
with parameter $\lambda$. We will now show this fact without using
the fact that the departure process is Poisson directly. Instead, we
will use it indirectly to induce PASTA for the reversed arrival
process to obtain that, following a departure, in steady-state, the
queue is empty with probability $1-\rho$ and non-empty with
probability $\rho$. If the queue is non-empty, the time until the
next departure is exponentially distributed with parameter $\mu$ --
this is the service time of the next customer. If the queue is
empty, we have to wait until the next customer arrival which is
exponentially distributed with parameter $\lambda$ and then we will
have to wait until the next departure which will take additional
time which is exponentially distributed. All together, if the queue
is empty, the time until the next departure is a sum of two
exponential random variables, one with parameter $\lambda$ and the
other with parameter $\mu$. Let $U_1$ and $U_2$ be two independent
exponential random variables with parameters $\lambda$ and $\mu$,
respectively. Define $U=U_1+U_2$, notice that $U$ is the convolution of $U_1$ and $U_2$,
and note that $U$ has hypoexponential
distribution.
Having the density $f_U(u)$, the density $f_D(t)$ of a
random variable $D$ representing the inter-departure time will be
given by
\begin{equation}
\label{interdeparture} f_D(t) = \rho \mu e^{-\mu t} + (1-\rho)
f_U(t).
\end{equation}
Knowing that $f_U(u)$ is a convolution of two exponentials, we
obtain
\begin{eqnarray*}
f_U(t) & = &   \int_{u=0}^t \lambda e^{-\lambda u} \mu e^{-\mu (t-u) } du\\
& = & \frac{\lambda \mu}{\mu - \lambda} \left( e^{-\lambda t} -
e^{-\mu t} \right).
\end{eqnarray*}
Then, by the latter and (\ref{interdeparture}), we obtain
\begin{equation}
\label{interdeparture1} f_D(t) = \rho \mu e^{-\mu t} + (1-\rho)
\frac{\lambda \mu}{\mu - \lambda} \left( e^{-\lambda t} - e^{-\mu t}
\right)
\end{equation}
which, after some algebra, gives
\begin{equation}
\label{interdeparture2} f_D(t) =  \lambda e^{-\lambda t}.
\end{equation}
This result is consistent with Burke's theorem.

\subsubsection*{Homework \ref{mm1queue}.\arabic{homework}}
 \addtocounter{homework}{1} \addtocounter{tothomework}{1} Complete all
the algebraic details in the derivation of equations
(\ref{interdeparture}) -- (\ref{interdeparture2}). $~~~\Box$

Another way to show consistency with Burke's theorem is the following.
Consider a stable $(\rho<1)$ M/M/1 queue.
Let $d_\epsilon$ be the unconditional number of departures, in steady state, that leave the
M/M/1 queue during a small interval of time of size $\epsilon$,  and
let $d_\epsilon(i)$ be the number of departures that leave the
M/M/1 queue during a small interval of time of size $\epsilon$ if
there are $i$ packets in our M/M/1 queue at the beginning of the interval.
Then, $P(d_\epsilon(i)>0)=o(\epsilon)$
if $i=0$, and $P(d_\epsilon (i)=1) = \epsilon \mu+ o(\epsilon)$ if $i>0$.
 Therefore,
in steady-state, $$P(d_\epsilon=1) = (1-\rho)0+ (\rho)\mu
\epsilon + o(\epsilon)=\epsilon \lambda + o(\epsilon),$$ which is a property consistent with the
assertion of Poisson output process  with parameter $\lambda$ in
steady-state.

\subsubsection*{Homework \ref{mm1queue}.\arabic{homework}}
 \addtocounter{homework}{1} \addtocounter{tothomework}{1}
So far, we have discussed the behavior of the M/M/1 departure
process in steady-state. You are now asked to demonstrate that the
M/M/1 departure process may not be Poisson with parameter $\lambda$
if we do not assume steady-state conditions. Consider an M/M/1 system
with arrival rate $\lambda$ and service rate $\mu$, assume that
$\rho=\lambda/\mu<1$ and that there are no customers in the system
at time $0$. Derive the distribution of the number of customers that
leave the system during the time interval $(0,t)$. Argue that this
distribution is, in most cases, not Poisson with parameter $\lambda
t$ and find a special case when it is.

\subsubsection*{Guide}
 Let $D(t)$ be a random variable representing the number
of customers that leave the system during the time interval $(0,t)$.
Let $X_p(\lambda t)$ be a Poisson random variable with parameter
$\lambda t$ and consider two cases: (a) the system is empty at time
$t$, and (b) the system is not empty at time $t$. In case (a),
$D(t)=X_p(\lambda t)$ (why?) and in case (b) $D(t)=X_p(\lambda
t)-Q(t)$ (why?) and use the notation used in Section \ref{ctmc}
$P_{00}(t)$ to denote the probability that in time $t$ the system is
empty, so the probability that the system is not empty at time $t$
is $1-P_{00}(t)$. Derive $P_{00}(t)$ using Eqs. (\ref{koback_bad})
and (\ref{koback_bad2}). Then, notice that $$D(t)=
P_{00}(t)X_p(\lambda t) + [1-P_{00}(t)][X_p(\lambda t) -Q(t)].
~~~\Box $$

Consider the limit
$$D_k(t) = \lim_{\Delta t \rightarrow 0} P[Q(t)=k \mid {\rm a~ departure ~occurs~
within} (t - \Delta t,t)].$$

Considering the fact that the reversed process is Poisson and
independence between departures before time $t$ and Q(t), we obtain
that
\begin{equation}
\label{pastaproofd} D_k(t) = P[Q(t)=k].
\end{equation}
Then, by taking the limit of both sides of (\ref{pastaproofd}), we
show that the queue size seen by a leaving customer is statistically
identical to the queue size seen by an independent observer.
$~~~\Box$

\subsubsection*{Homework \ref{mm1queue}.\arabic{homework}}
 \addtocounter{homework}{1} \addtocounter{tothomework}{1} Write a simulation of the M/M/1 queue by
 measuring queue size values in two ways: (1) just before arrivals and
 (2) just after departures. Verify that the results obtained for the
 mean queue size in steady-state are consistent. Use confidence intervals.
 Verify that the
 results are also consistent with analytical results. Repeat your
 simulations and computation for a wide range of parameter values
 (different $\rho$ values). Plot all the results in a graph
 including the confidence intervals (bars).
 $~~~\Box$

\subsection{Mean Busy Period and First Passage Time} \label{busyp}

The {\em busy period} of a single-server queueing system is defined
as the time between the point in time when the server starts being busy
and the point in time the server stops being busy. In other words,
it is the time elapsed from the moment a customer arrives at an
empty system until the first time the system is empty again.
Recalling the first passage time concept defined in Section
\ref{ctmc}, and that the M/M/1 system is in fact a continuous-time
Markov chain, the busy period is also the first passage time from
state 1 to state 0. The end of a busy period is the beginning of the
so-called {\em idle period} - a period during which the system is
empty. We know the mean of the idle period in an M/M/1 queue. It is
equal to $1/\lambda$ because it is the mean time until a new
customer arrives which is exponentially distributed with parameter
$\lambda$. A more interesting question is what is the mean busy
period. Let $T_B$ and $T_I$ be the busy and the idle periods,
respectively. Noticing that $E[T_B]/(E[T_B] + E[T_I])$ is the
proportion of time that the server is busy, thus it is equal to
$\rho$. Considering also that $E[T_I]=1/\lambda$, we obtain
\begin{equation}
\frac{E[T_B]}{E[T_B] + \frac{1}{\lambda}}=\rho.
\end{equation}
Therefore,
\begin{equation}
\label{busy} E[T_B]=\frac{1}{\mu-\lambda}, ~~~~~~~{\rm for } ~~ \mu > \lambda.
\end{equation}
Note that we are still assuming the stability condition of $\mu > \lambda$.

Interestingly, for the M/M/1 queue, the mean busy period is equal to
the mean delay of a single customer! This may seem counter-intuitive.
However, we can realize that there are many busy periods,
each of which is made of a single customer service time. It is likely that
for the majority of these busy periods (service times), their length is
shorter than the mean delay of a customer.
Furthermore, the fact that for the M/M/1 queue, the mean busy period is equal to
the mean delay of a single customer
can be proven by considering an M/M/1 queue
with a service policy of Last In First Out (LIFO). So far we
have considered only queues whose service policy is First In
First Out (FIFO). Let us consider an M/M/1 with LIFO with preemptive
priority. In such a queue the arrival and
service rates $\lambda$ and $\mu$, respectively, are the same as
those of the FIFO M/M/1, but in the LIFO
queue, the customer just arrived has preemptive
priority over all other customers that arrived before it, and
preempts the customers currently in service. (More information on LIFO queues is available in Section \ref{mg1lifo}.)

The two queues we consider, M/M/1 FIFO and
M/M/1 LIFO are both birth-and-death processes with the same
parameters so their respective queue size processes are
statistically the same. Then, by Little's formula, their respective
mean delays are also the same. Also, the delay of a customer in an
M/M/1 LIFO queue we consider is equal to the busy period in M/M/1
FIFO queue for the following reasons. Firstly, as the queue size processes in M/M/1 LIFO and FIFO queues are statistically the same, so are their respective busy periods. Secondly, under LIFO, the delay of an arriving message is independent of the queue size just before its arrival because all existing customers at the arrival time will have lower priority than the arriving customer. Thirdly, in the M/M/1 LIFO queue, the delay of a customer that arrives at an empty queue is equal to its own busy period. That is, its busy period starts at its arrival time and will end at its departure because all the customers that arrive during the time it spent in the system will complete their service before it, and when it completes its service, the system will be again empty at the moment of its departure. In conclusion, under M/M/1 LIFO, the mean delay of a customer is equal to the busy period, and therefore, the mean delay must be equal to the mean busy period in M/M/1 with the FIFO service policy.

\subsubsection*{Homework \ref{mm1queue}.\arabic{homework}}
 \addtocounter{homework}{1} \addtocounter{tothomework}{1} Derive an expression for the mean first passage
 time for M/M/1 from state
 $n$ to state 0 and from state 0 to state $n$, for $n\geq 3$. $~~~\Box$

\subsubsection*{Homework \ref{mm1queue}.\arabic{homework}}
 \addtocounter{homework}{1} \addtocounter{tothomework}{1} For a wide range of parameter values,
 simulate an
 M/M/1 system with FIFO service policy and an
 M/M/1 system with LIFO service policy
 with preemptive priority
 and compare their respective results for the
 mean delay, the variance of the delay, the mean queue
 size and the mean busy period.$~~~\Box$

\subsection{Dimensioning Based on Meeting Required Mean Delay}
\label{dim_mean}

A problem that often arises in practice is associated with resource allocation and dimensioning. Assume that we know the traffic demand; what is the minimal (least cost) service rate (or link capacity) such that a given delay requirement is met? A simple version of this problem is the following. We are given the required mean delay and the arrival rate, assume that M/M/1 conditions hold (SSQ, Poisson arrivals, and exponential service times), and we are asked to find the smallest value of $\mu$, called $\mu^*$, such that required mean delay ($E[D]_R$)is met.

To find $\mu^*$, we solve

 $$ E[D]_R = \frac{1}{\mu^* - \lambda}$$

 and we obtain

 \begin{equation}
\label{DIMMM1mu}\mu^* = \frac{1+ \lambda E[D]_R}{E[D]_R}.\end{equation}

 This will give the lowest (cheapest) service rate  ($\mu^*$) that meets the delay requirement ($E[D]_R$). Any service rate $\mu$, so that $\mu < \mu^*$ will violate the mean delay requirement, and any service rate $\mu$, so that $\mu > \mu^*$, will be unnecessarily too expensive. Therefore, $\mu^*$ is the optimal choice for service rate.

 A second dimensioning problem that arises is the following. We are given the service rate $\mu$ and the required mean  delay ($E[D]_R$). What is the highest possible arrival rate $\lambda^*$ that meets the delay requirement ($E[D]_R$).

 Now we need first to check if a feasible solution exists. That is, we need to check if
\begin{equation}
\label{DIMMM1}\frac{1}{\mu} \leq E[D]_R. \end{equation}
Notice that the delay ($D$) includes the service time, so if (\ref{DIMMM1}) does not hold, namely, the required mean delay is lower than the mean service time, no feasible solution exists, we cannot satisfy this requirement.

 On the other hand, if the condition (\ref{DIMMM1}) holds, then we simply solve the equation

  $$ E[D]_R = \frac{1}{\mu - \lambda^*}$$

  and we obtain

   $$\lambda^* = \frac{\mu E[D]_R -1 }{E[D]_R}.$$

\subsubsection*{Homework \ref{mm1queue}.\arabic{homework}}
 \addtocounter{homework}{1} \addtocounter{tothomework}{1}
 Show that $\mu^*$ obtained by (\ref{DIMMM1mu}) satisfies the condition $$ \frac{1}{\mu^*} \leq E[D]_R. $$
 $~~~\Box$

\subsection{Effect of Rising Internet Bit-rate on Link Efficiency and QoS}
\label{rising speed}

With the ever-increasing Internet bit-rate, it is interesting to observe the implication of Equation (\ref{EDMM1}) on the efficiency of links with the rising bit-rate of Internet links.
Observe that while the server utilization in M/M/1 is the ratio of $\lambda$ to $\mu$ (i.e., $\rho = \lambda/\mu$), the mean delay is a function of the difference between $\mu$ and $\lambda$. Therefore, we can maintain constant $E[D]$ by keeping the difference between $\mu$ and $\lambda$ constant and increasing $\lambda$ and $\mu$ towards infinity. Under such limiting conditions the utilization will approach unity. This can be shown by the following derivation. From Equation (\ref{EDMM1}), we have:

$$(\mu  -\lambda) E[D] = 1 $$
or
$$\mu - \lambda = \frac{1}{E[D]}. $$
Dividing both sides by $\mu$, we obtain
$$1 - \rho = \frac{1}{\mu E[D]}$$
or
$$ \rho = 1 - \frac{1}{\mu E[D]}.$$

Therefore, for constant $E[D]$, $\rho$ approaches one as $\mu$ (and also $\lambda$) approaches infinity.

Notice also that by Equation (\ref{EDMM1}), the mean delay in M/M/1 is given by $$E[D] = \frac{1}{\mu-\lambda}  =
\frac{1}{\mu(1-\rho)},$$ so for a fixed $\rho$ and arbitrarily large $\mu$, the delay is arbitrarily small.

\subsubsection*{Homework \ref{mm1queue}.\arabic{homework}}
\addtocounter{homework}{1} \addtocounter{tothomework}{1}
Packets destined for a given destination arrive at a router according to a Poisson process with a rate of 2000 packets per millisecond. The packet size follows an exponential distribution with a mean of 625 bytes. The router has a very large buffer and serves these packets by transmitting them through a single 10.5 Gb/s output link. The service policy is First Come First Served. Assume that the system is in steady state.

Compute the mean queue size (packets) and the mean total packet delay (including queueing time and service time). What do you observe from the answer?

\subsubsection*{Solution}

 Considering the Poisson arrivals, exponential service times and infinite buffer, M/M/1 is an appropriate model for this case.

 $\lambda	= 2000~~{\rm [packet/millisecond]} = 2 \times 10^{6} ~~{\rm [packet/s]}$

 $\mu = 10.5 ~~{\rm [Gb/s]} / (625 \times 8)~~{\rm [bits]} = 2.1 \times 10^{6} ~~{\rm [packet/s]}.$

Consider an M/M/1 queue for the case:
$\lambda   = 2 \times 10^{6}$ and
$\mu = 2.1 \times 10^{6}$. We obtain

$$\rho = \frac{\lambda}{\mu} = 0.952 ~~{\rm approx.}$$

$$E[Q] = \frac{\rho}{[1- \rho]} =   20 ~~{\rm [packets]} $$

$$E[D] = \frac{E[Q]}{\lambda}  = 10^{-5}~~{\rm [seconds]}. $$

The delay is very small even if the utilization is high because of the high bit-rate (service rate). ~~~$\Box$


\subsection{Multiplexing}
\label{muxmm1}
In telecommunications, the concept of multiplexing refers to a variety of schemes or techniques
that enable multiple traffic streams from possibly different sources to share a common transmission resource.
In certain situations such sharing of a resource can lead to a significant improvement in efficiency.
In this section, we use the M/M/1 queueing model where traffic from multiple sources that follow Poisson processes enters an infinite-length buffer and are served by a single server with exponential service time. This model enables us to gain insight into efficiency gain of multiplexing.

Note that this discussion of efficiency gain of multiplexing is very much related and complementary to our previous discussion on efficiency gain due to the rising bit-rate of Internet links in Section (\ref{rising speed}).

An important and interesting observation we can make by considering
the M/M/1 queueing performance results (\ref{mm1})--(\ref{vardelaymm1})
is that while the queue-size statistics are dependent only on $\rho$
(the ratio of the arrival rate to the service rate), the delay
statistics (mean and distribution) are a function of what we call
the {\em spare capacity} (or {\it mean net input}) which is the
difference between the service rate and the arrival rate.

Assume that our traffic model obeys the M/M/1 assumptions. Then, if
the arrival rate increases from $\lambda$ to $N\lambda$ and we
increase the service rate from $\mu$ to $N\mu$  (maintaining the
same $\rho$), the mean queue size and its distribution will remain
the same. However, in this scenario, the mean delay does not remain
the same. It reduces by $N$ times to $1/[N(\mu-\lambda)]$.

This is applicable to a situation where we have $N$ individual M/M/1 queues
each of which with arrival rate
$\lambda$ and service rate $\mu$. Then, we superpose (multiplex) all the arrival processes together
which results in a Poisson process of rate $N\lambda$. An interesting question is the following.
If we replace all the individual servers (each of which has service rate $\mu$) with one fast
server that serves the superposed Poisson stream of rate $N\lambda$, what service rate
this fast server should operate at.

If our QoS measure of interest is the mean delay or the probability that
the delay exceeds a certain value, and if for a given arrival rate
$\lambda$ there is a service rate $\mu$ such that our delay-related QoS measure
is just met, then if the arrival rate increases from
$\lambda$ to $N\lambda$, and we aim to find the service rate $\mu^*$
such that the delay-related QoS measure is just met, we will need to
make sure that the spare capacity is maintained, that is
\begin{equation}
\mu-\lambda = \mu^* - N\lambda
\end{equation}
or
\begin{equation}
\label{mustar}
\mu^* = \mu + (N-1)\lambda
\end{equation}
so by the latter and the stability condition of $\mu > \lambda$, which implies that the right-hand side of (\ref{mustar}) is less than $\mu + (N-1)\mu$ for $N>1$, we
must have that $\mu^* < N\mu$ for $N>1$. We can therefore define a measure for
multiplexing gain to be given by
\begin{equation}
M_{mg} = \frac{N\mu - \mu^*}{N\mu}
\end{equation}
so by (\ref{mustar}), we obtain
\begin{equation}
\label{multiplexingain}
M_{mg} = \frac{N-1}{N}   (1-\rho).
\end{equation}
Recalling the stability condition $\rho < 1$ and the fact that
$\pi_0=1-\rho$ is the proportion of time that the server is idle at
an individual queue, Eq.\@ (\ref{multiplexingain}) implies that
$(N-1)/N$ is the proportion of this idle time gained by
multiplexing. For example, consider the case $N=2$, that is, we
consider multiplexing of two M/M/1 queues each with parameters
$\lambda$ and $\mu$. In this case, half of the server idle time (or
efficiency wastage) in an individual queue can be gained back by
multiplexing the two streams to be served by a server that serves at
the rate of $\mu^* = \mu + (N-1)\lambda= \mu + \lambda$. The
following four messages follow from Eq.\@ (\ref{multiplexingain}).
\begin{enumerate}
\item The multiplexing gain is positive for all $N>1$.
\item The multiplexing gain increases with $N$.
\item The multiplexing gain is bounded above by $1-\rho$.
\item In the limiting condition as  $N\rightarrow \infty$,
the multiplexing gain approaches its bound $1-\rho$.
\end{enumerate}

The $1-\rho$ bound means also that if $\rho$
is very close to 1, then the multiplexing gain diminishes because in
this case the individual M/M/1 queues are already very efficient in terms of server utilization so
there is little room for improvement. On the other hand, if we have a case where
the QoS requirements are strict (requiring very low mean queueing
delay) such that the utilization $\rho$ is low, the potential
for multiplexing gain is high.

Let us now apply our general discussion on multiplexing to obtain
insight into performance comparison between two commonly used
multiple access techniques used in telecommunications. One such
technique is called {\it Time Division Multiple Access} (TDMA)
whereby each user is assigned one or more channels (in the form of
time slots) to access the network. Another approach, which we call
{\it full multiplexing} (FMUX), is to let all users to separately
send the data that they wish to transmit to a switch which then
forwards the data to the destination. That is, all the data is
stored in one buffer (in the switch) which is served by the entire
available link capacity.

To compare the two approaches, let us consider $N$ users
each transmitting packets at an average rate of $R_u$ [bits/second].
The average packet size denoted $S_u$ [bits] is assumed equal for
the different users. Let $\hat{\lambda}$ [packets/second] be the
packet rate
 generated by each of the users. Thus, $\hat{\lambda} = R_u/S_u$. Under TDMA,
 each of the users obtains a service rate of  $B_u$ [bits/sec]. Packet sizes are assumed to be exponentially
 distributed with mean $S_u$ [bits], so the service rate in packets/second denoted $\hat{\mu}$ is given by
$\hat{\mu} = B_u/S_u$. The packet service time is therefore
exponentially distributed with parameter $\hat{\mu}$. Letting
$\hat{\rho}=\hat{\lambda}/\hat{\mu}$, the mean queue size under
TDMA is given by
\begin{equation}
\label{meanmm1tdma} E[Q_{TDMA}]=\frac{\hat{\rho}}{1-\hat{\rho}},
\end{equation}
and the mean delay is
\begin{equation}
\label{meanDmm1tdma} E[D_{TDMA}]= \frac{1}{\hat{\mu}-\hat{\lambda}}.
\end{equation}
In the FMUX case, the total arrival rate is $N\hat{\lambda}$ and the
service rate is $N\hat{\mu}$, so in this case, the ratio between the
arrival and service rate remains the same, so the mean queue size
that only depends on this ratio remaining the same
\begin{equation}
\label{meanmm1fm} E[Q_{FMUX}]=\frac{\hat{\rho}}{1-\hat{\rho}} =
E[Q_{TDMA}].
\end{equation}
However, we can observe an $N$-fold reduction in the mean delay:
\begin{equation}
\label{meanDmm1fm} E[D_{FMUX}]=
\frac{1}{N\hat{\mu}-N\hat{\lambda}}=\frac{E[D_{TDMA}]}{N}.
\end{equation}
Consider a telecommunication provider that wishes to meet packet
delay requirement of its $N$ customers, assuming that the delay that
the customers experienced under TDMA was satisfactory, and assuming
that the M/M/1 assumptions hold, such a provider does not need a total
capacity of $N\hat{\mu}$ for the FMUX alternative. It is sufficient
to allocate $\hat{\mu} + (N-1)\hat{\lambda}$.

\subsubsection*{Homework \ref{mm1queue}.\arabic{homework}}
 \addtocounter{homework}{1} \addtocounter{tothomework}{1}
 Consider a telecommunication provider that aims to serve a network of 100 users
 each transmitting data at an overall average rate of 1 Mb/s. The mean packet size is
 1 kbit. Assume that packet lengths are exponentially distributed
 and that the process of packets generated by each user follows a Poisson process.
 Further assume that the mean packet delay requirement is 50 milliseconds.
 How much total capacity (bit-rate) is required to serve the 100
 users under TDMA and under FMUX.

 \subsubsection*{Guide}
 The arrival rate of each user is 1 Mb/s / 1 kbit = 1000 packets/s.
 For TDMA, use Eq.\@ (\ref{meanDmm1tdma}) and substitute $E[D_{TDMA}]=0.05$ and $\hat{\lambda}=1000$, to compute $\hat{\mu}$.
 This gives $\hat{\mu}=1020$ [packets/s] or bit-rate of 1.02 Mb/s per each user. For 100 users the required rate is 102,000
 packets/s or bit-rate of 102 Mb/s. For FMUX the required rate is $\hat{\mu} + (N-1)\hat{\lambda}$
 which is 100,020 packets/s or 100.02 Mb/s (calculate and verify it).
 The savings in using FMUX versus TDMA is therefore 1.98 Mb/s.
 $~~~\Box$

\subsection{Dimensioning Based on Delay Distribution}
\label{dimdelay}
In  Section \ref{dim_mean} we have considered dimensioning based on average delay. That is, the aim was to meet or maintain QoS measured by average delay. It may be, however, more practical to aim for a percentile of the delay distribution; e.g., to require that
no more than 1\% of the packets will experience over 100-millisecond delay.

In the context of the M/M/1 model, we define two dimensioning problems.
The first problem is: for a given $\lambda$, $t\geq 0$ and $\alpha$, find a minimal value for $\mu$, denoted $\mu^*$, such that
$$P(D>t) = e^{-(\mu^*-\lambda)t} \leq \alpha.$$

The solution is:

$$\mu^* = \lambda - \frac{\ln (\alpha)}{t}.$$

The second problem is: for a given $\mu$, $t\geq 0$ and $\alpha$, find maximal $\lambda$, denoted $\lambda^*$, such that

$$P(D>t) = e^{-(\mu-\lambda^*)t} \leq \alpha.$$

In this case, for a certain range, the solution is not feasible because the delay includes the service time and can never be less than the service
time. That is, for certain parameter values, even if the arrival rate is very low, the delay
requirements cannot be met, simply because the service time requirement exceeds the total delay requirement. To find the feasible range set $\lambda^*=0$, and obtain

$$\mu > \frac{-\ln (\alpha)}{t}.$$

This requires that the probability of the service time requirement exceeding $t$ is less than $\alpha$.

In other words, if this condition does not hold there is no feasible solution to the optimal dimensioning problem.

If a solution is feasible, the $\lambda^*$ is obtained by solving the equation

$$P(D>t) = e^{-(\mu-\lambda^*)t}= \alpha$$

and we obtain

$$\lambda^* = \frac{\ln (\alpha)}{t} + \mu.$$

\subsection{A Markov-chain Simulation of M/M/1}

A simulation of am M/M/1 queue can be made as a special case of G/G/1 as described before, or it can be
simplified by taking advantage of the M/M/1 Markov-chain structure
if we are not interested in performance measures that are associated with times (such as delay distribution).
If our aim is to evaluate queue size statistics or blocking probability, we can avoid tracking the time.
All we need to do is to collect the relevant information about the process at PASTA time points without even knowing what are the
times at these points. Generally speaking, using the random walk simulation approach, also called the {\it Random Walk simulation} approach, we simulate the evolution of the state of the process based on the transition
probability matrix and collect information on the values of interest at selective PASTA points without being concerned about the time.
We will now explain how these ideas may be applied to a few relevant examples.

If we wish to evaluate the mean queue size of an M/M/1 queue, we can write the following simulation.

Variables and input parameters: $Q$ = queue size; $\hat{E}(Q)$ = estimation for the mean queue size;
$N$ = number of $Q$-measurements taken so far which is also equal to the number of arrivals so far;
$MAXN$ = maximal number of $Q$-measurements taken;  $\mu$ = service rate;
$\lambda$ = arrival rate.

Define function: $I(Q) = 1$ if $Q > 0$; $I(Q) = 0$ if $Q = 0$.\\
Define function: $R(01)$ = a uniform $U(0,1)$ random variate. A new value for R(01) is generated every time it is called.

Initialization: $Q = 0$; $\hat{E}[Q] = 0$; $N = 0$.

1. If $R(01) \leq \lambda/(\lambda + I(Q)\mu)$, then $N = N + 1$, $\hat{E}(Q) = [(N-1) \hat{E}(Q) + Q]/N$, and $Q = Q + 1$;

else, $Q = Q - 1$.

2. If $N < MAXN$ go to 1; else, print $\hat{E}(Q)$.

This signifies the simplicity of the simulation. It has only two IF
statements: one to check if the next event is an arrival or a
departure according to Eq.\@ (\ref{Pwin}), and the second is merely
a stopping criterion.

\subsubsection*{Comments:}
\begin{enumerate}
\item The operation $Q = Q + 1$ is performed after the $Q$ measurement is
taken. This is done because we are interested in $Q$ values seen by arrivals just before they arrive.
If we include the arrivals after they arrive, we violate the PASTA principle. Notice that if we do that,
we never observe a $Q=0$ value which, of course, will not lead to an accurate estimation of $E[Q]$.
\item If the condition $R(01) \leq \lambda/(\lambda + I(Q)\mu)$ holds we have an
arrival. Otherwise, we have a departure.
This condition is true with probability $\lambda/(\lambda + I(Q)\mu)$. If $Q=0$, then $I(Q)=0$ in which case the next event is an arrival with probability 1. This is clearly
intuitive. If the system is empty, no departure can occur, so the next event must be an
arrival. If $Q>0$, the next event is an arrival with probability
$\lambda/(\lambda + \mu)$ and a departure with probability
$\mu/(\lambda + \mu)$. We have here a competition between two exponential
random variables: one (arrival) with parameter $\lambda$ and the other (departure) with
parameter $\mu$. According to the discussion in Section \ref{exponential} and as
mentioned in the introduction to this section, the probability that the
arrival ``wins'' is $\lambda/(\lambda + \mu)$, and
the probability that the departure ``wins'' is $\mu/(\lambda + \mu)$.
\item In the case of a departure, all we do is decrement the queue size; namely,
$Q = Q - 1$. We do not record the queue size at these points because, according to PASTA, arrivals
see time averages. (Notice that due to reversibility, if we measure the queue size
immediately {\bf after} departure points, we will also see time averages.)
 \end{enumerate}

\subsubsection*{Homework \ref{mm1queue}.\arabic{homework}}
\addtocounter{homework}{1} \addtocounter{tothomework}{1} Simulate an M/M/1 queue using a Markov-chain simulation
to evaluate the
mean queue size for the cases of Section \ref{gg1simul}. Compare
the results with the results obtained analytically and with those
obtained using the G/G/1 simulation principles. In your comparison,
consider accuracy (closeness to the analytical results), the length
of the confidence intervals, and running times. $~~~\Box$

\newpage
\section{M/M/$\infty$}
\label{secmminf}

\setcounter{homework}{1} 

The next queueing system we consider is the M/M/$\infty$ queueing system where the number of servers is infinite. Because the
number of servers is infinite, the buffer capacity is unlimited and
arrivals are never blocked. We assume that the arrival process is
Poisson with parameter $\lambda$ and each server renders a service
which its time is exponentially distributed with parameters $\mu$. As in the
case of M/M/1, we assume that the service times are independent and
independent of the arrival process.

Based on the discussion in Section \ref{offered_carried}, for M/M/$\infty$ we have the following results:

$$A =~ {\rm offered ~traffic} = {\rm carried ~traffic}~ = E[Q] $$

because all the traffic is admitted as there are no losses in  M/M/$\infty$, and all the offered traffic is admitted and carried by the system.

The last equality can be explained in two ways. First, one of the definitions of carried traffic is the long run (steady state) average number of busy servers, which in M/M/$\infty$ equal to the long run mean number of customers in the system, and also by Little's formula,
the mean number of
customers in the system is equal to the arrival rate ($\lambda$)
times the mean time a customer spends in the system, which is equal to $1/\mu$ in the case of M/M/$\infty$. Because there are no losses in M/M/$\infty$, all the arriving traffic enters the service system, so  we obtain
\begin{equation} \label{eqeqA} E[Q]=\lambda (1/\mu) = A. \end{equation}

In practice, the number of servers (channels or circuits) is limited, and the
offered traffic is higher than the carried traffic because
some of the calls are blocked due to call congestion when all
circuits are busy. A queueing model which describes this more realistic case is the M/M/$k$/$k$ queueing model discussed in the next chapter.

Notice (again) that we use
 the notation $A$ here for the offered traffic, which is equal to the ratio
$\lambda/\mu$ while we used the notation $\rho$ for this ratio in
the M/M/1 case as $\rho$ is defined, in Section \ref{littles}, by the ratio $A/k$ and in an SSQ $k=1$ so $\rho=A$. Clearly, both $A$ and $\rho$ represent the offered traffic by
definition, but $\rho$ is  used for SSQs and $A$ for multi-server queues. In addition, $\rho$ and $A$ represent
 the mean number of busy servers in the M/M/1 and M/M/$\infty$ cases, respectively. Recall that we used the notation $E[N_s]$ for the mean number of busy servers (bounded above by 1) for M/M/1 and G/G/1. We have shown in (\ref{ENsGG1}) that this is true for a G/G/1 queue (and therefore also for M/M/1), and in (\ref{eqeqA}) we showed
 that it is true for an  M/M/$\infty$ queue.


In M/M/1, we must have that $\rho$ cannot exceed unity for stability.
In M/M/1, $\rho$ also represents the server utilization which cannot exceed unity.
However, in   M/M/$\infty$, $A$ can take any non-negative value and we often have $A>1$.
M/M/$\infty$ is stable
for any $A\geq 0$. Notice that in M/M/$\infty$,  the service rate increases with the number of busy servers, and
when we reach a situation where the number of busy servers $j$ is higher than $A$ (namely $j>A=\lambda/\mu$),
we will have that the system service rate is higher than the arrival rate (namely $j\mu > \lambda$).

As discussed in Chapter \ref{general},   offered traffic is measured in {\em erlangs}.
One erlang represents the traffic load of one arrival, on average, per mean service time.

\subsection{Steady-State Equations}

As for M/M/1, the queue-size process of an M/M/$\infty$ system can also be viewed as a continuous-time Markov chain with the
state being the queue size (the number of customers in the system).
As for M/M/1, since in M/M/$\infty$ queue-size process is stationary and its
transitions can only occur
upwards by one or downwards by one, the queue-size process is a birth-and-death
process and therefore it is reversible.

As in M/M/1, the arrival rate is independent of changes in the
queue-size. However, unlike M/M/1, in M/M/$\infty$, the service rate
does change with the queue size. When there are $n$ customers in the system, and at the same time, $n$ servers are busy, the service rate
is $n\mu$, and the time until the next event is exponentially
distributed with parameter $\lambda+ n\mu $, because it is a
competition between $n+1$ exponential random variables: $n$ with
parameter $\mu$ and one with parameter $\lambda$.

Considering a birth-and-death process that represents the queue evolution
of an M/M/$\infty$ queueing system, and its reversibility property,
the steady-state probabilities $\pi_i$ (for $i=0, 1, 2, \ldots $)
of having $i$ customers in the system satisfy
the following detailed balance (steady-state) equations:
 \\
$\pi_0 \lambda = \pi_1 \mu$\\
$\pi_1 \lambda = \pi_2 2\mu $\\
\ldots~\\
and in general:
\begin{equation}
\label{mminfsse}
\pi_n \lambda = \pi_{n+1} (n+1)\mu, ~{\rm for} ~n=0,~1,~2,~\ldots~.
\end{equation}

The sum of the steady-state probabilities must be equal to one,
so we again have the additional normalizing equation
\begin{equation}
\label{sumto1mminf}
\sum_{j=0}^{\infty} \pi_j =1.
\end{equation}

The infinitesimal generator of M/M/$\infty$ is given by

$Q_{i,i+1}=\lambda $ for $i$=0, 1, 2, 3, ~ \ldots\\
$Q_{i,i-1}=i\mu $ for $i$= 1, 2, 3, 4, ~ \ldots\\
$Q_{0,0}=-\lambda $ \\
$Q_{i,i}=-\lambda -i\mu $ for $i$=1, 2, 3,~\ldots~.

That is,

$${\bf Q} =
\begin{bmatrix}
    -\lambda & \lambda  & 0 & 0 & 0 & \ldots \\
    \mu & -(\lambda + \mu) & \lambda  & 0 & 0 & \ddots \\
    0 & 2\mu & -(\lambda + 2\mu) & \lambda  & 0 &  \ddots \\
    0 & 0 & 3\mu & -(\lambda + 3\mu) & \lambda  &  \ddots \\
    \vdots & \hspace{-8 mm} \ddots & \hspace{-16 mm}  \ddots & \hspace{-16 mm}  \ddots & \ddots &  \ddots
    \end{bmatrix}.
$$

\subsection{Solving the Steady-State Equations}

Using the $A$ notation we obtain

$\pi_1 = A \pi_0 $\\
$\pi_2 = A \pi_1/2 = A^2 \pi_0$/2\\
$\pi_3 = A \pi_2/3 = A^3 \pi_0$/(3!)

and in general:

\begin{equation}
\label{mminf}
\pi_n = \frac{A^n \pi_0}{n!}  ~{\rm for} ~n=0,~1,~2,~\ldots~ .
\end{equation}

To obtain $\pi_0$, we sum up both sides of Eq.\@ (\ref{mminf}), and
because the sum of the $\pi_n$s equals one, we obtain
\begin{equation}
\label{mminfsum}
1 = \sum_{n=0}^\infty \frac{A^n \pi_0}{n!}.
\end{equation}
By the definition of Poisson random variable, see Eq.\@
(\ref{poissonrv}), we obtain
\begin{equation}
1=\sum_{i=0}^{\infty} e^{-\lambda}\frac{\lambda^i} {i!}.
\end{equation}
Thus,
$$e^{\lambda} = \sum_{i=0}^{\infty} \frac{\lambda^i} {i!}$$
which is also the well-known Maclaurin series expansion of
$e^{\lambda}$. Therefore, Eq.\@ (\ref{mminfsum}) reduces to
\begin{equation}
1=\pi_0 e^{A},
\end{equation}
or
\begin{equation}
\pi_0=e^{-A}.
\end{equation}
Substituting the latter in Eq.\@ (\ref{mminf}), we obtain
\begin{equation}
\label{mminfin}
\pi_n = \frac{e^{-A} A^n }{n! } ~{\rm for} ~n=0,~1,~2,~\ldots~.
\end{equation}
By Eq.\@ (\ref{mminfin}) we observe that the distribution of the
number of busy channels (simultaneous calls or customers) in an
M/M/$\infty$ system is Poisson with parameter $A$.

\subsection{State Transition Diagram of M/M/$\infty$}

The state transition diagram of M/M/$\infty$ is similar to that of M/M/1 except that the rate downwards from state $n$ ($n=1,2,3, \ldots $) is
$n\mu$ rather than $\mu$ reflecting the fact at state $n$ there are $n$ servers serving the $n$ customers. The state transition diagram of M/M/$\infty$ is:

$$
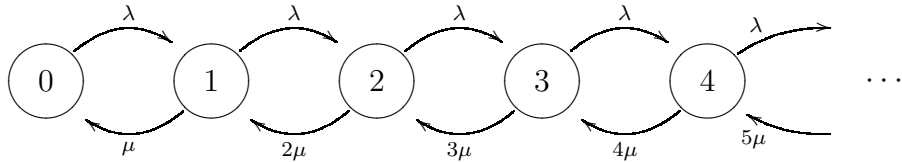
\xymatrix@C=32pt{
*++++[o][F]{0}\ar@/^{20pt}/[r]^{\lambda} &
*++++[o][F]{1}\ar@/^{20pt}/[l]^{\mu}\ar@/^{20pt}/[r]^{\lambda} &
*++++[o][F]{2}\ar@/^{20pt}/[l]^{2\mu}\ar@/^{20pt}/[r]^{\lambda} &
*++++[o][F]{3}\ar@/^{20pt}/[l]^{3\mu}\ar@/^{20pt}/[r]^{\lambda} &
*++++[o][F]{4}\ar@/^{20pt}/[l]^{4\mu}\ar@/^{6pt}/[r]!<-20pt,20pt>^{\lambda} &
*++++[]{\cdots}\ar@/^{6pt}/[]!<-20pt,-20pt>;[l]^{5\mu}
}
$$

We observe that the rates of
transitions between the states in the state transition diagram of M/M/$\infty$ are consistent with the rates in the
detailed balance equations of M/M/$\infty$ (\ref{mminfsse}).

\subsection{Insensitivity}
\label{insensitiveMMinf}

The above results for $\pi_i$, $i=0,1,2 ~\ldots$ and for the mean
number of busy servers are insensitive to the shape of the service
time (holding time) distribution
\cite{benes57,newell66,riordan51,Sevas57}; all we need to know is the mean
service time and the results are insensitive to higher
moments. In other words, the above results apply
to an M/G/$\infty$ system. This is important because it makes the
model far more robust which allows us to use its analytical results
for many applications where the service time is not exponential.
This insensitivity property is valid also for the M/G/$k$/$k$ system
\cite{gross08,ross96,Sevas57,Taylor2010}.

To explain the insensitivity property of M/G/$\infty$ with respect
to the mean occupancy, consider an arbitrarily long period of time
$L$ and also consider the queue size process that
represents the number of busy servers at any point in time between 0
and $L$. The average number of busy servers is obtained by the area
under the queue size process function divided by $L$. This area is
closely approximated by the number of arrivals during $L$ which is
$\lambda L$ times the mean holding (service) time of each arrival
($1/\mu$). Therefore, the mean number of busy servers, which is also
equal to the mean number of customers in the system (queue size), is
equal to $A = \lambda/\mu$ (notice that the $L$ is canceled out
here). Since all the traffic load enters the system ($A$) is also
the carried traffic load.

The words ``closely approximated'' are used here because there are
some customers that arrive before $L$ and receive service after $L$
and there are other customers that arrive before time 0 and are
still in the system after time 0. However, because we can choose $L$
to be arbitrarily long, their effect is negligible.

In the above discussion, we do not use moments higher than the
mean of the holding time, the mean number of busy servers (or mean
queue size) is insensitive to the shape of the holding-time
distribution and it is only sensitive to its mean.

Moreover, the distribution of the number of busy servers in
M/G/$\infty$ is also insensitive to the holding time distribution.
This can be explained as follows. We know that the arrivals follow a
Poisson process. Poisson process normally occurs in nature by having
a very large number of independent sources each of which generates
occasional events (arrivals) \cite{newell66} - for example, a large
population of customers making phone calls. These customers are
independent of each other. In M/G/$\infty$, each one of the arrivals
generated by these customers is able to find a server and its
arrival time, service time and departure time is independent of all
other arrivals (calls). Therefore, the event that a customer
occupies a server at an arbitrary point in time in steady-state is
also independent of the event that any other customer occupies a
server at that point in time. Therefore, the server occupancy events
are also due to many sources generating occasional events. This
explains the Poisson distribution of the server occupancy. From the
above discussion, we know that the mean number of servers is equal
to $A$, so we always have, in M/G/$\infty$,  in steady-state, a
Poisson distributed number of servers with parameter $A$ which is
independent of the shape of the service-time distribution.

\subsection{Applications}
\label{mminfmaccess}

\subsubsection{A Multi-access Model}
\label{multiaccessexample}

An interesting application of the M/M/$\infty$ system is the
following multi-access problem (see Problem 3.8 in \cite{BG92}).
Consider a stream of packets whose arrival times follow a
Poisson process with parameter $\lambda$. If the inter-arrival
times of any pair of packets (not necessarily a consecutive pair)
is less than the transmission time of the packet that arrived
earlier out of the two, these two packets are said to collide.
Assume that packets have independent exponentially distributed
transmission times with parameter $\mu$. What is the probability of no
collision?

Notice that a packet can collide with any one or more of the
packets that arrived before it. In other words, it is possible
that it may not collide with its immediate predecessor, but it may
collide with a packet that arrived earlier. However, if it does
not collide with its immediate successor, it will not collide with
any of the packets that arrive after the immediate successor.

Therefore, the probability that an arriving packet will not
collide on arrival can be obtained to be the probability of an
M/M/$\infty$ system to be empty, that is, $e^{-A}$. While the
probability that its immediate successor will not arrive during
its transmission time is $\mu/(\lambda +\mu)$. The product of the
two, namely $e^{-A}\mu/(\lambda +\mu)$, is the probability of no
collision.

For more information on multi-access modeling and analyzes see Chapter \ref{multiaccess}.

\subsubsection{Birth Rate Evaluation}

Another application of the M/M/$\infty$ system (or M/G/$\infty$
system) is to the following problem. Consider a city with
a population of 3,000,000, and assume that (1) there is no immigration
in and out of the city, (2) the birth rate $\lambda$ is constant (time-independent), and
(3)  life-time expectancy $\mu^{-1}$ in the city is
constant. It is also given that the average lifetime of
people in this city is 78 years.
How to compute the birth rate?

Using the M/M/$\infty$ model (or
actually, the M/G/$\infty$ as human lifetime is not exponentially
distributed) with $E[Q]=3,000,000$ and $\mu^{-1}=78$, realizing
that $E[Q]=A=\lambda/\mu$, we obtain, $\lambda = \mu E[Q] =
3,000,000/78=38461$ new births per year or 105 new births per day.

\subsubsection*{Homework \ref{secmminf}.\arabic{homework}}
\addtocounter{homework}{1} \addtocounter{tothomework}{1}

Consider an M/M/$\infty$ queue, with $\lambda=120$ [call/s], and $\mu=3$ [call/s]. Find the steady-state probability that there are 120 calls in the system. This should be done by a computer. Use ideas presented in Section \ref{poissonrvs}
to compute the required probabilities.

\newpage
\section{M/M/$k$/$k$ and Extensions}
\label{secerlang}

\setcounter{homework}{1} 
We begin this chapter with the M/M/$k$/$k$ queueing
system where the number of servers is $k$ assuming that the
arrival process is Poisson with parameter $\lambda$ and that each server
renders a service that is exponentially distributed with parameter
$\mu$. Later in the chapter, we extend the model to cases where the arrival process is non-Poisson.

In the M/M/$k$/$k$ model, as in the other M/M/$\ldots$ cases, we assume that the
service times are mutually independent and are independent of the arrival
process. We will now discuss Erlang's derivation of the loss
probability of an M/M/$k$/$k$ system that leads to the well known
Erlang's Loss Formula, also known as Erlang B Formula.

\subsection{M/M/$k$/$k$: Offered, Carried and Overflow Traffic}
\label{overflow}
The offered traffic under M/M/$k$/$k$ is the same as under M/M/$\infty$ it is equal to $$A=\lambda/\mu.$$
However, because some of the traffic is blocked, the offered traffic is not equal to the carried traffic.
To obtain the carried traffic given a certain blocking probability $P_b$, we recall that the carried traffic is equal to
the mean number of busy servers. To derive the latter we again invoke Little's formula. We notice that the arrival rate into the service system is
equal to $(1-P_b)\lambda$ and that the mean time each customer (or call) spends in the system is $1/\mu$. The mean queue size (which is also the mean
number of busy servers in the case of the M/M/$k$/$k$ queue) is obtained \cite{Heyman80} to be given by
\begin{equation} \label{meanqmmkk} E[Q]=\frac{(1-P_b)\lambda}{\mu}=(1-P_b)A. \end{equation}
Therefore the carried traffic is equal to $(1-P_b)A$. Notice that since $P_b>0$ in M/M/$k$/$k$, the carried traffic here is lower than the
corresponding carried traffic for M/M/$\infty$, which is equal to $A$.




The {\em overflow traffic} (in the context of M/M/$k$/$k$ it is also called: {\it lost traffic}) is defined as the difference between the two. Namely,
$${\rm overflow ~traffic = offered ~ traffic - carried ~traffic}.$$
Therefore, for M/M/$k$/$k$, the overflow traffic is
$$A-(1-P_b)A=P_bA.$$


\subsection{The Steady-State Equations and Their Solution}
\label{dimErlang}
The steady-state equations for M/M/$k$/$k$ are the same as the first $k$
steady-state equations for M/M/$\infty$.

As for M/M/$\infty$, the queue-size process of an M/M/$k$/$k$
is a birth-and-death process where queue-size transitions can only occur
upwards by one or downwards by one. Therefore, the M/M/$k$/$k$ queue-size process is also reversible, which means that solving its
detailed balance equations and the normalizing equation yields
the steady-state probabilities of the queue size.

The difference between the two systems is
that the queue size of M/M/$k$/$k$
can never exceed $k$ while for M/M/$\infty$ it is unlimited.
In other words, the M/M/$k$/$k$ system is a truncated version of the
M/M/$\infty$ with losses occur at state $k$, but the ratio
$\pi_j/\pi_0$ for any $0 \leq j \leq k$ is the same in both
M/M/$k$/$k$ and M/M/$\infty$ systems.

Another difference between the two is associated with the physical interpretation of
reversibility. Although both systems are reversible, while in  M/M/$\infty$ the physical
interpretation is that in the reversed process, the input point process is the point process of the call completion times in reverse, in M/M/$k$/$k$ the reversed process is the superposition of call completion times in reverse and call blocking times in reverse. Notice that both types of points in time represent customer departures. At the call completion times customers depart because they completed their service, and at call blocking times, they depart immediately upon arrivals because no server is available to serve them.

The reversibility property of M/M/$k$/$k$ and the truncation at $k$
  imply that the detailed balance
equations for M/M/$k$/$k$ are the same as the first $k$
detailed balance (steady-state) equations for M/M/$\infty$. Namely,
these balance equations are:

\begin{equation}
\label{mmkksse}
\pi_n \lambda = \pi_{n+1} (n+1)\mu, ~{\rm for} ~n=0,~1,~2,~\ldots~k-1.
\end{equation}

The infinitesimal generator of M/M/$k$/$k$ is given by

$Q_{i,i+1}=\lambda $ for $i=0, 1, 2, 3, ~ \ldots, ~ k-1$\\
$Q_{i,i-1}=i\mu $ for $i= 1, 2, 3, 4, ~ \ldots, ~ k$\\
$Q_{0,0}=-\lambda $ \\
$Q_{i,i}=-\lambda -i\mu $ for $i=1, 2, 3,~ \ldots, ~ k-1$\\
$Q_{k,k}=-k\mu$.

That is,

$${\bf Q} =\begin{bmatrix}
-\lambda & \lambda & 0&0&0&0&0& \ldots &0 \\
\mu & -(\mu+\lambda) & \lambda & 0&0&0&0& \ldots &0 \\
0&2\mu & -(2\mu+\lambda) & \lambda  & 0&0&0& \ldots &0 \\
0&0&3\mu & -(3\mu+\lambda) & \lambda  & 0&0& \ldots &0 \\
0&0&0&\ddots&\ddots&\ddots  & 0 & \ldots &0 \\
\vdots&\vdots&\vdots&\ddots&\ddots&\ddots&~~~~~~\ddots &~~~~~~~ \ddots &\vdots\\
0&0&0&0&0&0&(k-1)\mu & -[(k-1)\mu+\lambda] & \lambda \\
0&0&0&0&0&0&0&k\mu & -k\mu
\end{bmatrix}
$$

The balance equations can also be described by the following state transition diagram of M/M/$k$/$k$:

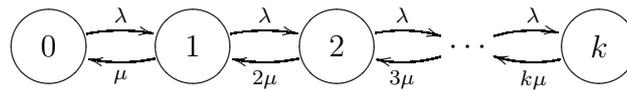
$$\xymatrix{
*++++[o][F]{0}\ar@/^{.5pc}/[r]^{\lambda} &
*++++[o][F]{1}\ar@/^{.5pc}/[l]^{\mu}\ar@/^{.5pc}/[r]^{\lambda} &
*++++[o][F]{2}\ar@/^{.5pc}/[l]^{2\mu}\ar@/^{.5pc}/[r]^{\lambda} &
{\cdots}\ar@/^{.5pc}/[l]^{3\mu}\ar@/^{.5pc}/[r]^{\lambda} &
*++++[o][F]{k}\ar@/^{.5pc}/[l]^{k\mu}}$$

The sum of the steady-state probabilities must be equal to one,
so we again have the additional normalizing equation
\begin{equation}
\label{sumto1mmkk}
\sum_{j=0}^{k} \pi_j =1.
\end{equation}

Accordingly, we obtain for
M/M/$k$/$k$:

\begin{equation}
\label{mmkk} \pi_n = \frac{A^n \pi_0}{n!}  ~{\rm for}
~n=0,~1,~2,~\ldots, k.
\end{equation}
To obtain $\pi_0$, we again sum up both sides of the latter. This
leads to
\begin{equation}
\label{mmkkp0} \pi_0 = \frac{1}{\sum_{n=0}^k \frac{A^n}{n!} }.
\end{equation}
Substituting Eq.\@ (\ref{mmkkp0}) in Eq.\@ (\ref{mmkk}), we obtain
\begin{equation}
\label{mmkkfin} \pi_n = \frac{ \frac{A^n}{n!} } {\sum_{n=0}^k
\frac{A^n}{n!} }  ~{\rm for} ~n=0,~1,~2,~\ldots, k.
\end{equation}

The relationship between (\ref{mmkkfin}) and (\ref{mminf}) is
now clearer. Observing (\ref{mmkkfin}) that gives
the distribution of the number of customers in an M/M/$k$/$k$ model,
it is apparent that it is a truncated version of (\ref{mminf}).
Since (\ref{mminf}) is merely the Poisson distribution,
(\ref{mmkkfin}) is the truncated Poisson distribution. Accordingly,
to obtain (\ref{mmkkfin}), we can simply consider (\ref{mminf}), and
firstly set $\pi_j=0$ for all $\pi_j$ with $j>k$. Then, for $0\leq
j\leq k$ we set the $\pi_j$ for the M/M/$k$/$k$ values by dividing
the $\pi_j$ values of (\ref{mminf}) by the sum ${\sum_{j=0}^k \pi_j
}$ of the $\pi_j$ values in the M/M/$\infty$ model. This is
equivalent to considering the M/M/$\infty$ model and deriving the
conditional probability of the process being in state $j$ for
$j=0,1,2, \ldots, k$, conditioning on the process being within the
states $j=0,1,2, \ldots, k$. This conditional probability is exactly
the steady-state probabilities $\pi_j$ of the M/M/$k$/$k$ model.

The most important quantity out of the values obtained by Eq.\@
(\ref{mmkkfin}) is $\pi_k$. It is the probability that all $k$
circuits are busy, so it is the proportion of time that no new calls
can enter the system; namely, they are blocked. It is therefore
called {\it time congestion}. The quantity $\pi_k$ for an
M/M/$k$/$k$ system loaded by offered traffic $A$ is usually denoted
by $E_k(A)$ and is given by:
\begin{equation}
\label{eka} E_k(A) = \frac { \frac{A^k}{k!} } {\sum_{n=0}^k
\frac{A^n}{n!} }.
\end{equation}
Eq.\@ (\ref{eka}) is known as Erlang's loss Formula or Erlang B
Formula, published first by A. K. Erlang in 1917 \cite{Erlang17}.

Due to the special properties of the Poisson process, in addition of
being the proportion of time during which the calls are blocked, $E_k(A)$
also gives the proportion of calls blocked due to congestion;
namely, it is the {\em call congestion} or {\em blocking probability}.
A simple way to explain that for an M/M/$k$/$k$ system, the call congestion (blocking probability) is equal to the time congestion is the following.
Let $L$ be an arbitrarily long period of time. The proportion of time during $L$
when all servers are busy, and every arrival is blocked is $\pi_k=E_k(A)$, so the time during $L$ when new arrivals are blocked is $\pi_k L$.
The mean number of blocked arrivals during $L$ is, therefore, equal to  $\lambda \pi_k L$.
The mean total number of arrivals during $L$ is $\lambda L$.
The blocking probability (call congestion) $P_b$ is the ratio between the two. Therefore:
$$P_b = \frac{\lambda \pi_k L}{\lambda L}=\pi_k=E_k(A).$$
Eq.\@ (\ref{eka}) has many applications in telecommunications
network design. Given its importance, it is necessary to be able to
compute Eq.\@ (\ref{eka}) quickly and exactly for large values of
$k$. This will enable us to answer a {\bf dimensioning problem}, for example: how many circuits are required so that the blocking probability is
no more than 1\% given offered traffic of $A=1000$?

\subsection{Recursion and Jagerman Formula}

Observing Eq.\@ (\ref{eka}), we notice the factorial terms which may
hinder such computation for a large $k$. We shall now present an
analysis which leads to a recursive relation between $E_m(A)$ and
$E_{m-1}(A)$ that gives rise to a simple and scalable algorithm for
the blocking probability. By Eq.\@ (\ref{eka}), we obtain
\begin{equation}
\label{recur_deriv}
\frac{E_m(A)}{ E_{m-1}(A)}
= \frac{  \frac { \frac{A^m}{m!} } {\sum_{j=0}^m
\frac{A^j}{j!} }  }{\frac { \frac{A^{m-1}}{(m-1)!} } {\sum_{j=0}^{m-1}
\frac{A^j}{j!} } }
= \frac{  \frac { \frac{A^m}{m!} } {\sum_{j=0}^m
\frac{A^j}{j!} }  }{\frac { \frac{A^{m-1}}{(m-1)!} } {\sum_{j=0}^{m}
\frac{A^j}{j!} - \frac{A^m}{m!}} }
= \frac {A}{m}(1-E_m(A) ).
\end{equation}
Isolating $E_m(A)$, this leads to
\begin{equation}
\label{recur} E_m(A) = \frac {A E_{m-1}(A) }{m  +A E_{m-1}(A) }
~{\rm for} ~m=~1,~2,~\ldots, k.
\end{equation}

\subsubsection*{Homework \ref{secerlang}.\arabic{homework}}
\addtocounter{homework}{1} \addtocounter{tothomework}{1} Complete
all the details in the derivation of Eq.\@ (\ref{recur}). $~~~\Box$

When $m=0$, there are no servers (circuits) available, and therefore all
customers (calls) are blocked, namely,
\begin{equation}
\label{recur0}
E_0(A) = 1.
\end{equation}
The above two equations give rise to a simple recursive algorithm
by which the blocking probability can be calculated for a large
$k$. An even more computationally stable way to compute $E_m(A)$
for large values of $A$ and $m$ is to use the inverse
\cite{Iver15}
\begin{equation}
\label{invers} I_m(A) = \frac{1}{E_m(A)}
\end{equation}
and the recursion
\begin{equation}
\label{Irecur} I_m(A) = 1+\frac {m}{A} I_{m-1}(A) ~{\rm for}
~m=~1,~2,~\ldots, k.
\end{equation}
with the initial condition $I_0(A)=1$.

A useful formula for $I_m(A)$ due to Jagerman \cite{Jag74} is:

\begin{equation}
\label{jag} I_m(A) = A\int_{0}^\infty e^{-Ay}(1+y)^m dy.
\end{equation}

\subsubsection*{Homework \ref{secerlang}.\arabic{homework}}
\addtocounter{homework}{1} \addtocounter{tothomework}{1}
Long long ago in a far-away land, John, an employee of a telecommunication provider company, was asked to derive the blocking probability of a switching system loaded by a Poisson arrival process of calls  where the offered load is given by $A = 180$.
These calls are served by $k = 200$ circuits.

The objective was to meet a requirement of no more than 1\% blocking probability. The company has been operating with $k = 200$ circuits for some time and there was a concern that the blocking probability exceeds the 1\% limit.

John was working for a while on the calculation of this blocking probability, but when he was very close to a solution, he won the lottery, resigned, and went to the Bahamas. Mary, another employee of the same company, was given the task of solving the problem. She found some of John's papers where it was revealed that for an M/M/196/196 model and $A = 180$, the blocking probability is approximately 0.016. Mary completed the solution in a few minutes. Assuming that John's calculations were correct, what is the solution for the blocking probability of M/M/200/200 with $A=180$? If the blocking probability in the case of $k = 200$  is more than 1\%, what is the smallest number of circuits that should be added to meet the 1\% requirement?

\subsubsection*{Solution}

Using the Erlang B recursion (\ref{recur})
and knowing that $E_{196}(180) = 0.016$, we obtain

$E_{197}(180) \approx 0.0145$\\
$E_{198}(180) \approx 0.013$\\
$E_{199}(180) \approx 0.0116$\\
$E_{200}(180) \approx 0.0103$.

One more circuit should be added to achieve:

$E_{201}(180) \approx 0.0092$
$~~~\Box$

\subsection{The Special Case: M/M/1/1}

\subsubsection*{Homework \ref{secerlang}.\arabic{homework}}
\addtocounter{homework}{1} \addtocounter{tothomework}{1}  Derive a formula for the blocking
probability of M/M/1/1 in four ways: (1) by Erlang B Formula
(\ref{eka}), (2) by the recursion (\ref{recur}), (3) by the
recursion (\ref{Irecur}), and (4) by Jagerman Formula (\ref{jag}).
$~~~\Box$

The reader may observe a fifth direct way to obtain a formula for the blocking probability of
M/M/1/1 using Little's formula. The M/M/1/1 system can have at most one customer in it. Therefore, its mean queue size is given by
$E[Q] = 0\pi_0 + 1\pi_1=\pi_1$ which is also its blocking probability. Noticing also that the arrival rate into the system
(made only of successful arrivals) is equal to $\lambda(1-E[Q])$, the mean time a customer stays in the system is $1/\mu$, and revoking Little's formula, we obtain
\begin{equation}
\label{mm11litt}
\frac{\lambda(1-E[Q])}{\mu} = E[Q].
\end{equation}
Isolating $E[Q]$, the blocking probability is given by
\begin{equation}
\label{mm11bp}
\pi_1 = E[Q] = \frac{A}{1+A}.
\end{equation}

\subsection{Lower bound of $E_k(A)$}
By (\ref{meanqmmkk}), and since $E[Q]$  satisfies  $E[Q]< k$ for a finite $k$, we obtain
$$ A(1-E_k(A)) < k, $$
or following simple algebraic operations
\begin{equation}
\label{eka_bound1} E_k(A) > 1 - \frac{k}{A}. \end{equation}
Clearly,  the latter is relevant only when  $k<A$; otherwise, a tighter bound is
obtained by $E_k(A)\geq 0$.  Thus, for the full range of $A\geq 0$ and $k=1,2,3\ldots $ values, $E_k(A)$  satisfies \cite{Whitt02}
\begin{equation}
\label{eka_bound2}
E_k(A) \geq {\rm max} \left(0,1 - \frac{k}{A} \right).
\end{equation}

\subsection{Monotonicity of $E_k(A)$}
 Intuitively, the blocking probability  $E_k(A)$ is expected to increase as the offered traffic $A$ increases and to decrease as the number of servers $k$ increases.

 To show that $E_k(A)$ is monotonic in $A$, we use induction on $k$ and the Erlang B recursion (\ref{recur}) \cite{Whitt02}.

  \subsubsection*{Homework \ref{secerlang}.\arabic{homework}}
\addtocounter{homework}{1}  \addtocounter{tothomework}{1} Complete all the details in the proof that $E_k(A)$ is monotonic in $A$. $~~~\Box$

  Then, to show that $E_k(A)$ is monotonic in $k$, recall that by  (\ref{recur_deriv}), for  $k=1,2,3, \ldots $ we have

$$
\frac{E_k(A)}{ E_{k-1}(A)}
= \frac {A}{k}[1-E_k(A) ].
$$
Then, by   (\ref{eka_bound1})  we have

$$1-E_k(A) < \frac {k}{A}.$$

Thus,

$$\frac{E_k(A)}{ E_{k-1}(A)} < \frac {A}{k} \frac {k}{A} = 1,$$

so
$${E_k(A)} < { E_{k-1}(A)}.$$

\subsection{The Limit $k \rightarrow \infty$ with $A/k$ Fixed}
As traffic and capacity (number of servers) increase, there is an
interest in understanding the blocking behavior in the limit $k$ and
$A$ both approach infinity with $A/k$ Fixed.
\subsubsection*{Homework \ref{secerlang}.\arabic{homework}}
\addtocounter{homework}{1} \addtocounter{tothomework}{1} Prove that
the blocking probability approaches zero for the case $A/k \leq 1$
and that it approaches $1 - k/A$ in the case $A/k > 1$.
\subsubsection*{Guide (by Guo Jun based on \cite{syski86})}
By (\ref{jag}), we have.
\begin{displaymath}
\frac{1}{ \mathbf{E}(A,k)}=\int_0^\infty e^{-t} \left(1+\frac{t}{A}
\right)^k\; dt.
\end{displaymath}
Consider the case where $k$ increases in such a way that $A/k$ is
constant. Then,
\begin{eqnarray*}
\lim_{k \to \infty} \frac{1}{ \mathbf{E}(A,k)}
&= \lim_{k \to
\infty} \int_0^\infty e^{-t} \left(1+\frac{t}{A}
\right)^k\; dt \\
&= \int_0^\infty e^{-t} \cdot e^{\frac{t}{A/k}}\; dt \\
&= \int_0^\infty e^{-(1-\frac{1}{A/k})t} \; dt.
\end{eqnarray*}
Then, observe that
\begin{equation}
\lim_{k \to \infty} \frac{1}{ \mathbf{E}(A,k)}
=\left\{\begin{array}{ll}
\infty & \mbox{if $A/k\leq 1$}\\
\frac{1}{1-\frac{1}{A/k}} & \mbox{if $A/k>1.$}
\end{array}
\right.
\end{equation}
 And the desired result follows.
 $~~~\Box$

 \subsubsection*{Homework \ref{secerlang}.\arabic{homework} (Jiongze Chen)}
\addtocounter{homework}{1} \addtocounter{tothomework}{1} Provide an
alternative  proof to the results of the previous homework using the
Erlang B recursion.
\subsubsection*{Guide} Set $a=A/k$. Notice that in the limit $E_k(ak) \cong
E_{k+1}(a(k+1))\cong E_{k}(A)$ and provide a proof for the cases
$k<A$ and $k=A$. Then, for the case $k>A$, first show that the
blocking probability decreases as $k$ increases for a fixed $A$
(using the Erlang B recursion), and then argue that if the blocking
probability already approaches zero for $A=k$, it must also approach
zero for $k>A$.
 $~~~\Box$

\subsubsection*{Homework \ref{secerlang}.\arabic{homework}}
\addtocounter{homework}{1} \addtocounter{tothomework}{1} Provide
an intuitive explanation of the results of the previous homework.
\subsubsection*{Guide}
Due to the insensitivity property, M/M/$k$/$k$ and M/D/$k$/$k$
experience the same blocking probability if the offered traffic in
both systems is the same. Observe that since the arrivals follow a
Poisson process, the variance is equal to the mean. Also, notice that
as the arrival rate $\lambda$ increases, the Poisson process
approaches a Gaussian process. Having the variance equal to the
mean, the standard deviation becomes negligible relative to the mean
for a very large $\lambda$. With negligible variability, M/D/$k$/$k$
behaves like D/D/$k$/$k$, and the results follow.
 $~~~\Box$

\subsection{M/M/$k$/$k$: Dimensioning and Utilization}
\label{dim2Erlang}
Taking advantage of the monotonicity of Erlang formula, we can
also solve the dimensioning problem. We simply keep incrementing
the number of circuits and calculate, in each case, the blocking
probability. When the desired blocking probability (e.g., 1\%) is
reached, we have our answer.

\subsubsection*{Homework \ref{secerlang}.\arabic{homework}}
\addtocounter{homework}{1} \addtocounter{tothomework}{1} In relation to the above discussed monotonicity of $E_k(A)$, prove that if $A>A'$, then
$E_n(A)>E_n(A')$ for $n=1,2,3 \ldots $. \\ {\bf Hint:} Consider the Erlang B recursion and use induction on $n$.$~~~\Box$

We have already derived the mean number of busy circuits in an M/M/$k$/$k$ system fed by $A$
erlangs in (\ref{meanqmmkk}) using Little's formula. Substituting $\pi_k$ for $P_b$ in (\ref{meanqmmkk}), we obtain
$$E[Q]=(1-\pi_k)A.$$
Note that it can also be computed by the weighted sum
$$E[Q]=\sum_{i=0}^k i\pi_i.$$

Accordingly, the utilization of an M/M/$k$/$k$ system is given by
\begin{equation}
\label{utilmmkk}
\hat{U}=\frac{(1-\pi_k)A}{k}. \end{equation}

\subsubsection*{Homework \ref{secerlang}.\arabic{homework}}
\addtocounter{homework}{1} \addtocounter{tothomework}{1} Prove that $\sum_{i=0}^k
i\pi_i=(1-\pi_k)A.$ $~~~\Box$

In the following Table, we present the minimal values of $k$
obtained for various values of $A$ such that the blocking
probability is no more than 1\%, and the utilization obtained in
each case. It is clearly observed that the utilization increases
with the traffic.

Note that normally it is impossible to find $k$ that gives exactly 1\% blocking probability, so we conservatively choose $k$ such that
$E_k(A) <1\%$, but we must make sure that $E_{k-1}(A) > 1\%$. Namely, if we further reduce $k$, the blocking probability requirement to be no more than 1\% is violated.

\begin{center}
\renewcommand{\arraystretch}{1.4}
\begin{tabular}{|c|c|c|c| }  \hline
$A$ & $k$ & $E_k(A)$ & Utilization  \\  \hline 20 & 30 & 0.0085
 &
66.10\%
 \\  \hline
 100 & 117 & 0.0098  & 84.63\%
 \\
\hline 500 & 527 & 0.0095 & 93.97\%  \\
\hline 1000  &  1029   & 0.0099 & 96.22\%  \\
\hline 5000   & 5010  &  0.0100 & 98.81\%  \\
\hline 10000  & 9970  &  0.0099 & 99.30\% \\
\hline
\end{tabular}
\end{center}

\subsubsection*{Homework \ref{secerlang}.\arabic{homework}}
\addtocounter{homework}{1}  \addtocounter{tothomework}{1} Reproduce the results of the above Table. $~~~\Box$

We also notice that for the case of $A=10,000$ erlangs, to maintain
no more than 1\% blocking, $k$ value less than $A$ is required.
Notice, however, that the carried traffic is not $A$ but
$A(1-E_k(A))$. This means that for $A \geq 10,000$, dimensioning
simply by $k=A$ will mean no more than 1\% blocking and no less than
99\% utilization - not bad for such a simple rule of thumb! This
also implies that if the system capacity is much larger than
individual service requirement, very high efficiency (utilization)
can be achieved without a significant compromise on the quality of
service.  Let us now further examine the case $k=A$.

\subsection{M/M/$k$/$k$ under Critical Loading}
\label{critload}

A system where the offered traffic load is equal to the system capacity is called {\it critically loaded} \cite{Hunt89}.
Accordingly, in a critically loaded Erlang B System, we have $k=A$.
From the table below, it is clear that if we maintain  $k=A$ and we increase them both, the blocking probability decreases, the utilization increases, and interestingly, the product $E_k(A)\sqrt{A}$ approaches a constant, which we denote $\tilde{C}$, that does not depend on $A$ or $k$. This implies that in the limit, the blocking probability decays at the rate of $1/\sqrt{k}$. That is, for a critically loaded Erlang B system, we obtain
\begin{equation}
\label{limiterlang}
\lim_{k \rightarrow \infty} E_k(A)= \frac{\tilde{C}}{\sqrt{k}}.
\end{equation}

\begin{center}
\renewcommand{\arraystretch}{1.4}
\begin{tabular}{|c|c|c|c|c| }  \hline
$A$ & $k$ & $E_k(A)$ & Utilization  & $E_k(A)\sqrt{A}$ \\
\hline 10 & 10 & 0.215 &    78.5 \%     & 0.679
 \\  \hline
 50 & 50 & 0.105    & 89.5\% & 0.741
  \\
 \hline 100 & 100 & 0.076   & 92.4\% &  0.757
\\
\hline 500 & 500 & 0.035 & 96.5\% & 0.779
 \\
\hline 1000  &  1000   &
0.025   & 97.5\%    & 0.785
\\
\hline 5000   & 5000  &
0.011   &  98.9\% &     0.792
 \\
\hline 10000  & 10000  &
0.008   &  99.2\%   &  0.79365
\\
\hline 20000 &  20000 & 0.00562 &   99.438\%     & 0.79489

\\
\hline 30000 &  30000 & 0.00459 &   99.541\%     & 0.79544
\\
\hline 40000 &  40000 &
0.00398 & 99.602\%  & 0.79576
\\
\hline 50000 &  50000 &
0.00356 & 99.644\%  & 0.79599
\\
\hline
\end{tabular}
\end{center}
To explain the low blocking probability in a critically loaded large system, we refer back to our homework problem related to an Erlang B system with large capacity where the ratio $A/k$ is maintained
constant. In such a case, the standard deviation of the traffic is
very small relative to the mean, so the traffic behaves close to
deterministic. If 100 liters per second of water are offered,  at a
constant rate, to a pipe that has a capacity of 100 liters per second,
then the pipe can handle the offered load with very small losses.

\subsection{Insensitivity and Many Classes of Customers}
\label{manyclasses}

We have discussed in Section \ref{insensitiveMMinf}, the
distribution, and the mean of the number of busy servers is
insensitive to the shape of the service time distribution (although it is still sensitive to the mean of the service time) in the cases
of M/G/$\infty$ and M/G/$k$/$k$. For M/G/$k$/$k$, also the blocking
probability is insensitive to the shape of the service time
distribution \cite{gross08,ross96,Sevas57}.

However, we must make it very clear that the insensitivity property
does not extend to the arrival process. We still require a Poisson arrival
process for the Erlang B formula to apply. If we have a more burtsy
arrival process (e.g., Poisson batch arrivals), we will have more
losses than predicted by the Erlang B formula, and if we have a smoother
arrival process than Poisson, we will have fewer losses than
predicted by the Erlang B formula. To demonstrate it, let us compare an
M/M/1/1 system with a D/D/1/1 system. Suppose that each of these two
systems is fed by $A$ erlangs and that $A<1$.

Arrivals into the D/D/1/1 system with $A<1$ will never experience losses
because the inter-arrivals are longer than the service times, so
the service of a customer is always completed before the arrival
of the next customer. Accordingly, by Little's formula: $E[Q]=A$,
and since $E[Q]=0\times\pi_0 + 1\times\pi_1$, we have that
$\pi_1=A$ and $\pi_0=1-A$. In this case, the blocking probability
$P_b$ is equal to zero and not to $\pi_1$. The utilization will be given by $\hat{U}=\pi_1=A$.

By contrast, for the M/M/1/1 system, $P_b=E_1(A)=E[Q]=\pi_1=A/(1+A)$,
so $\pi_0=1-\pi_1=1/(1+A)$. To obtain the utilization we can
either use the fact that it is the proportion of time our single
server is busy, namely it is equal to $\pi_1=A/(1+A)$, or we can
use Equation (\ref{utilmmkk}) for $\hat{U}$ in M/M/$k$/$k$ system and obtain
\begin{equation}
\hat{U}=(1-\pi_1)A=[1-A/(1+A)]A=A/(1+A).
\end{equation}

This comparison is summarized in the following table:
\begin{center}
\renewcommand{\arraystretch}{1.4}
\begin{tabular}{|c|c|c| }  \hline
 & M/M/1/1 & D/D/1/1   \\  \hline
$\pi_0$ & $1/(1+A)$ & $1-A$  \\  \hline $\pi_1$  & $A/(1+A)$ &  $A$ \\
\hline
$\hat{U}$ & $A/(1+A)$ & $A$  \\
\hline
$P_b$ & $A/(1+A)$ & 0  \\
\hline
$E[Q]$ & $A/(1+A)$ & $A$ \\
\hline
\end{tabular}
\end{center}

Clearly, the steady-state equations (\ref{mmkk}) will not apply to
a D/D/1/1 system.

We have already mentioned that for M/G/$k$/$k$ the distribution of
the number of busy servers and therefore also the blocking
probability is insensitive to the shape of the service time
distribution (moments higher than the first). All we need is to know
that the arrival process is Poisson and the ratio of the arrival
rate to the service rate of a single server, and we can obtain the
blocking probability using the Erlang B formula. Let us now consider
the following problem.

Consider two classes of customers (packets). Class $i$ customers, $i=1,2$
arrive following an independent Poisson process at a rate of $\lambda_i$ each of which requires independent and exponentially
distributed service time with parameter $\mu_i$, for $i=1,2$. There are
$k$ independent servers without waiting room (without additional buffer). The
aim is to derive the blocking probability of customers of each type.

The combined arrival process of all the customers is a Poisson
process with parameter $\lambda=\lambda_1+\lambda_2$. Because the
probability of an arbitrary customer belonging to the first class is
$$p=\frac{\lambda_1}{\lambda_1+\lambda_2}=\frac{\lambda_1}{\lambda},$$
the service time of an arbitrary customer has hyperexponential
distribution because with probability $p$ it is exponentially
distributed with parameter $\mu_1$, and with probability $1-p$, it
is exponentially distributed with parameter $\mu_2$.

Therefore, by the Law of iterated expectation, the mean service time (holding time) is given by
$$ E[S] = \frac{p}{\mu_1} + \frac{1-p}{\mu_2}$$
so $A=\lambda E[S]$, and Erlang B can then be used to obtain the blocking probability $E_k(A)$.

Furthermore, let $$A_i=\frac{\lambda_i}{\mu_i} ~~~~ i=1,2$$
and observe that
$$E[S] = \left(\frac{\lambda_1}{\lambda}\right) \left(\frac{1}{\mu_1}\right) + \left(\frac{\lambda_2}{\lambda}\right) \left(\frac{1}{\mu_2}\right) = \frac{A_1+A_2}{\lambda}.$$
Then, $$A=\lambda E[S]=A_1 + A_2.$$

\subsubsection*{Homework \ref{secerlang}.\arabic{homework} \cite{BG92}}
\addtocounter{homework}{1} \addtocounter{tothomework}{1}
Consider $n$ classes of calls that are offered to a loss system with $k$ servers and without a waiting room. The arrival processes of the different classes follow independent Poisson processes.
The arrival rate of calls of class $i$ is $\lambda_i$ and their mean holding time is $1/\mu_i$, $i=1,2, \ldots, n$.
Find the blocking probability.
\subsubsection*{Guide} $$\lambda = \sum_{i=1}^n \lambda_i $$
$$E[S] = \sum_{i=1}^n \left(\frac{\lambda_i}{\lambda}\right) \left(\frac{1}{\mu_i}\right)$$
$$A=\lambda E[S]=\sum_{i=1}^n A_i$$
where
$$A_i=\frac{\lambda_i}{\mu_i} ~~~~ i=1,2, \ldots n.$$
Then, the blocking probability $E_k(A)$ is obtained by Erlang B formula.  $~~~\Box$

\subsubsection*{Homework \ref{secerlang}.\arabic{homework} }
\addtocounter{homework}{1} \addtocounter{tothomework}{1}
\begin{enumerate} \item Consider an M/M/$\infty$ queueing system
with the following twist. The servers are numbered 1, 2, $\ldots$
and an arriving customer always chooses the server numbered lowest
among all the free servers it finds. Find the proportion of time
that each of the servers is busy \cite{BG92}.

\subsubsection*{Guide} Notice that the input (arrival) rate into the system
comprises servers $n+1, n+2, n+3 \ldots$ is equal to $\lambda
E_n(A)$. Then, using Little's formula notice that the mean number of
busy servers among $n+1, n+2, n+2 \ldots$  is equal to $AE_n(A)$.
Repeat the procedure for the system comprises servers $n+2, n+3, n+4
\ldots$, you will observe that the mean number of busy servers in
this system is  $AE_{n+1}(A)$. Then, considering the difference
between these two mean values, you will find that the mean number
of busy servers in a system comprises only of the $n+1$ server is
$$A[E_n(A)-E_{n+1}(A)].$$ Recall that the mean queue size (mean
number of busy servers) of a system that comprises only the single server is (probability of the server being busy) times 1 + (probability of the server being idle) times 0, which is equal to the probability that the
server is busy, we obtain that $A[E_n(A)-E_{n+1}(A)]$ is the probability that the server is busy.

An alternative way to look at this problem is the following.
Consider the system made only of the $n+1$ server. The offered
traffic into this system is $AE_{n}(A)$, and the rejected traffic of
this system is $AE_{n+1}(A)$. Therefore, the carried traffic of this
system is $A[E_n(A)-E_{n+1}(A)]$. This means that the arrival rate
of customers that actually enters this single-server system is

$$\lambda_{{\rm enters}(n+1)} = \lambda[E_n(A)-E_{n+1}(A)]$$
and since the mean time spent in this system is $1/\mu$, we have
that the mean queue size in this single server system is
$$\lambda_{{\rm enters}(n+1)} \frac{1}{\mu}= A[E_n(A)-E_{n+1}(A)]$$
which is the carried load.
 Based on the arguments above it is equal to the proportion of time
 the $n+1$ server is busy.
\item Show that if the number of servers is finite $k$, the  proportion of time that server $n+1$ is busy is
$$A\left(1-\frac{E_k(A)}{E_n(A)}\right) E_n(A)-A\left(1-\frac{E_k(A)}{E_{n+1}(A)}\right)E_{n+1}(A)
= A[E_n(A)-E_{n+1}(A)]$$
and provide intuitive arguments to why the result is the same as in the infinite server case.
\item Verify the results by discrete-event and Markov-chain simulations. $~~~\Box$ \end{enumerate}

\subsubsection*{Homework \ref{secerlang}.\arabic{homework} }
\addtocounter{homework}{1} \addtocounter{tothomework}{1}
Consider an M/M/$k$/$k$  queue with a given arrival rate $\lambda$ and mean holding
time $1/\mu$. Let $A= \lambda/\mu$. Let $E_k(A)$ be the blocking probability. An
independent Poisson inspector inspects the M/M/$k$/$k$ queue at times $t_1, t_2,
t_3, \ldots$~. What is the probability that the first arrival after an inspection
is blocked?

{\bf Answer:} $$\frac{E_k(A)\lambda}{k\mu +\lambda}. ~~~\Box $$

\subsubsection*{Homework \ref{secerlang}.\arabic{homework}}
\addtocounter{homework}{1} \addtocounter{tothomework}{1}
Bursts of data of exponential lengths with mean $1/\mu$ that arrive according to a Poisson process are transmitted through a bufferless optical switch. All arriving bursts compete for $k$ wavelength channels at a particular output trunk of the switch. If a burst arrives and all $k$ wavelength channels are busy, the burst is dumped at the wavelength bit rate. While it is being dumped,
if one of the wavelength channels becomes free, the remaining portion of the burst is successfully transmitted through the wavelength channel.
\begin{enumerate}
\item Show that the mean loss rate of data $E[Loss]$ is given by $$E[Loss]=1P(X=k+1)+2P(X=k+2) + 3P(X=k+3)+\ldots$$ where $X$ is a Poisson random variable with parameter $A=\lambda/\mu$.
    \item Prove that $$E[Loss]=\frac{A  \gamma(k,A)}{\Gamma(k)} - \frac{k \gamma(k+1,A)}{\Gamma(k+1)}$$
    where $\Gamma(k)$ is the Gamma function and $\gamma(k,A)$ is the lower incomplete Gamma function.  \end{enumerate}
\subsubsection*{Background information and guide}
    The Gamma function is defined by
\begin{equation}
\Gamma(a)=\int_{0}^{\infty} t^{a-1} e^{-t} dt.
\end{equation}
The lower incomplete Gamma function is defined by
\begin{equation}
\gamma(a,x)=\int_{0}^{x} t^{a-1} e^{-t} dt.
\end{equation}
The upper incomplete Gamma function is defined by
\begin{equation}
\Gamma(a,x)=\int_{x}^{\infty} t^{a-1} e^{-t} dt.
\end{equation}
Accordingly,
$$\gamma(a,x) + \Gamma(a,x) = \Gamma(a).$$
For an integer $k$, we have
\begin{equation}
\Gamma(k)=(k-1)!
\end{equation}
\begin{equation}
\Gamma(k,x)=(k-1)!e^{-x}\sum_{m=0}^{k-1} \frac{x^m}{m!}.
\end{equation}
Therefore,
\begin{equation}
e^{-A}\sum_{m=0}^{k-1} \frac{A^i}{i!} = \frac{\Gamma(k,A)}{\Gamma(k)}
\end{equation}
so
\begin{equation}
1- e^{-A}\sum_{m=0}^{k-1} \frac{A^i}{i!} = 1-\frac{\Gamma(k,A)}{\Gamma(k)}=\frac{\Gamma(k)-\Gamma(k,A)}{\Gamma(k)}=\frac{\gamma(k,A)}{\Gamma(k)}.
\end{equation}
Now notice that
\begin{eqnarray*}
E[Loss] &
= &   1 \times P(X=k+1)+2 \times P(X=k+2) + 3 \times P(X=k+3)+\ldots\\ &
= &   \sum_{i=k+1}^\infty (i-k)A^i \frac{e^{-A}}{i!}\\ &
= & A\sum_{i=k+1}^\infty A^{i-1} \frac{e^{-A}}{(i-1)!} - k\sum_{i=k+1}^\infty A^i \frac{e^{-A}}{i!}\\ &
= & A\sum_{i=k}^\infty A^i \frac{e^{-A}}{i!} - k\sum_{i=k+1}^\infty A^i \frac{e^{-A}}{i!}\\ &
= & A\left[1-e^{-A}\sum_{i=0}^{k-1}\frac{A^i}{i!}\right] -k\left[ 1-e^{-A}\sum_{i=0}^{k}\frac{A^i}{i!}  \right] \\ &
= & \frac{A  \gamma(k,A)}{\Gamma(k)} - \frac{k \gamma(k+1,A)}{\Gamma(k+1)}. ~~~~~~\Box
\end{eqnarray*}

\subsection{First Passage Time and Time Between Events in M/M/$k$/$k$}

 As discussed, the term first passage time between  states $i$ and $j$ is the time it takes for a process in state $i$ to enter state $j$ for the first time.  To enhance the understanding of Markov chains in general and M/M/$k$/$k$ in particular, we consider here several examples associated with the first passage time,  as well as problems associated with times between events in the context of M/M/$k$/$k$. These examples will enhance understanding of Markov chains and provide preparation for the section on Markov-chain simulation of M/M/$k$/$k$. It is also important to understand that derivations of first passage times and times between events are applicable to a more general class of Markov chains and Markovian queueing models and they are clearly not limited to M/M/$k$/$k$.

\subsubsection*{Homework \ref{secerlang}.\arabic{homework}}
\addtocounter{homework}{1} \addtocounter{tothomework}{1}

Consider an M/M/4/4 system with an arrival rate $\lambda = 1$ per minute and a holding time of 180 seconds. Assume that at a given point in time $t$, there are 3 servers busy.  In each of the questions below, write the formula first and then substitute correctly the numerical values and compute the final answer in numerical value.

\begin{enumerate}
 \item Find the mean time from $t$ until the next arrival.

{\bf Solution}

The mean time from $t$ until the next arrival is $$\frac{1}{\lambda} = \frac{1}{1} =1.$$

\item Find the mean time from $t$ until the next event (either arrival or departure). (3 marks)

{\bf Solution}

$$\mu = \frac{1}{3} [{\rm min.}^{-1}].$$

The mean time from $t$ until the next event is

$$\frac{1}{\lambda + 3\mu} = \frac{1}{1 + 3 \times (1/3)} =0.5 [{\rm min.}].$$

\item Find the probability that the next event will be an arrival.

{\bf Solution}

$$\frac{\lambda}{\lambda + 3\mu} = \frac{1}{1 + 3 \times (1/3)} =0.5.$$

\item Find the probability that the next event will be a departure.

{\bf Solution}

$$\frac{3\mu}{\lambda + 3\mu} = \frac{3\times (1/3)}{1 + 3 \times (1/3)} =0.5.$$

\item Find the mean time from $t$ until the next departure.

{\bf Solution}

Let $X$ be the time from $t$ until the next departure. Then, $X$ can be divided into two parts: $X_1$ and $X_2$, where $X_1$ is the time until the next event (which can be an arrival or a departure) and $X_2$ is the time from the next event until the first departure.

Then, $E[X_1]$ was already obtained in Part 2 of this question, and we have:

$$E[X_1]=\frac{1}{\lambda + 3\mu} = \frac{1}{1 + 3 \times (1/3)} =0.5 [{\rm min.}].  $$

And $E[X_2]$ is obtained by the Law of Iterated Expectations as follows:

$$E[X_2] = \frac{\lambda}{\lambda + 3\mu} \times \frac{1}{4\mu} +  \frac{3\mu}{\lambda + 3\mu} \times 0 = 0.5 \times 0.75 = 0.375[{\rm min.}].$$

Then, we obtain

$$E[X] = E[X_1] + E[X_2]= 0.5 +0.375 = 0.875 [{\rm min.}]. $$

\end{enumerate}

\subsubsection*{Homework \ref{secerlang}.\arabic{homework}}
\addtocounter{homework}{1} \addtocounter{tothomework}{1}

Consider an M/M/2/2 queueing system with $\lambda$ being the arrival rate and $1/\mu$ being the mean holding time. Assume the case $\lambda = 9$ and $\mu =5$. At time $t$, there are two customers in the system. What is the mean time from time $t$ until the first time that the system is empty? (In other words, what is the mean first passage time from state 2 to state 0.)

\subsubsection*{Solution}

Let $m_{2,0}$ be the mean first passage time from state 2 to state 0.

Let $m_{1,0}$ be the mean first passage time from state 1 to state 0.

Then, $$m_{2,0} = \frac{1}{2\mu} + m_{1,0}$$
and
$$m_{1,0} = \frac{1}{\lambda +\mu} + \frac{\lambda}{\lambda +\mu}m_{2,0}.  $$

Substituting the above values for $\lambda$ and $\mu$, we obtain

$$m_{2,0} = \frac{1}{10} + m_{1,0}$$
and
$$m_{1,0} = \frac{1}{14} + \frac{9}{14}m_{2,0}$$

or

$$m_{2,0} = \frac{1}{10} +  \frac{1}{14} + \frac{9}{14}m_{2,0}$$

from which we obtain

$$m_{2,0} = \frac{14}{50} +  \frac{1}{5} = \frac{24}{50}=0.48.$$

\subsection{A Markov-chain Simulation of M/M/$k$/$k$}

We have described a Markov-chain simulation in the context of the M/M/1 queue. In a similar way,
we can use a Markov-chain simulation to evaluate the blocking probability of an M/M/$k$/$k$ queue, as follows.

Variables and input parameters: \\$k$ = number of servers;\\$Q$ = number of customers in the
system (queue size);\\ $B_p$ = estimation for the blocking
probability; \\$N_a$ = number of customer arrivals counted so far;\\
$N_b$ = number of blocked customers counted so far;\\ $MAXN_a$ =
maximal number of customer arrivals (it is used for the stopping
condition);
\\$\mu$ = service rate;\\ $\lambda$ = arrival rate.

Define function: $R(01)$ = a uniform $U(0,1)$ random variate. A new
value for $R(01)$ is generated every time it is called.

Initialization: $Q = 0$; $N_a=0$, $N_b=0$.

1. If $R(01) \leq \lambda/(\lambda + Q\mu)$, then $N_a = N_a + 1$; if $Q=k$, then $N_b=N_b+1$, else $Q = Q + 1$;

else, Q = Q - 1.

2. If $N_a < MAXN_a$ go to 1; else, print $B_p=N_b/N_a$.

Again, it is a very simple program of two IF statements: one to
check if the next event is an arrival or a departure, and the other a stopping criterion.

\subsubsection*{Homework \ref{secerlang}.\arabic{homework}}
\addtocounter{homework}{1} \addtocounter{tothomework}{1} Simulate an M/M/$k$/$k$ queue based on the Markov-chain
 simulation principles to evaluate the blocking probability for
a wide range of parameter values. Compare the results you obtained
with equivalent results obtained analytically using the Erlang B Formula
and with equivalent M/M/$k$/$k$ queue blocking probability results
obtained based on discrete-event simulations as discussed in Section
\ref{gg1simul}. In your comparison, consider accuracy (closeness to
the analytical results), the length of the confidence intervals and
running times. $~~~\Box$

\subsubsection*{Homework \ref{secerlang}.\arabic{homework}}
\addtocounter{homework}{1} \addtocounter{tothomework}{1} Use the principles presented in Section
\ref{gg1simul} for discrete-event simulation of a
G/G/1 queue to write three computer simulations for an M/G/$k$/$k$
queue where G is the service time modeled by
\begin{enumerate}
\item exponential distribution
\item uniform distribution (continuous)
\item deterministic (constant) value -- it always takes the same value
\item Pareto distribution
\item non-parametric distribution,
\end{enumerate}
where the mean service times in all cases are the same. Choose your own parameters for the different distributions (or values for the non-parametric distribution), but you must make sure to keep the mean service time the same in all cases.

Demonstrate the insensitivity property of  M/G/$k$/$k$ by showing that the blocking probability is the same in all cases. Obtain the blocking probability also by Erlang B formula and by a Markov-chain simulation for the case of exponentially distributed service times.
Repeat your simulation experiments for different cases using a wide range of parameter values. $~~~\Box$

\subsubsection*{Homework \ref{secerlang}.\arabic{homework}}
\addtocounter{homework}{1} \addtocounter{tothomework}{1} Simulate equivalent U/U/$k$/$k$,  M/U/$k$/$k$ ($U$ denotes here a uniform random variable)
and M/M/$k$/$k$ models. (You may use programs you have written in previous assignments.
Run these simulations
 for a wide range of parameter values and compare them
numerically. Compare them also with equivalent results obtained
analytically using the Erlang B Formula. Again, in your comparison, consider accuracy (closeness to
the analytical results), the length of the confidence intervals, and
running times. While in the previous assignment, you learned the effect of the method used on accuracy and running time, this time, try also to learn how the different models affect the accuracy and running times. $~~~\Box$

\subsubsection*{Homework \ref{secerlang}.\arabic{homework}}
\addtocounter{homework}{1} \addtocounter{tothomework}{1}
Use the M/M/$k$/$k$ model to compare the utilization of an optical switch with full wavelength conversion and without wavelength conversion.
\subsubsection*{Background information and guide}
Consider a switch with $T_I$ input trunks and $T_O$ output trunks. Each trunk comprises $F$ optical fibers, each of which comprises $W$ wavelengths. Consider a particular output trunk and assume that the traffic directed to it
follows a Poisson process with parameter $\lambda$ and that any packet is
of exponential length with parameter $\mu$. Let $A=\lambda/\mu$. In the case of full wavelength
conversion, every packet from any wavelength can be converted to any other wavelength, so the Poisson traffic with parameter $\lambda$ can all be directed to the output trunk and can use any of the $k=FW$ links. In the case of
no wavelength conversion, if a packet arrives on a given wavelength at an input port, it must continue on the same wavelength at the output port, so now consider $W$ separate systems, each has only $F$ links per trunk.
Compare the efficiency that can be achieved for both alternatives,
if the blocking probability is set limited to 0.001.
In other words, in the wavelength conversion case, you have an M/M/$k$/$k$ system with $k=FW$, and in the non-wavelength conversion case, you have $k=F$. Compute the traffic $A$ that gives a blocking probability of 0.001 in each case and compare efficiency. Realistic ranges are $40 \leq W \leq 100$ and $10 \leq F \leq 100$.
$~~~\Box$

\subsection{Preemptive Priorities}

So far, we have considered single-class traffic without any priority given to some
calls over others. Let us now consider an extension of the M/M/$k$/$k$ loss system, where some calls (customers) have preemptive priorities over
other lower-priority calls. The first to address this problem was Katzschner \cite{Katzschner70} in 1970, and later the problem was revisited in \cite{Vu02,Yang14,Yang24}. In the following, we provide the analysis of \cite{Vu02}.
In this case, when a higher priority arrives, and none of the $k$ servers is free, but some of the servers serve lower priority calls, the higher priority call preempts one of the lower priority calls in service and enters service instead of the preempted call.
We consider arriving calls to be of $m$ priority types. Where priority 1 represents the highest priority and priority $m$ represents the lowest priority. In general, if $i<j$, then priority $i$ is higher than priority $j$, so a priority $i$ arrival may preempt a priority $j$ customer upon its arrival.

The arrival process of priority  $i$ customers follows a Poisson process with rate $\lambda_i$, for $i=1,2,3, \ldots, m.$ The service time of all the customers is exponentially distributed with parameter $\mu$. The offered traffic of priority $i$ customers is given by  $$A_i = \frac{\lambda_i}{\mu}, ~~ i=1,2,3, \ldots, m.$$

Let $P_b(i)$ be the blocking probability of priority  $i$ customers. Because the priority 1 traffic access the system regardless of low-priority loading, for the case $i=1$, we have
$$P_b(1)= E_k (A_1).$$
To obtain $P_b(i)$ for $i>1$, we first observe that the blocking probability of all the traffic of priority $i$ and higher priorities, namely, the traffic generated by priorities $1, 2, \ldots, i$, is given by
$$E_k (A_1+A_2+\ldots, A_i).$$

To explain this observation, we first need to notice that we are applying results from the case with multiple classes of traffic without priorities to the case with priorities. Interestingly, because of the memoryless property of exponential distribution, they result in the same blocking probability. Observe that in a priority system, when a new high-priority call arrives in a full system where there are some low-priority calls in progress, the high-priority call will preempt a low-priority call and will cause one call loss. Then, in the non-priority system, the new arriving call will be blocked, and again one call loss will be recorded. Because of the exponential distribution of call holding time, and because they have the same mean, both systems, after the arrival of the new call, have the same number of customers in service and their remaining service times are exponentially distributed with parameter $\mu$. The potential inefficiency caused in the priority system by serving the low-priority call that has been preempted by the high priority has no effect relative to the non-priority system where such preemptions and potential waste of resources do not exist. In both systems, namely, the priority system and the non-priority system, the remaining service time of the call in progress is the same, and one blocked call has been recorded.

Next, we observe that the lost traffic of priority $i$, $i=1,2,3, \ldots, m$, is given by the lost traffic of priorities  $1,2,3, \ldots, i$ minus the lost traffic of priorities  $1,2,3, \ldots, i-1$, namely,
$$(A_1+A_2+\ldots, A_i)E_k (A_1+A_2+\ldots, A_i) - (A_1+A_2+\ldots, A_{i-1})E_k (A_1+A_2+\ldots, A_{i-1}).$$

Therefore, the value of $P_b(i)$ for $i>1$ can be obtained as the ratio of the lost traffic of priority $i$ to the offered traffic of priority $i$, that is,

$$P_b(i)=\frac{\left(\sum_{j=1}^i A_j \right)E_k \left(\sum_{j=1}^i A_j \right) - \left(\sum_{j=1}^{i-1} A_j \right)E_k \left(\sum_{j=1}^{i-1} A_j \right)}{A_i}.$$

Recall that the Erlang B formula, which provides the blocking probability for the M/M/$k$/$k$ queueing system, is insensitive to the distribution of the holding time beyond its mean. One important question is whether or not this insensitivity property is also applicable to the above solution  for M/M/$k$/$k$ with strict multiple priorities.  In other words, we question whether or not the above solution holds in general when the distribution of the holding times is not restricted to the exponential distribution, but it applies to cases where the holding times have other distributions with mean $1/\mu$.

Based on the described-above construction of the solution for M/M/$k$/$k$ with strict multiple priorities and the known insensitivity property of M/M/$k$/$k$ with a single priority, the reader may reach an incorrect conclusion that the blocking probabilities in an M/M/$k$/$k$ system with strict multiple priorities are also insensitive to the relevant holding times distributions except for their means. However, as described in \cite{Yang24}, the insensitivity property is only applicable for the top priority traffic because the service of the top priority customers is not affected at all by the lower priority customers due to the preemptive nature of the priority regime. Therefore, they can be considered as if they are alone in the $k$ server loss system. On the other hand, this insensitivity property does not apply to the lower priority customers.

\subsubsection*{Homework \ref{secerlang}.\arabic{homework}}
\addtocounter{homework}{1} \addtocounter{tothomework}{1}
Assume that the traffic offered to a 10-circuit system is composed of two priority traffic: high and low. The arrival rate of the high-priority traffic is 5 calls/minute, and that of the low-priority traffic is 4 calls/minute. The call holding time of both traffic classes is exponentially distributed with a mean of 3 minutes. Find the blocking probability of each of the priority classes.
$~~~\Box$

\subsubsection*{Homework \ref{secerlang}.\arabic{homework}}
\addtocounter{homework}{1} \addtocounter{tothomework}{1}

Consider the above-described loss system with Poisson arrivals and strict priority regime, and obtain the above results for the blocking probabilities of the different priority customers directly from the steady-state equations. See \cite{Yang24}  for a guide.

\subsubsection*{Homework \ref{secerlang}.\arabic{homework}}
\addtocounter{homework}{1} \addtocounter{tothomework}{1}

Demonstrate by simulations that for the case of loss system with Poisson arrivals and strict priority regime, the insensitivity property does not apply to customers of priorities $2, 3, \ldots, m$. See \cite{Yang24}  for a guide.

\subsection{Overflow Traffic of M/M/$k$/$k$}

In many practical situations, traffic that cannot be admitted to a $k$ server group overflows to another server group.
In such a case, the overflow traffic is not Poisson, but it is more bursty than a Poisson process.
In other words, the variance of the number of arrivals in an interval is higher than the mean number of arrivals in that interval.

It is, therefore, important to characterize such overflow traffic by its variance and its mean. In particular, consider an M/M/$k$/$k$
queueing system with input offered traffic $A$ and let $M$ [Erlangs] be the
traffic overflowed from this $k$-server system. As discussed in Section  \ref{overflow},
\begin{equation} \label{Mequal} M=AE_{k}(A).\end{equation}

Let $V$ be the variance of the overflow traffic. Namely, $V$ is the variance of the number of busy servers
in an infinite server system to which the traffic overflowed  from our M/M/$k$/$k$ is offered.

The following result, known as Riordan Formula, can be used to obtain $V$.

\begin{equation} \label{Vequal} V = M \left( 1 - M + \frac{A}{k+1+M-A}\right). \end{equation}

The derivation of this result is provided in Appendix I
authored by J. Riordan in \cite{Wilkinson56}.

Note that $M$ and $V$ of the overflow traffic are completely
determined by $k$ and $A$.

The variance to mean ratio of a traffic stream is called {\it Peakedness}. In our case, the peakedness of the overflow traffic is denoted
$Z$ and is given by  $$Z=\frac{V}{M},$$ and it normally satisfies $Z>1$.

\subsection{Multi-server Loss Systems with Non-Poisson Input}
\label{msls}

Consider a generalization of an M/M/$k$/$k$ system to the case where
the arrival process is not Poisson. As mentioned in the previous section, one way non-Poisson traffic occurs is when we consider a secondary server group to which traffic overflows from a primary server group. If we know the offered traffic to the primary server group (say $A$) and the numbers of servers in the primary and secondary groups are given by $k_1$ and $k_2$, respectively, then the blocking probability of the secondary server group is obtained by

\begin{equation} \label{BP2nd}
 P_b (\mbox{secondary})=\frac{E_{k_1+k_2}(A)}{E_{k_1}(A)}.
\end{equation}
In this case, we also know the mean and variance of the traffic offered to the secondary server group which is readily obtained by Equations (\ref{Mequal}) and (\ref{Vequal}).

However, the more challenging problem is the following: given a multi-server loss system with $k_2$ servers loaded by non-Poisson offered traffic with mean $M$ and variance $V$, find the blocking probability. First, let us clarify the definitions of $M$ and $V$ by saying that if this traffic is offered to a system with an infinite number of servers (instead of $k_2$), the mean and variance of the number of busy servers will be $M$ and  $V$, respectively.
This offered traffic could have come from or overflowed from various sources, and unlike the previous problem, here we do not know anything about the original offered traffic streams or the characteristics of any primary systems. All we know are the values of $M$ and $V$.

This problem does not have an exact solution, but reasonable approximations are available. We will now present two approximations:
\begin{enumerate}
\item Equivalent Random Method (ERM)
\item Hayward Approximation.
\end{enumerate}

\subsubsection*{Equivalent Random Method (ERM)}

We wish to estimate the blocking probability for a system with $k_2$ servers loaded by
non-pure chance offered traffic with mean $M$ and variance $V.$

Under the ERM, due to \cite{Wilkinson56}, we model the system as if the traffic
was the overflow traffic from a primary system with $N_{eq}$ circuits and offered traffic $A_{eq}$ that follows Poisson process. If we find such $A_{eq}$ and $N_{eq}$, then by Eq. (\ref{BP2nd}),
the blocking probability in our $k_2$-server system, denoted $PB_{k_2}(N_{eq}, A_{eq})$, will be estimated by:
$$PB_{k_2}(N_{eq}, A_{eq})= \frac{E_{N_{eq}+k_2}(A_{eq})}{E_{N_{eq}}(A_{eq})}.$$

To approximate $A_{eq}$  and $N_{eq}$, we use the following:

\begin{equation}
\label{aeq}
A_{eq} = V + 3 Z (Z -1);
\end{equation}

\begin{equation}
\label{neq}
N_{eq} = \frac{A_{eq}(M+Z)} {M+Z-1} - M - 1.
\end{equation}
Note that Equation (\ref{aeq}) is an approximation, but Equation (\ref{neq}) is exact and it results in an approximation only because Equation (\ref{aeq}) is an approximation.

One important benefit of ERM is that it can also provide an estimation of the variance of the overflow traffic of the secondary system not only its mean, by using Riordan Formula for a  system of $k_1 + k_2$ servers of both the primary and the secondary.

\subsubsection*{Hayward Approximation}

The Hayward approximation \cite{Fredericks80} is based on the following result. A multi-server system with $k$ servers fed by traffic with mean $M$ and variance $V$ has a similar blocking probability to that of an M/M/$k$/$k$ system with offered load $\frac{M}{Z}$ and $\frac{k}{Z}$ servers, hence Erlang B formula that gives $E_{\frac{k}{Z}}\left(\frac{M}{Z}\right)$ can be used if $\frac{k}{Z}$ is rounded to an integer.

The Hayward Approximation is simpler than ERM, but it cannot provide an estimation of the variance of the overflow traffic.

\subsubsection*{Homework \ref{secerlang}.\arabic{homework}}
\addtocounter{homework}{1} \addtocounter{tothomework}{1}

The XYZ Corporation has measured their offered traffic to a set of 24 circuits and found that it has a mean of $M = 21$ and a variance of $V = 31.5$. Arriving calls cannot be queued and delayed before they are served.  This means that if a call arrives and all the circuits (servers) are busy, the call is blocked and cleared from the system.

Use the Hayward Approximation as well as the Equivalent Random Method to estimate the blocking probability.
\subsubsection*{Solution}

$M = 21$;~~                   $ V = 31.5$;

Peakedness:        $ Z = \frac{31.5}{21} = 1.5.$

\subsubsection*{Hayward Approximation}
The mean offered traffic in the equivalent system = 21 / 1.5 = 14.

The number of servers in the equivalent system = 24/1.5 = 16.

Blocking probability = $E_{16}(14) = 0.1145.$

\subsubsection*{Equivalent Random Method}
$$A_{eq} = V + 3Z(Z-1) = 31.5 + 3 \time 1.5 \times 0.5 = 33.75.$$
$$N_{eq} = \frac{ A_{eq} (M + Z)}{M+Z-1} – M - 1
= \frac{33.75(21+1.5)}{21 +1.5 -1} - 21 - 1 = 13.32.$$

We will use $N_{eq} = 13$ conservatively.

From Erlang B:

$E_{13}(33.75) = 0.631$.

$E_{13+24}(33.75) = 0.075$.

Then, the blocking probability is obtained by  $$\frac{E_{13+24}(33.75)}{ E_{13}(33.75)} = \frac{0.07536}{ 0.631} = 0.1194. ~~~\Box $$

\subsubsection*{Homework \ref{secerlang}.\arabic{homework}}
\addtocounter{homework}{1} \addtocounter{tothomework}{1}
Assume that non-Poisson traffic with mean = 65 Erlangs and variance = 78 is offered to a loss system.
Use both Hayward Approximation
and the Equivalent Random Method to estimate the minimal number of circuits required to guarantee
that the blocking probability is not more than 1\%.

\subsubsection*{Solution}

Let us use the notation $N^*$ to represent the minimal number of circuits required to guarantee
that the blocking probability is not more than 1\%.
Previously, we use $k_2$ to represent the {\it given} number of servers in the secondary system. Now we use the notation $N^*$ to represent the
desired number of servers in the system.

Given, $M=65$ and $V=78$, the peakedness is given by
$$Z = \frac{78}{65}=1.2.$$

\subsubsection*{Equivalent Random Method}

By (\ref{aeq}): $A_{eq} = 78 + 3 \times 1.2 \times 0.2 = 78.72$.

By (\ref{neq}): $N_{eq} = \frac{78.72 (65+1.2)}{65 +1.2 -1} - 65 -1 =13.92736 = 14 ~~{\rm approx.}$

A conservative rounding would be to round it down to $N_{eq}=13$. This will result in a more conservative dimensioning because a lower value for $N_{eq}$ will imply a higher loss in the primary system, which leads to higher overflow traffic from the primary to the secondary system. This in turn will require more servers in the secondary system to meet the required blocking probability level.

In the present case, because the result is 13.92736 (so close to 14), we round it up to  $N_{eq}=14.$ In any case, we need to be aware of the implications of our choice. We will repeat the calculation using the more conservative choice of $N_{eq}=13.$

The blocking probability in the primary equivalent system is given by,

$$E_{14}(78.72) = 0.825.$$

Next, find the minimal value for $N^*+14$ such that
$$\frac{E_{N^*+14}(78.72)}{0.825} \leq 0.01,$$

or,

$$E_{N^*+14}(78.72) \leq 0.00825.$$

By Erlang B formula: $N^* + 14 = 96$, so the number of required servers is estimated by
$N^* = 82$.

Notice that the choice of $N^* + 14 = 96$ is made because $$E_{96}(78.72) = 0.0071$$ satisfies the requirement
and
$$E_{95}(78.72) = 0.0087,$$
does not satisfy the requirement.

Now we consider $N_{eq}=13$.

The blocking probability in the primary equivalent system is given by,

$$E_{13}(78.72) = 0.8373.$$

Next, find the minimal value for $N^*+13$ such that
$$\frac{E_{N^*+13}(78.72)}{0.8373} \leq 0.01,$$

or,

$$E_{N^*+13}(78.72) \leq 0.008373.$$

By Erlang B formula: we choose $N^* + 13 = 96$ because as above $E_{96}(78.72) = 0.0071$ satisfies the requirement
and $E_{95}(78.72) = 0.0087$ does not.

Because of our conservative choice of $N_{eq}=13$, the more conservative choice for the
number of required servers is $N^* = 83$.

\subsubsection*{Hayward Approximation}

Mean offered traffic in the equivalent system = 65/1.2 = 54.16667.

By the Erlang B formula, the number of required servers in the equivalent system for 1\% blocking is 68.

Then, $68  \times 1.2 = 81.6.$ Rounding up conservatively, we obtain that 82 servers are required.

Now the designer will need to decide between 82 servers based on the Hayward approximation
and based on a reasonable but less conservative approach according to the Equivalent Random Method,
or 83 according to the very conservative approach based on the Equivalent Random Method. $~~~\Box$

\subsubsection*{Homework \ref{secerlang}.\arabic{homework}}
\addtocounter{homework}{1} \addtocounter{tothomework}{1}

Consider again the two loss systems the primary and the secondary and use them to compare numerically
between:
\begin{enumerate}
\item the exact solution;
\item
 the Hayward approximation;
 \item the Equivalent Random Method approximation;
 \item a 3rd approximation that is based on the assumption that the arrival process
 into the secondary system follows a Poisson process. For this approximation assume that the traffic lost in the primary system is offered to the secondary system following a Poisson process.
 \end{enumerate}
\subsubsection*{Guide}
For the comparison, at first assume that you know $A$, $k_1$, and $k_2$ and compute $M$ and $V$, i.e., the mean and variance of the offered load to the secondary system as well as the blocking probability of traffic in the secondary system using (\ref{BP2nd}).

Next, assume that $A$ and $k_1$ are not known but $k_2$ is known; also known are $M$ and $V$, i.e., the mean and variance of the offered load to the secondary system that you computed previously. And evaluate the blocking probability  using both Hayward, the Equivalent Random Method, and the Poisson approximations.

Compare the results for a wide range of parameters. $~~~\Box$

\newpage
\section{{ {M/M/$k$}}}
\label{mmk}

\setcounter{homework}{1} 

The M/M/$k$  queue is a generalization of the M/M/1 queue to the case
of $k$ servers. As in M/M/1, for an M/M/$k$ queue, the buffer is
infinite and the arrival process is Poisson with rate $\lambda$. The
service time of each of the $k$ servers is exponentially distributed
with parameter $\mu$. As in the case of M/M/1, we assume that the
service times are independent and independent of the arrival
process.

\subsection{Steady-State Equations and Their Solution}

Letting $A=\lambda/\mu$, and assuming the stability
condition $\lambda<k\mu$, or $A<k$, the M/M/$k$  queue gives rise to
the following steady-state equations:

$\pi_1 = A \pi_0 $\\
$\pi_2 = A \pi_1/2 = A^2 \pi_0$/2\\
$\pi_3 = A \pi_2/3 = A^3 \pi_0$/(3!)\\
~\ldots~ \\
$\pi_k = A \pi_{k-1}/k = A^k \pi_0/(k!)$ \\
$\pi_{k+1} = A \pi_{k}/k = A^{k+1} \pi_0/(k!k)$ \\
$\pi_{k+2} = A \pi_{k+1}/k = A^{k+2} \pi_0/(k!k^2)$ \\
~\ldots~ \\
$\pi_{k+j} = A \pi_{k+j-1}/k = A^{k+j} \pi_0/(k!k^j) ~~{\rm for}~ j=1,~ 2, ~3,~\ldots$\\

and in general:

\begin{equation}
\label{mmk1} \pi_n = \frac{A^n \pi_0}{n!}  ~{\rm for}
~n=0,~1,~2,~\ldots,~k-1
\end{equation}
and
\begin{equation}
\label{mmk2} \pi_n = \frac{A^n \pi_0}{k!k^{n-k}} ~{\rm for}
~n=k,~k+1,~k+2,~\ldots~.
\end{equation}

These balance equations can also be described by the following state transition diagram of M/M/$k$:

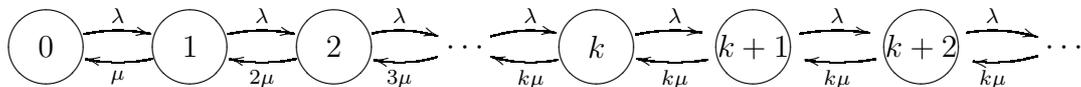
$$\xymatrix{
*++++[o][F]{0}\ar@/^{.5pc}/[r]^{\lambda} &
*++++[o][F]{1}\ar@/^{.5pc}/[l]^{\mu}\ar@/^{.5pc}/[r]^{\lambda} &
*++++[o][F]{2}\ar@/^{.5pc}/[l]^{2\mu}\ar@/^{.5pc}/[r]^{\lambda} &
{\cdots}\ar@/^{.5pc}/[l]^{3\mu}\ar@/^{.5pc}/[r]^{\lambda} &
*++++[o][F]{k}\ar@/^{.5pc}/[l]^{k\mu}\ar@/^{.5pc}/[r]^{\lambda} &
*++[o][F]{k+1}\ar@/^{.5pc}/[l]^{k\mu}\ar@/^{.5pc}/[r]^{\lambda} &
*++[o][F]{k+2}\ar@/^{.5pc}/[l]^{k\mu}\ar@/^{.5pc}/[r]^{\lambda} &
{\cdots}\ar@/^{.5pc}/[l]^{k\mu}
}$$

To obtain $\pi_0$, we sum up both sides of Eqs. (\ref{mmk1}) and
(\ref{mmk2}), and  because the sum of the $\pi_n$s equals one, we
obtain an equation for $\pi_0$, which its solution is
\begin{equation}
\pi_0=\left(\sum_{n=0}^{k-1}
\frac{A^n}{n!}+\frac{A^k}{k!}\frac{k}{(k-A)}\right)^{-1}.
\end{equation}
Substituting the latter in Eqs. (\ref{mmk1}) and (\ref{mmk2}), we
obtain the steady-state probabilities $\pi_n,~~ n=0,~1,~2,~\ldots~.$

\subsection{Erlang C Formula}

Of special interest is the so-called Erlang C formula. It
represents the proportion of time that all $k$ servers are busy
and is given by:

\begin{equation}
\label{ErC} C_k(A) = \sum_{n=k}^\infty \pi_n =
\frac{A^k}{k!}\frac{k}{(k-A)}\pi_0=\frac{\frac{A^k}{k!}\frac{k}{(k-A)}}{\sum_{n=0}^{k-1}
\frac{A^n}{n!} + \frac{A^k}{k!}\frac{k}{(k-A)}}.
\end{equation}
\subsubsection*{Homework \ref{mmk}.\arabic{homework}}
\addtocounter{homework}{1} \addtocounter{tothomework}{1} Derive
Eq.\@ (\ref{ErC}). $~~~\Box$

By  Eqs. (\ref{eka}) and (\ref{ErC}) we obtain the following
relationship:

\begin{equation}
\label{ErCfromErB} C_k(A) = \frac{kE_k(A)}{k-A[1-E_k(A)]}.
\end{equation}

\subsubsection*{Homework \ref{mmk}.\arabic{homework}}
\addtocounter{homework}{1} \addtocounter{tothomework}{1}
\begin{enumerate}
\item  Derive Eq.\@ (\ref{ErCfromErB});
\item  Show that $C_k(A) \geq E_k(A)$. $~~~\Box$
\end{enumerate}

An elegant result for $C_k(A)$ is the following \begin{equation}
\label{elegant} \frac{1}{C_k(A)} = \frac{1}{E_k(A)} -
\frac{1}{E_{k-1}(A)}. \end{equation}

\subsubsection*{Homework \ref{mmk}.\arabic{homework}}
\addtocounter{homework}{1} \addtocounter{tothomework}{1}
Prove Eq.
(\ref{elegant}).

In the following table, we add the corresponding $C_k(A)$ values
to the table of the previous section. We can observe the
significant difference between $C_k(A)$ and $E_k(A)$ as the ratio
$A/k$ increases. Clearly, when $A/k>1$, the M/M/$k$ queue is
unstable.

\begin{center}
\renewcommand{\arraystretch}{1.4}
\begin{tabular}{|c|c|c|c| }  \hline
$A$ & $k$ & $E_k(A)$ & $C_k(A)$  \\  \hline 20 & 30 & 0.0085
 &
0.025
 \\  \hline
 100 & 117 & 0.0098  & 0.064
 \\
\hline 500 & 527 & 0.0095 & 0.158  \\
\hline 1000  &  1029   & 0.0099 & 0.262  \\
\hline 5000   & 5010  &  0.0100 & 0.835  \\
\hline 10000  & 9970  &  0.0099 & unstable \\
\hline
\end{tabular}
\end{center}

\subsubsection*{Homework \ref{mmk}.\arabic{homework}}
\addtocounter{homework}{1} \addtocounter{tothomework}{1} Reproduce the results of the above table.
$~~~\Box$

\subsection{Mean Queue Size, Delay, Waiting Time and Delay Factor}
Let us reuse the following notations:\\
$Q$ = a random variable representing the total number of customers in the system (waiting in the queue
and being served);\\
$N_Q$ = a random variable representing the total number of customers waiting in the
queue (this does not include those customers
being served);\\
$N_s$ = a random variable representing the total number of customers that are being served;\\
$D$ = a random variable representing the total delay in the system (this
includes the time a customer waits in the queue
and in service);\\
$W_Q$ = a random variable representing the time a customer waits in the queue
(this excludes the time a customer spends in service);\\
$S$ = a random variable representing the service time.\\
 $\hat{D}$ = The delay of a delayed customer including the service time;\\
 $\hat{W_Q}$ = The delay of a delayed customer in the queue excluding the service time.

Using the above notation, we have
\begin{equation}
\label{mmkqsize}
E[Q] = E[N_Q] + E[N_s]
\end{equation}
and
\begin{equation}
\label{mmkqdel}
E[D] = E[W_Q] + E[S].
\end{equation}
Clearly, $$E[S]=\frac{1}{\mu}.$$
To obtain $E[N_s]$ for the M/M/$k$ queue, we use Little's formula for the system made of servers. If we consider the system
of servers (without considering the waiting room outside the
servers), we notice that since there are no losses, the arrival
rate into this system is $\lambda$ and the mean waiting time of each
customer in this system is $E[S]=1/\mu$.
Therefore, by Little's formula, the mean number of busy servers is given by
\begin{equation}
\label{ENsmmk}
E[N_s]=\frac{\lambda}{\mu}=A.
\end{equation}
To obtain $E[N_Q]$, we consider two mutually exclusive and exhaustive events:
$\{Q\geq k\}$, and $\{Q< k\}$.
Recalling the Law of Iterated Expectation (\ref{meancondind}), we have
\begin{equation}
\label{meanNQmmk} E[N_Q] = E[N_Q \mid Q\geq k] P(Q \geq  k) + E[N_Q \mid Q< k] P(Q< k).
\end{equation}
To derive $E[N_Q \mid Q\geq k]$, we notice that the evolution of the M/M/$k$
queue during the time when $Q\geq k$ is equivalent to that of an M/M/1 queue with arrival rate $\lambda$
and service rate $k\mu$. This can also be seen by observing the state transition diagram from state $k$ onwards.
Then, under the condition $Q\geq k$, the mean queue size of such M/M/1 queue is equal to
$\rho/(1-\rho)$ where $\rho = \lambda / (k\mu) = A/k$. Thus,
$$E[N_Q \mid Q\geq k]=\frac{A/k}{1-A/k}=\frac{A}{k-A}.$$
Therefore, since $E[N_Q \mid Q< k] =0$ and $P(Q \geq  k) = C_k(A)$, we obtain by  (\ref{meanNQmmk}) that
\begin{equation}
\label{meanmmk} E[N_Q]=C_k(A)\frac{A}{k-A}.
\end{equation}

\subsubsection*{Homework \ref{mmk}.\arabic{homework}}
\addtocounter{homework}{1} \addtocounter{tothomework}{1} Derive
Eq.\@ (\ref{meanmmk}) by a direct approach using
$E[N_Q]=\sum_{n=k}^\infty (n-k)\pi_n$.

\subsubsection*{Guide}
By (\ref{mmk2}),
$$E[N_Q]=\sum_{n=k}^\infty (n-k)\pi_n= \sum_{n=k}^\infty (n-k) \frac{A^n \pi_0}{k!k^{n-k}} $$
Set $i=n-k$, to obtain
$$E[N_Q]= \sum_{i=0}^\infty i \frac{A^{i+k} \pi_0}{k!k^{i}} = \frac{\pi_0 A^k}{k!}
\sum_{i=0}^\infty i \left( \frac{A}{k}\right)^i=C_k(A)\frac{k-A}{k} \frac{A/k}{(1-A/k)^2},$$
and (\ref{meanmmk}) follows.
$~~~\Box$

\subsubsection*{Homework \ref{mmk}.\arabic{homework}}
\addtocounter{homework}{1} \addtocounter{tothomework}{1}
Confirm the
consistence between (\ref{meanmmk}) and (\ref{meanmm1}). $~~~\Box$

By (\ref{mmkqsize}), (\ref{ENsmmk}) and (\ref{meanmmk}), we obtain
\begin{equation}
\label{mmkqsizef}
E[Q] = C_k(A)\frac{A}{k-A} + A.
\end{equation}
By (\ref{meanmmk}), invoking Little's formula, we obtain
\begin{equation}
\label{mmkEWQ}
E[W_Q] = \frac{C_k(A)\frac{A}{k-A}}{\lambda} = \frac{C_k(A)}{\mu k-\lambda}.
\end{equation}
Notice that $E[W_Q]$ is the ratio between the probability of having all
servers busy and the spare capacity of the system.

\subsubsection*{Homework \ref{mmk}.\arabic{homework}}
\addtocounter{homework}{1} \addtocounter{tothomework}{1}
Show the result of (\ref{mmkEWQ}) by using the law of iterated expectations directly without using (\ref{meanmmk}).
 $~~~\Box$

The mean delay is readily obtained by adding the mean service time to $E[W_Q]$. Thus,
\begin{equation}
\label{mmkED}
E[D] = \frac{C_k(A)}{\mu k-\lambda} + \frac{1}{\mu}.
\end{equation}

Another useful measure is the so-called {\it delay factor}
\cite{ciscovoip}. It is defined as the ratio of the mean waiting
time in the queue to the mean service time. Namely, it is given by
\begin{equation}
\label{Dfac} D_F = \frac{E[W_Q]}{1/\mu}= \frac{\frac{C_k(A)}{\mu
k-\lambda}}{\frac{1}{\mu}} = \frac{C_k(A)}{k-A}.
\end{equation}
The rationale for using the delay factor is that in some applications users
that require long service time may be willing to wait a long time in
the queue in direct proportion to the service time.

\subsection{Mean Delay of Delayed Customers}
\label{meandelayedmmk}
In Section \ref{meandelayed}, we have shown how to derive, for the case of M/M/1,
$E[\hat{D}]$ and  $E[\hat{W_Q}]$, namely, the mean delay of a delayed customer including
the service time and excluding the service time, respectively. We now extend the same ideas to the case of M/M/$k$.
As in the previous case, to obtain $E[\hat{W_Q}]$, we use Little's formula where we consider the queue (without the servers) as the system
and the arrival rate of the delayed customers which in the present case is $\lambda C_k(A)$.

Therefore, $$E[\hat{W_Q}]= \frac{AC_k(A)}{\lambda C_k(A) (k-A)}=\frac{1}{k\mu - \lambda},  $$
and $$E[\hat{D}] = E[\hat{W_Q}] + \frac{1}{\mu}= \frac{1}{k\mu-\lambda} + \frac{1}{\mu}.$$

As in Section \ref{meandelayed}, we can check the latter using the Law of Iterated Expectation as follows:
\begin{eqnarray*}
E[D] & = &  (1-C_k(A)) E[S]  + C_k(A) E[\hat{D}]\\
& = & (1-C_k(A))\frac{1}{\mu} + C_k(A) \left( \frac{1}{k\mu-\lambda} + \frac{1}{\mu}\right) =  \frac{C_k(A)}{\mu k-\lambda} + \frac{1}{\mu},
\end{eqnarray*}
and again we observe that consistency is achieved and note that this consistency check is an alternative way to obtain $E[\hat{D}]$.

\subsection{Dimensioning}
\label{dim_mmk} One could solve the dimensioning problem of finding,
for a given $A$, the smallest $k$ such that $C_k(A)$ or the mean
delay is lower than a given value. Using Eq. (\ref{ErCfromErB}),
and realizing that the value of $C_k(A)$ decreases as $k$ increases,
the dimensioning problem with respect to $C_k(A)$ can be solved in
an analogous way to the M/M/$k$/$k$ dimensioning problem. Having
$C_k(A)$  for value of $k$, one can also obtain the
minimal $k$ such that the mean delay is bounded by a given value using Eq.\@
(\ref{mmkED}). A similar procedure can be used to find the minimal
$k$ such that
 delay factor requirement is met.

\subsection{Utilization}

The utilization of an M/M/$k$ queue is the ratio of the mean number of
busy servers to $k$, therefore the utilization
of an M/M/$k$ queue is obtained by
\begin{equation} \hat{U}=\frac{E[N_s]}{k}=\frac{A}{k}.
\end{equation}

\subsubsection*{Homework \ref{mmk}.\arabic{homework}}
\addtocounter{homework}{1} \addtocounter{tothomework}{1}
Write a computer program that computes the minimal $k$, denoted $k^*$,
subject to a bound on $E[D]$. Run the program for a wide range of parameter
values and plot the results. Try to consider meaningful relationships, e.g.,
plot the spare capacity $k^*\mu-\lambda$ and utilization as a function of
various parameters and discuss implications. $~~~\Box$

\subsubsection*{Homework \ref{mmk}.\arabic{homework}}
\addtocounter{homework}{1} \addtocounter{tothomework}{1} Consider the M/M/$2$ queue with arrival rate $\lambda$ and service rate $\mu$ of each server. \begin{enumerate}
\item Show that $$\pi_0=\frac{2-A}{2+A}.$$
\item Derive formulae for $\pi_i ~~{\rm for} ~ i=1,2,3,~\ldots~.$
\item Show that $$C_2(A)=\frac{A^2}{2+A}.$$ Note that for $k=2$, it is convenient to use $C_2(A)=1-\pi_0-\pi_1$.
\item Derive a formula for $E[N_s]$ using the sum: $\pi_1+2C_2(A)$ and show that $$E[N_s]=\pi_1+2C_2(A)=A.$$
\item Derive $E[Q]$ in two ways, one using the sum $\sum_{i=0}^\infty i\pi_i$ and the other using Eqs. (\ref{meanNQmmk}) -- (\ref{mmkqsizef}), and show that in both ways you obtain $$E[Q]=\frac{4A}{4-A^2}. ~~~\Box $$
\end{enumerate}

\subsubsection*{Homework \ref{mmk}.\arabic{homework}}
\addtocounter{homework}{1} \addtocounter{tothomework}{1} Queueing theory is a useful tool for decisions on hospital resource
allocation \cite{Butterfield07,chan18,Green2001,Green2006,Weber2006}.
 In particular,
the M/M/$k$ model has been considered \cite{Green2001,Green2006}.
Consider the following example from \cite{Weber2006}. Assume that a patient stays at an Intensive Care Unit (ICU) for an exponentially
distributed period of time with an average time of 2.5 days.
Consider two hospitals. Patients' arrivals at each of the hospitals
follow a Poisson process. They arrive at Hospital 1 at the rate of
one patient per day and the arrival rate at Hospital 2 is 2 patients per day.
Assume that Hospital 2 has 10 ICU beds. Then, the management of
Hospital 1 which has never studied queueing theory believes that they
need only 5 beds because they think that if they have half the
traffic load they need half the number of beds. Your job is to
evaluate and criticize their decision. Assuming an M/M/$k$ model,
calculate the mean delay and the probability of having all servers
busy for each of the two systems. Which one performs better? If you
set the probability of having all servers busy in Hospital 2 as the
desired QoS standard, how many beds Hospital 1 will need? Maintaining the same QoS standard, provide a table with the number of beds needed in Hospital 1 if it has traffic arrival rates $\lambda_1 = 4, 8, 16, 32, 64$ patients per day. For each of the $\lambda_1$ values, provide a simple rule to estimate the number
of beds $n$ in Hospital 1, maintaining the same QoS standard.
Provide rigorous arguments to justify this rule for large values of
$\lambda_1$ and $n$.

{\bf Hint:} Observe that as $n$ grows with $\lambda_1$,
$n-\lambda_1$, approaches $C\sqrt{n}$ for some constant $C$ (find
that constant!). For rigorous arguments, study \cite{Halfin81}.
$~~~\Box$


\newpage
\section{Engset Loss Formula}
\label{secengset}

\setcounter{homework}{1} 

The Engset loss formula applies to telephony situations where the
number of customers is not too large relative to the number of available
circuits. Such situations include  an exchange in a small rural
community, PABX, or a lucrative satellite service to a small
number of customers. Let the call holding times be IID exponentially
distributed with mean $1/\mu$ and the time until an idle source
attempts to make a call is also exponentially distributed but with mean $1/\hat{\lambda}$.
We also assume that there is no dependence between the holding
times and the idle periods of the sources. Let the number of
customers (sources of traffic) be $M$,
the number of  circuits $k$, and the
blocking probability $P_b$.

This gives rise to a finite source loss model as opposed to the M/M/$k$/$k$ model where the arrival process follows a Poisson process which can describe traffic generated by an infinite number of sources (customers). In this book we only consider finite source loss models that do not involve queueing. For finite source models that involve queueing, the reader is referred to \cite{sztrik}.

The reader will recall that in M/M/1, the arrival rate, as well as the service rate, are independent of the state of the system, and in M/M/$\infty$, the arrival rate is also independent of the
number of customers in the system, but the service rate is state
dependent. In the present case, when the number of customers is
limited, we have a case where both the arrival and the service rates are state-dependent.

As in M/M/$k$/$k$, in our Engset finite source loss model, the service rate is $n\mu$ when there are $n$
busy circuits (namely $n$ customers are making phone calls).
However, unlike M/M/$k$/$k$, where the arrival rate is state-independent, in the present case, the arrival rate depends on the number of customers that are being served.
Under the Engset finite source model, busy customers
do not make new phone calls thus they do not contribute to the
arrival rate. Therefore, if $n$ circuits are busy, the arrival
rate is $(M-n)\hat{\lambda}$. As a result, considering both arrival and
service processes, at any point in time, given that there are
 $n$ customers in the system, and at the
same time, $n$ servers/circuits are busy, and the time until the next
event is exponentially distributed with parameter $(M-n)\hat{\lambda}+
n\mu $, because it is a competition between $M$ exponential random
variables: $n$ with parameter $\mu$ and $M-n$ with parameter
$\hat{\lambda}$.

An important question we must always answer in any Markov-chain analysis is how many states we have. If $M>k$, then the number of states is $k+1$, as in  M/M/$k$/$k$. However, if $M<k$, the number of states is $M+1$ because no more than $M$ calls can be in progress at the same time. Therefore, the number of states is ${\rm min}\{M,k\}+1$.

\subsection{Steady-State Equations and Their Solution}
Considering a finite state birth-and-death process that represents the queue evolution
of the above-described queueing system with
$M$ customers (sources) and $K$ servers, we obtain the following steady-state equations:\\
$\pi_0 M\hat{\lambda} = \pi_1 \mu$\\
$\pi_1 (M-1)\hat{\lambda} = \pi_2 2\mu $\\
$\pi_2 (M-2)\hat{\lambda} = \pi_3 3\mu $\\
\ldots\\
and in general:

\begin{equation}
\label{engsetsse1} \pi_n (M-n)\hat{\lambda} = \pi_{n+1} (n+1)\mu, ~{\rm
for} ~n=0,~1,~2,~ \ldots, {\rm min}\{M,k\}-1.
\end{equation}

Therefore, after standard algebraic manipulations of
(\ref{engsetsse1}), we
can write $\pi_n$, for \\$n=0,~1,~2,~ \ldots, {\rm min}\{M,k\}$, in
terms of $\pi_0$ as follows:

\begin{equation}
\label{engsetsse2} \pi_n  =  \left(\begin{array}{c}M\\n
 \end{array}\right)\left(\frac{\hat{\lambda}}{\mu}
\right)^n
 \pi_{0}, ~{\rm for}
~n=0,~1,~2,~ \ldots, {\rm min}\{M,k\},
\end{equation}

These steady-state equations are described by the following state transition diagram for the case $M>k$:

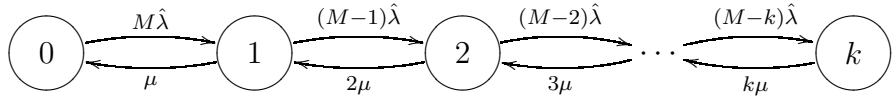
$$\xymatrix@C=48pt{
*++++[o][F]{0}\ar@/^{.5pc}/[r]^{M\hat{\lambda}} &
*++++[o][F]{1}\ar@/^{.5pc}/[l]^{\mu}\ar@/^{.5pc}/[r]^{\left(M-1\right)\hat{\lambda}} &
*++++[o][F]{2}\ar@/^{.5pc}/[l]^{2\mu}\ar@/^{.5pc}/[r]^{\left(M-2\right)\hat{\lambda}} &
{\cdots}\ar@/^{.5pc}/[l]^{3\mu}\ar@/^{.5pc}/[r]^{\left(M-k\right)\hat{\lambda}} &
*++++[o][F]{k}\ar@/^{.5pc}/[l]^{k\mu}}$$

Using the notation $\hat{\rho} = \hat{\lambda}/\mu$, we obtain

\begin{equation}
\label{engsetsse3} \pi_n  =  \left(\begin{array}{c}M\\n
 \end{array}\right)\hat{\rho}^n
 \pi_{0}, ~{\rm for}
~n=0,~1,~2,~ \ldots, {\rm min}\{M,k\}.
\end{equation}


Of course, the sum of the steady-state probabilities must be equal
to one, so we again have the additional normalizing equation
\begin{equation}
\label{sumto1engset}
\sum_{j=0}^{{\rm min}\{M,k\}} \pi_j =1.
\end{equation}

By  (\ref{engsetsse3}) together with the normalizing Eq.\@
(\ref{sumto1engset}), we obtain

$$\pi_0  =  \frac{ 1}{  \sum_{j=0}^{\min \{M,k\} }
 \left(\begin{array}{c}M\\j
 \end{array}\right)\hat{\rho}^j }.$$

Therefore, by  (\ref{engsetsse3}), we obtain

\begin{equation}
\label{engsetpis} \pi_n  =  \frac{ \left(\begin{array}{c}M\\n
 \end{array}\right)\hat{\rho}^n }{\sum_{j=0}^{\min \{M,k\} }
 \left(\begin{array}{c}M\\j
 \end{array}\right)\hat{\rho}^j }, ~{\rm for} ~n=0,~1,~2,~ \ldots, {\rm min}\{M,k\}.
\end{equation}

\subsection{Blocking Probability}

Now, what is the blocking probability $P_b$? When $k \geq M$,
clearly, $P_b=0$, as there is never a shortage of circuits.

To derive the blocking probability for the case when $k < M$, we
first realize that unlike in the case of Erlang Formula, $\pi_{k}$
does not give the blocking probability. Still, $\pi_{k}$ is the
probability of having $k$ busy circuits, or the proportion of time
that all circuits are busy which is the so-called {\em
time-congestion}, but it is not the probability that a call is
blocked -- the so-called {\em call-congestion}. Unlike the case of
Erlang B Formula, here, call congestion is not equal to time
congestion. This is because, in the Engset model, the arrival process does not follow a Poisson process.
In fact, the arrival rate is dependent on the state of the system. When the system is full the arrival rate is lower and could be much lower, than when the system is empty.

In particular, when $i$ circuits are busy, the arrival rate is
$\hat{\lambda}(M-i)$, therefore to find the proportion of calls blocked, or the blocking probability denoted $P_b$, we compute the ratio between calls arrive when there are $k$ circuits busy and the
total calls arrive. This gives

\begin{equation}
\label{engset1} P_b=\frac{\hat{\lambda} (M-k)\pi_k }{\hat{\lambda} \sum_{i=0}^k
(M-i)\pi_i}.
\end{equation}

Substituting (\ref{engsetsse2}) and (\ref{engsetsse3}) in (\ref{engset1}) and performing few algebraic manipulations,
we obtain the Engset loss formula that gives the blocking
probability for the case $M > k$ as follows.
\begin{equation}
\label{engset} P_b=\frac{\left(\begin{array}{c}M-1\\k
 \end{array}\right)\hat{\rho}^k}{\sum_{i=0}^k
 \left(\begin{array}{c}M-1\\i
 \end{array}\right)\hat{\rho}^i}.
\end{equation}
Notice that  $\hat{\rho}$, defined above by $\hat{\rho} = \hat{\lambda}/\mu$, is the traffic generated by
a {\bf free} customer. An interesting interpretation
of (\ref{engset}) is that the call congestion, or the blocking
probability, when there are $M$ sources is equal to the time
congestion when there are $M-1$ sources. This can be intuitively
explained as follows. Consider an arbitrary tagged source (or
customer). For this particular customer, the proportion of time it
cannot access is equal to the proportion of time the $k$ circuits
are all busy by the other $M-1$ customers. During the rest of the
time our tagged source can successfully access a circuit.

\subsubsection*{Homework \ref{secengset}.\arabic{homework}}
\addtocounter{homework}{1} \addtocounter{tothomework}{1} Perform the
derivations that lead to Eq.\@ (\ref{engset}). $~~~\Box$

\subsection{Obtaining the Blocking Probability by a Recursion}

Letting $B_i$ be the blocking probability given that the number of
 circuits (servers) is $i$, the Engset loss formula can be solved
numerically by the following recursion:

\begin{equation}
\label{recurengset}
B_{i}=\frac{\hat{\rho}(M-i)B_{i-1}}{i+\hat{\rho}(M-i)B_{i-1}} ~~~~~~~~~
i=1,2,3, ~\ldots,~ k
\end{equation}
with the initial condition
\begin{equation}
\label{initial}
 B_0=1.
\end{equation}

\subsubsection*{Homework \ref{secengset}.\arabic{homework}}
\addtocounter{homework}{1} \addtocounter{tothomework}{1} Derive Eqs. (\ref{recurengset}) and (\ref{initial}).
\subsubsection*{Guide}
By (\ref{engset}) and the definition of $B_i$ we have
$$B_i=\frac{\left(\begin{array}{c}M-1\\i
 \end{array}\right)\hat{\rho}^i}{\sum_{j=0}^i
 \left(\begin{array}{c}M-1\\j
 \end{array}\right)\hat{\rho}^j}
$$
and
$$B_{i-1}=\frac{\left(\begin{array}{c}M-1\\i-1
 \end{array}\right)\hat{\rho}^{i-1}}{\sum_{j=0}^{i-1}
 \left(\begin{array}{c}M-1\\j
 \end{array}\right)\hat{\rho}^j}.
$$
Consider the ratio $B_i/B_{i-1}$ and after some algebraic manipulations
(that are somewhat similar to the derivations of the Erlang B recursion)
you will obtain
$$\frac{B_i}{B_{i-1}} = \frac{\rho(M-i)}{i}\left(1-B_i\right)$$ which leads to (\ref{recurengset}).
Notice that $B_0=1$ is equivalent to the statement that if there are no circuits (servers) (and $M>0,\hat{\rho}>0 $) the blocking probability is equal to one.
$~~~\Box$

\subsection{Insensitivity}

In his original work \cite{Engset1918}, Engset assumed that the
idle time, as well as the holding time, are exponentially
distributed. These assumptions have been relaxed in \cite{cohen57}
and now it is known that the Engset formula applies also to arbitrary
idle and holding time distributions (see also \cite{Hui90}).

\subsection{Load Classifications and Definitions}

An important feature of the Engset model is that a customer already
engaged in a conversation does not originate calls. This leads to an interesting peculiarity. If we fix the number of customers
(assuming $M>k$) and reduce $k$, the offered traffic increases
because reduction in $k$ leads to increase in $P_b$ and reduction in
the average number of busy customers which in turn leads to an increase
in idle customers each of which offers more calls, so the offered
load increases.

Let us now discuss the concept of the so-called {\em intended}
offered load \cite{Akimaru02} under the Engset settings. We know that
$1/\hat{\lambda}$ is the mean time until a free customer makes a
call (will attempt to seize a circuit). Also, $1/\mu$ is the mean
holding time of a call. If a customer is never blocked, it is
behaving like an on/off source, alternating between on and off
states, being on for an exponentially distributed period of time
with mean $1/\mu$, and being off for an exponentially distributed
period of time with mean $1/\hat{\lambda}$. For each cycle of
average length $1/\hat{\lambda}+1/\mu$, a source will be busy, on
average, for a period of $1/\mu$. Therefore, in steady-state, the
proportion of time a source is busy is $\hat{\lambda}/(\hat{\lambda}
+ \mu)$, and since we have $M$ sources, the {\em intended} offered
load is given by

\begin{equation}
\label{intend}  T = M\frac{\hat{\lambda}}{\hat{\lambda} + \mu} =
\frac{\hat{\rho} M}{(1+\hat{\rho})}.
\end{equation}

This {\em intended} offered load is equal to the offered traffic
load and the carried traffic load if $M\leq k$, namely, when
$P_b=0$. However, when $M > k$ (thus $P_b>0$), the offered traffic
load and the carried traffic load are not equal. Let $T_c$ and
$T_o$ be the {\em carried} and the {\em offered} traffic load,
respectively. The carried traffic is the mean number of busy
circuits and it is given by
\begin{equation}
\label{TC0} T_c = \sum_{i=0}^k i \pi_i.
\end{equation}
The offered traffic is obtained by averaging the intensities of the free customers
weighted by the corresponding probabilities of their numbers, as follows.
\begin{equation}
\label{TO0} T_o = \sum_{i=0}^k \hat{\rho} (M-i) \pi_i.
\end{equation}
To compute the values for $T_c$ and $T_o$ in terms of the blocking
probability $P_b$, we first realize that
\begin{equation}
\label{TC1} T_c = T_o (1-P_b),
\end{equation}
and also,
\begin{equation}
\label{TO1} T_o = \sum_{i=0}^k \hat{\rho} (M-i) \pi_i = \hat{\rho} M - \hat{\rho}
\sum_{i=0}^k i \pi_i =\hat{\rho}(M-T_c)
\end{equation}
and by (\ref{TC1}) -- (\ref{TO1}) we obtain
\begin{equation}
\label{TC} T_c = \frac{\hat{\rho} M(1-P_b)}{[1+\hat{\rho} (1-P_b)]}
\end{equation}
and
\begin{equation}
\label{TO} T_o = \frac{\hat{\rho} M}{[1+\hat{\rho} (1-P_b)]}.
\end{equation}
Notice that when  $P_b=0$, we have
\begin{equation}
\label{eqclass} T_o=T=T_c,
\end{equation}
 and when $P_b>0$, we obtain by (\ref{intend}), (\ref{TC}) and (\ref{TO}) that
\begin{equation}
\label{inequal} T_o>T>T_c.
\end{equation}
\subsubsection*{Homework \ref{secengset}.\arabic{homework}}
\addtocounter{homework}{1} \addtocounter{tothomework}{1} Using
(\ref{intend}), (\ref{TC}) and (\ref{TO}), show (\ref{eqclass}) and
(\ref{inequal}). $~~~\Box$

Notice also that
the above three measures may be divided by $k$ to obtain the
relevant traffic load per server.

\subsection{The Small Idle Period Limit}
Let $\hat{\lambda}$ approach infinity, while $M$, $k$ and $\mu$ stay fixed and assume $M>k$. Considering the steady-state equations (\ref{engsetsse1}),
their solution at the limit is $\pi_i=0$ for $i=0, ~1, ~2, \ldots, k-1$ and  $\pi_k=1$. To see that consider  the equation
$$\pi_0 M\hat{\lambda} = \pi_1 \mu.$$
Assume $\pi_0>0$, so as $\hat{\lambda} \rightarrow \infty$, we must have $\pi_1>1$ which leads to contradiction; thus, $\pi_0=0$, repeating the same
argument for the steady-state equations (\ref{engsetsse1}) for $n=1, ~2, \ldots, k-1$, we obtain that $\pi_i=0$ for $i=0, ~1, ~2, \ldots, k-1$.
Then, because  $\sum_{i=0}^k \pi_i =1$, we must have $\pi_k=1$.
Therefore by (\ref{TC0}), $$T_c=k.$$
and by (\ref{TO0}),
$$T_o = \hat{\rho} (M-k) \rightarrow \infty. $$
Intuitively, this implies that as $k$ channels (circuits) are constantly busy serving $k$ customers, the remaining $M-k$ sources (customers)
reattempt to make calls at an infinite rate.
In this case, by (\ref{intend}), the intended traffic load is
$$T =   \frac{\hat{\rho} M}{(1+\hat{\rho})}  \rightarrow M. $$
Note also that under this condition, by (\ref{TC1}), noticing that $T_o \rightarrow \infty$, and that $T_c=k$, we obtain that the blocking probability approaches 1.

\subsection{The Many Sources Limit}
Let $M$ approach infinity and $\hat{\lambda}$ approaches zero in a way that maintains the intended offered load constant.
In this case, since $\hat{\lambda} + \mu \rightarrow \mu$,  the limit of the intended load will take the form
\begin{equation}
 \label{rho} \lim T = M\frac{\hat{\lambda}}{\mu} =
{\hat{\rho} M}.
\end{equation}
Furthermore, under this limiting condition, the terms ${\hat{\rho}(M-i)}$, $i=1,2,3, \ldots, k$, in (\ref{recurengset}) can be substituted by
${\hat{\rho}M}$ which is the limit of the intended traffic load. It is interesting to observe that if we substitute $A={\hat{\rho} M}$ for the
${\hat{\rho} (M-i)}$ terms in  (\ref{recurengset}), equations (\ref{recur}) and (\ref{recurengset}) are equivalent.
This means that if the number of sources increases and the arrival rate of each source decreases in a way that the intended
load stays fixed, the blocking probability obtained by Engset loss formula approaches that of Erlang B formula.

\subsection{Obtaining the Blocking Probability by Successive Iterations}

In many cases, $\hat{\rho}$ is not available and instead the offered load $T_o$ is available. Then, it is convenient to
 obtain the blocking probability $P_b$ in terms of $T_o$. By Eq.\@ (\ref{TO}) we obtain,
\begin{equation}
\label{rhoTO} \hat{\rho} = \frac{T_o}{M-T_o(1-P_b)}.
\end{equation}
The latter can be used together with Eq.\@ (\ref{engset}) or
(\ref{recurengset}) to obtain $P_b$ by an iterative process. One
begins by setting an initial estimate value to $P_b$ (e.g.
$P_b=0.1$). Then, this initial estimate is substituted into Eq.\@
(\ref{rhoTO}) to obtain an estimate for $\hat{\rho}$, and then the value
you obtain for $\hat{\rho}$ is substituted in Eq.\@ (\ref{engset}),
or use the recursion  (\ref{recurengset}), to obtain another value
for $P_b$ which is then substituted in Eq.\@ (\ref{rhoTO}) to obtain
another estimate for $\hat{\rho}$. This iterative process continues
until the difference between two successive estimations of $P_b$ is
arbitrarily small.

\subsubsection*{Homework \ref{secengset}.\arabic{homework}}
\addtocounter{homework}{1} \addtocounter{tothomework}{1} Consider
the case $M=20$, $k=10$, $\hat{\lambda}=2$, $\mu=1$. Compute $P_b$
using the recursion of Eq.\@ (\ref{recurengset}). Then, compute $T_o$
and assuming $\rho$ is unknown, compute $P_b$ using the iterative
processes starting with various initial estimations. Compare the
results and the running time of your program. $~~~\Box$

\newpage
\section{State Dependent SSQ} \label{sdssq}

\setcounter{homework}{1} 

In the queueing model discussed in the previous chapter, the arrival and service rates vary based on the state of the system. In this section,
we consider a general model of a Markovian queue where the arrival and
service rates depend on the number of customers in the
system. Having this general model, we can apply it to
many systems whereby capacity is added (service rate increases)
and/or traffic is throttled back as queue size increases.

In particular, we study a model of a single-server queue in which
the arrival process is a state-dependent Poisson process. This is a
Poisson process that its rate $\lambda_i$ fluctuates based on the
queue size $i$. The service rate $\mu_i$ also fluctuates based on
$i$. That is, when there are $i$ customers in the system, the
service is exponentially distributed with parameter $\mu_i$. If
during service, before the service is complete, the number of
customers changes from $i$ to $j$ ($j$ could be either $i+1$ or
$i-1$), then the remaining service time changes to exponentially
distributed with parameter $\mu_j$. We assume that the number of
customers in the queue is limited by $k$.

This model gives rise to a birth-and-death process described in
Section \ref{ctmc}. The state-dependent arrival and service rates
$\lambda_i$ and $\mu_i$ are equivalent to the birth-and-death rates
$a_i$ and $b_i$, respectively.

Following the birth-and-death model of Section \ref{ctmc} the
infinitesimal generator for our Markovian state-dependent queue-size
process is given by

$Q_{i,i+1}=\lambda_i $ for $i=0, 1, 2, 3, \ldots, k$\\
$Q_{i,i-1}=\mu_i $ for $i= 1, 2, 3, 4, \ldots, k$\\
$Q_{0,0}=-\lambda_0 $ \\
$Q_{i,i}=-\lambda_i -\mu_i $ for $i=1, 2, 3, \ldots, k-1$\\
$Q_{k,k}=-\mu_k$.

Then, the steady-state equations $0=\mathbf{\Pi Q}$, can be written
as:
\begin{equation}
\label{1stsse} 0=-\pi_0\lambda_0  + \pi_1 \mu_1
\end{equation}
and
\begin{equation}
\label{Beqs} 0=\pi_{i-1} \lambda_{i-1} - \pi_i (\lambda_i + \mu_i) +
\pi_{i+1} \mu_{i+1} ~~{\rm for}~~ i=1,2,3, \ldots, k-1.
\end{equation}
There is an additional last dependent equation
\begin{equation}
\label{lasteq}0=\pi_{k-1}\lambda_{k-1} - \pi_k (\mu_k)
\end{equation}
which is redundant.

These balance equations can also be described by the following state transition diagram:

$$
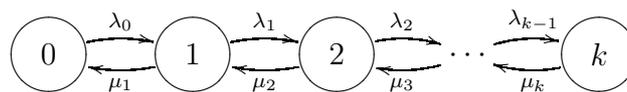
\xymatrix{
*++++[o][F]{0}\ar@/^{.5pc}/[r]^{\lambda_0} &
*++++[o][F]{1}\ar@/^{.5pc}/[l]^{\mu_1}\ar@/^{.5pc}/[r]^{\lambda_1} &
*++++[o][F]{2}\ar@/^{.5pc}/[l]^{\mu_2}\ar@/^{.5pc}/[r]^{\lambda_2} &
{\cdots}\ar@/^{.5pc}/[l]^{\mu_3}\ar@/^{.5pc}/[r]^{\lambda_{k-1}} &
*++++[o][F]{k}\ar@/^{.5pc}/[l]^{\mu_k}}$$

The normalizing equation
\begin{equation}
\label{normeq} \sum_{i=0}^k \pi_i =1 \end{equation} must also be
satisfied.

Notice that the equation
$$ 0=-\pi_0\lambda_0  + \pi_1 \mu_1$$
and the first equation of the set (\ref{Beqs}), namely,
$$ 0=\pi_{0} \lambda_{0} - \pi_1 (\lambda_1 + \mu_1) +
\pi_{2} \mu_{2}$$
gives
$$ 0= - \pi_1 \lambda_1 + \pi_2 \mu_{2}$$
which together with the second equation of the set (\ref{Beqs}), namely,
$$ 0=\pi_{1} \lambda_{1} - \pi_2 (\lambda_2 + \mu_2) +
\pi_{3} \mu_{3}$$
gives
$$ 0= - \pi_2 \lambda_2 + \pi_3 \mu_{3}$$
and in general, we obtain the set of $k$ equations:
$$ 0= - \pi_{i-1} \lambda_{i-1} + \pi_i \mu_i~~~~~~~~~i=1,2,3, \ldots, k$$
or the recursive equations:
\begin{equation} \label{sdrecur} \pi_i = \rho_i \pi_{i-1} ~~{\rm for}~~ i=1,2,3, \ldots, k \end{equation}
where $$\rho_i = \frac{\lambda_{i-1}}{\mu_i}~~ {\rm for}~~i=1,2,3, \ldots k.$$

Defining also $\rho_0 \equiv 1$, by (\ref{sdrecur}), we obtain
\begin{equation} \label{sdallthepis} \pi_i = \rho_i \rho_{i-1} \rho_{i-2} \ldots, \rho_1 \pi_0   ~~{\rm for}~~ i=0,1,2,3, \ldots, k. \end{equation}

\subsubsection*{Homework \ref{sdssq}.\arabic{homework}}
\addtocounter{homework}{1} \addtocounter{tothomework}{1}
Drive $\pi_i$ for $i=0, 1, 2 \ldots, k$.

\subsubsection*{Guide}
Summing up equations (\ref{sdrecur}) will give an equation with $1-\pi_0$ on its left-hand side and a constant times $\pi_0$ on its right-hand side.
This linear equation for $\pi_0$ can be readily solved for $\pi_0$. Having $\pi_0$, all the other $\pi_i$ can be  obtained by (\ref{sdrecur}).
 $~~~\Box$

Having obtained the $\pi_i$ values, let us derive the blocking probability.
As in the case of M/M/$k$/$k$, the proportion of time that the buffer is full is
given by $\pi_k$. However, the proportion of time that the buffer is full is not the blocking probability.
This can be easily seen in the case $\lambda_k=0$. In this case, no packets arrive when the buffer is full, so no losses occur, but we may still have
$\rho_i>0 ~~{\rm for}~~ i=1,2,3, \ldots, k,$ so $\pi_k>0$.

As in the case of the Engset model, the blocking probability is the ratio of the
number of arrivals during the time that the buffer is full to the total number of arrivals. Therefore,
\begin{equation}
\label{sdblockprob} P_b = \frac{\lambda_k \pi_k}{\sum_{i=0}^k \lambda_i \pi_i}.
\end{equation}
Notice that, as in the Engset model, the PASTA principle does not apply here since the arrival process is not Poisson.  However, if the arrival rates do not depend on the state of the
system, even if the service rates do, the arrival process becomes Poisson and the blocking probability is
equal to $\pi_k$. To see this simply set $\lambda_i=\lambda$ for all
$i$ in Eq.\@ (\ref{sdblockprob}) and we obtain $P_b =\pi_k$.

\subsubsection*{Homework \ref{sdssq}.\arabic{homework}}
\addtocounter{homework}{1} \addtocounter{tothomework}{1} Consider a single-server Markovian queue
with state-dependent arrivals and service times. You are free to choose
the $\lambda_i$ and $\mu_i$ rates, but make sure they are
different for different $i$ values. Set the buffer size at $k=200$.
Solve the steady-state equations using the method of successive substitutions (see Section \ref{succsub})
and using a standard method (e.g. using Cramer's rule). Compare the results and the
computation time. Then, obtain the blocking probability by simulation
and compare the results and the running times with the equivalent blocking probability results and running times obtained by solving the
state equations. Repeat the results for a wide range of parameters
by using various $\lambda_i$ vectors. $~~~\Box$


\newpage
\section{Queueing Models with Finite Buffers}
\label{finitebuffer}

\setcounter{homework}{1} 

We have already encountered several examples of queueing systems
where the number of customers/packets in the system is limited.
Examples include the M/M/$k$/$k$ system, the Engset system, and the
state-dependent SSQ described in the previous chapter. Given that in
real life all queueing systems have a limited capacity, it is
important to understand the performance behavior of such queues. A
distinctive characteristic of a queue with finite capacity is the
possibility of a blocking event. In practice, blocking probability
evaluation is an important performance measure of a queueing system.
For example, depending of the type of service and protocol used,
packets lost in the Internet due to buffer overflow are either
retransmitted which increases delay, or never arrive at their
destination which may adversely affect QoS perceived by users.

In a loss system such as M/M/$k$/$k$  an arriving customer may either experience blocking or is immediately admitted to service. Then, in a delay system such as the M/M/1 queue, an arriving customer  may either experience delay or is immediately admitted to service. By comparison, in a system with a finite queue discussed here, where the total number of buffer places $N$ is larger than the number of servers $k$ (which is equal to the number of buffer places at the service system), an arriving customer may encounter either one of the following three possible events: \begin{enumerate}
\item It may  be admitted to service immediately.
\item It may be placed in the queue until a server is available.
\item It may be blocked because all servers are busy and all buffer places are occupied.
\end{enumerate}

We begin the chapter by considering two extreme SSQ systems with finite
buffers. The first is a D/D/1/$N$ system where the blocking
probability is equal to zero as long as the arrival rate is not
higher than the service rate and the second one is a model where a
single large burst (SLB) arrives at time zero. We call it an
SLB/D/1/$N$ queue. In such a queue, for an arbitrarily small arrival
rate, the blocking probability approaches unity. These two extreme
examples signify the importance of using the right traffic model,
otherwise, the blocking probability estimation can be very inaccurate.
These two extreme cases will be followed by four other
cases of Markovian queues with finite buffers: the  M/M/1/$N$,
the M/M/$k$/$N$ for $N>k$, the
MMPP(2)/M/1/$N$ and the M/${\rm E_m}$/1/$N$ Queues.

\subsection{D/D/1/$N$}
As in our discussion on deterministic queues,
we assume that if an arrival and a departure occur at the same point in time, the departure
occurs before the arrival. We consider the case $N>1$ because the case of D/D/1/1 was discussed under the D/D/$k$/$k$ queue in Section \ref{ddkk}.
For the case of $\rho=\lambda/\mu < 1$, the evolution of the D/D/1/$N$, with $N>1$, is the same as that of a D/D/1 queue.
In such a case, there is never more than one customer/packet in the system, thus no losses occur and an arriving customer/packet always finds the server available to serve it. That is, out of the three possible events mentioned above, only the first event occurs in this particular case.
Let us now consider the case $\rho=\lambda/\mu > 1$. In this case, the queue reaches a persistent congestion state
where the queue size fluctuates between $N$ and $N-1$.
Recall that we assume $N>1$. In this case, whenever
a packet completes its service, there is always another customer/packet queued which enters service immediately after the previous one left the system.
Therefore, the server generates output at a constant rate of $\mu$. We also know that the arrival rate is $\lambda$. Therefore, the loss rate (namely, the number of arriving customers/packets that are lost/blocked per time unit) is
$\lambda$ - $\mu$ so the blocking probability is given by

\begin{equation}
\label{dd1kbp} P_B = \frac{\lambda-\mu}{\lambda}.
\end{equation}

In this case, where $\lambda>\mu$, and $N>1$, of the three cases mentioned above only the second and third events occur, and no customer will enter service immediately upon its arrival.

\subsection{Single Large Burst/D/1/$N$}
In this case, we have an arbitrarily large burst $L_B >> N$ [packets] arrives at time 0, and no further packets ever arrive. For this case
the blocking probability is
\begin{equation}
\label{slbd1kbp} P_B = \frac{L_B-N}{L_B}
\end{equation}
and since $L_B >> N$, we have that $P_B\approx 1$.
Notice that in this case, $L_B$ packets arrive during a period of time $T$, with $T \rightarrow \infty$, so the arrival rate approaches zero.
This case demonstrates that we can have an arbitrarily small arrival rate with a very high blocking probability.

\subsection{M/M/1/$N$}

As in the M/M/1 case, the M/M/1/$N$ queue-size process increases by
only one and decreases by only one, so it is also a birth-and-death process.
However, unlike the case of the M/M/1 birth-and-death process where the
state-space is infinite, in the case of the M/M/1/$N$ birth-and-death
process, the state space is finite and limited by the buffer size.

The M/M/1/$N$ queue is a special case of the state-dependent SSQ considered in the previous section. If we
set $\lambda_i=\lambda$ for all $i=0,1,2 ~\ldots, N-1$, $\lambda_i=0$ for all $i\geq N$
and $\mu_i=\mu$, for all $i=1,2 ~\ldots, N$, in the model of the previous section, that model
is reduced to  M/M/1/$N$.

As
$N$ is the buffer size, the infinitesimal generator for the M/M/1/$N$ queue-size process is given by

$Q_{i,i+1}=\lambda $ for $i=0, 1, 2, 3,~ \ldots, ~N-1$\\
$Q_{i,i-1}=\mu $ for $i= 1, 2, 3, 4, ~\ldots, ~N$\\
$Q_{0,0}=-\lambda $ \\
$Q_{i,i}=-\lambda -\mu $ for $i=1, 2, 3, ~\ldots,~ N-1$\\
$Q_{k,k}= -\mu.$

Substituting this infinitesimal generator in Eq.\@ (\ref{ssqij}) and
performing some simple algebraic operations, we obtain the following
steady-state
equations for the M/M/1/$N$ queue.\\
$\pi_0 \lambda = \pi_1 \mu$\\
$\pi_1 \lambda = \pi_2 \mu $\\
~\ldots~ \\
and in general:
\begin{equation}
\label{sseqmm1k} \pi_i \lambda = \pi_{i+1} \mu, ~{\rm for}
~i=0,~1,~2,~\ldots, ~N-1.
\end{equation}

These balance equations can also be described by the following state transition diagram of M/M/1/$N$:

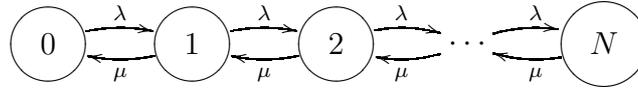
$$\xymatrix{
*++++[o][F]{0}\ar@/^{.5pc}/[r]^{\lambda} &
*++++[o][F]{1}\ar@/^{.5pc}/[l]^{\mu}\ar@/^{.5pc}/[r]^{\lambda} &
*++++[o][F]{2}\ar@/^{.5pc}/[l]^{\mu}\ar@/^{.5pc}/[r]^{\lambda} &
{\cdots}\ar@/^{.5pc}/[l]^{\mu}\ar@/^{.5pc}/[r]^{\lambda} &
*++++[o][F]{N}\ar@/^{.5pc}/[l]^{\mu}
}$$

The normalizing equation is:
\begin{equation}
\label{sumto1mm1k} \sum_{j=0}^{N} \pi_j =1.
\end{equation}

Setting $\rho=\lambda/\mu$, so we obtain,

$\pi_1 = \rho \pi_0 $\\
$\pi_2 = \rho \pi_1 = \rho^2 \pi_0$\\
$\pi_3 = \rho \pi_2 = \rho^3 \pi_0$

and in general:

\begin{equation}
\label{rhoandpiis} \pi_i = \rho^i \pi_0  ~{\rm for}
~i=0,~1,~2,~\dots N.
\end{equation}

Summing up both sides of (\ref{rhoandpiis}), we obtain (for the case $\rho \neq 1$)

\begin{equation}
\label{rhoandpiis2} 1 = \sum_{i=0}^N \rho^i \pi_0  =
\pi_0\frac{1-\rho^{N+1}}{1-\rho}.
\end{equation}
Therefore,
\begin{equation}
\label{rhoandpi0} \pi_0 = \frac{1-\rho}{1-\rho^{N+1}}.
\end{equation}
Substituting the latter in (\ref{rhoandpiis}), we obtain (for the case $\rho \neq 1$)
\begin{equation}
\label{rhoandpiisnew} \pi_i = \rho^i \frac{1-\rho}{1-\rho^{N+1}}
~{\rm for} ~i=0,~1,~2,~\dots N.
\end{equation}
Of particular interest is the blocking probability $\pi_N$ given by
\begin{equation}
\label{pkmm1k} \pi_N = \rho^N \frac{1-\rho}{1-\rho^{N+1}} =
\frac{\rho^N-\rho^{N+1}}{1-\rho^{N+1}}=\frac{\rho^N(1-\rho)}{1-\rho^{N+1}}.
\end{equation}
Notice that since M/M/1/$N$ has a finite state space, stability is
assured even if $\rho>1$.

Considering the probabilities of the three events encountered by an arriving customer, the steady state probability of the first and the third events are given by $\pi_0$ and $\pi_N$, respectively, and the probability of the second event, also called the {\it delay probability} is given simply by $1-\pi_0-\pi_N$.

Clearly, if $N=1$, the delay probability is equal to zero, because for $N=1$: $1-\pi_0-\pi_N = 1-\pi_0-\pi_1=0$.

\subsubsection*{Homework \ref{finitebuffer}.\arabic{homework}}
\addtocounter{homework}{1} \addtocounter{tothomework}{1}
Complete the above derivations for the case $\rho=1$, noticing that equation (\ref{rhoandpiis2}) for this case is:
$$ 1 = \sum_{i=0}^N \rho^i \pi_0  =
\pi_0 (N+1).$$ Alternatively, use the L'Hopital rule to obtain the limit:
$$\lim_{\rho \rightarrow 1} \frac{1-\rho^{N+1}}{1-\rho}.$$ Make sure that the results are consistent. $~~~\Box$

\subsubsection*{Homework \ref{finitebuffer}.\arabic{homework}}
\addtocounter{homework}{1} \addtocounter{tothomework}{1}
Consider an M/M/1/5 queueing model with $\lambda$  being the arrival rate and  $\mu$ being the service rate of the server. Assume that steady-state conditions hold.

\begin{enumerate}
\item Write the steady state equations.

Solve parts 2-7 for the case  $\lambda = \mu = 1$, and provide numerical answers.

\item  Solve the steady-state equations.
\item  Find the blocking probability.
\item  Find the probability that the system is empty.
\item  Find the mean number of customers in the system (including those in the queue and in service).
\item  Find the mean delay of the customers that are not blocked.
\item  Find the mean delay of delayed customers that are not blocked.
\end{enumerate}

\subsubsection*{Solution}
1.
$$\pi_0 \lambda = \pi_1 \mu$$
$$\pi_1 \lambda = \pi_2 \mu$$
$$\pi_2 \lambda = \pi_3 \mu$$
$$\pi_3 \lambda = \pi_4 \mu$$
$$\pi_4 \lambda = \pi_5 \mu$$
and the normalizing equation
$$\sum_{i=0}^5 \pi_i = 1$$

2.
$$\pi_1  = \pi_0 $$
$$\pi_2 = \pi_1 =  \pi_0 $$
$$\pi_3  = \pi_2 = \pi_1 =  \pi_0  $$
$$ \pi_4 = \pi_3  = \pi_2 = \pi_1 =  \pi_0  $$
$$ \pi_5 = \pi_4 =\pi_3  = \pi_2 = \pi_1 =  \pi_0  $$
$$\pi_5 + \pi_4 + \pi_3  + \pi_2 + \pi_1 + \pi_0 =1 $$
$$6\pi_0 = 1$$
$$\pi_0 = \frac{1}{6}$$
$$\pi_5 = \pi_4 = \pi_3  = \pi_2 = \pi_1 =  \pi_0 = \frac{1}{6}. $$

3. The blocking probability is
$$ \pi_5 = \frac{1}{6}. $$

4. The probability that the system is empty is
$$ \pi_0  = \frac{1}{6}.  $$
5. The mean number of customers in the system is
$$E[Q] = \sum_{i=0}^5 i\pi_i = (1+2+3+4+5)\times  \frac{1}{6}  = \frac{15}{6}  =  \frac{5}{2} = 2.5.  $$
6. Let $E[D]$ be the mean delay of the customers that are not blocked.
$$E[D]=\frac{E[Q]}{\lambda(1-\pi_5)}  = \frac{2.5}{1\times (1-\frac{1}{6}  )}=3. $$
7. Let $E[\hat{D}]$ be the mean delay of delayed customers that are not blocked.
By the law of iterated expectations:
$$E[D] = \frac{\pi_0 }{1-\pi_5}\times \frac{1}{\mu} + \frac{ \pi_1 + \pi_2 + \pi_3 + \pi_4 }{1-\pi_5}E[\hat{D}].$$
$$3 = \frac{1/6}{1-1/6}\times \frac{1}{1} + \frac{4/6}{1-1/6}E[\hat{D}].$$
$$3 = \frac{1}{5} + \frac{4}{5}E[\hat{D}].$$
$$15 = 1 + 4E[\hat{D}].$$
$$ E[\hat{D}] = \frac{14}{4}=3.5. ~~~\Box $$

\subsubsection*{Homework \ref{finitebuffer}.\arabic{homework}}
\addtocounter{homework}{1} \addtocounter{tothomework}{1} Consider an M/M/1/2 queue. \begin{enumerate}
\item Show that the blocking probability is equal to $E[N_Q]$ i.e., the mean number of customers in the queue (excluding the one in service).
\item Derive $E[N_s]$ (the mean number of customers at the server), $E[Q]$ (the mean number of customers in the system) and show that $E[N_Q]+E[N_s]=E[Q]$.
\end{enumerate}
\subsubsection*{Guide}
First notice that $E[N_Q] = 0(\pi_0 + \pi_1) + 1(\pi_2) = \pi_2$ which is the blocking probability.

Next,

$$\pi_0 = \frac{1}{1+\rho + \rho^2} ~~~~~~~~~~ \pi_1 = \frac{\rho}{1+\rho + \rho^2} ~~~~~~~~~~~~ \pi_2 = \frac{\rho^2}{1+\rho + \rho^2}. $$

Therefore,

$$E[Q] = 0 \pi_0 + 1 \pi_1 + 2 \pi_2 = \frac{\rho+2\rho^2}{1+\rho + \rho^2}.$$

By Little's formula,
$$E[N_s] = \lambda (1-\pi_2) \frac{1}{\mu}= \rho(1-\pi_2) = \frac{\rho(1 + \rho)}{1+\rho + \rho^2}  = \frac{\rho + \rho^2}{1+\rho + \rho^2}.$$

Clearly,  $E[N_Q]+E[N_s]=E[Q]$. $~~~\Box$

\subsubsection*{A numerical solution for the M/M/1/$N$ queue steady-state
probabilities}

A numerical solution for the M/M/1/$N$ queue steady-state probabilities follows. Set an initial value for $\pi_0$ denoted
$\hat{\pi}_0$ at an arbitrary value. For example, $\hat{\pi}_0=1$;
then compute the initial value for $\pi_1$ denoted $\hat{\pi}_1$,
using the equation $\hat{\pi}_0 \lambda = \hat{\pi}_1 \mu$,
substituting  $\hat{\pi}_0=1$. Then, use your result for
$\hat{\pi}_1$ to compute the initial value for $\pi_2$ denoted
$\hat{\pi}_2$ using $\hat{\pi}_1 \lambda = \hat{\pi}_2 \mu $, etc.
until all the initial values $\hat{\pi}_N$  are obtained. To obtain
the corresponding ${\pi}_N$ values, we normalize the $\hat{\pi}_N$
values as follows.


\begin{equation}
\label{pinormalization} \pi_N = \frac{\hat{\pi}_N}{\sum_{i=0}^N \hat{\pi}_i}.
\end{equation}

\subsubsection*{Homework \ref{finitebuffer}.\arabic{homework}}
\addtocounter{homework}{1} \addtocounter{tothomework}{1} Consider an M/M/1/$N$ queue with
$N=\rho=1000$, estimate the blocking probability. {\bf Answer:} 0.999.$~~~\Box$

\subsubsection*{Homework \ref{finitebuffer}.\arabic{homework}}
\addtocounter{homework}{1} \addtocounter{tothomework}{1} A well known approximate formula
that links TCP's flow rate $R_{TCP}$ [packets/sec], its round trip time (RTT), denoted $RTT$, and
TCP packet loss rate $L_{TCP}$ is \cite{Mat97}:
\begin{equation}
\label{tcp}
R_{TCP}=\frac{1.22}{RTT \sqrt{L_{TCP}}}.
\end{equation}
Consider a model of TCP over an M/M/1/$N$. That is, consider many TCP connections with a given RTT all passing  through
a bottleneck modeled as an M/M/1/$N$ queue. Assuming that packet sizes are exponentially distributed,
estimate TCP throughput, using Equations (\ref{pkmm1k}) and (\ref{tcp}) for a given RTT, mean packet size and service rate of the M/M/1/$N$ queue.
Compare your results with those obtained by
ns2 simulations \cite{ns2}.

\subsubsection*{Guide} Use the method of iterative fixed-point solution. See \cite{fes2006} and \cite{Gar03}.
 $~~~\Box$

 \subsubsection*{Homework \ref{finitebuffer}.\arabic{homework}}
\addtocounter{homework}{1} \addtocounter{tothomework}{1} Consider a state dependent Markovian SSQ described as follows.\\
$\lambda_i=\lambda ~~{\rm for}~~ i=0,1,2, \ldots, N-1$\\
$\lambda_N=\alpha\lambda $ where $0\leq \alpha \leq 1$\\
$\mu_i=\mu ~~{\rm for}~~ i=1,2,3 \ldots, N.$\\
This represents a congestion control system (like TCP) that reacts to congestion by reducing the arrival rate.
Derive the blocking probability and compare it with that of an M/M/1/$N$ SSQ with
arrival rate of $\lambda$ and service rate of $\mu$.
$~~~\Box$

\subsection{M/M/$k$/$N$}

This Markovian queue can also be viewed as a special case of the state dependent SSQ considered in the previous section
setting $\lambda_i=\lambda$ for all $i=0,1,2 ~\ldots, N-1$, $\lambda_i=0$ for all $i\geq N$, and setting $\mu_i=i\mu$, for all $i=1,2 ~\ldots, k$, and
$\mu_i=k\mu$, for all $i=k+1,k+2 ~\ldots, N$.

We can observe that the rates
between states $0, 1, 2, \ldots, k, k+1, k+2, \ldots, N$ are the same as in the M/M/$k$ queue. Therefore, the detailed balance equations of M/M/$k$/$N$ are the same as the first $N-1$ equations of the M/M/$k$ queue. They are:

\begin{eqnarray}
          \lambda \pi_0 &=& \mu \pi_1             \nonumber \\
  \lambda \pi_1 &=& 2\mu \pi_2  \nonumber \\
        ....                     ....         \nonumber \\
  \lambda \pi_{i-1} &=& i \mu \pi_i ~~~ i \leq k  \nonumber \\
    \lambda \pi_{i-1} &=& k \mu \pi_i ~~~ i = k+1, k+2, \ldots, N.  \nonumber
\end{eqnarray}

The state transition diagram that describes these detailed balance steady-state equations is:

$$
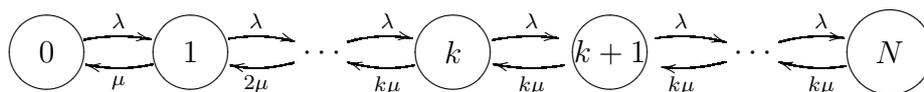
\xymatrix{
*++++[o][F]{0}\ar@/^{.5pc}/[r]^{\lambda} &
*++++[o][F]{1}\ar@/^{.5pc}/[l]^{\mu}\ar@/^{.5pc}/[r]^{\lambda} &
{\cdots}\ar@/^{.5pc}/[l]^{2\mu}\ar@/^{.5pc}/[r]^{\lambda} &
*++++[o][F]{k}\ar@/^{.5pc}/[l]^{k\mu}\ar@/^{.5pc}/[r]^{\lambda} &
*++[o][F]{k+1}\ar@/^{.5pc}/[l]^{k\mu}\ar@/^{.5pc}/[r]^{\lambda} &
{\cdots}\ar@/^{.5pc}/[l]^{k\mu}\ar@/^{.5pc}/[r]^{\lambda} &
*++++[o][F]{N}\ar@/^{.5pc}/[l]^{k\mu}
}$$

The normalizing equation is:
\begin{equation}
\label{sumto1mmkn} \sum_{j=0}^{N} \pi_j =1.
\end{equation}

Consider the notation used previously:  $$A=\frac{\lambda}{\mu},$$
and
 $$\rho=\frac{A}{k}.$$
From the detailed balance equations, we obtain
     \[  \pi_i = \left\{ \begin{array}{ll}
           \frac{A^{i}}{i!}\pi_0       &  {\rm for} ~~0 \leq i \leq k\\
\ \\
           \frac{A^{k}}{k!}\left(\frac{A}{k}\right)^{i-k}\pi_0 & {\rm for} ~~k <
i \leq N. \end{array}
                 \right.                  \]
                Notice that $\pi_k$ is the probability that there are $k$ customers in the system, namely all servers are busy and the queue is empty. This probability is given by:

 \begin{equation}
 \label{piksimple}
 \pi_k = \frac{A^{k}}{k!}\pi_0.
 \end{equation}
 Also, $\pi_N$, the probability that an arriving customer is blocked, is given by

  \begin{equation}
 \label{piNsimple}
 \pi_N=\pi_0\frac{A^{k}}{k!}\left(\frac{A}{k}\right)^{N-k} =  \pi_k \left(\frac{A}{k}\right)^{N-k}.
  \end{equation}

Summing up the $\pi_i$s, using the normalising equation, and isolating $\pi_0$, we obtain:
\begin{equation}
\label{mmknpi0}
\pi_0 = \left(\sum_{j=0}^{k-1} \frac{A^{j}}{j!} +
                    \frac{A^{k}}{k!}
\sum_{j=k}^{N}\left(\frac{A}{k}\right)^{j-k}\right)^{-1}.
\end{equation}

Summing up the second sum (the geometrical series) in (\ref{mmknpi0}), we obtain for the case $\rho \neq 1$ the following:

\begin{eqnarray*}
\sum_{j=k}^{N}\left(\frac{A}{k}\right)^{j-k} & = & \sum_{j=k}^{N} \rho ^{j-k} \\
                 & = & 1 +\rho + \rho^2 + \ldots, + \rho^{N-k} \\
                 & = & \frac{(1-\rho^{N-k+1})}{1-\rho}.
            \end{eqnarray*}
This leads, in the case of $\rho \neq 1$, to the following result for $\pi_0$.
\begin{equation}
\label{mmknpi0_2} \pi_0 = \left(\sum_{j=0}^{k-1} \frac{A^{j}}{j!} +
                    \frac{A^{k}}{k!} \frac{(1-\rho^{N-k+1})}{1-\rho} \right)^{-1}. \end{equation}

For the case $\rho = 1$, $\pi_0$ can be derived by observing that the second sum in Eq. (\ref{mmknpi0}) can be simplified, namely,
\begin{eqnarray*}
\sum_{j=k}^{N}\left(\frac{A}{k}\right)^{j-k} & = & \sum_{j=k}^{N} \rho ^{j-k} \\
                 & = & 1 +\rho + \rho^2 + \ldots, + \rho^{N-k} \\
                 & = & N-k+1.
            \end{eqnarray*}
Therefore, for $\rho=1$,

\begin{equation}
\label{mmknpi0_3} \pi_0 = \left(\sum_{j=0}^{k-1} \frac{A^{j}}{j!} +
                    \frac{A^{k}}{k!}
                    (N-k+1)\right)^{-1}.
\end{equation}

Notice also that using L'Hopital law we obtain
$$\lim_{\rho \rightarrow 1} \frac{(1-\rho^{N-k+1})}{1-\rho} = \frac{-(N-k+1)}{-1}=N-k+1$$
which is consistent with Eq. (\ref{mmknpi0_3}).

 Noticing that by Eq. (\ref{mmknpi0}), the expression for $\pi_0$ can be rewritten as

 \begin{equation}
\label{mmknpi01}
\pi_0 = \left(\sum_{j=0}^{k} \frac{A^{j}}{j!} +
                    \frac{A^{k}}{k!}
\sum_{j=k+1}^{N}\left(\frac{A}{k}\right)^{j-k}\right)^{-1}.
\end{equation}

 Then, by (\ref{mmknpi01}) and (\ref{piksimple}), in
 the case $\rho \neq 1$, we obtained

 $$\pi_k = \left(E_k^{-1}(A) + \rho \frac{(1-\rho^{N-k})}{1-\rho}\right)^{-1} $$

 where $E_k(A)=$ the Erlang B blocking probability for an M/M/$k$/$k$ system with offered traffic $A$.

For the case $\rho = 1$, we obtain

 $$\pi_k = \left(E_{k}^{-1}(A) + N-k)\right)^{-1}. $$



A call is delayed if it finds all servers
busy and there is a free place in the queue.
Notice that in our discussion on the M/M/$k$ queue, a call is delayed if it arrives and finds all servers busy and the probability of an arriving call being delayed is, for M/M/$k$, by the PASTA principle, the proportion of time all servers are busy. For the M/M/$k$/$N$ queue, there is the additional condition that the queue is not full as in such a case, an arriving call will be blocked.
Therefore, the probability that an arriving call is delayed (the delay probability) is:

\begin{eqnarray*}
 P({\rm delay}) & = & \sum_{j=0}^{N-k-1} \pi_{j+k} \\
                 & = & \pi_k \sum_{j=0}^{N-k-1} \rho^{j}  \\
                 & = & \pi_k \frac{1-\rho^{N-k}}{1-\rho}.
            \end{eqnarray*}

 We can observe that under the condition $\rho<1$, and $N \rightarrow \infty$ the M/M/$k$/$N$ reduces to  M/M/$k$. We can also observe that the M/M/1/$N$ and M/M/$k$/$k$ are also special cases of  M/M/$k$/$N$, in the instances of $k=1$ and $N=k$, respectively.

 \subsubsection*{Homework \ref{finitebuffer}.\arabic{homework}}
\addtocounter{homework}{1} \addtocounter{tothomework}{1} Show that the results of M/M/$k$, M/M/1/$N$, and M/M/$k$/$k$ for $\pi_0$ and $\pi_k$ are consistent with the results obtained of M/M/$k$/$N$. $~~~\Box$

Next we derive the mean number of customers waiting in the queue $E[N_Q]$.

\begin{eqnarray*}
 E[N_Q] & = & \sum_{j=k+1}^{N} (j-k)\pi_{j} \\
        & = & \pi_k \sum_{j=k+1}^{N} (j-k)\rho^{j-k}\\
        & = & \frac{A^{k}}{k!}\pi_0 \sum_{j=k+1}^{N} (j-k)\rho^{j-k}\\
        & = & \frac{A^{k}\rho}{k!}\pi_0 \sum_{j=k+1}^{N} (j-k)\rho^{j-k-1}\\
        & = & \frac{A^{k}\rho}{k!}\pi_0 \sum_{i=1}^{N-k} i\rho^{i-1}.
            \end{eqnarray*}
            Now, as before, we consider two cases: $\rho=1$ and $\rho \neq 1$. In the case of $\rho = 1$, we have:
            $$\sum_{i=1}^{N-k} i\rho^{i-1} = 1+2+\ldots, + N-k= \frac{(1+N-k)(N-k)}{2}.$$
            Therefore,
            $$E[N_Q]_{\rho=1}=\frac{\pi_0 A^{k}\rho(1+N-k)(N-k)}{2k!}.$$
            In the case of $\rho \neq 1$, the mean number of customers in the queue is derived as follows:
            \begin{eqnarray*}
 E[N_Q]_{\rho\neq 1} & = & \frac{\pi_0A^{k}\rho}{k!} \frac{d}{d\rho}\left(\sum_{i=0}^{N-k}\rho^i\right)\\
 & = & \frac{\pi_0A^k\rho}{k!} \frac{d}{d\rho}\left(\frac{1-\rho^{N-k+1}}{1-\rho}\right)\\
 & = & \frac{\pi_0A^k\rho[1-\rho^{N-k+1}-(1-\rho)(N-k+1)\rho^{N-k}]}{k!(1-\rho)^2}.
            \end{eqnarray*}

            As in our previous discussion on the M/M/$k$ queue, we have
\begin{equation}
\label{mmkNqsize}
E[Q] = E[N_Q] + E[N_s]
\end{equation}
and
\begin{equation}
\label{mmkNqdel}
E[D] = E[W_Q] + E[S].
\end{equation}
We know that, $$E[S]=\frac{1}{\mu}.$$
To obtain $E[N_s]$ for the M/M/$k$/$N$ queue, we again use Little's formula for the system made of servers. Recall that in the case of the M/M/$k$ queue, the arrival rate into this system was $\lambda$, but now the arrival rate should exclude the blocked customers, so now in the case of the M/M/$k$/$N$ queue the arrival rate of customers that actually access the system of servers is $\lambda(1-\pi_N)$. The
mean waiting time of each customer in that system is $E[S]=1/\mu$ (as in M/M/$k$).
Therefore, by Little's formula, the mean number of busy servers is given by
\begin{equation}
\label{ENsmmkN}
E[N_s]=\frac{\lambda(1-\pi_N)}{\mu}=A(1-\pi_N).
\end{equation}
Having $E[N_Q]$ and $E[N_s]$ we can obtain the mean number of customers in the system
$E[Q]$ by Eq. (\ref{mmkNqsize}).

Then, by Little's formula, we obtain
$$E[D] = \frac{E[Q]}{\lambda(1-\pi_N)}$$
and
$$E[W_Q] = \frac{E[N_Q]}{\lambda(1-\pi_N)}.  $$
Also, since $$E[S]=\frac{1}{\mu},$$
by (\ref{mmkNqdel}), we also have the relationship:
$$E[D] = E[W_Q] + \frac{1}{\mu}.$$

If we are interested in the mean delay of the delayed customers, denoted $E[\hat{D}] = E[D|Delayed]$, we notice that.

$$E[\hat{D}] = E[D|Delayed] = E[W_Q|Delayed] +\frac{1}{\mu},$$
where $E[W_Q|Delayed]$ is the mean waiting time in the queue of delayed customers. To obtain $E[W_Q|Delayed]$, we again use Little's formula considering the system that includes only the delayed customers, as follows:

$$E[W_Q|Delayed] = \frac{E[N_Q]}{\lambda P({\rm delay})},$$

where $P({\rm delay}) = \pi_k+\pi_{k+1}+ \ldots, + \pi_{N-1}.$

An alternative approach is to use the law of iterated expectations by solving the following equation for $E[\hat{D}]$.

$$E[D] = \frac{\pi_0 +\pi_1 + \ldots, + \pi_{k-1}}{1-\pi_N}\times \frac{1}{\mu} + \frac{\pi_k +\pi_{k+1}+ \ldots, +\pi_{N-1}}{1-\pi_{N}}\times E[\hat{D}].$$

 \subsubsection*{Homework \ref{finitebuffer}.\arabic{homework}}
Consider an M/M/3/4 queueing model with $\lambda$  being the arrival rate and $\mu$  being the service rate of each server. Assume that steady-state conditions hold. Let $A = \lambda/\mu$.
\begin{enumerate}
\item  Write the steady-state equations.

Solve items 2-7, for the case:  $\lambda = 3$ and  $\mu  = 1$.

\item Solve the steady-state equations.
\item Find the blocking probability.
\item Find the probability that the system is empty.
\item Find the mean number of customers in the system (including those in the queue and in service).
\item Find the mean number of customers in the queue (excluding the customer in service).
\item Find the mean delay of the customers that are not blocked.
\end{enumerate}

 \subsubsection*{Solutions}
1. The steady-state equations are:
$$\pi_0 \lambda = \pi_1 \mu$$
$$\pi_1 \lambda = \pi_2 2\mu$$
$$\pi_2 \lambda = \pi_3 3\mu$$
$$\pi_3 \lambda = \pi_4 3\mu.   $$
And the normalizing equation is:
$$\sum_{i=0}^4 \pi_i = 1.  $$
2.
$$(1+3+4.5+4.5+4.5)\pi_0 = 1$$
$$\pi_0 = \frac{2}{35}=0.057$$
$$\pi_1 = \frac{6}{35}=0.171$$
$$\pi_2 = \frac{9}{35}=0.257$$
$$\pi_3 = \frac{9}{35}=0.257$$
$$\pi_4 = \frac{9}{35}=0.257.   $$
3. The blocking probability is given by
$$\pi_4 = \frac{9}{35}= 0.257142857.$$
4. The probability that the system is empty is
$$\pi_0 = \frac{2}{35}=0.057142857.$$
5. The mean number of customers in the system is obtained as follows.
$$E[Q] = \sum_{i=0}^4 i\pi_i = 0\times  \frac{2}{35} ~+ ~ 1\times  \frac{6}{35} ~+~  2\times  \frac{9}{35} ~+~ 3\times  \frac{9}{35} ~+~ 4\times  \frac{9}{35} = \frac{87}{35}= 2.485714286.  $$

6.
$$E[N_Q] = 0\times (\pi_0+\pi_1+\pi_2+\pi_3) + 1\times \pi_4 = \pi_4 = 0.257142857. $$
To check the result, we use an alternative way to solve the problem.
The mean number of customers at the server is obtained by
$$E[N_s] = 0\times \pi_0 + 1\times \pi_1 +2\times \pi_2  + 3\times (\pi_3+\pi_4) $$
 $$ =  0\times  \frac{2}{35} ~+ ~ 1\times  \frac{6}{35} ~+~  2\times  \frac{9}{35} ~+~ 3\times \left( \frac{9}{35} ~+~   \frac{9}{35}\right) = 2.228571429.$$
Then,
$$E[N_Q] = E[Q] - E[N_s] = 2.485714286 - 2.228571429= 0.257142857.  $$
The results are consistent. $~~~\Box$

\subsubsection*{Homework \ref{finitebuffer}.\arabic{homework}}
\addtocounter{homework}{1} \addtocounter{tothomework}{1} Show that the lower bound for the M/M/$k$/$k$ blocking probability   obtained in (\ref{eka_bound2}), namely

$$\max \left(0,1 - \frac{k}{A} \right)$$

is also applicable to the
 M/M/$k$/$N$ queue. $~~~\Box$

\subsection{MMPP(2)/M/1/$N$}
In Section \ref{mmpp}, we described the MMPP and its two-state special case -- the MMPP(2).
Here we consider an SSQ where the MMPP(2) is the arrival process.

The MMPP(2)/M/1/$N$ Queue is an SSQ with buffer size $N$ characterized
by an MMPP(2) arrival process with parameters  $\lambda_0$,
$\lambda_1$, $\delta_0$, and $\delta_1$, and exponentially
distributed service time with parameter $\mu$. The service times are
mutually independent and are independent of the arrival process. Unlike the Poisson arrival process,
the inter-arrival times in the case of the MMPP(2) process are not independent.
As will be discussed, such dependency affects
queueing performance, packet loss and utilization.

The MMPP(2)/M/1 queue process is a continuous-time Markov chain, but
its states are two-dimensional vectors and not scalars. Each state
is characterized by two scalars: the mode $m$ of the arrival process
that can be either $m=0$ or $m=1$ and the queue size. Notice that
all the other queueing systems we considered so far were based on a
single-dimensional state-space.

Let $\pi_{im}$ for $i=0,1,2 \ldots, N$ be the probability that the
arrival process is in mode $m$ and that there are $i$ packets in
the system. After we obtain the $\pi_{im}$ values, the steady-state
queue size probabilities can then be obtained by
$$\pi_i=\pi_{i0}+\pi_{i1} ~~{\rm for} ~~i=0,1,2, \ldots, N.$$

As in the homework on IPP in Section \ref{ipphw}, the mode process itself is a two-state continuous-time
Markov chain, so the probabilities of the arrival mode being in
state $j$, denoted $P(m=j), ~~ {\rm for}~~ j=0,1$, can be solved
using the following equations:
$$ P(m=0)\delta_0 = P(m=1)\delta_1 $$
and the normalizing equation
$$ P(m=0) + P(m=1)=1. $$
Solving these two equations gives the steady-state probabilities $P(m=0)$
and $P(m=1)$ as functions of the mode duration parameters $\delta_0$ and $\delta_1$, as follows:
\begin{equation}
\label{pm0} P(m=0) =  \frac{\delta_1}{\delta_0+\delta_1}
\end{equation}
\begin{equation}
\label{pm1}
 P(m=1) =  \frac{\delta_0}{\delta_0+\delta_1}.
\end{equation}
Because the probability of the arrival process to be in mode $m$
(for $m=0,1$) is equal to $\sum_{i=0}^N \pi_{im}$, we obtain by
(\ref{pm0}) and (\ref{pm1})
\begin{equation}
\label{pm02sum}\sum_{i=0}^N \pi_{im} =
\frac{\delta_{1-m}}{\delta_{1-m}+\delta_m} ~~{\rm for} ~~m=0,1.
\end{equation}

The average arrival rate, denoted $\lambda_{av}$, is given by
\begin{equation}
\label{lambdaav} \lambda_{av} = P(m=0)\lambda_0 + P(m=1)\lambda_1 =
\frac{\delta_1}{\delta_0+\delta_1}\lambda_0 +
\frac{\delta_0}{\delta_0+\delta_1}\lambda_1.
\end{equation}

Denote $$\rho = \frac{\lambda_{av}}{\mu}.$$ The MMPP(2)/M/1/$N$ queueing
process is a stable, irreducible, and aperiodic continuous-time
Markov chain with finite state-space (because the buffer size $N$ is
finite).  We again remind the reader that the condition $\rho < 1$
is not required for stability in a finite buffer queueing system, or
more generally, in any case of a continuous-time Markov chain with
finite state space. Such a system is stable even if $\rho > 1$.

An important performance factor in queues with MMPP(2) input is the
actual time the queue stays in each mode. Even if the apportionment
of time between the modes stays fixed, the actual time can make a big
difference. This is especially true for the case
$\rho_1=\lambda_1/\mu>1$ and $\rho_2=\lambda_1/\mu<1$, or vice versa.
In such a case, if the actual time of staying in each mode is long,
there will be a long period of overload when a long queue is built
up and/or many packets are lost, followed by long periods of light
traffic during which the queues are cleared. In such a case we say
that the traffic is {\it bursty} or strongly correlated. (As mentioned above, here inter-arrival
times are not independent.) On the other hand, if the time of
staying in each mode is short; i.e., the mode process exhibits
frequent fluctuations, the overall traffic process is smoothed out
and normally long queues are avoided. To see this numerically one
could set initially $\delta_0=\delta_0^*$ $\delta_1=\delta_1^*$
where, for example, $\delta_0=1$ and $\delta_1^*=2$, or
$\delta_0^*=\delta_1^*=1$, and then set $\delta_m=\psi\delta_m^*$
for $m=0,1$. Letting $\psi$ move towards zero will mean infrequent fluctuations of the mode process that may lead to bursty traffic (long stay in each mode), and letting $\psi$ move towards infinity means frequent fluctuations of the mode process. In the exercises
below the reader is asked to run simulations and numerical
computations to obtain blocking probability and other measures for a
wide range of parameter values. Varying $\psi$ is one good way to
gain insight into performance/burstiness effects.

Therefore, the $\pi_{im}$ values can be obtain by solving the
following finite set of steady-state equations:
\begin{equation} \label{mmppsseq} 0= \mathbf{\Pi}\mathbf{Q}
\end{equation}
where $\mathbf{\Pi}$  $=$ $[\pi_{00}, \pi_{01}, \pi_{10}, \pi_{11},
\pi_{20}, \pi_{21}, \ldots, \pi_{N-1,0}, \pi_{N-1,1}, \pi_{N0},
\pi_{N1}] $, and the infinitesimal generator $2N\times 2N$ matrix is
$\mathbf{Q}$ $=$ $[Q_{\mathbf{i},\mathbf{j}}] $, where  $\mathbf{i}$
and $\mathbf{j}$ are two-dimensional vectors. Its non-zero entries
are:

\vspace{1 mm}

 $Q_{00,00} = -\lambda_0 - \delta_0$; ~~~~$Q_{00,01} =\delta_0$;~~~~ $Q_{00,10} = \lambda_0$;

\vspace{1 mm}

 $Q_{01,00} = \delta_1$;~~~~ $Q_{01,01} =-\lambda_1 - \delta_1$;~~~~ $Q_{01,11} = \lambda_1$;

\vspace{1 mm}

For $N>i>0$, the non-zero entries are:

\vspace{1 mm}

 $Q_{i0,i0} = -\lambda_0 -\delta_0-\mu$;~~~~ $Q_{i0,i1} = \delta_0$;~~~~ $Q_{i0,(i+1,0)} =\lambda_0$;

\vspace{1 mm}

$Q_{i1,i0} = \delta_1$;~~~~ $Q_{01,01} = -\lambda_1
-\delta_1-\mu$;~~~~ $Q_{i1,(i+1,1)} = \lambda_1$;

\vspace{1 mm}

and

\vspace{1 mm}

$Q_{N0,(N-1,0)} = \mu$;~~~~ $Q_{N0,N0} = -\delta_0-\mu$;~~~~
$Q_{N0,N1} = \delta_0$;

\vspace{1 mm}

$Q_{N1,(N-1,1)} = \mu$;~~~~ $Q_{N1,N1} = -\delta_1-\mu$;~~~~
$Q_{N1,N0} = \delta_1$.

\vspace{1 mm}

 In addition, we have the normalizing equation
\begin{equation}
\label{normmmpp} \sum_{i=0}^N \sum_{m=0}^1 \pi_{im} =1.
\end{equation}


The state transition diagram for the MMPP(2)/M/1/$N$ queue is:

$$
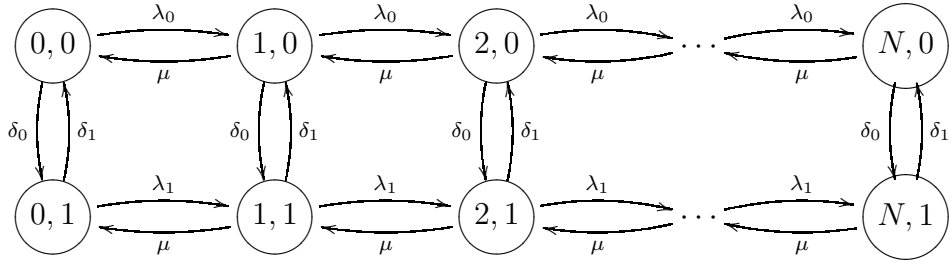
\xymatrix@C=48pt@R=36pt{
*+++[o][F]{0,0}\ar@/^{.5pc}/[r]^{\lambda_0}\ar@/_{.5pc}/[d]_{\delta_0} &
*+++[o][F]{1,0}\ar@/^{.5pc}/[l]^{\mu}\ar@/^{.5pc}/[r]^{\lambda_0}\ar@/_{.5pc}/[d]_{\delta_0} &
*+++[o][F]{2,0}\ar@/^{.5pc}/[l]^{\mu}\ar@/^{.5pc}/[r]^{\lambda_0}\ar@/_{.5pc}/[d]_{\delta_0} &
{\cdots}\ar@/^{.5pc}/[l]^{\mu}\ar@/^{.5pc}/[r]^{\lambda_0} &
*+++[o][F]{N,0}\ar@/^{.5pc}/[l]^{\mu}\ar@/_{.5pc}/[d]_{\delta_0} \\
*+++[o][F]{0,1}\ar@/^{.5pc}/[r]^{\lambda_1}\ar@/_{.5pc}/[u]_{\delta_1} &
*+++[o][F]{1,1}\ar@/^{.5pc}/[l]^{\mu}\ar@/^{.5pc}/[r]^{\lambda_1}\ar@/_{.5pc}/[u]_{\delta_1} &
*+++[o][F]{2,1}\ar@/^{.5pc}/[l]^{\mu}\ar@/^{.5pc}/[r]^{\lambda_1}\ar@/_{.5pc}/[u]_{\delta_1} &
{\cdots}\ar@/^{.5pc}/[l]^{\mu}\ar@/^{.5pc}/[r]^{\lambda_1} &
*+++[o][F]{N,1}\ar@/^{.5pc}/[l]^{\mu}\ar@/_{.5pc}/[u]_{\delta_1}
}$$

An efficient way, that normally works well for solving the steady-state equations of the\\ MMPP(2)/M/1/$N$ queue is any one of the
methods of successive substitutions (Jacobi,
Gauss-Seidel, and over-relaxation)
 discussed in Section \ref{succsub}.

Since the arrival process is not Poisson, the PASTA principle does not apply here. Therefore, to obtain the blocking probability $P_b$ we again notice that $\pi_N=\pi_{N0}+\pi_{N1}$ is the
proportion of time that the buffer is full. The proportion of packets that are lost is therefore the ratio of the number of packets
that arrive during the time that the buffer is full to the total number of packets that arrive. Therefore,
\begin{equation}
\label{mmpppb} P_b=\frac{\lambda_0\pi_{N0} +\lambda_1 \pi_{N1}}{\lambda_{av}}.
\end{equation}

As an example, we hereby provide the infinitesimal generator for
$N=2$:

\renewcommand{\arraystretch}{1.4}
\begin{tabular}{| c| c |c |c |c |c|c | }
  \hline
   & \textbf{00} & \textbf{01} & \textbf{10} & \textbf{11} & \textbf{20}  & \textbf{21}  \\ \hline
  \textbf{00} & $-\lambda_0-\delta_0$ & $\delta_0$ &  $\lambda_0$  & 0 & 0 & 0
  \\ \hline
  \textbf{01} & $\delta_1$ & $-\lambda_1-\delta_1$ & 0 & $\lambda_1$ & 0 & 0
  \\   \hline
    \textbf{10} & $\mu$ & 0  &  $-\lambda_0-\delta_0-\mu$ & $\delta_0$ &  $\lambda_0$
    & 0   \\   \hline
       \textbf{11} & 0 & $\mu$ & $\delta_1$   & $-\delta_1-\mu$ & 0 &  $\lambda_1$
     \\   \hline
        \textbf{20} & 0 & 0  & $\mu$ &0& $-\lambda_0-\delta_0-\mu$ & $\delta_0$   \\   \hline
           \textbf{21} & 0 & 0& 0& $\mu$ & $\delta_1$   & $-\delta_1-\mu$
     \\   \hline
\end{tabular}


\subsubsection*{Homework \ref{finitebuffer}.\arabic{homework}}
\addtocounter{homework}{1} \addtocounter{tothomework}{1} Consider an MMPP(2)/M/1/1 queue with
$\lambda_0=\delta_0=1$ and $\lambda_1=\delta_1=2$ and $\mu=2$.
\begin{enumerate} \item Without using a computer solve the steady-state
equations by standard methods to obtain $\pi_{00}, \pi_{01},
\pi_{10}, \pi_{11}$ and verify that (\ref{pm02sum}) holds. \item
Obtain the blocking Probability.
\item Find the proportion of time that the server is idle. \item
Derive an expression and a numerical value for the utilization.
\item Find the mean queue size. $~~~\Box$
\end{enumerate}

\subsubsection*{Homework \ref{finitebuffer}.\arabic{homework}}
\addtocounter{homework}{1} \addtocounter{tothomework}{1} Consider an
MMPP(2)/M/1/200 queue with
$\lambda_0=1$, $\delta_0=10^{-3}$, $\lambda_1=2$,
$\delta_1=2\times 10^{-3}$ and $\mu=1.9$.
\begin{enumerate} \item Solve the steady-state
equations by successive substitutions (see Section \ref{succsub}) to obtain the $\pi_{im}$ values
and verify that (\ref{pm02sum}) holds.
\item Obtain the blocking Probability.
\item Find the proportion of time that the server is idle. \item
Obtain numerical value for the utilization. \item Find the mean
queue size. \item Compare the results obtained with those obtained
before for the case $N=1$ and discuss the differences. $~~~\Box$
\end{enumerate}
\subsubsection*{Homework \ref{finitebuffer}.\arabic{homework}}
\addtocounter{homework}{1} \addtocounter{tothomework}{1} Consider again the
MMPP(2)/M/1/200 queue. Using
successive substitutions, obtain the mean queue size for a wide
range of parameter values and discuss differences. Confirm your
results by simulations with confidence intervals. Compare the
results with those obtained by successive substitutions and
simulations of an equivalent M/M/1/200 queue that has the same
service rate and its arrival rate is equal to $\lambda_av$ of the
MMPP(2)/M/1/200. Provide interpretations and explanations for all
your results. $~~~\Box$
\subsubsection*{Homework \ref{finitebuffer}.\arabic{homework}}
\addtocounter{homework}{1} \addtocounter{tothomework}{1} Consider again the
MMPP(2)/M/1/200 queue
and its M/M/1/200 equivalence. For a wide range of parameter
values, compute the minimal service rate $\mu$ obtained such that the
blocking probability is no higher than $10^{-4}$ and observe the
utilization. Plot the utilization as a function of the
parameter $\psi$ to observe the effect of burstiness on the
utilization. Confirm your results, obtained by successive
substitutions, by simulations using confidence intervals.
Demonstrate that as $\psi\rightarrow \infty$ the performance
(blocking probability and utilization) achieved approaches that of
the M/M/1/200 equivalence. Discuss and explain all the results you
obtained. $~~~\Box$

\subsection{M/${\rm E_m}$/1/$N$}

We consider here an M/${\rm E_m}$/1/$N$ SSQ model characterized by a
Poisson arrival process with parameter $\lambda$, buffer size of
$N$, and service time that has Erlang distribution with $m$ phases
(${\rm E_m}$) with mean $1/(\mu)$. Such a service time model arises
in situations when the standard deviation to mean ratio of the
service time is lower than one (recall that for the exponential
random variable this ratio is equal to one).

\subsubsection*{Homework \ref{finitebuffer}.\arabic{homework}}
\addtocounter{homework}{1} \addtocounter{tothomework}{1} Derive and plot the standard deviation to mean ratio as a function of $m$ for an ${\rm E_m}$ random variable.  $~~~\Box$

This queueing system can be analyzed using a two-dimensional state-space representing the number of customers and the number of phases still remaining to be served for the customer in service. However, it is simpler if we are able to
represent the system by a single dimension state-space. In the present case this can be done by considering the total number of phases
as the state, where each of the items (phases) in the queue is served at the rate of $m\mu$ and an arrival adds $m$ items to the queue.
The total number of items (phases) is limited to $m\times N$ because the queue size is limited to $N$ customers each of which required $m$ service phases.  Notice the one-to-one correspondence between the single dimension vector $(0, 1, 2, \ldots, m\times N)$ and the ordered set $(0, 11, 12, \ldots, 1m, 21, 22, \dots 2m, 31, \ldots, Nm)$ where the first element is 0, and the others are 2-tuples where the first is the number of customers and the second is the number of phases remains for the customer in service.

Let $\pi_i$ be the probability that there are $i$ items (phases) in the queue  $i=0,1,2,3, \ldots, m\times N$. For clarity of presentation also define $$\pi_i=0 ~~~{\rm for} ~~i<0.$$

This model is a continuous-time Markov chain, so the steady-state probabilities $\pi_i$, $i=1,2,3, \ldots, m\times N$ satisfy the following local-balance steady-state equations:

\begin{eqnarray*}
\lambda\pi_{0}  & = & m\mu\pi_{1}\\
(\lambda +m\mu)\pi_{i} & = & m\mu\pi_{i+1} + \lambda\pi_{i-m} ~~{\rm for}~~  i= 2,3, \ldots, m\times N-1.
\end{eqnarray*}
The first equation equates the probability flux of leaving state 0 (to state $m$) with the probability flux of entering state 0
only from state $1$ - where there is only one customer in the system who is in its last service phase (one item).
The second equation equates the probability flux of leaving state i (either by an arrival or by completion of the service phase) with the probability flux of
entering state i (again either by an arrival, i.e., a transition from below from state $i-m$, or from above by phase service completion from state i+1).

The probability of having $i$ customers in the system, denoted $P_i$, is
obtained by $$ P_i = \sum_{j=1}^m \pi_{(i-1)m+j}.$$ The blocking
probability is the probability that the buffer is full namely
$P_N$. The mean queue size is obtained by $$E[Q] = \sum_{i=1}^N
i\pi_{i}.$$ The mean delay is obtained by Little's formula: $$E[D] =
\frac{E[Q]}{\lambda}.$$

\subsubsection*{Homework \ref{finitebuffer}.\arabic{homework}}
\addtocounter{homework}{1} \addtocounter{tothomework}{1} Plot the state transition diagram for the M/${\rm E_m}$/1/$N$ considering the number of phases as the state.  $~~~\Box$

\subsubsection*{Homework \ref{finitebuffer}.\arabic{homework}}
\addtocounter{homework}{1} \addtocounter{tothomework}{1} Consider an M/${\rm E_m}$/1/$N$ queue.
For a wide range of parameter values (varying $\lambda,\mu,m,N$) using
successive substitutions, obtain the mean queue size, mean delay, and blocking probability and
discuss the differences. Confirm your
results by simulations using confidence intervals. Provide
interpretations and explanations to all
your results. $~~~\Box$

\subsection{Saturated Queues}
Saturated queues are characterized by having all the servers busy all the time (or almost all the time).
In such a case it is easy to estimate the blocking probability for queues with finite buffers, by simply considering the so-called {\it fluid flow model}.
Let us consider, for example, An M/M/$k$/$N$ queue, and assume that either the arrival rate  $\lambda$ is much higher
than the total service rate of all $k$ servers
$k\mu$, i.e., $\lambda>>k\mu$, or that  $\lambda>k\mu$ and $N>>0$. Such conditions will guarantee that the servers will be busy all (or most of) the time.
Since all $k$ servers are busy all the time, the output rate of the system is $k\mu$ packets/s and since the
input is $\lambda$ packets/s during a very long period of time $L$, there will be $\lambda L$ arrivals and $k\mu L$ departures.
Allowing $L$ to be arbitrarily large so that the initial transient period during which the buffer is filled can be ignored, the blocking probability can be evaluated by
\begin{equation}
\label{saturpb} P_b=\frac{\lambda L-k\mu L}{\lambda L}=\frac{\lambda-k\mu}{\lambda}=\frac{A-k}{A},
\end{equation}
where $A=\lambda/\mu$.

Another way to see (\ref{saturpb}) is by recalling that the overflow traffic is equal to the offered traffic minus the carried traffic. The offered traffic is $A$, the carried traffic in a saturated M/M/$k$/$N$ queue is equal to $k$ because all $k$ servers are continuously busy so the mean number of busy servers is equal to $k$ and the overflow traffic is equal to $AP_b$. Thus, $$A-k=AP_b$$ and (\ref{saturpb}) follows.

\subsubsection*{Homework \ref{finitebuffer}.\arabic{homework}}
\addtocounter{homework}{1} \addtocounter{tothomework}{1} Consider an
M/M/$k$/$N$ queue. Write and solve the steady-state equations to
obtain an exact solution for the blocking probability. A numerical
solution is acceptable. Validate your results by both Markov-chain
and discrete-event simulations using confidence intervals. Then, demonstrate that as $\lambda$ increases the blocking probability
approaches the result of (\ref{saturpb}). Present your results for a
wide range of parameter values (varying $\lambda,\mu,N,k$). Provide
an interpretation of your results. $~~~\Box$

\subsubsection*{Homework \ref{finitebuffer}.\arabic{homework}}
\addtocounter{homework}{1} \addtocounter{tothomework}{1} Consider again
an M/M/1/$N$ queue with
$N=\rho=1000$ and estimate the blocking probability, but this time use
the saturated queue approach. {\bf Answer:} 0.999. $~~~\Box$

\subsubsection*{Homework \ref{finitebuffer}.\arabic{homework}}
\addtocounter{homework}{1} \addtocounter{tothomework}{1}

Consider an M/M/4/6 queueing system with $\lambda$ being the arrival rate and $\mu$ being the service rate of each server. At time $T^*$, there are six customers in the system, where Steve is the customer that has arrived last and his position is last in the queue. What is the probability that among the six customers in the system at time $T^*$, Steve is the last to complete his service and to leave the entire system? In other words, you are asked to provide the probability that the other five customers will complete their service before Steve. You must provide a clear explanation for your answer.

\subsubsection*{Solution}

After the first two customers leave the system, Steve will join the service and then there are four customers in service served by the four servers. Because of the memoryless property of the exponential distribution, each of them has an equal chance to be last, so Steve will be last with probability 1/4.

In more detail, at the moment Steve joins the service, he is being served together with the three other remaining customers of the original six. The probability that he is first out of the remaining three customers (out of the original five) is 1/4. The probability he is second is (3/4)(1/3) = 1/4, The probability he is third is (3/4)(2/3)(1/2) = 1/4, so the probability that he is last is
1-(1/4)-(1/4)-(1/4) = 1/4. $~~~\Box$

\newpage \section{Processor Sharing}
\label{processorsharing}

\setcounter{homework}{1} 

In a processor sharing (PS) queueing system the server capacity is shared equally among all the customers that are present in the system. This model is applicable to a time-shared computer system where a central processor serves all the jobs present in the system simultaneously at an equal service rate. Another important application of PS is for a multiplicity of TCP connections that share a common bottleneck.
The Internet router at the bottleneck simultaneously switches
(serves) the flows generated by the users, while TCP congestion
control mechanism guarantees that the service rate obtained by the
different flows are equal. As any of the other models considered in
this book, the PS model is only an approximation for the various
real-life scenarios. It does not consider overheads and waste
associated with various real-life operations
of computer systems, and therefore it may be expected to
underestimate queueing delay.
If the server capacity to render service is $\mu$ [customers per
time-unit] and there are $i$ customers in the system, each of the
customers is served at the rate of $\mu/i$. As soon as
a customer arrives,  its service starts.

\subsection{The M/M/1-PS queue}
\label{mm1ps}

The M/M/1-PS queue is characterized by Poisson arrivals and
exponentially distributed service-time requirement, as the ordinary (FIFO) M/M/1
queue), but its service regime is assumed to be processor sharing.
In particular, we assume that the process of the
number of customers $i$ in the system is a continuous time
Markov chain, where customers arrive according to a Poisson process
with parameter $\lambda$ [customers per time-unit] and that the
service time required by an arriving customer is exponentially
distributed with parameter $\mu$.
 We also assume the stability condition of $\lambda<\mu$.

Let us now consider the transition rates of the continuous-time
Markov chain for the number of customers in the system associated
with the M/M/1-PS model. Firstly, we observe that the transition
rates from state $i$ to state $i+1$ is $\lambda$ as in the M/M/1
model. We also observe that the rates from state $i$ to state $i+j$
for $j>1$ and from state $i$ to state $i-j$ for $j>1$ are all equal
to zero (again, as in M/M/1). The latter is due to the fact that the
probability of having more than one event, arrival or departure,
occurring at the same time is equal to zero. To derive the rates from
state $i$ to state $i-1$ for $i\geq 1$ notice that at state $i$,
assuming that no arrivals occur, the time until a given customer
completes its service is exponentially distributed with rate
$\mu/i$. Therefore, the time until the first customer, out of the $i$ customers, that completes its service is the minimum of $i$
exponential random variables each of which with rate $\mu/i$, which
is exponentially distributed with rate $i(\mu/i)=\mu$. Therefore,
the transition rate from state $i$ to state $i-1$ is equal to $\mu$
(again, as in M/M/1). These imply that the process of the number of
customers in the system associated with the M/M/1-PS model is
statistically the same as the continuous-time Markov chain that
describes the M/M/1 (FIFO) queue. Therefore, the queue size
steady-state distribution $\{\pi_i\}$ and the mean queue-size $E[Q]$
given by equations (\ref{mm1}) and (\ref{meanmm1}), respectively,
 are also applied to the M/M/1-PS model. That is,

$$ \pi_i = \rho^i (1-\rho)  ~{\rm for} ~i=0,~1,~2,~\ldots~ $$

and

\begin{equation}  \label{EQPS} E[Q]=\frac{\rho}{1-\rho}. \end{equation}

 By Little's formula
the result obtained for the mean delay $E[D]$ in Eq.
(\ref{meanDmm1}) is also applicable to the M/M/1-PS model:

\begin{equation} \label{E_D}
E[D]= \frac{1}{(1-\rho)\mu}=\frac{1}{\mu-\lambda}. \end{equation}

However, the delay distribution of M/M/1 given by Eq.
(\ref{delaymm1}) does not apply to M/M/1-PS.

Having obtained the mean delay for a customer in the M/M/1-PS queue,
an interesting question is: what is the mean delay of a customer that
requires an amount of service $x$? Here $x$ represents the
time that the customer spends in the system to complete its service
assuming that there are no other customers being served and all the
server capacity can be dedicated to it. Let $E[D|x]$ be the mean delay of such a customer conditional on its service requirement $x$.

For an M/M/1 queue under the FIFO discipline, denoted in this section as M/M/1-FIFO, the
time a customer waits in the queue is not a function of its service requirement $x$ because
it depends only on the service requirements of other customers. Only
after the customer completes its waiting time in the queue, $x$ will
affect its total delay simply by being added to the waiting time in
the queue.

In particular, for M/M/1-FIFO,

\begin{equation}
\label{EDMM1s} E[D|x] = E[W_Q] + x =  \frac{\rho}{\mu - \lambda}+x,
\end{equation}
where the last equality is obtained by substituting the value of $E[W_Q]$ from (\ref{EWMM1}).

By comparison, in the case of the M/M/1-PS queue, for a customer that requires service time of $x$ from a dedicated server, its mean delay in the system from the moment it arrives until its service is complete
 has a linear relationship with $x$ \cite{kleinrock,kleinrock2}. That is,

\begin{equation} \label{linear} E[D|x] = cx, \end{equation}
for some constant $c$.
That is, under PS, the mean delay of a message has a linear relationship with its required service time.

The implication of (\ref{linear}) is that, on average, if a customer requires twice
as much service than another customer, its mean delay will be twice
that of the mean delay of the other customer.

Note that
as soon as a message arrives, its service starts. It does not need to
wait in a queue to the start of a service. However, the service rate of an individual message changes with the load.
We know that under our stability assumption, the process of the number
of customers in the system is a stable and stationary continuous-time Markov chain. In fact, it is a birth-and-death process because
the transitions are only up by one or down by one. Therefore, the
infinitesimal service rate obtained by a test customer will also
follow a stable and stationary continuous-time Markov chain. This provides an intuitive explanation of
the linear relationship of (\ref{linear}).

To obtain the parameter $c$, which is required to obtain $E[D|x]$ for M/M/1-PS, we first obtain the mean delay of an average size message  and
then invoke the linearity of (\ref{linear}). This is done as follows.

Taking the mean with respect to $x$ on both sides of (\ref{linear}) and invoking the law of iterated expectations,
we obtain

$$ E[D] = c\frac{1}{\mu}, $$

and by (\ref{E_D}) this leads to

$$ \frac{1}{(1-\rho)\mu} = c\frac{1}{\mu}.$$

Thus,

$$ c = \frac{1}{1-\rho},$$
so by the latter and (\ref{linear}), for M/M/1-PS, we obtain

\begin{equation} \label{linearx} E[D|x] = \frac{x}{1-\rho}.  \end{equation}

Comparing the two approaches PS and FIFO, PS is better for small messages and FIFO is better for large messages. Normally, the customer/message does not choose the system. But if the system administrator wants a system that does not allow jobs or customers that have very large service demands to cause extreme congestion and large delays for small ones, they use PS.

\subsubsection*{Homework \ref{processorsharing}.\arabic{homework}}
\addtocounter{homework}{1} \addtocounter{tothomework}{1}
Use equations (\ref{EDMM1s}) and  (\ref{linearx}) to plot graphs for $E[D|x]$ for M/M/1-FIFO and M/M/1-PS for a wide range of parameter values for $\rho$, $\mu$ and $x$ to compare and illustrate the effects you observe on the behavior of $E[D|x]$ for the two queueing models. $~~~\Box$

As mentioned above, while the queue-size steady-state probabilities for M/M/1-PS  and for M/M/1-FIFO are the same, which implies that the mean queue size and the mean delay (by Little's formula) for both are also equal, the delay distributions for the two queues are different. Also the variances for the two queueing models are different.

The variance for the delay for M/M/1-PS queue \cite{Ott84,Virtamo2018} is given by

\begin{equation}
\label{vardelaymm1ps}
Var[D] = \frac{1}{(\mu-\lambda)^2}  \left(\frac{2+\rho}{2-\rho}\right) = \frac{1}{\mu^2(1-\rho)^2} \left(\frac{2+\rho}{2-\rho}\right).
\end{equation}

Comparing equations (\ref{vardelaymm1}) and (\ref{vardelaymm1ps}), we observe that the variance of the delay for the M/M/1-PS queue is larger than that for M/M/1-FIFO, for $1 > \rho >0$, by the factor of
 $$ \left(\frac{2+\rho}{2-\rho}\right), $$

 which increases from 1 to 3 as $\rho$ increases from 0 to 1.

 The latter can be intuitively explained by the fact that under PS, the delay of the very long messages (or customers that require very long service time) is longer than under FIFO. This is because under PS, a very long message not only needs to wait until the messages that it finds in the system at the time of its  arrival complete their service, which is the case under FIFO; it also has to wait until many of the messages that arrive {\it after} it arrives complete their service. This effect is more pronounced as $\rho$ increases.
 Since the range of message delay values is between nearly zero (experienced by very short messages that arrive at an empty queue) and the longest delays experienced by the longest messages, the variance of the delay under M/M/1-PS is longer than under M/M/1-FIFO. This intuitive explanation can also be supported by equations (\ref{EDMM1s}) and  (\ref{linearx}). For large $x$ values, $x$ becomes the dominant term in (\ref{EDMM1s}), which is smaller than $x/(1-\rho)$ in (\ref{linearx}) in the range $0 < \rho <1$, and the difference between the two increases with the value of $\rho$ in this range. Finally, it remains to explain why the factor between the delay results of equations (\ref{vardelaymm1}) and (\ref{vardelaymm1ps}) is bounded (by 3). This is explained by noticing that $x$ is exponentially distributed, so very large $x$ values are very rare. Therefore, their effect on the variance is limited.

\subsection{Insensitivity}
\label{mm1psinsesitivity}

One important property of a processor-sharing queue is that the mean
number of customers in the system $E[Q]$ the mean delay of a customer
$E[D]$, and the mean delay of a customer with service requirement $x$, $E[D(x)]$,
given by Eqs. (\ref{EQPS}) and (\ref{E_D}), and (\ref{linearx}), respectively,
are insensitive to the shape
of the distribution of the
service-time requirements of the customers. In other words, these
results
apply also to the M/G/1-PS model characterized by Poisson arrivals,
generally distributed service-time requirements and a processor
sharing service policy. The M/G/1-PS model is a generalization of the M/M/1-PS model where we relax the exponential distribution of the service time requirements of the M/M/1-PS model, but retain the other characteristics of the M/M/1-PS model, namely, Poisson arrivals and
processor sharing service discipline.

Furthermore, the insensitivity property applies also to the
 distribution of the number of customers in the system, but not to the delay distribution.
 This means that the geometric
 distribution of the steady-state number of customers in the system of M/M/1
 applies also to the M/G/1-PS model and it is insensitive to the
 shape of the distribution of the service time requirement.
Notice that these M/M/1 results extend to the M/M/1-PS and
 M/G/1-PS models, but do not extend to the M/G/1 model. See
 discussion on the M/G/1 queue in Chapter \ref{mg1}.

Although the insensitivity applies to the distribution of the number
of customers in the M/G/1-PS model, it does not apply to the delay
distribution of M/G/1-PS.

Finally, notice the similarity between the M/G/1-PS and the M/G/$\infty$ models.
They are both insensitive to the shape of the distribution of the service time
requirement in terms of mean delay and mean number of customers in
the system. In both, the insensitivity applies to the
distribution of the number of customers in the system but does not
apply to the
 delay distribution.

\subsubsection*{Homework \ref{processorsharing}.\arabic{homework}}
\addtocounter{homework}{1} \addtocounter{tothomework}{1}

Consider packets arriving at a multiplexer where the service discipline is based on processor sharing. Assume that the service rate of the multiplexer is 2.5 Gb/s. The mean packet size is 1250 bytes. The packet arrival process is assumed to follow a Poisson process with a rate of 200,000 [packet/sec] and the packet size is exponentially distributed.
\begin{enumerate}
\item	Find the mean number of packets in the multiplexer.
\item	Find the mean delay of a packet.
\item	Find the mean delay of a packet of size 5 kbytes.
\end{enumerate}

\subsubsection*{Solution}

\underline{1. Find the mean number of packets in the multiplexer.}

mean packet size = $1250 \times 8 = 10,000$ bits

$$\mu = \frac{2,500,000,000}{10,000} = 250,000 ~~{\rm packet/s}$$

$$\rho = \frac{200,000}{250,000} = 0.8$$

$$E[Q] = \frac{0.8}{1-0.8} = 4.$$

\underline{2. Find the mean delay of a packet.}

By Little's formula

$$E[D]=\frac{E[Q]}{\lambda}=\frac{4}{200,000}=0.00002 ~~{\rm sec.} = 20 ~~{\rm microseconds}. $$

\underline{3. Find the mean delay of a packet of size 5 kbytes.}

Let $x$ be the time that the 5 kbytes packet is delayed if it is the only one in the system.

$$x = \frac{5000\times 8}{2.5 \times 1,000,000,000} = 16 ~~
{\rm microseconds}. $$

Now we will use the time units to be microseconds.

$$E[D(x)] = \frac{x}{1-\rho}= \frac{16}{1-0.8}=80 ~~{\rm microseconds}. $$

A packet four times larger than an average-sized packet will be delayed four times longer.
 $~~~\Box$

 \subsubsection*{Homework \ref{processorsharing}.\arabic{homework}}
\addtocounter{homework}{1} \addtocounter{tothomework}{1}

Assume that packets that arrive at a processor sharing system are classified into  $n$ classes of traffic, where the $i$th class is characterized by Poisson arrivals with arrival rate $\lambda_i, ~~i=1,2,~\ldots, n$, and required mean holding time (assuming a packet is alone in the system) $h_i, ~~i=1,2,~\ldots, n$. The server rate is $\mu$. \begin{enumerate}
\item Find the  mean delay of a packet.
\item Find the  mean delay of a packet that requires service time $x$.
\item Find the mean number of packets in the system.
\item Find the mean number of class $i$, $i=1, 2, ~\ldots,~ n$, packets in the system.
\item Show that the mean number of packets in the system is equal to the sum of the means obtained for classes $i$, $i=1, 2, ~\ldots,~ n$.
\end{enumerate}

\subsubsection*{Guide}
The arrival rate of all packets is given by
$\lambda = \sum_{i=1}^n \lambda_i. $
The mean holding time of a packet is given by $$h=\frac{\sum_{i=1}^n  \lambda_i h_i}{\lambda}. $$
Then, $$\rho = \lambda h = \sum_{i=1}^n  \lambda_i h_i.$$
Invoke insensitivity and use equations (\ref{E_D}) and (\ref{linearx}).

To find the mean number of packets in the system, you can use either Little's formula, or
the M/M/1 model.

Next, considering the mean delay of a packet that requires $h_i$ time if it is the only one in the system, obtain the mean system time of class $i$ packets. Then, having the arrival rate $\lambda_i$ and the mean system time of class $i$ packets, by Little's formula, obtain the mean number of class $i$ customers in the system.

Finally, with the help of some algebra, you can also show that the mean number of packets in the system is equal to the sum of the means obtained for classes $i$, $i=1, 2, ~\ldots,~ n$.
 $~~~\Box$

 \subsubsection*{Homework \ref{processorsharing}.\arabic{homework}}
 \addtocounter{homework}{1} \addtocounter{tothomework}{1}
 Consider a multiplexer with one output link that serves data messages at a rate of 100 Gb/s (Giga-bit/second). The multiplexer has a very large buffer, which is sufficient to store all data messages that arrive. There are two types of data messages.
 \begin{itemize}
\item The sizes of the messages of the first type are exponentially distributed with a mean of $10^8$ bits and they arrive at the multiplexer following a Poisson process with a rate of 576 messages per second.
\item The sizes of the messages of the second type are always exactly $10^8$ bits. This means that their standard deviation is equal to zero and they arrive at the multiplexer following a Poisson process with a rate of 384 messages per second.
 \end{itemize}
Assume that steady-state conditions hold.

Consider two service disciplines: First In First Out (FIFO) and Processor sharing (PS).  For each of these two cases (FIFO and PS), find the mean message delay from the moment a message arrives until it leaves the system (this includes queueing time and service time). Explain the differences between the results obtained for FIFO versus those obtained for PS.

 \subsubsection*{Solution}
 The mean message size is $10^8$ bits because this is the mean message size of both types.

$$\mu = \frac{10^{11}}{10^8} = 1000 ~{\rm [messages/second]}.$$

$$\lambda  = 576 + 384 = 960 ~{\rm [messages/second]}. $$

$$\rho=\frac{960}{1000}=0.96. $$
By EVVE, the variance of the length of a randomly chosen message is obtained by
$$
\frac{576}{960} \times (10^8)^2 + 0 + 0 = 0.6 \times 10^{16} ~ [{\rm bit}^2],
$$
and its standard deviation is
$$ 0.7746 \times 10^8 ~[{\rm bit}]. $$
Therefore, the standard deviation of the service time is given by
$$
\sigma_s = 0.7746 \times \frac{10^8  }{10^{11}} = 0.7746 \times 10^{-3} ~ [{\rm second}].
$$
Note that the mean message service time is one millisecond, and since the message sizes of one stream are exponentially distributed and of a second stream of messages are deterministic, and since the mean
message length of both streams is equal to each other, we can expect that the standard deviation of the service time will be lower than the mean.
Then, the variance of the service time is given by
$$
{\sigma_s}^2 = (0.7746 \times 10^{-3})^2 = 0.6 \times 10^{-6} ~ [{\rm second}^2].
$$
For the FIFO case based on the M/G/1 queueing model, the mean queue size is evaluated, by substituting  $\lambda = 960, \rho=0.96$ and $\sigma_s^2= 0.6\times 10^{-6}$  in the PK formula to obtain
$$
E[Q]=\rho + \frac{\rho^2 +
\lambda^2\sigma_s^2}{2(1-\rho)} =0.96 + \frac{(0.96)^2 +
960^2\times 0.6 \times 10^{-6} }{2(1-0.96)} = 19.39.
$$
Then, the mean delay under FIFO is obtained by Little’s formula:
$$
E[D] = E[Q]/\lambda = 19.39/960 = 0.02~ {\rm [second]}.
$$
The mean delay under M/G/1-PS is the same as the mean delay under M/M/1-PS because of the insensitivity property, which is equal to the mean delay under M/M/1.
$$E[D]=\frac{1}{\mu-\lambda}=\frac{1}{1000-960}=0.025  ~ {\rm [second]}.$$
We observe that the mean delay under processor sharing service discipline is slightly higher. This is explained by the fact that, as discussed, the standard deviation of the message size is somewhat below its mean. This indicates that the variance is lower than that of the exponentially distributed message size. In such a case, the M/G/1 model gives a lower delay than the M/M/1 (or M/M/1-PS). Intuitively, with lower variance, there are fewer cases where long messages delay many small ones under FIFO, while these small messages can be served and leave the system quicker under FIFO than under PS.
  $~~~\Box$

  \subsubsection*{Homework \ref{processorsharing}.\arabic{homework}}
  \addtocounter{homework}{1} \addtocounter{tothomework}{1}
  Consider a multiplexer with one output link that serves data messages at a rate of 9.8 Gb/s (Giga-bit/second). The multiplexer has a very large buffer which is sufficient to store all data messages that arrive. There are two types of data messages.
  \begin{itemize}
\item The sizes of the messages of the first type have a mean of 8 Mbytes and a standard deviation of 5 Mbytes and they arrive at the multiplexer following a Poisson process with rate of 20 messages per second.
\item The sizes of the messages of the second type have a mean of 16 Mbytes and a standard deviation of 15 Mbytes and they arrive at the multiplexer following a Poisson process with a rate of 60 messages per second.
  \end{itemize}
Consider the system to be in a steady state and consider two types of service disciplines: First In First Out (FIFO) and Processor sharing (PS).
\begin{enumerate}
\item For each case: FIFO and PS, what queueing model you will choose for the performance evaluation? Provide justifications.
\item Based on the model you chose for each case (FIFO and PS), find the mean message delay from the moment a message arrives until it leaves the system. You must show all steps. Numerical answers are required. Explain the differences between the results obtained for FIFO versus those obtained for PS.
\item For the case of FIFO in Part 2 compare your result to the mean message delay in the case where the arrival process is as in Part 2, but the message length is deterministic (all messages service times are equal which is equal to the mean service time obtained in Part 2).
\item Consider an individual message with a length of 5 Mbytes together with the above information.
\end{enumerate}
a.	What is its mean delay under FIFO?\\
b.	What is its mean delay under PS?
Explain the differences between the results obtained for FIFO versus those obtained for PS.

   \subsubsection*{Solution}

1.

   For FIFO, the model M/G/1 is chosen because it models Poisson arrivals, with FIFO service discipline, and general service times (non-exponential and non-deterministic).
For PS, the model M/G/1-PS is chosen because it models Poisson arrivals, with processor-sharing service discipline, and general service times (non-exponential and non-deterministic).

2.

The Mean message size is obtained by $$8\times \frac{20}{80} + 16 \times \frac{60}{80} = 14~ {\rm  Mbytes = 112~ Mbits.}$$

$\mu = 9800/112 = 87.5$ messages/second.\\ $\lambda  = 20 + 60 = 80$ messages/second. \\ $\rho=\frac{80}{87.5}=0.914. $

To obtain the variance of a randomly chosen message, we first obtain the variance of the message size by EVVE:

$$5^2\times \frac{20}{80}+15^2 \times \frac{60}{80}+(16-8)^2\times 0.25\times 0.75=187 ~{\rm   Mbyte}^2. $$

The first two terms give the mean of the conditional variances and the third is the variance of the conditional means. The conditional mean can be written as $8+8X$  where $X$  is a Bernoulli random variable with parameter $p=0.75$, and its variance is obtained by the third term.
The standard deviation of the size of a randomly chosen message is

$\sqrt{187}=13.67 $  Mbyte = 109.4 Mbits.
The standard deviation of the service time of a randomly chosen message is

$\frac{109.4}{9800}=0.011 $   sec.
so

$\sigma_s^2 = (0.011)^2 = 0.00012~~ $
For the FIFO case based on the M/G/1 queueing model, the mean queue size is evaluated, by substituting $\lambda = 80$, $\rho=0.914$, and $\sigma_s = 4 \times 10^{-5}$  in the PK formula to obtain

$$
E[Q]=\rho + \frac{\rho^2 +
\lambda^2\sigma_s^2}{2(1-\rho)} =0.914 + \frac{(0.914)^2 +
80^2\times 0.00012}{2(1-0.914)} = 10.4.
$$
Then, the mean delay under FIFO is obtained by Little’s formula:

$$
E[D] = \frac{E[Q]}{\lambda} = \frac{10.4}{80} = 0.13 {\rm   [sec.]}.
$$

The mean delay under M/G/1-PS is the same as the mean delay under M/M/1-PS because of the insensitivity property.
$$E[Q]=\frac{\rho}{1-\rho}=\frac{0.914}{1-0.914}=10.667.$$

The mean message delay is obtained by Little’s formula:
$$E[D] = \frac{E[Q]}{\lambda}=\frac{10.667}{80}=0.133~{\rm seconds}$$

As in the previous homework, we observe that the mean delay under processor sharing is slightly higher. This is, again, because the standard deviation of the message size (109.4 Mbits.) is slightly below the mean message size (112 Mbits). This indicates that the variance is slightly lower than that of the exponentially distributed message size with the same mean in an equivalent M/M/1 (or M/M/1-PS) model which leads in turn to lower mean queue size and delay.

3.

The relevant model of this case is now M/D/1, for which we obtain

$$
E[Q]=\rho + \frac{\rho^2}{2(1-\rho)} =0.914 + \frac{(0.914)^2}{2(1-0.914)} = 5.79.
$$

And the mean delay is obtained by

$$E[D] = \frac{E[Q]}{\lambda}=\frac{5.79}{80}=0.072~{\rm seconds}.   $$

The variance of the message size in the case of M/D/1 is even smaller than in the FIFO case of 6.2 as for M/D/1 this variance is equal to zero, so both the mean queue size and the mean delay are lower.

4.

(a) FIFO

$$
E[W_Q] = E[D] - \frac{1}{\mu} = 0.133 - \frac{1}{87.5} = 0.133 - 0.011 = 0.122 ~~{\rm seconds.}
$$

What is the mean delay of a message of size 5 Mbytes?

$$x=\frac{\rm  5~Mbytes}{\rm 9.8~Gb/s}=\frac{40}{9800}=0.004~ {\rm sec.} ~~~~~~~~~~E[D|x] = 0.004 + 0.122 = 0.126~ {\rm sec.}$$

(b) PS

What is the mean delay of a message of size 5 Mbytes?

$$~~~~~~E[D|x] = \frac{x}{1-\rho}=\frac{0.004}{1-0.914}=0.046 ~{\rm sec.}$$

We observe that for this short message of 5 Mbytes, the delay is shorter than that under FIFO. Despite the higher mean delay under PS, the shorter message benefits from the property of PS that discriminates in favor of short messages.   $~~~\Box$

\newpage
\section{Multi-service Loss Model}
\label{multiservice}
\setcounter{homework}{1} 

We have discussed in Section \ref{manyclasses}, a case of a Markovian
multi-server loss system ($k$ servers without additional waiting
room), involving different classes (types) of customers where customers
belong to different classes (types) may be characterized by different
arrival rates and holding times. There we assumed that each admitted
arrival is always served by a single server. We will now extend
the model to the case where customers of some classes may require
service by more than one server simultaneously. This is applicable
to a telecommunications network designed to meet heterogeneous
service requirements of different applications. For example, it is
clear that a voice call will require a lower service rate than a movie
download. In such a case, a movie download will belong to a class
that requires more servers/channels than that of the voice call. By
comparison, the M/M/$k$/$k$ system is a multi-server single-service
loss model, while here we consider a multi-server multi-service loss
model and we are interested in the blocking probability of each
class of traffic.

As the case is with the Erlang Loss System, the blocking probability is also
an important performance measure in the more general multi-service system with a finite number of servers.
However, unlike the case in the M/M/$k$/$k$ system where all customers experience the same
blocking probability, in the case of the present multi-service system, customers belonging to
different classes experience different blocking probabilities. This is intuitively clear.
Consider a system with 10 servers and assume that seven out of the
10 servers are busy. If a customer that requires one server arrives,
it will not be blocked,
but if a new arrival, that requires five servers, will be blocked. Therefore,
in many cases, customers that belong to a class that requires more servers will
experience a higher blocking probability. However, there are cases, where customers
of different classes experience the same blocking probability.
See the relevant homework question below.

This chapter covers certain key issues on multi-service models, but
it provides intuitive explanations rather than rigorous proofs. For
more extensive coverage and rigorous treatments, the reader is
referred to \cite{Iver15} and \cite{ross95} and to earlier
publications on the topic
\cite{Fortet64,Iver87,kauffman81,roberts81,ronnblum59,ross89,ross89b,ross90}.

\subsection{Model Description}

Consider a set of $k$ servers that serve arriving customers that
belong to $I$ classes. Customers from class $i$ require simultaneous
$s_i$ servers and their holding times are assumed exponentially
distributed with mean $1/\mu_i$. (As the case is for the M/M/$k$/$k$
system, the results of the analysis presented here are insensitive
to the shape of the distribution of the holding time, but since we
use a continuous-time Markov-chain modeling, this exponential
assumption is made for now.)  Class $i$ customers arrive according
to an independent Poisson process with arrival rate $\lambda_i$. The
holding times are independent of each other, of the arrival
processes, and of the state of the system.

Define
$$A_i = \frac{\lambda_i}{\mu_i}.$$

As discussed, an admitted class-$i$ customer will use $s_i$ servers
for the duration of its holding time which has a mean of $1/\mu_i$.
After its service time is complete, all these $s_i$ servers are
released and they can serve other customers. When a class-$i$
customer arrives, and cannot find $s_i$ free servers, its service is
denied and it is blocked and cleared from the
system. An important measure is the probability that an arriving
class-$i$ customer is blocked. This is called the class-$i$ customer
blocking probability denoted $B(i)$.

The following is a state transition diagram for a multi-service loss model for the case: $I=2$, $k=5$, $s_1=2$, and $s_2=1$:
$$
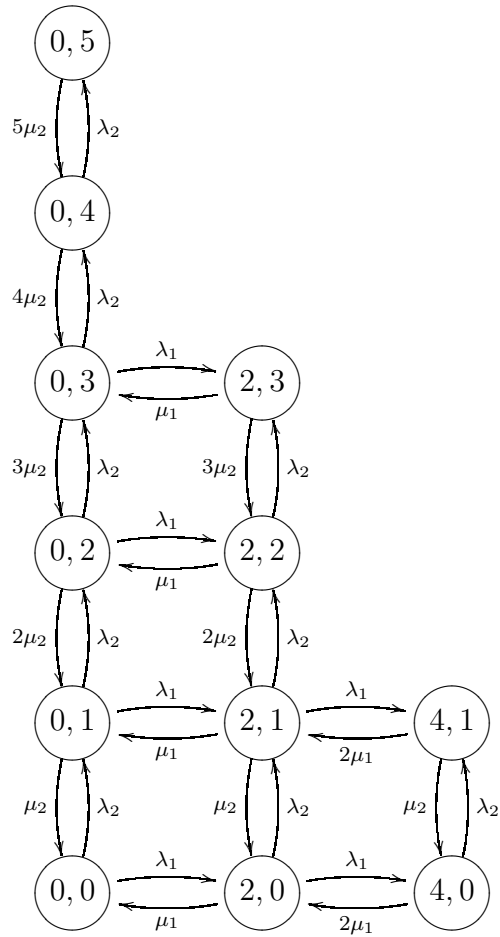
\xymatrix@=36pt{
*+++[o][F]{0,5}\ar@/_{.5pc}/[d]_{5\mu_2}\\
*+++[o][F]{0,4}\ar@/_{.5pc}/[d]_{4\mu_2}\ar@/_{.5pc}/[u]_{\lambda_2} \\
*+++[o][F]{0,3}\ar@/_{.5pc}/[d]_{3\mu_2}\ar@/_{.5pc}/[u]_{\lambda_2}\ar@/^{.5pc}/[r]^{\lambda_1} &
*+++[o][F]{2,3}\ar@/_{.5pc}/[d]_{3\mu_2}\ar@/^{.5pc}/[l]^{\mu_1}\\
*+++[o][F]{0,2}\ar@/_{.5pc}/[d]_{2\mu_2}\ar@/_{.5pc}/[u]_{\lambda_2}\ar@/^{.5pc}/[r]^{\lambda_1} &
*+++[o][F]{2,2}\ar@/_{.5pc}/[d]_{2\mu_2}\ar@/_{.5pc}/[u]_{\lambda_2}\ar@/^{.5pc}/[l]^{\mu_1}\\
*+++[o][F]{0,1}\ar@/_{.5pc}/[d]_{\mu_2}\ar@/_{.5pc}/[u]_{\lambda_2}\ar@/^{.5pc}/[r]^{\lambda_1} &
*+++[o][F]{2,1}\ar@/_{.5pc}/[d]_{\mu_2}\ar@/_{.5pc}/[u]_{\lambda_2}\ar@/^{.5pc}/[r]^{\lambda_1}\ar@/^{.5pc}/[l]^{\mu_1} &
*+++[o][F]{4,1}\ar@/_{.5pc}/[d]_{\mu_2}\ar@/^{.5pc}/[l]^{2\mu_1}\\
*+++[o][F]{0,0}\ar@/_{.5pc}/[u]_{\lambda_2}\ar@/^{.5pc}/[r]^{\lambda_1} &
*+++[o][F]{2,0}\ar@/_{.5pc}/[u]_{\lambda_2}\ar@/^{.5pc}/[r]^{\lambda_1}\ar@/^{.5pc}/[l]^{\mu_1} &
*+++[o][F]{4,0}\ar@/_{.5pc}/[u]_{\lambda_2}\ar@/^{.5pc}/[l]^{2\mu_1}
}$$

\subsection{Attributes of the Multi-service Loss Model}

The multi-service system model as defined above has the following
important attributes. \begin{enumerate} \item {\bf Accuracy and
Scalability:} The process of the number of customers in the system
of the various classes is \underline{reversible}. This property
implies that the detailed balance equations hold and together with
the normalizing equation lead to \underline{exact solution}. It is
far easier to solve the detailed balance equations than the global
balance equations and therefore the exact solution is scalable to
problems of realistic size. \item {\bf Robustness -- Insensitivity:}
The blocking probabilities depend on the customers' holding time
distributions only through their means. That is, they are
insensitive to the shape of the holding time distributions. This
\underline{insensitivity property} implies that holding time (packet
size, or flow size) can have any distribution. All that we need to know about the holding times of the calls of the various services are their means,
and the exact blocking probability for each service type can be obtained using the detailed balance equations as if the holding times follow exponential distributions.
It is known that the Internet flows follow a
heavy-tailed distribution such as Pareto. Due to this insensitivity
property, the model is robust enough to be exact even for
heavy-tailed holding time distributions. This makes the analyzes and
results of multi-service systems very relevant for real-life
telecommunications systems and networks. \item {\bf Applicability:}
Given the wide diversity of bandwidth requirements of the Internet
services and the limited capacity of communications links, there is a
clear need for a model that will provide performance evaluation in
terms of blocking probability. The M/M/$k$/$k$ which is a special
case of this model (for the case of a single service class) has been
a cornerstone in telephony used by engineers to design and dimension
telephone networks for almost a century due to its accuracy,
scalability and robustness. In telephony, we have had one type of service, i.e.,
phone calls all requiring the same link capacity. As we have entered
the Internet age, the multi-service model, given its accuracy,
scalability and robustness can play an important role. As
discussed, the insensitivity and scalability properties of the
M/M/$k$/$k$ system extends to the multi-service system model and
makes it applicable to practical scenarios. For example, a
transmission trunk or lightpath \cite{grooming} has limited capacity
which can be subdivided into many wavelength channels based on
wavelength division multiplexing (WDM) and each wavelength channel
is further subdivided into TDM sub-channels. Although the assumption
of Poisson arrivals of Internet flows during a busy hour that demand
capacity from a given trunk or a lightpath may be justified because
they are generated by a large number of sources, the actual demand
generated by the different flows/connections vary significantly from
a short SMS or email, through voice calls, to large movie downloads,
and far larger data bursts transmitted between data centers or
experimental data generated, for example, by the Large Hadron
Collider (LHC). These significant variations imply a large variety
of capacity allocated to the various flows/connections and also
large variety in their holding times, so that the restrictive
exponentially distributed holding time assumption may not be
relevant. Therefore, the insensitivity property of the multi-service
loss model is key to the applicability of the multi-service model.
\end{enumerate}

\subsection{A Simple Example with $I=2$ and $k=2$}

Consider a multi-service system with two classes of services (voice and video).
Both traffic streams of voice and video calls follow a Poisson
process and their holding times are exponentially distributed.
The arrival rate of the voice service is $\lambda_1 = 0.3$ calls per minute and the average voice service time $1/\mu_1$  is 3 minutes.
The arrival rate of the video service is $\lambda_2 = 0.2$ calls per minute and
the average video service-time $1/\mu_2$ is 5 minutes.
The system has two channels (servers).

We now aim to calculate the blocking probability of the arriving voice calls and the
arriving video calls in the case where the voice service requires one
channel per call and video service  requires two channels per call.
The system has two channels (servers), i.e., $k=2$.

Let $j_i$ be the number of channels used to serve class-$i$ customers for $i
= 1, 2$. Then, the state space is all feasible pairs $\{j_1,j_2\}$,
namely: (0,0), (1,0), (2,0), (0,2).

\subsubsection*{Homework \ref{multiservice}.\arabic{homework}}
\addtocounter{homework}{1} \addtocounter{tothomework}{1}
Plot the state transition diagram for this case with
$I=2$ and $k=2$.
 $~~~\Box$

Let $\pi_{j_1,j_2}$ be the
steady-state probability of being in state $(j_1,j_2)$. Then, we
obtain the following global balance equations.

 \begin{eqnarray*}
(\lambda_1 + \lambda_2)\pi_{0,0} &=& \mu_1 \pi_{1,0} + \mu_2 \pi_{0,2} \\
(\mu_1 + \lambda_1)\pi_{1,0} &=& \lambda_1 \pi_{0,0} +  2\mu_1 \pi_{2,0} \\
2\mu_1 \pi_{2,0} &=& \lambda_1 \pi_{1,0}  \\
  \mu_2 \pi_{0,2} &=&  \lambda_2\pi_{0,0}.
  \end{eqnarray*}

Each of these equations focuses on one state and represents the
balance of the total probability flux out and into the state. The
first equation focuses on the state (0,0), the second on (1,0), the
third on (2,0), and the fourth on (0,2).

By the first and the fourth equations we can obtain a fifth
equation:

$$  \mu_1 \pi_{1,0} =  \lambda_1\pi_{0,0}. $$

The same result is obtained by the second and third equations.

The third, fourth, and fifth equations are a complete set of detailed balance equations representing the balance of probability flux between each pair of neighboring states. These three detailed balance equations together with the normalizing equation

$$\pi_{0,0} + \pi_{1,0} + \pi_{2,0}  + \pi_{0,2}=1$$

yield a unique solution for the steady-state probabilities:
$\pi_{0,0}, \pi_{1,0}, \pi_{2,0},$ and $\pi_{0,2}$.

This shows that this 4-state multi-service system is reversible. As
the case is with the M/M/$k$/$k$ system, the physical interpretation
of the reversibility property includes
 the lost calls. For the system in the forward direction, we have multiple of Poisson processes for different types (classes) of calls, and for the system in the reversed direction, we will also have the same processes if we include as output (input in reverse) the lost calls.

The reversibility property applies also
to the general case of a multi-service system, so it is
sufficient to solve the detailed balance equations together with the
normalizing equation to obtain the steady-state probabilities of the
process.

Having obtained the steady-state probability, we can obtain the blocking probability
for the voice and for the video calls. Notice that the voice calls are only blocked
when the system is completely full. Therefore, the voice-blocking probability is:

$$ \pi_{2,0}  + \pi_{0,2}.$$

However, video calls are blocked also when there is only one channel free.
Therefore, the video-blocking probability is

$$\pi_{1,0} + \pi_{2,0}  + \pi_{0,2}.$$

Actually, in our example, the video calls can only access in state (0,0), so
the video-blocking probability is also given by

$$1 - \pi_{0,0}.$$

\subsubsection*{Homework \ref{multiservice}.\arabic{homework}}
\addtocounter{homework}{1} \addtocounter{tothomework}{1}
Compute the blocking probability of the voice calls and of the video calls for the
above small example with $I=2$ and $k=2$.
\subsubsection*{Answer}
Voice blocking probability = 0.425.\\
Video blocking probability = 0.7.
 $~~~\Box$

 \subsection{Other Reversibility Criteria for Markov chains}

 We have realized the importance of the reversibility property in simplifying
 steady-state
 solutions of Markov chains, where we can solve the simpler detailed balance equations
 and avoid the complexity of solving
the global balance equations. It is therefore important to know ways
that we can identify if a
continuous-time Markov chain is reversible.
We  use here the opportunity that is being presented
by considering the multi-service model which is an example
of the more general multi-dimensional Markov chain to discuss useful
reversibility criteria for general Markov chains that  have a wider scope of applicability that goes beyond multi-service systems.

Note that all the discussion on continuous-time Markov chains has an analogy in discrete-time Markov chains. However, we focus here on
stationary, irreducible, and aperiodic continuous-time Markov chains which is the model used for multi-service systems and is also relevant to many other models in this book. Accordingly,
whenever we mention a continuous-time Markov chain in this chapter,
we assume that it is stationary, irreducible, and aperiodic.

We already know that if the detailed balance equations together with the normalizing equation have a unique solution for the steady-state probabilities, the process is reversible. Here we describe other ways to identify if a process is reversible. However, before discussing specific reversibility criteria, we shall introduce several relevant graph theory concepts (as in \cite{kelly79}) to help visualize
relationship associated with probability flux balances.

Consider a graph $G=G(V,E)$ where $V$ is the set on vertices and $E$ is the set of
edges. We associate $G$
with a continuous-time Markov chain as follows. Let the set $V$
represent the set of states in the continuous-time Markov chain.
The graph $G$ will have an edge between nodes ${\bf x}$ and ${\bf
y}$ in $G$ if there is an edge between the two
states in the corresponding
continuous-time Markov chain, i.e., if there is
positive rate either from  the state ${\bf x}$ to state ${\bf y}$, and/or from state ${\bf y}$ to state ${\bf x}$, in the corresponding continuous-time Markov chain. We consider only cases where the continuous-time Markov chain is irreducible; therefore, the
corresponding graph must be connected \cite{kelly79}. We define a {\it cut} in the graph $G$ as a division of $G$ into two mutually exclusive sets of nodes $A$ and $\bar{A}$ such that  $A \cap \bar{A}=G$.

From the global balance equations, it can be shown (see
\cite{kelly79} for details) that for stationary
continuous-time Markov chain the probability flux across a cut
in one way is equal to the probability flux in the opposite way.

Now it is easy to show that a continuous-time Markov chain
that its graph is a tree must be reversible. Because in this case,
every edge in $G$ is a cut and therefore the probability flux across
any edge must be balanced. As a result, the detailed balanced equations hold.

Notice that all the reversible processes that we have discussed so
far, including single-dimension birth-and-death processes,
 such as the queue size processes of M/M/1,
M/M/$\infty$, and M/M/$k$/$k$, and the process associated with the above
discussed
multi-service example with $I=2$ and $k=2$
are all trees. We can therefore appreciate that
the tree criterion of reversibility is
applicable to many useful processes.
However, there are many reversible continuous-time Markov chains
that are not trees and there is a need
for further criteria to identify reversibility.

One important class of reversible processes is
the general multi-service problem with any finite $I$ and $k$.
We have already demonstrated the reversibility property for the
 small example with $I=2$ and $k=2$ that its associated graph is a tree. Let us now
 consider a slightly larger example
 where $k$ is increased from 2 to 3. All other
 parameter values are as before: $I=2$, $\lambda_1 = 0.3$, $1/\mu_1 = 3$,
 $\lambda_2 = 0.2$, and $1/\mu_2 = 5$. The associated graph of this multi-service problem is no longer
 a tree, but we already know that it is reversible because the general
 queue size process(es)
 of the multi-service model is reversible.

The detailed balance equations of this multi-service problem are:

 \begin{eqnarray*}
(i+1)\mu_1 \pi_{i+1,0} &=& \lambda_1 \pi_{i,0}, \qquad i=0,1,2.  \\
\mu_1 \pi_{1,2} &=& \lambda_1 \pi_{0,2}\\
  \mu_2 \pi_{0,2} &=&  \lambda_2\pi_{0,0}\\
    \mu_2 \pi_{1,2} &=&  \lambda_2\pi_{1,0}.
  \end{eqnarray*}

  Because the reversibility property applies
to the general case of a multi-service system,  it is
sufficient to solve the detailed balance equations together with the
normalizing equation

  $$\pi_{0,0} + \pi_{1,0} + \pi_{2,0} + \pi_{3,0}  + \pi_{0,2} + \pi_{1,2}=1.$$
This yields a unique solution for the steady-state probabilities:
$\pi_{0,0}, \pi_{1,0}, \pi_{2,0}, \pi_{3,0}, \pi_{0,2}$ and $\pi_{1,2}.$

Having obtained the steady-state probability, we can obtain the blocking probability
for the voice and for the video calls. As in the previous case,
the voice calls are only blocked
when the system is completely full. Therefore, the voice-blocking probability is:

$$ \pi_{3,0}  + \pi_{1,2}$$

and as the video calls are blocked also when there is only one channel free, the blocking probability of the video calls is

$$\pi_{3,0} + \pi_{2,0}  + \pi_{0,2} + \pi_{1,2}.$$

The associated graph of the continuous-time Markov chain that represents this multi-service problem with $k=3$ is not a tree, but this Markov chain is still reversible. Using this example with $k=3$,
we will now illustrate another criterion for reversibility called {\bf Kolmogorov criterion} that applies to a general continuous-time Markov chain and not only to those whose associated graphs are trees.
Graphs that are not trees, by definition, have cycles
and this criterion is based on conditions that apply to every
cycle in the graph that represents the Markov chain.
Furthermore, the Kolmogorov criterion has the desired feature that establishes the reversibility property directly from the given transition rates without the need to compute other results, such as steady-state probabilities.

To establish the Kolmogorov criterion, let $i$ and $j$ be two
neighboring states in a  continuous-time Markov chain and define $R(i,j)$ as the transition rate from state $i$ to state $j$. The following Theorem is known as the Kolmogorov criterion.

A stationary continuous-time Markov chain is reversible if and only if
for any cycle defined by the following
finite sequence of states $i_1, i_2, i_3, \ldots, i_n, i_1$
its transition rates satisfy:
 \begin{eqnarray}
R(i_1, i_2)R(i_2, i_3) &\ldots, & R(i_{n-1}, i_k)R(i_n, i_1)
\nonumber\\
 &=& R(i_1, i_n)R(i_n, i_{n-1}) \ldots, R(i_3, i_2)R(i_2,
 i_1).
  \end{eqnarray}

  The Kolmogorov criterion essentially says that a sufficient and necessary condition for a continuous-time Markov chain to be
  reversible is that for every cycle in the graph associated
  with the Markov chain, the product of the rates in one direction of
  the cycle starting in a given state and ending up in the same state is equal to the product of the rates in the opposite direction.

To illustrate the Kolmogorov Criterion, consider in our example with $k=3$, the circle composed of the states (0,0), (0,2), (1,2), and (1,0). According to the above-mentioned detailed balance equations, we
obtain the following rates in one direction:

 \begin{eqnarray*}
 R([0,0],[0,2]) &=& \lambda_2  \\
 R([0,2],[1,2]) &=& \lambda_1  \\
 R([1,2],[1,0]) &=& \mu_2  \\
 R([1,0],[0,0]) &=& \mu_1  \\
  \end{eqnarray*}

  and in the opposite direction:

   \begin{eqnarray*}
R([0,0],[1,0]) &=& \lambda_1  \\
R([1,0],[1,2]) &=& \lambda_2  \\
R([1,2],[0,2]) &=& \mu_1  \\
R([0,2],[0,0]) &=& \mu_2.  \\
  \end{eqnarray*}

  We can see that the product of the rates in one direction (which is $\lambda_1\lambda_2\mu_1\mu_2$) is equal to the product of the rates in the opposite direction.

\subsection{Computation}
One simple method to compute the steady-state probabilities is to
set an arbitrary initial value to one of them, to use the detailed balance
equations to obtain values for the neighbors, the neighbors' neighbors, etc. until they all have values that satisfy the detailed balance equations.
Finally, normalize all the values.

Having the steady-state probabilities, the blocking probability of all
classes can be found by adding up, for each class $i$
the steady-state probabilities of all the states where the server occupancy
is higher than $k-s_i$.

Let $\pi_{{\bf i}}$ be the steady-state probability of the being in state
$i$ after the normalization and $\hat{\pi}_{{\bf i}}$
the steady-state probability of the being in state $i$ before the normalization. Let ${\bf \Psi}$ be the set of all states. Therefore,

\begin{equation}
\label{normalization}
\pi_{\bf i} = \frac{\hat{\pi_{\bf i}}}{\sum_{\bf i \in {\bf \Psi}} \hat{\pi}_{\bf
i}}.
\end{equation}

To illustrate this approach, let again consider the above example
with $I=2$, $k=3$, $\lambda_1 = 0.3$, $1/\mu_1 = 3$,
 $\lambda_2 = 0.2$, and $1/\mu_2 = 5$.

Set $\hat{\pi_{0,0}}=1$, then
  \begin{eqnarray*}\hat{\pi_{1,0}}&=&\hat{\pi}_{0,0}\frac{\lambda_1}{\mu_1}\\
  &=& 1 \times \frac{0.3}{1/3}=0.9.
    \end{eqnarray*}
    Next,
      \begin{eqnarray*}\hat{\pi_{2,0}}&=&\hat{\pi}_{1,0}\frac{\lambda_1}{2\mu_1}\\
  &=& 0.9 \times \frac{0.3}{2/3}=0.45
    \end{eqnarray*}
    and
      \begin{eqnarray*}\hat{\pi_{3,0}}&=&\hat{\pi}_{2,0}\frac{\lambda_1}{3\mu_1}\\
  &=& 0.45 \times \frac{0.3}{3/3}=0.3.
     \end{eqnarray*}
      Moving on to the states (0,2) and (1,2), we obtain:
      \begin{eqnarray*}\hat{\pi_{0,2}}&=&\hat{\pi}_{0,0}\frac{\lambda_2}{\mu_2}\\
  &=& 1 \times \frac{0.2}{1/5}=1
    \end{eqnarray*}
and
   \begin{eqnarray*}\hat{\pi_{1,2}}&=&\hat{\pi}_{0,2}\frac{\lambda_1}{\mu_1}\\
  &=& 1 \times \frac{0.3}{1/3}= 0.9.
    \end{eqnarray*}

    To normalize we compute
  $$  {\sum_{\bf i \in {\bf \Psi}} \hat{\pi}_{\bf i}}= 1+0.9+0.45+0.3+1+0.9=4.55.   $$
  Therefore,
     \begin{eqnarray*}
\pi_{0,0} &=& \frac{1}{4.55}= 0.21978022  \\
\pi_{1,0} &=& \frac{0.9}{4.55}= 0.197802198  \\
\pi_{2,0} &=& \frac{0.45}{4.55}= 0.098901099  \\
\pi_{3,0} &=& \frac{0.3}{4.55}= 0.065934066  \\
\pi_{0,2} &=& \frac{1}{4.55}= 0.21978022  \\
\pi_{1,2} &=& \frac{0.9}{4.55}=0.197802198.
  \end{eqnarray*}

Therefore, the voice-blocking probability is:

 \begin{eqnarray*}B_{voice}&=&\pi_{3,0}+\pi_{1,2}\\
  &=& 0.065934066 + 0.197802198 = 0.263736264
    \end{eqnarray*}

    and the video-blocking probability is

     \begin{eqnarray*}B_{video}&=&\pi_{2,0}+\pi_{3,0}+\pi_{0,2}+\pi_{1,2}\\
  &=& 0.098901099+0.065934066 + 0.21978022+0.197802198 =
  0.582417582.
    \end{eqnarray*}

    We can see that reversibility makes it easier to solve for
    steady-state probabilities. However, if we consider a
    multi-service system where $I$ and $k$ are very large, it may be challenging to solve the problem in a reasonable time.

    There are two methods to improve the efficiency of the
    computation.
    \begin{enumerate}
    \item {\bf The Kaufman Roberts Algorithm:} This algorithm is based on the recursion of the number of busy
    servers. For details on this algorithm see
\cite{Fortet64,Iver15,kauffman81,roberts81,ross95}.
    \item  {\bf The Convolution Algorithm:} This algorithm is based on the aggregation of traffic
    streams. In other words, if one is interested in the blocking
    probability of traffic type $i$, the algorithm successively
    aggregates by convolution all other traffic types until we have
    a problem with $I=2$, namely, traffic type $i$ and all other
    types together. Then, the problem can be easily solved.
    For details on this algorithm, see \cite{Iver87,Iver15,ross95,ross90}.
       \end{enumerate}

\subsection{A General Treatment}
So far, we discussed the properties of multi-service systems through
simple examples. Now we present general definitions and concepts.
For notation convenience, let us consider a slightly different
Markov chain than the one we considered above. Previously we
considered the state space to represent the number of busy servers
(channels) occupied by each of the services (traffic types).
Now we will consider a continuous-time Markov chain where the state space
represents the number of customers (calls) of each traffic type
rather than the number of busy servers. The two approaches are
equivalent because at any point in time, the number of channels for service $i$
in the system is a
multiplication by a factor of $s_i$ of the number of customers
(calls) of service $i$ in the system.

Let $j_i$ be the number of class-$i$ customers in the system for $i
= 1, 2, \ldots, I$. Let $$\overrightarrow{j}=(j_1,j_2,  \ldots,
j_I)$$ and $$\overrightarrow{s}=(s_1,s_2,  \ldots, s_I).$$
 Then, $$\overrightarrow{j} \overrightarrow{s} = \sum_{i=1}^I j_i s_i$$ is the number of busy servers.
Now we consider an $I$-dimensional continuous-time Markov chain
where the state space is defined by all feasible vectors
$\overrightarrow{j}$ each of which represents a multi-dimensional possible state of the system.
In particular, we say that the state
$\overrightarrow{j}$ is feasible if
$$\overrightarrow{j} \overrightarrow{s} = \sum_{i=1}^I j_i s_i \leq k.$$
Let {\bf F} be a set of all feasible  vectors
$\overrightarrow{j}$.

A special case of this multi-service model is the M/M/$k$/$k$ model
where $I=1$. If we consider the M/M/$k$/$k$ model and let $k
\to \infty$, we obtain the M/M/$\infty$ model described in Section
\ref{secmminf}. Accordingly, the M/M/$\infty$ model is the special
case ($I=1$) of the multi-service model with $k=\infty$. In our
discussion in Section \ref{dimErlang}, the distribution of the
number of customers in an M/M/$k$/$k$ model, given by
(\ref{mmkkfin}) is a truncated version of the distribution of the
number of customers in an M/M/$\infty$ model. As we explain there,
the former distribution can be obtained using the latter by
truncation.

In a similar way, we begin by describing a multi-service system with
an infinite number of servers. Then, using truncation, we derive the
distribution of the number of customers of each class for a case
where the number of servers is finite.

\subsubsection{Infinite Number of Servers}

For the case $k=\infty$, every arrival of any class $i$ customer can always find $s_i$
free servers; therefore, this case can be viewed as $I$ independent
uni-dimensional continuous-time Markov chains, where $X_i(t), i=1,2,
\ldots, I$, represents the evolution of the number of  class-$i$
customers in the system and characterized by the birth rate
$\lambda_i$ and the death-rate $j_i\mu_i$.

Let $\pi_i(j_i)$ be the steady-state probability of the process
$X_i(t), i=1,2, \ldots, I$ being in state $j_i$. Then, $\pi_i(j_i)$
satisfy the following steady-state equations:

$$\lambda_i \pi_i(0) = j_i \mu_i \pi_i(1)$$
$$\lambda_i \pi_i(j_i) = j_i \mu_i \pi_i(j_i+1) {\rm for} ~~{j_i=1,2,3, \ldots} $$
and the normalizing equation
$$\sum_{j_i=0}^\infty \pi_i(j_i) =1.$$
These equations are equivalent to the equations that represent the
steady-state equations of the M/M/$\infty$ model. Replacing $n$ for
$j_i$, $\lambda$ for $\lambda_i$, and $\mu$ for $j_i \mu_i$ in the
above equations, we obtain the M/M/$\infty$ steady-state equations.
This equivalence has also a physical interpretation. Simply consider
a group of $s_i$ servers as a single server serving each class-$i$
customer. Following the derivations is Section \ref{secmminf} for
the M/M/$\infty$ model, we obtain:
\begin{equation}
\pi_i(j_i) = \frac{e^{-A_i} A_i^{j_i} }{j_i! } ~{\rm for}
~j_i=0,~1,~2,~\ldots~.
\end{equation}

Since the processes $X_i(t),~ i=1,2, \ldots, I$, are independent, the
probability $p(\overrightarrow{j})=p(j_1,j_2, \ldots, j_I)$ that in steady-state
$X_1(t)=j_1$, $X_2(t)=j_2$,  $ \ldots$, $X_I(t)=j_I$, is given by
\begin{equation}
\label{infinitek} p(\overrightarrow{j}) = p(j_1,j_2, \ldots, j_I) =
\prod_{i=1}^{I} \frac{e^{-A_i} A_i^{j_i }}{j_i! } e^{-A_i}.
\end{equation}

The solution for the steady-state joint probability distribution of
a multi-dimensional process, where it is obtained as a product of
steady-state distribution of the individual single-dimensional
processes, such as the one given by (\ref{infinitek}), is called a
{\it product-form solution}.

A simple example to illustrate the product-form solution is to
consider a two-dimensional multi-service loss system with $k=\infty$, and to observe
that to satisfy the detailed balance equations,
the steady-state probability of the state $(i,j)$ $\pi_{ij}$ is
the product of $$\pi_{0j} = \pi_{00} \frac{A_2^j}{j!}$$
and $$\frac{A_1^i}{i!}.$$

Then, realizing that $$\pi_{00} = \pi_{0}(1)\pi_{0}(2)$$
where $\pi_{0}(1)$ and $\pi_{0}(2)$ are the probabilities that
the independent systems of services 1 and 2 are empty, respectively.
Thus,
 \begin{eqnarray*}\pi_{ij}&=&\pi_{00} \frac{A_1^i}{i!}\frac{A_2^j}{j!}\\
  &=& \pi_{0}(1)\pi_{0}(2) \frac{A_1^i}{i!}\frac{A_2^j}{j!} \\
   &=& \left(\pi_{0}(1) \frac{A_1^i}{i!}\right) \left(\pi_{0}(2)\frac{A_2^j}{j!} \right)
    \end{eqnarray*}
and the product form has directly been obtained from the detailed balanced equations.
This illustrates the relationship between reversibility and product form
solution.

Next, we consider a system with a finite number of servers,
observe  that for such a system, the detailed balance equations also
gives a product form solution because the equation
$$\pi_{ij}=\pi_{00} \frac{A_1^i}{i!}\frac{A_2^j}{j!}$$
which results directly from the detailed balance equation still
holds. Note that $\pi_{00}$ is not the same in the infinite and
finite $k$ cases, and it is normally different for different $k$ values.

\subsubsection*{Homework \ref{multiservice}.\arabic{homework}}
\addtocounter{homework}{1} \addtocounter{tothomework}{1} Provide an
example where $\pi_{00}$ is the same for different $k$ values. $~~~\Box$

\subsubsection{Finite Number of Servers}

Consider a multi-service system model where the number of
servers is limited to $k$. We are interested in the probability $B(m)$
that a class $m$ customer is blocked. We begin by deriving the state probability vector
$p(\overrightarrow{j})$ for all $\overrightarrow{j} \in {\bf F}$. By the definition of
conditional probability, $p(\overrightarrow{j})$ conditional on $\overrightarrow{j} \in {\bf F}$
is given by
\begin{equation}
\label{finitek}
p(\overrightarrow{j}) = p(j_1,j_2, \ldots, j_I) = \frac{1}{C} \prod_{i=1}^{I} \frac{e^{-A_i} A_i^{j_i
}}{j_i! }~~ \overrightarrow{j} \in {\bf F} \end{equation}
where
$$ C = \sum_{\overrightarrow{j} \in {\bf F}}  \prod_{i=1}^{I} \frac{e^{-A_i} A_i^{j_i
}}{j_i! }.$$

\subsubsection*{Homework \ref{multiservice}.\arabic{homework}}
\addtocounter{homework}{1} \addtocounter{tothomework}{1} Derive
(\ref{finitek}) by truncating (\ref{infinitek}).
\subsubsection*{Guide}
Consider the steady-state probability distribution of
$\overrightarrow{j}$ for the case $k=\infty$ give by
(\ref{infinitek}). Then, set $p(\overrightarrow{j})=0$ for all
$\overrightarrow{j}$ not in  {\bf F}, and normalize the
probabilities $\overrightarrow{j} \in {\bf F}$ by dividing by them
by the probability that the infinite server process is in a feasible
state considering that the number of servers $k$ is finite. Then, cancel out the exponentials and obtain (\ref{finitek}). $~~~\Box$

Let ${\bf F}(m)$ be the subset of the states in which an arriving class $m$
customer will not be blocked. That is
\begin{equation}
\label{admitting_states}
{\bf F}(m) = \{ \overrightarrow{j} \in {\bf F} ~{\rm such ~that} ~\sum_{i=1}^I s_i j_i \leq k - s_m \}.
\end{equation}
Then, \begin{equation}
\label{bkblocking}
B(m) = 1 - \sum_{\overrightarrow{j}\in {\bf F}(m) } p(\overrightarrow{j}), ~m = 1,2, \ldots, I.
\end{equation}
Therefore, by (\ref{finitek}), we obtain
\begin{equation}
\label{bkblocking_explicit}
B(m) = 1 -  \frac{\displaystyle  \sum_{\overrightarrow{j}\in {\bf F}(m) } \prod_{i=1}^{I} \frac{ A_i^{j_i
}}{j_i! }}{\displaystyle \sum_{\overrightarrow{j}\in {\bf F} } \prod_{i=1}^{I} \frac{ A_i^{j_i
}}{j_i! }}.
\end{equation}

\subsubsection*{Homework \ref{multiservice}.\arabic{homework}}
\addtocounter{homework}{1} \addtocounter{tothomework}{1}
Consider the case with $k=3$, $s_1=1$, $s_2=2$, $\lambda_1 = \lambda_2 = 1$,
and $\mu_1=\mu_2=1$. Find the Blocking probabilities $B(1)$ and $B(2)$.
\subsubsection*{Guide}
Let $(i,j)$ be the state in which there are $i$ class 1 and $j$ class 2 customers in the system.

The Set {\bf F} in this example is given by $${\bf F}=\{(0,0),~
(0,1),~ (1,0),~ (1,1),~ (2,0),~ (3,0)\}.$$ Write and solve the
steady-state equations for the steady-state probabilities of the
states in the set  {\bf F}. Alternatively, you can use
(\ref{finitek}).

Then, $${\bf F}(1)=\{(0,0),~ (0,1),~ (1,0),~ (2,0)\}.$$
and
$${\bf F}(2)=\{(0,0),~ (1,0)\}.$$

Use (\ref{bkblocking_explicit}) to obtain the blocking probability. $~~~\Box$

\subsection{Critical Loading}
\label{critloadMS}
As discussed in Section \ref{critload}, a critically loaded system is a one
where the offered traffic load is equal to the system capacity.
Accordingly, in a critically loaded multi-service loss system, the following
condition holds
\begin{equation}
\sum_{i=1}^I A_i = k.
\end{equation}
Given the tremendous increase in the capacity of telecommunications networks and
systems and in the number of human and non-human users of the Internet, the
case of large $k$ is of particular interest. As we have learned in the case of M/M/$k$/$k$, when the total capacity of the system is very large relative to the capacity required by any individual user, critical loading is an efficient dimensioning rule. The result for the asymptotic behavior of the blocking probability under the critical loading condition can be extended to the case of a multi-service loss system as follows.

\begin{equation} \label{critloadedMS} \lim_{k \rightarrow \infty}
B(i) = \frac{s_i {C_{MS}}}{\sqrt{k}}, ~~ i=1,2, ~\ldots, ~I
\end{equation} where ${C_{MS}}$ is a constant independent of $i$ and
$k$. Notice that if there is only one class ($I=1$) and $s_1=1$,
this asymptotic result reduces to (\ref{limiterlang}) by setting
${C_{MS}}=\tilde{C}$. Notice that as in the special case of
M/M/$k$/$k$, the asymptotic blocking probability decays at the rate
of $1/{\sqrt{k}}$, and also notice that the asymptotic class $i$
blocking probability is linear with $s_i$. This means that in the
limit, if each of the class 1 customers requires one server and each of the class 2 customers requires two servers, then a class 2 customer
will experience twice the blocking probability experienced by a
class 1 customer. Recall that, in this case, a class 1 customer
requires only one server to be idle for it to be able to access a
server and to obtain service, while a class 2 customer requires two
idle servers to obtain service otherwise, according to our
multi-service loss model, it is blocked and cleared from the system.

\subsubsection*{Homework \ref{multiservice}.\arabic{homework}}
\addtocounter{homework}{1} \addtocounter{tothomework}{1} Consider
the case $\lambda_1 = 1$, $s_1=1$, $\mu_1=1$,  $\lambda_2 = 2$,
$s_1=2$, $\mu_2=2$, $k = 4$. Obtain the blocking probability of each
class in two ways: (1) by a discrete event simulation, (2) by
solving the steady-state equations or (\ref{finitek}) and using
Eq.\@ (\ref{bkblocking}), and (3) by using the recursive algorithm.
$~~~\Box$
\subsubsection*{Homework \ref{multiservice}.\arabic{homework}}
\addtocounter{homework}{1} \addtocounter{tothomework}{1} Provide examples where customers
that belong to different classes experience the same blocking probability. Verify the equal blocking probability
using (\ref{finitek}), by the recursive algorithm.  and by simulations. \subsubsection*{Guide}
One example is with $k=6$, and two classes of customers $s_1=6$ and $s_2=5$.
Provide other examples and verify the equal blocking probability
using the analysis that leads to (\ref{bkblocking}) and simulations.
$~~~\Box$
\subsubsection*{Homework \ref{multiservice}.\arabic{homework}}
\addtocounter{homework}{1} \addtocounter{tothomework}{1} Demonstrate by simulations
the robustness of the multi-service loss model
to the shape of the holding time distribution. \subsubsection*{Guide}
Simulate various multi-service loss systems with exponential holding time versus equivalent
systems where the holding times distributions are
hyper-exponential (the variance is larger than exponential), deterministic (where the variance is
equal to zero), and Pareto (choose cases where the valiance is finite). Demonstrate that the
blocking probability for each class is the same when the mean holding time is the same
regardless of the choice of the holding time distribution.
$~~~\Box$
\subsubsection*{Homework \ref{multiservice}.\arabic{homework}}
\addtocounter{homework}{1} \addtocounter{tothomework}{1}
Study and program  the convolution algorithm described in \cite{Iver87,Iver15,ross95}.
Also write a program for the recursion algorithm and for the method based on (\ref{finitek}).
For a given (reasonably large) problem, compute the blocking probability for each class.
Make sure it is the same for all three alternatives. Then,  compare for a wide
range of parameter values the running times of the various algorithms and explain the differences.  $~~~\Box$

\subsubsection*{Homework \ref{multiservice}.\arabic{homework}}
\addtocounter{homework}{1} \addtocounter{tothomework}{1}
Provide an example of a continuous-time Markov chain that represents
a queueing model that is not reversible.
\subsubsection*{Guide}
Consider the MMPP(2)/M/1/1 queue, and show cases where the continuous-time Markov chain is not reversible.
$~~~\Box$

\newpage
\section{Discrete-Time Queue}
\label{discrete}
\setcounter{homework}{1} 

To complement the considerable attention we have given to continuous-time queues, we will now provide an example of a discrete-time
queueing system. Discrete-time models are very popular studies of
computer and telecommunications systems because in some cases, time
is divided into fixed length intervals (time slots) and packets of
information called cells are of fixed length, such that exactly one
cell can be transmitted during a time slot. Examples of such cases
include technologies, such as ATM and the IEEE 802.6 Metropolitan
Area Network (MAN) standard.

Let the number of cells that join the queue at different time slots
be an IID random variable. Let $a_i$ be the probability of $i$ cells
joining the queue at the beginning of any time slot. Assume that at
any time slot, if there are cells in the queue, one cell is served,
namely, removed from the queue. Further, assuming that arrivals occur at the
beginning of a time slot means that if a cell arrives during a
time slot it can be served in the same time slot.

In this case, the queue size process follows a discrete-time Markov chain with state-space $\Theta$ composed of all the nonnegative
integers, and a Transition Probability Matrix ${\bf P}=[{P_{ij}}]$
given by

\begin{equation}
\label{ssdt1} P_{i,i-1} = a_0 ~~{\rm for} ~~i \geq 1
\end{equation}

and

$$P_{0,0} = a_0 +  a_1 $$

$$P_{i,i} = a_1 ~~{\rm for}~~ i \geq 1$$

$$P_{i,i+1} = a_2 ~~{\rm for}~~ i \geq 0$$

and in general

\begin{equation}
\label{ssdt} P_{i,i+k} = a_{k+1} ~~{\rm for}~~ i \geq 0, k \geq 1.
\end{equation}

Defining the steady-state probability vector by
 $ \mathbf{\Pi}=[\pi_0, \pi_1, \pi_2, \ldots ] $, it can be obtained
 by solving the steady-state equations:

 $$\mathbf{\Pi}=\mathbf{\Pi} \textbf{P}.$$

 together with the normalizing equation
$$\sum_{i=0}^\infty \pi_i =1.$$

To solve for the $\pi_i$s, we will begin by writing down the steady-state equations as follows

$\pi_0 = \pi_0 P_{00} +  \pi_1 P_{10}$

$\pi_1 = \pi_0 P_{01} +  \pi_1 P_{11} + \pi_2 P_{21}$

$\pi_2 = \pi_0 P_{02} +  \pi_1 P_{12} + \pi_2 P_{22} + \pi_3 P_{32}$

$\pi_3 = \pi_0 P_{03} +  \pi_1 P_{13} + \pi_2 P_{23} + \pi_3 P_{33}
+ \pi_4 P_{43}$

and in general

$\pi_n = \sum_{i=0}^{n+1}  \pi_i P_{i,n}$ for $n \geq 0$.

Substituting (\ref{ssdt1}) and (\ref{ssdt}) in the latter, we obtain

\begin{equation}
\label{ssdt1st} \pi_0 = \pi_0[a_0 + a_1] + \pi_1 a_0
\end{equation}
\begin{equation}
\label{ssdt2nd} \pi_1 = \pi_0 a_2 + \pi_1 a_1 + \pi_2 a_0
\end{equation}
\begin{equation}
\label{ssdt3rd} \pi_2 = \pi_0 a_3 + \pi_1 a_2 + \pi_2 a_1 + \pi_3
a_0
\end{equation}
and in general
\begin{equation}
\label{ssdt2}  \pi_n =  \sum_{i=0}^{n+1} \pi_i a_{n+1-i}~~{\rm
for}~~ n \geq 1.
\end{equation}

Defining $\Pi(z)$ the Z-transform of the $\mathbf{\Pi}$ vector and
$A(z)$ as the Z-Transform of $[a_0, a_1, a_2, \ldots]$, multiplying
the $n$th equation of the set (\ref{ssdt1st}) -- (\ref{ssdt2}) by
$z^n$, and summing up, we obtain after some algebraic operations
\begin{equation}
\label{PiZ}  \Pi(z) = \pi_0 a_0 - \pi_0 z^{-1} a_0 + z^{-1} A(z)\Pi (z)
\end{equation}
which leads to
\begin{equation}
\label{PiZ1}   \Pi(z)= \frac{\pi_0 a_0 (1-z^{-1})}{ 1- z^{-1}A(z)}.
\end{equation}
Then, deriving the limit of $\Pi(z)$ as $z\rightarrow 1$ by applying L'Hopital rule,  denoting $A^\prime(1) = \lim_{z \rightarrow 1} A^\prime(z)$, and noticing that $\lim_{z\rightarrow 1} \Pi(z) = 1$ and $\lim_{z\rightarrow 1} A(z) = 1$, we obtain,
\begin{equation}
\label{Pi0spec}  \pi_0 = \frac{1- A^\prime(1) } {a_0}.
\end{equation}
This equation is somewhat puzzling. We already know that the
proportion of time the server is idle must be equal to one minus the
utilization. We also know that $A^\prime(1)$ is the mean arrival
rate of the number of arrivals per time slot and since the service
rate is equal to one, $A^\prime(1)$ is also the utilization; so what
is wrong with Eq.\@ (\ref{Pi0spec})? The answer is that nothing
wrong with it. What we call $\pi_0$ here is not the proportion of
time the server is idle. It is the probability that the queue is
empty at the slot boundary. There may have been one cell served in
the previous slot, and there may be an arrival (or more) in the next
slot which keeps the server busy.

The proportion of time the server is idle is in fact $\pi_0 a_0$, which is the probability of an empty queue at the slot boundary times
the probability of no arrivals in the next slot, and the consistency
of Eq.\@ (\ref{Pi0spec}) follows.

\subsubsection*{Homework \ref{discrete}.\arabic{homework}}
\addtocounter{homework}{1} \addtocounter{tothomework}{1}  Provide in detail
all the algebraic operations and the application of the L'Hopital rule
to derive equations (\ref{PiZ}), (\ref{PiZ1}) and (\ref{Pi0spec}).
\subsubsection*{Guide}
Multiplying
the $n$th equation of the set (\ref{ssdt1st}) -- (\ref{ssdt2}) by
$z^n$ and summing up, we obtain an equation for $\Pi(z)$ by focussing first on terms involving $\pi_0$ and then on the remaining terms.
For the remaining terms, it is convenient to focus first on terms involving $a_0$ and then on those involving $a_1$, etc. Notice in the following that all the remaining terms can be presented by a double summation.
\begin{eqnarray*}
\Pi(z) &
= &   \pi_0 a_0 z^0 + \pi_0\sum_{i=1}^\infty a_i z^{i-1} + \sum_{j=0}^\infty \left[ a_j \sum_{i=1}^\infty \pi_i z^{i-(1-j)} \right] \\ &
= &   \pi_0 a_0 + \pi_0 z^{-1}[A(z)-a_0] + z^{-1} A(z) [\Pi(z) - \pi_0 ] \\ &
= &  \pi_0 a_0 - \pi_0 z^{-1} a_0 + z^{-1} A(z)\Pi (z)
\end{eqnarray*}
and (\ref{PiZ1}) follows.

L'Hopital rule says that if functions $a(x)$ and $b(x)$ satisfy
$\lim_{x\rightarrow l^*} a(x)=0$ and $\lim_{x\rightarrow l^*}
b(x)=0$, then $$\lim_{x\rightarrow l^*}\frac{a(x)}{b(x)} =
\frac{\lim_{x\rightarrow l^*} a(x)}{\lim_{x\rightarrow l^*} b(x)}.$$
Therefore, from (\ref{PiZ1}) we obtain \begin{eqnarray*}
\lim_{z\rightarrow 1} \Pi(z) & = &   \lim_{z\rightarrow 1}
\frac{\pi_0 a_0 (1-z^{-1})}{ 1- z^{-1}A(z)}                  \\ & =
&    \lim_{z\rightarrow 1} \frac{\pi_0 a_0 z^{-2}   }{ z^{-2}A(z) -
z^{-1}A^\prime(z)  }. \end{eqnarray*} Substituting
$\lim_{z\rightarrow 1} \Pi(z) = 1$ and $\lim_{z\rightarrow 1} A(z) =
1$, we obtain, $$1=\frac{\pi_0 a_0}{1-A^\prime(z) }$$ and
(\ref{Pi0spec}) follows. $~~~\Box$ \subsubsection*{Homework
\ref{discrete}.\arabic{homework}} \addtocounter{homework}{1}
\addtocounter{tothomework}{1} Derive the mean and variance of the
queue size using the Z-transform method and verify your results by
simulations over a wide range of parameter values using confidence
intervals. $~~~\Box$

\newpage
\section{M/G/1 and Extensions}
\label{mg1}

\setcounter{homework}{1} 

The M/G/1 queue is a generalization of the M/M/1 queue where the
service time is no longer exponential. We now assume that the
service times are IID with mean $1/\mu$ and standard deviation
$\sigma_s$. The arrival process is assumed to be Poisson with rate
$\lambda$ and we will use the previously defined notation:
$\rho=\lambda/\mu$. As in the case of M/M/1, we assume that the
service times are independent and are independent of the arrival
process. In addition to M/M/1, another commonly used special case of
the M/G/1 queue is the M/D/1 queue where the service time is
deterministic.

The generalization from M/M/1 to M/G/1 brings with it a
significant increase in complexity. No longer can we use the Markov-chain structure that was so useful in the previous analyzes
where both service and inter-arrival times are memoryless. Without
the convenient Markov-chain structure, we will use different
methodologies as described in this section.

This Chapter provides queueing analyzes and results for the M/G/1 queue with FIFO and LIFO service disciplines as well as for a model of the M/G/1 queue with priorities.

\subsection{Pollaczek Khinchine Formula: Residual Service Approach \cite{BG92}}

The waiting time in the queue of an arriving customer to an M/G/1 queue is the remaining service time of the customer in service plus the sum of the service times of all the customers in the queue ahead of the arriving customer. Therefore, the mean waiting time in the queue
is given by
\begin{equation}
\label{WR} E[W_Q]=E[R]+\frac{E[N_Q]}{\mu}
\end{equation}
where $E[R]$ denotes the mean (unconditional) residual service time. For M/M/1, by the law of iterated expectations, $E[R]=\rho/\mu$, which is the probability of having one customer in service, which is equal to $\rho$, times the mean residual service time of that customer, which is equal to $1/\mu$ due to
the memoryless property of the exponential distribution, plus the
probability of having no customer in service (the system is
empty), which is $1-\rho$, times the mean residual service time if
there is no customer in service, which is equal to zero.

If we consider the residual service time conditioning on the event that the server is busy, namely that there is a customer in service that is being served, then the mean of the conditional residual service time, for the case of M/M/1, is equal to $1/\mu$.

\subsubsection*{Homework \ref{mg1}.\arabic{homework}}
\addtocounter{homework}{1} \addtocounter{tothomework}{1} Verify that
Eq.\@ (\ref{WR}) holds for M/M/1. $~~~\Box$

Observe that while  Equation (\ref{WR}) is based on considerations
at time of arrival, Little's formula
$$E[N_Q] = \lambda E[W_Q]$$
could be explained based on considerations related to a point in time when a customer
leaves the queue and enters the server. Recall the intuitive explanation of Little's formula in Section \ref{general}
which can be  applied to a system composed of the queue excluding the server.
Consider a customer that just left the queue leaving behind on average $E[N_Q]$
customers that have arrived during
the customer's time in the system which is on average $\lambda E[W_Q]$.

 By Little's formula and (\ref{WR}), we obtain,
\begin{equation}
\label{WR1} E[W_Q]=\frac{E[R]}{1-\rho}.
\end{equation}
It remains to obtain $E[R]$ to obtain results for the mean values
of waiting time and queue size.

Now that the service time is generally distributed, we
encounter certain interesting effects. Let us ask ourselves the
following question. If we randomly inspect an M/G/1 queue, will
the mean remaining (residual) service time of the customer
in service be longer or shorter than the mean service time? A
hasty response may be: shorter. Well, let us consider the
following example. There are two types of customers. Each of the
customers of the first type requires $10^6$ service units, while
each of the customers of the second type requires $10^{-6}$
service units. Assume that the proportion of the customers of the
first type is $10^{-7}$, so the proportion of the customers of
the second type is $1-10^{-7}$. Assume that the capacity of the
server to render service is one service unit per time unit and
that the mean arrival rate is one customer per time unit. As the
mean service time is of the order of $10^{-1}$, and the arrival
rate is one, although the server is idle 90\% of the time, when it
is busy it is much more likely to be busy serving a customer of
the first type despite the fact that these are very rare, so the
residual service time, in this case, is approximately $0.1\times
10^6/2=50,000$ which is much longer than the $10^{-1}$ mean
service time. Intuitively, we may conclude that the
residual service time is affected significantly by the variance of
the service time.

Notice that what we have computed above is the unconditional mean
residual service time which is our $E[R]$. Conditioning on the event that the server is busy, the mean residual service time will
be 10 times longer. We know that if the service time is
exponentially distributed, the conditional residual service time
of the customer in service has the same distribution as the
service time due to the memoryless property of the exponential
distribution. Intuitively, we may expect that if the variance of
the service time is greater than its exponential equivalence (an exponential random variable with the same mean), then the mean
residual service time (conditional) will be longer than the mean service time.
Otherwise, it will be shorter. For example, if the service time is
deterministic of length $d$, the conditional mean residual service
time is $d/2$, half the size of its exponential equivalence.

To understand why if the variance of the service time is sufficiently large, the residual service time can be longer than the mean service time, assume you have a queue with many customers who require very short (nearly zero) service time, but occasionally you have customers that require very long service times and these demanding customers (although they are a small minority of the customer population) actually occupy most of the time of the server. At a random point in time, you are likely to find one of these demanding customers in service (because they occupy most of the server time), and its residual service time is very long, but most customers require nearly zero service time, so the average service time can be lower than the residual service time.

To compute the (unconditional) mean residual service time $E[R]$,
consider the process $\{R(t), t\geq 0\}$ where $R(t)$ is the
residual service time of the customer in service at time $t$. And
consider a very long time interval $[0,T]$. Then,
\begin{equation} \label{ER} E[R]=\frac{1}{T}\int_0^T R(t) dt.
\end{equation} Following \cite{BG92}, let $S(T)$ be the number of
service completions by time $T$ and $S_i$ the $i$th service time.
Notice that the function $R(t)$ takes the value zero during times
when there is no customer in service and jumps to the value of $S_i$
at the point of time when the $i$th service time commences. During a
service time, it linearly decreases with a rate of one and reaches zero
at the end of a service time. Therefore, the area under the curve
$R(t)$ is equal to the sum of the areas of $S(T)$ isosceles right
triangles where the side of the $i$th triangle is $S_i$. Therefore,
for large $T$, we can ignore the last possibly incomplete triangle,
so we obtain \begin{equation} \label{ER1}
E[R]=\frac{1}{T}\sum_{i=1}^{S(T)} \frac{1}{2}
S_i^2=\frac{1}{2}\frac{S(T)}{T}\frac{1}{S(T)}\sum_{i=1}^{S(T)}
 S_i^2.
\end{equation}
Letting $T$ approach infinity, the latter gives
\begin{equation}
\label{ER2} E[R]=\frac{1}{2}\lambda \overline{S^2}
\end{equation}
where $\overline{S^2}$ is the second moment of the service time.

By (\ref{WR1}) and (\ref{ER2}), we obtain
\begin{equation}
\label{PK1} E[W_Q]=\frac{\lambda \overline{S^2}}{2(1-\rho)}.
\end{equation}
Thus, considering (\ref{EDEW}), we obtain
\begin{equation}
\label{PK2}  E[D]=\frac{\lambda \overline{S^2}}{2(1-\rho)}+ \frac{1}{\mu}.
\end{equation}

Recalling that $\sigma_s^2=\overline{S^2}-(1/\mu)^2$, Eq.\@ (\ref{PK2}) leads to
the well-known Pollaczek Khinchine formula for the mean delay in an M/G/1 system:

\begin{equation}
\label{PK_Delay}  E[D]= \frac{\lambda (\sigma_s^2 + \mu^{-2}) }{2(1-\rho)}+  \frac{1}{\mu} = \frac{\rho + \lambda\mu\sigma_s^2}{2(\mu-\lambda)} + \frac{1}{\mu}.
\end{equation}

Using Little's formula, we obtain the Pollaczek Khinchine formula for the mean number of
customers in an M/G/1 system:
\begin{equation}
\label{PK}  E[Q]=\rho + \frac{\rho^2 +
\lambda^2\sigma_s^2}{2(1-\rho)}.
\end{equation}
Observe that according to the Pollaczek Khinchine formula, if we have two M/G/1 queueing systems, where they both have the same arrival and service rates, but for one the variance of the service time is higher than that of the other, the one with the higher variance will experience higher mean queue size and delay.

Note that although there are two versions of the Pollaczek Khinchine formula: one for the mean delay and one for the mean queue size, we use the term Pollaczek Khinchine formula in singular to mean either one of them. This is common in the literature. Note also that it is straightforward to derive formulae for $E[N_Q]$,  the mean number only in the queue (excluding those in service), and $E[W_Q]$, the mean queueing waiting time excluding the time in service, similar to the way it was done in Section \ref{perf_measure_mean}.


\subsection{Pollaczek-Khinchine Formula: by Kendall's Recursion \cite{kendall51}}

Let us now derive (\ref{PK}) in a different way.
Letting $q_i$ be the number of customers in the system immediately
following the departure of the $i$th customer, the following recursive relation,
is obtained.
\begin{equation}
\label{recurMG1}  q_{i+1} = q_i + a_{i+1} - I(q_i)
\end{equation}
where $a_{i}$ is the number of arrivals during the service time of the $i$th customer,
and $I(x)$ is a function defined for $x\geq 0$, taking the value 1 if $x>0$, and the value 0 if $x=0$.
This recursion was first introduced by Kendall \cite{kendall51}, so we will call it Kendall's Recursion.
Some call it a ``Lindley's type Recursion''  in reference to an equivalent recursion for the G/G/1 waiting time in \cite{Lindley52}.
Along with Little's and Erlang B formulae, and the Pollaczek-Khinchine equation, the Kendall's and Lindley's
recursions are key foundations of queueing theory.

To understand the recursion (\ref{recurMG1}), notice that there are
two possibilities here: either $q_i=0$ or $q_i>0$.

If $q_i=0$,  then the $i+1$th customer arrives in an empty system. In this
case $I(q_i)=0$ and the number of customers in the system when the
$i+1$th customer leaves must be equal to the number of customers that arrives
during the service time of the $i+1$th customer.

If $q_i>0$, then the $i+1$th customer arrives in a nonempty system.
It starts its service when the $i$th customer leaves. When it starts
its service there are $q_i$ customers in the system. Then, during its
service time, additional $a_{i+1}$ customers arrive. And when it leaves
the system there must be $q_i + a_{i+1} -  1$ (where the `-1'
represents the departure of the $i+1$th customer).

Squaring both sides of (\ref{recurMG1}) and taking expectations, we obtain
\begin{equation}
\label{eq2}
E[q_{i+1}^2] = E [q_i^2] + E[I(q_i)^2]+ E[a_{i+1}^2] - 2E[q_iI(q_i)] + 2E[q_ia_{i+1}] - 2 E[I(q_i)a_{i+1}].
\end{equation}
Notice that in steady-state $E[q_{i+1}^2]=E [q_i^2]$, $I(q_i)^2=I(q_i)$, $E[I(q_i)^2]=E[I(q_i)]=\rho$ , and that
for any $x\geq 0$, $xI(x)=x$, so $q_iI(q_i)=q_i$.
Also notice that because of the independence between $a_{i+1}$ and $q_i$, and because (by (\ref{meancondind})) the mean number of arrivals during service time in M/G/1 is equal to $\rho$,
we obtain in steady-state that $E[I(q_i)a_{i+1}]=\rho^2$ and $E[q_ia_{i+1}]=E[q_i]\rho$.
Therefore, considering (\ref{eq2}), and setting the steady-state notation $E[a]=E[a_i]$
and $E[Q]=E[q_i]$, we obtain after some algebra
\begin{equation}
\label{eqpk}
E[Q] = \frac{\rho+E[a^2]-2\rho^2}{2(1-\rho)}.
\end{equation}
To obtain $E[a^2]$, we notice that by EVVE,
\begin{equation}
\label{evvevara}
Var[a]=E[Var[a\mid S]] + Var[E [ a\mid S]]=\lambda E[S] +\lambda^2 \sigma_s^2 =\rho + \lambda^2 \sigma_s^2
\end{equation}
recalling that $S$ is the service time and that $\sigma_s^2$ is its
variance. Also recall that $Var[a]=E[a^2]-(E[a])^2$ and since
$E[a]=\rho$, we have by Eq.\@ (\ref{evvevara}) that
$$E[a^2]=Var[a]+\rho^2=\rho + \lambda^2 \sigma_s^2+\rho^2.$$
Therefore,
\begin{equation}
\label{eqpka}
E[Q] = \frac{2\rho + \lambda^2 \sigma_s^2 -\rho^2}{2(1-\rho)}
\end{equation}
or
\begin{equation}
\label{PKnew}  E[Q]=\rho + \frac{\rho^2 +
\lambda^2\sigma_s^2}{2(1-\rho)}
\end{equation}
which is identical to (\ref{PK}) - the Pollaczek-Khinchine Formula.

\subsubsection*{Homework \ref{mg1}.\arabic{homework}}
\addtocounter{homework}{1} \addtocounter{tothomework}{1} Re-derive the Pollaczek-Khinchine Formula
in the two ways presented above with attention to all the details
(some of which are skipped in the above derivations). $~~~\Box$


\subsection{Special Cases: M/M/1 and M/D/1}

Now let us consider the special case of exponential service time.
Namely, the M/M/1 case. To obtain $E[Q]$ for M/M/1, we substitute
$\sigma_s^2=1/\mu^2$ in (\ref{PK}), and after some algebra, we
obtain
\begin{equation}
\label{PKMM1}  E[Q]=\frac{\rho}{1-\rho}
\end{equation}
which is consistent with (\ref{meanmm1}).

Another interesting case is the M/D/1 queue in which case we have: $\sigma_s^2=0$.
Substituting the latter in (\ref{PK}), we obtain after some
algebra
\begin{equation}
\label{PKMD1}  E[Q]=\frac{\rho}{1-\rho} \times \frac{2-\rho}{2}.
\end{equation}
Because the second factor of (\ref{PKMD1}), namely $(2-\rho)/2$,
is less than one for the range $0<\rho<1$, we clearly see that the mean number of customers in an M/M/1 queue is higher than that of an
M/D/1 queue with the same arrival and service rates.

Furthermore, the inverse of this factor, given by

$$\frac{2}{2-\rho}$$

is the ratio of the M/M/1 mean queue size to that of M/D/1 queue size assuming both have the same traffic intensity $\rho$. Noting that it is monotonically decreasing in $\rho$ within  the relevant range of $0<\rho <1$, and that it takes the values 1 and 2, for $\rho=0$ and $\rho=1$, respectively. Accordingly, this ratio is bounded between 1 and 2. This implies that due to the  service variability in M/M/1, its mean queue size is always higher (but it is never higher by more than 100\%) than that of the equivalent M/D/1.

Considering (\ref{PKMD1}), and invoking Little's formula, the mean delay for an  M/D/1 queue is obtained by

\begin{equation}
\label{PKMD1delay}  E[D]=\frac{1}{2\mu} \times \frac{2-\rho}{1-\rho} = \frac{1}{2\mu} \times \frac{2\mu -\lambda}{\mu-\lambda}.
\end{equation}

\subsubsection*{Homework \ref{mg1}.\arabic{homework}}
\addtocounter{homework}{1} \addtocounter{tothomework}{1}
Derive (\ref{PKMD1delay}) from (\ref{PKMD1}) invoking Little's formula.

\subsubsection*{Homework \ref{mg1}.\arabic{homework}}
\addtocounter{homework}{1}
Show that $E[W_Q]$ (the time spent in the queue but not in service) for M/D/1 is half that of its M/M/1 counterpart assuming that the mean service times in both systems are the same.

\subsubsection*{Guide}
This can be done in several ways. Here is one way to show it.

First recall that $$E[W_Q] = E[D] - 1/\mu = E[Q]/\lambda -1/\mu.$$

For M/M/1, by (\ref{PKMM1}),

$$E[W_Q] = \frac{\rho}{(1-\rho)\lambda} - \frac{1}{\mu} = \frac{1}{\mu-\lambda} - \frac{1}{\mu} = \frac{\rho}{\mu-\lambda}.$$

For M/D/1, by (\ref{PKMD1delay}),

$$E[W_Q] = \frac{1}{2\mu} \times \frac{2-\rho}{1-\rho} - \frac{1}{\mu} = \frac{\rho}{2(\mu-\lambda)}.$$

Clearly,  $E[W_Q]$  for M/D/1 is half  that of its M/M/1 counterpart.

Complete all steps. $~~~\Box$

\subsubsection*{Homework \ref{mg1}.\arabic{homework}}
\addtocounter{homework}{1}

Consider a multiplexer with one output link that serves data messages based on the first come first served principle. The multiplexer has a very large buffer. Data messages arrive at the multiplexer following a Poisson process at a rate of 100 messages per second. The message sizes have a mean of 25 Mbytes (Mbytes = Mega-bytes) (one byte = 8 bits) and the standard deviation of the message size is also 25 Mbytes. It is required that the mean message delay in the steady state will not exceed 10 ms (ms = millisecond = 0.001 second). What should be the minimal bit rate (in Gb/s = Gigabits per second) of the multiplexor to satisfy this delay requirement?  Show all the steps and explain them clearly. As a first step, choose an appropriate queueing model for solving this problem and justify your choice.  Then, describe and justify the remaining steps until you reach the solution.

\subsubsection*{Solution}
The appropriate model is M/M/1.

Justification:

First, notice the following four points:
\begin{enumerate}
\item The very large buffer justifies the infinite buffer model.
\item The Poisson arrivals are given
\item 	The single output link justifies the single server assumption.
\item 	The service policy is given as first come first served.
\end{enumerate}

Then, the remaining issue is why we do not choose the more complex M/G/1 model, as it was not given that the message size or service time is exponentially distributed. Indeed, M/G/1 can be used to solve the problem correctly. However, using M/M/1 will give exactly the same answer and the solution is far simpler. Notice that we dimension the output rate only based on the mean delay, and since the Pollaczek–Khinchine formula is only based on the first two moments of the service time distribution, and since here the mean service time is equal to the standard deviation of the service time (because it is given that the mean message length is equal to the standard deviation of the service time) which is a property of the exponential distribution, M/M/1 will give the correct result here.

Mean message size = 25 Mbytes = $25 \times 8$ Mbits = 200 Mbits.

Let $X$ [Mb/s] be the Multiplexer bit rate.

$$\mu =\frac{\rm X }{\rm 200~Mbits} ~[{\rm messages/sec.}]$$  \\ $$\lambda=100 ~ {\rm messages/sec}$$~~~\\   $$\rho=\frac{\lambda}{\mu}=\frac{2\times 10^4}{X}.$$

\subsubsection*{Homework \ref{mg1}.\arabic{homework}}
\addtocounter{homework}{1}

Use Eq. (\ref{ER2}) to show that the value of the mean (unconditional) residual service time $E[R]$ for M/D/1 is half of its value for M/M/1.

\subsubsection*{Solution}
For M/M/1, we have:

$$E[R]=\frac{1}{2}\lambda\overline{S^2}=\frac{1}{2}\lambda\left(Var\left[S\right]+E\left[S\right]^2\right)=\frac{1}{2}\lambda\left(\left(\frac{1}{\mu}\right)^2+\left(\frac{1}{\mu}\right)^2\right)=\frac{1}{2}\lambda\left(2\left(\frac{1}{\mu}\right)^2\right)=\rho\frac{1}{\mu}.$$

For M/D/1, we have:

$$E[R]=\frac{1}{2}\lambda\overline{S^2}=\frac{1}{2}\lambda\left(Var\left[S\right]+E\left[S\right]^2\right)=\frac{1}{2}\lambda\frac{1}{\mu^2}=\frac{1}{2}\rho\frac{1}{\mu}.$$

Note that we commented above that the conditional mean residual service time of deterministic service is half that of its exponential counterpart. Then, we can expect that this also applies to the unconditional means, because the  unconditional mean residual service time is equal to the conditional mean residual service time times $\rho$.

Note also that we have observed that not only the mean residual service time (conditional and unconditional) of deterministic service is half that of its exponential counterpart, but also the mean waiting time of M/D/1 is half its M/M/1 counterpart.

\subsection{Busy Period}

We have defined and discussed the concept of the busy period in Section
\ref{busyp} in the context of the M/M/1 queue. The same analysis
applies to the case of the M/G/1 system, and we obtain:
\begin{equation}
\label{busymg1} E[T_B]=\frac{1}{\mu-\lambda}.
\end{equation}
What we learn from this is that the mean busy period is insensitive
to the shape of the service time distribution. In other words, the
mean busy periods of M/M/1 and M/G/1 systems are the same if the
mean arrival and service rates are the same.

\subsubsection*{Homework \ref{mg1}.\arabic{homework}}
\addtocounter{homework}{1} \addtocounter{tothomework}{1}
\begin{enumerate} \item
Prove that $$\frac{E[T_B]}{E[T_B] + E[T_I]}$$ is the
proportion of time
that the server is busy.
\item Show that Equation (\ref{busymg1}) also applies to an M/G/1 queue. $~~~\Box$
\end{enumerate}
\subsubsection*{Homework \ref{mg1}.\arabic{homework}}
\addtocounter{homework}{1} \addtocounter{tothomework}{1}
Consider an M/G/1 queueing system with the following twist. When a new customer arrives at an empty system, the server is not available immediately. The customer then rings a bell and the server arrives an exponentially distributed  amount of time with parameter $\zeta$ later. As in M/G/1,
customers arrive in accordance with a Poisson process with rate $\lambda$ and the mean service time is $1/\mu$.
Service times are mutually independent and independent of the inter-arrival times. Find the mean busy period defined as a continuous period during which the server is busy.
\subsubsection*{Guide}
Explain and solve the following two equations:
$$\frac{E[T_B]}{E[T_B] + E[T_I]}=\rho=\frac{\lambda}{\mu}$$
and
$$E[T_I] = \frac{1}{\lambda} + \frac{1}{\zeta}. ~~~\Box $$

\subsubsection*{Homework \ref{mg1}.\arabic{homework}}
\addtocounter{homework}{1} \addtocounter{tothomework}{1}
Consider a multiplexer with one output link that serves data at a rate of 2.4 Gb/s. The multiplexer has a very large buffer. Data messages arrive at the multiplexer following a Poisson process at a rate of 12 messages per second. The message sizes have a mean of 20 Mbytes and a standard deviation of 15 Mbytes. The arriving messages are processed and transmitted by the output link based on the first come first served principle. Consider the system to be in a steady state. 	Choose an appropriate queueing model for this multiplexer system, and clearly justify your choice.

Based on the model you choose, answer the following.
\begin{itemize}
\item Find the mean number of messages in the entire multiplexer system (in service and in the queue).
\item Find the mean delay of a message (including service time).
\end{itemize}

\subsubsection*{Solution}
The appropriate model is M/G/1.

{\bf Justification:}
\begin{enumerate}
\item	The very large buffer justifies the infinite buffer model.
\item The Poisson arrivals are given.
\item	Service time is neither exponential nor deterministic, so the general service time assumption is justified.
\item	The service policy is given as first come first served.
\item	The single output link justifies the single server assumption.
\end{enumerate}

Next,

2.4 Gb/s = 2,400 Mb/s.\\
Mean message size = 20 Mbytes = $20 \times 8 Mbits$ = 160 Mbits.\\
Standard deviation of message size = 15 Mbytes = 120 Mbits\\
Standard deviation of service time = 120 / 2,400 = 0.05 [sec.]

$$\mu =\frac{\rm 2,400 }{\rm 160~Mbits}=15~{\rm messages/sec.}~~~~~~\lambda=12 ~ {\rm messages/sec}~~~~~~~ \rho=\frac{\lambda}{\mu}=\frac{12}{15}=0.8$$
Then, the mean number of messages in the entire multiplexer system is obtained by:
$$E[Q]=\rho + \frac{\rho^2 +
\lambda^2\sigma_s^2}{2(1-\rho)}=0.8+\frac{0.8^2 +
12^2\times 0.05^2}{2(1-0.8)}=3.3$$
and the mean delay is obtained by Little's formula as:
$$E[D] = \frac{E[Q]}{\lambda} = \frac{3.3}{12} = 0.275~\rm{[sec.]}. ~~~\Box $$

\subsubsection*{Homework \ref{mg1}.\arabic{homework}}
\addtocounter{homework}{1} \addtocounter{tothomework}{1}
Consider a single server queue with an infinite-size buffer. There are two traffic streams of messages that arrive at the queue for processing. Each of the two traffic streams follows a Poisson process: one with arrival rate $ \lambda_1= 0.5 $  message per millisecond (ms) and the other with arrival rate  $ \lambda_2 = 1 $   messages per ms. All messages that arrive are served based on the first come first served principle.
The messages of the first stream require exponentially distributed service times with a mean of 0.6 ms and the messages of the second stream require exponentially distributed service times with a mean of 0.3 ms.
Find the mean number of messages (regardless of the stream they belong to) in the system (this includes messages waiting in the queue and in service).

\subsubsection*{Solution}
The total arrival rate is obtained by $\lambda = \lambda_1+\lambda_2 = 0.5+1 = 1.5 ~~{\rm messages/ms.} $

The mean service time is obtained by the Law of iterated Expectation:

\[
E[S] = \frac{\lambda_1}{\lambda_1+\lambda_2}0.6 + \frac{\lambda_2}{\lambda_1+\lambda_2}0.3 =\frac{0.5}{1.5}\times 0.6 +\frac{1}{1.5}\times 0.3=0.4~ {\rm ms.}
\]

Another way is:

\[
E[S]=E[0.3 + 0.3 X]
\]
where $X$ is a Bernoulli random variable with parameter $p=0.5/1.5=1/3$.

Thus,

$$E[S]=E[0.3 + 0.3 X] = 0.3 +0.3 E[X] = 0.3 + 0.3\times \frac{1}{3}=0.4  ~ {\rm ms.}$$

The variance of the message time can be obtained by the Law of Total Variance (or EVVE), as follows

\[
\sigma_S^2 = \frac{1}{3}\times (0.6)^2 + \frac{2}{3}\times  (0.3)^2  + Var[0.3 +0.3X]
\]
where $X$ is a Bernoulli random variable with parameter $p=1/3$.

Therefore,

\[
Var[0.3  +0.3X] = (0.3)^2Var[X] = 0.09 \times p(1-p) = 0.09 \times \frac{1}{3}\times\frac{2}{3}=0.09\times \frac{2}{9}=0.02
\]
so
$$\sigma_S^2 = 0.2.$$

Now we use the Pollaczek-Khinchine Formula for the mean number of messages in an M/G/1 system:

$$\rho=\lambda E[S] = 1.5\times 0.4 = 0.6.$$
\[
E[Q]=\rho + \frac{\rho^2 +
\lambda^2\sigma_s^2}{2(1-\rho)}=0.6 +\frac{0.6^2+(1.5)^2\times 0.2}{2(1-0.6)}=1.6125~~ {\rm messages}. ~~~\Box
\]



\subsection{M/G/1-LIFO}
\label{mg1lifo}
The M/G/1-LIFO queue possesses similar properties to the M/G/1-PS queue that we discussed in Section \ref{mm1psinsesitivity}.
They are both insensitive to the shape of the service time
distribution.

We have already mentioned in Section \ref{busyp} that the queue size process of
M/M/1 is the same as that of its M/M/1-LIFO equivalence. Therefore,
they also have the same mean queue size and delay. Due to the
insensitivity of M/G/1-LIFO, the M/M/1 results for the mean queue
size, mean delay, and queue size distribution are applicable also to
M/G/1-LIFO.

Specifically, if we are given an M/G/1-LIFO queue with an arrival rate
$\lambda$ and mean service time $1/\mu$, denote $\rho=\lambda/\mu$, then the queue size
distribution is given by:

\begin{equation}  \pi_i = \rho^i (1-\rho)  ~{\rm for}
~i=0,~1,~2,~\ldots~. \end{equation}

The mean queue size is given by
\begin{equation}
E[Q]=\frac{\rho}{1-\rho} \end{equation}
and the mean delay is given by
\begin{equation}  E[D]=
\frac{1}{\mu-\lambda}. \end{equation}

To show why M/G/1-LIFO queue is insensitive and satisfies the above equations, notice that
an arriving customer that upon its arrival finds $i$ customers in the
system will be served only during the time when the system is in state $i+1$. Furthermore,
all the customers served when the system is in state $i+1$ will be customers that have arrived when the system is in state $i$. Therefore,
the time the system spends in the state $i+1$ comprises exactly the
 service times of the customers that arrive when the system is in  state $i$.
Now consider a long interval of time $T$. As we denote by $\pi_i$ the proportion of time that the system is in state $i$, then during a long time interval $T$, the mean number of arrivals in state $i$ is $\lambda \pi_i T$ and their total service time is equal to $\lambda \pi_i T(1/\mu)=\rho \pi_i T$. Accordingly, the proportion of time the system is in state $i+1$ is given by
$$\pi_{i+1}=\frac{\rho \pi_i T}{T}=\rho \pi_i.$$ Since the latter holds for $i=0, 1, 2, \ldots$, we observe that the queue-size distribution of M/G/1-LIFO obeys the steady-state equations of M/M/1 regardless of the shape of the service time distribution.

\subsection{M/G/1 with $m$ Priority Classes}

Let us consider an M/G/1 queueing system with $m$ priority classes.
Let $\lambda_j$ and $\mu_j$ be the arrival and service rate of
customers belonging to the $j$th priority class for $j=1, 2, 3,
\ldots, m$. The mean service time of customers belonging to the
$j$th priority class is, therefore, equal to $1/\mu_j$. The second
moment of the service time of customers belonging to the $j$th
priority class  is denoted $\overline{S^2(j)}$. We assume that priority class $j$ has a higher priority than priority class $j+1$, so
Class 1 represents the highest priority class and Class $m$ is the
lowest. For each class $j$, the arrival process is assumed to be
Poisson with parameter $\lambda_j$, and the service times are assumed
mutually independent and independent of any other service times  of
customers belonging to the other classes and are also independent
of any inter-arrival times. Let $\rho_j=\lambda_j/\mu_j$. We assume
that $\sum_{j=1}^m \rho_j<1$. We will consider two priority
policies: {\em nonpreemptive} and {\em preemptive resume}.

\subsection{Nonpreemptive}
Under this regime, a customer in service will complete its service even if a customer of a higher priority class arrives while it is being served.
Let  $E[N_Q(j)]$ and $E[W_Q(j)]$ represent the mean number of class $j$ customers in the queue
excluding the customer in service and the mean waiting time of a class $j$
customer in the queue (excluding its service time), respectively. Further, let $R$ be the residual service time (of all customers of all priority classes).
In similar way we derived (\ref{ER2}), we obtain:
\begin{equation}
\label{ER2np} E[R]=\frac{1}{2} \sum_{j=1}^m \lambda_j \overline{S^2(j)}.
\end{equation}
\subsubsection*{Homework \ref{mg1}.\arabic{homework}}
\addtocounter{homework}{1} \addtocounter{tothomework}{1} Derive
Eq.\@ (\ref{ER2np}). $~~~\Box$

As in Eq.\@ (\ref{WR}), we have for the highest priority,
\begin{equation}
\label{WRpr} E[W_Q(1)]=E[R]+\frac{E[N_Q(1)]}{\mu_1},
\end{equation}
and similarly to (\ref{WR}), we obtain
\begin{equation}
\label{WR1pr} E[W_Q(1)]=\frac{E[R]}{1-\rho_1}.
\end{equation}
Regarding the second priority, $E[W_Q(2)]$ is the sum of the mean residual service time $E[R]$,
the mean time it takes to serve the Class 1 customers in the queue  $E[N_Q(1)]/\mu_1$, the mean time it takes to serve the Class 2 customers in the queue  $E[N_Q(2)]/\mu_2$, and the mean time it takes to serve all the Class 1 customers that arrive during the waiting time in the queue for the Class 2 customer $E[W_Q(2)]\lambda_1/\mu_1=E[W_Q(2)]\rho_1$. Putting it together
\begin{equation}
\label{WR1pr2} E[W_Q(2)]=E[R] + \frac{E[N_Q(1)]}{\mu_1} + \frac{E[N_Q(2)]}{\mu_2}+ E[W_Q(2)]\rho_1.
\end{equation}
By the latter and Little's formula for Class 2 customers, namely,
$$E[N_Q(2)]=\lambda_2E[W_Q(2)],$$
we obtain
\begin{equation}
\label{WR12pr2} E[W_Q(2)]=\frac{E[R]+\rho_1 E[W_Q(1)]}{1-\rho_1-\rho_2}.
\end{equation}
By Eqs. (\ref{WR12pr2}) and (\ref{WR1pr}), we obtain
\begin{equation}
\label{WR12pr22} E[W_Q(2)]=\frac{E[R]}{(1-\rho_1)(1-\rho_1-\rho_2)}.
\end{equation}
\subsubsection*{Homework \ref{mg1}.\arabic{homework}}
\addtocounter{homework}{1} \addtocounter{tothomework}{1}
Show that for $m=3$,
\begin{equation}
\label{WR12pr23} E[W_Q(3)]=\frac{E[R]}{(1-\rho_1-\rho_2)(1-\rho_1-\rho_2-\rho_3)}.
\end{equation}
and that in general
\begin{equation}
\label{WR12pr2m} E[W_Q(j)]=\frac{E[R]}{(1-\sum_{i=1}^{j-1} \rho_i)(1-\sum_{i=1}^j \rho_i)}.
\end{equation}
 $~~~~~~~~~~~~~~~~~~~~~~~~~~~~~~~~~~~~~~~~~~~~~~~~~~~~~~~~~~~~~~~~~~~~~~~~~~~~~~~~~~~~~~~~~~~\Box$

 The mean delay for a $j$th priority class customer, denoted $E[D(j)]$, is  given by
\begin{equation}
\label{Dprj} E[D(j)]=E[W_Q(j)] + \frac{1}{\mu_j} {\rm ~for~} j=1,2,3, \ldots, m.
\end{equation}
\subsubsection*{Homework \ref{mg1}.\arabic{homework}}
\addtocounter{homework}{1} \addtocounter{tothomework}{1}
Consider the case of $m=2$, $\lambda_1=\lambda_2=0.5$  with $\mu_1=2$ and $\mu_2=1$. Compute the average delay for each class and the overall average delay. Then, consider the case of $m=2$, $\lambda_1=\lambda_2=0.5$
with $\mu_1=1$ and $\mu_2=2$ and compute the average delay for each class
and the overall average delay. Explain the difference between the two cases and draw conclusions.
Can you generalize your conclusions?
 $~~~\Box$
\subsection{Preemptive Resume}
In this case, an arriving customer of priority $j$ never waits for a customer
of a lower priority class (of Class $i$ for $i>j$) to complete its service.
Therefore, when we derive the delay of a customer of priority $j$,
we can ignore all customers of class $i$ for all $i>j$. Accordingly, the mean delay
of a priority $j$ customer satisfies the following equation
\begin{equation}\label{delpreres}
    E[D(j)]=\frac{1}{\mu_j} + \frac{R(j)}{1-\sum_{i=1}^j \rho_i } +  E[D(j)] \sum_{i=1}^{j-1}\rho_i
\end{equation} where $R(j)$ is the mean residual time of all
customers of classes $i=1,2, \ldots, j$ given by $$R(j) =
\frac{1}{2} \sum_{i=1}^j \lambda_i \overline{S^2(i)}.$$ The first
term of Eq.\@ (\ref{delpreres}) is simply the mean service time of a
$j$th priority customer. The second term is the mean time it takes
to clear all the customers of priority $j$ or higher that are
already in the system when a customer of Class $j$ arrives. It is
merely Eq.\@ (\ref{WR1}) that gives the mean time of waiting in the
queue in an M/G/1 queueing system where we replace $\rho$ of
(\ref{WR1}) by $\sum_{i=1}^j \rho_i$, which is the total traffic load
offered by customers of priority $j$ or higher. From the point of view of the $j$th priority customer, the order of the customers ahead
of it will not affect its mean delay, so we can ``mix'' all these
customers up and consider the system as M/G/1. The last term of
Eq.\@ (\ref{delpreres}) is the mean total work introduced to the
system by customers of priorities higher than $j$ that arrive during
the delay time of our $j$ priority customer. Notice that we use the
$\rho_i$s there because $\rho_i=\lambda_i (1/\mu_i)$ is
the product of the mean rate of customer arrivals and the mean work
they bring to the system for each priority class $i$.

Eq.\@ (\ref{delpreres}) leads to

\begin{equation}\label{delpreres1}
    E[D(1)]=\frac{(1/\mu_1)(1-\rho_1)+R(1)}{1-\rho_1},
\end{equation}
and
\begin{equation}\label{delpreresj}
    E[D(j)]=\frac{(1/\mu_j)\left(1-\sum_{i=1}^j \rho_i \right)+R(j)}{ \left(1-\sum_{i=1}^{j-1} \rho_i \right)\left(1-\sum_{i=1}^{j} \rho_i\right) }.
\end{equation}

\subsubsection*{Homework \ref{mg1}.\arabic{homework}}
\addtocounter{homework}{1} \addtocounter{tothomework}{1}
Derive Eqs. (\ref{delpreres1}) and (\ref{delpreresj}).
 $~~~\Box$

\subsubsection*{Homework \ref{mg1}.\arabic{homework}}
\addtocounter{homework}{1} \addtocounter{tothomework}{1}
Consider a single server queue with two classes of customers: Class 1 and Class 2, where Class 1 customers have preemptive resume priority over Class 2 customers. Class $i$ customer arrivals follow a Poisson process with parameter $\lambda_i$, and their service times are exponentially distributed with mean $1/\mu_i$, $i=1,2$. \begin{enumerate}
\item Derive formulae for the mean delay (including service time) of each of the classes. \item Assume $\mu=\mu_1=\mu_2$, let $\rho_i=\lambda_i/\mu$, $i=1,2$,
    and assume $\rho_1+\rho_2<1$. Maintain $\rho_1$ and $\rho_2$
 fixed  and let $\mu$ approach infinity, show that under these conditions, the mean delay of either traffic class approaches zero.
 \item Now assume the conditions $\rho_1<1$, but $\rho_1+\rho_2>1$, again let $\mu=\mu_1=\mu_2$ approach infinity and show that under these conditions, the mean delay of traffic Class 1 approaches zero.
  \end{enumerate}

     \subsubsection*{Guide}
     For exponentially distributed service times with mean $1/\mu$, we have
     $$ R(1) = \frac{1}{2} \left( \lambda_1 \frac{2}{\mu^2} \right) =  \frac{\rho_1}{\mu}.$$
       $$ R(2) = \frac{1}{2} \left( \lambda_1 \frac{2}{\mu^2}  + \lambda_2 \frac{2}{\mu^2}\right) =\frac{\rho_1+\rho_2}{\mu}.$$
       $$E[D(1)]=\frac{(1/\mu)(1-\rho_1)+\rho_1/\mu}{1-\rho_1}=\frac{1}{\mu(1-\rho_1)}.   $$
       This is not a surprise. It is the mean delay obtained by M/M/1 if all the traffic is of class 1 customers. We can observe clearly that if $\rho_1$ stays fixed and $\mu$ approaches infinity, the mean delay approaches zero. This applies to both 2 and 3.
       $$E[D(2)]=\frac{(1/\mu)(1-\rho_1-\rho_2)+R(2)}{ (1-\rho_1)(1-\rho_1-\rho_2) }.$$
       Substituting $R(2)$, we obtain,
       $$E[D(2)]=\frac{(1/\mu)(1-\rho_1-\rho_2)+(1/\mu)(\rho_1+\rho_2)}{ (1-\rho_1)(1-\rho_1-\rho_2) }= \frac{1}{\mu (1-\rho_1)(1-\rho_1-\rho_2) }.$$
       Now we can observe that if $\rho_1+\rho_2<1$, as $\mu$ approaches infinity, the mean delay also of priority 2, approaches zero. $~~~\Box$

 The last homework problem solution implies the following.
 For the M/M/1 with priorities model, if the queues of all priorities are stable and if the service rate is arbitrarily high, then the mean delay is arbitrarily low, regardless of the utilization.  Then, in such a case, there is not much benefit in implementing priorities. However, if, for example, $\rho_1+\rho_2>1$ but $\rho_1<1$, then priority 1 customers clearly benefit from having priority even if the service rate (and also arrival rate) is arbitrarily large.
 Notice that we have observed similar results for M/M/1 without priorities. Also, notice that we consider here a scenario where the service rate is arbitrarily high, and the utilization is fixed, which means that the arrival rate is also arbitrarily high.

\newpage
\section{Kingman's Formula}
 \label{kingmanF}
\setcounter{homework}{1} 
In this chapter, we consider an infinite-buffer SSQ where the inter-arrivals and service times are IID, and the service times are independent of the arrival process. This is a further generalization of the M/G/1 queue, discussed in the previous chapter, where the arrival process was Poisson, so the inter-arrival times were exponentially distributed IID random variables.

In the next chapter, we will further generalize the SSQ by not requiring independence between inter-arrival times or between service times and arrivals.

For our infinite-buffer SSQ with IID inter-arrivals and service times and no dependencies between service times and the arrival process, We further consider that the system is under heavy traffic, namely, $\rho \rightarrow 1$. For such an SSQ, we can obtain the following approximation, known as Kingman's formula \cite{kingman61}.

  \begin{equation}
 \label{kingman}
 E[W_Q] \approx \left( \frac{\rho}{1-\rho}   \right)  \left(  \frac   {c_a^2 +c_s^2}{2}         \right) E[S].
 \end{equation}

  where $c_a$ is the coefficient of variation of the inter-arrival times  which is the ratio of the standard deviation of the inter-arrival times to the mean of the inter-arrival time, and $c_s$ is the coefficient of variation for service times which is the ratio of the standard deviation of the service times to the mean of the service time.

\subsubsection*{Homework \ref{kingmanF}.\arabic{homework}}
\addtocounter{homework}{1} \addtocounter{tothomework}{1}

Show that Kingman's formula is consistent with the Pollaczek Khinchine Formula in the case where the arrival process is Poisson (Credit: Ivanovich). $~~~\Box$

\subsubsection*{Homework \ref{kingmanF}.\arabic{homework}}
\addtocounter{homework}{1} \addtocounter{tothomework}{1}

Demonstrate by simulations the accuracy of Kingman's for an infinite-buffer SSQ with IID inter-arrivals and service times and no dependencies between service times and the arrival process as $\rho$ approaches 1. You may consider, for example, different uniform distributions for the inter-arrivals and service times.
$~~~\Box$

\newpage
\section{Queues with General Input}
\label{gg1q}
\setcounter{homework}{1} 

In many situations where there is a non-zero correlation
between inter-arrival times, the Poisson assumption for the arrival process, which makes queueing models amenable to analysis, does not apply.
In this case, we consider
more general single-server queues, such as G/GI/1 and G/G/1, or their finite buffer equivalent G/GI/1/$k$ and G/G/1/$k$.
In fact, the performance of a queue can be very different if we no longer assume that inter-arrival times are IID. Consider for example
the blocking probability of an M/M/1/$N$ queue as a function of $\rho=\lambda/\mu$,
then the blocking probability will
gradually increase with $\rho$ and approaches one as $\rho \rightarrow \infty$.
However, we recall our discussion of the SLB/D/1/$N$ where we demonstrate that we can construct an
example of a finite buffer queue where the blocking
probability approaches one for an arbitrarily low value of $\rho=\lambda/\mu$.

Note that we have already covered some results applicable to G/G/1.
We already know that for G/G/1, the utilization $\hat{U}$ representing the proportion of time the server is busy
satisfies $\hat{U}=\lambda/\mu.$ We know that G/G/1 is work conservative, and we also know that Little's formula
$$
 E[Q]=\lambda E[D]
$$
is applicable to G/G/1.

\subsection{Reich's Formula}

We would like to introduce here a new and important concept the
{\it virtual waiting time}, and a formula of wide applicability in the study of
G/G/1 queues known as {\it Reich's formula} \cite{benes63,cohen82,reich58}.

The virtual waiting time, denoted $W_q(t)$, is the time that a packet has to wait in the queue (not including its own service)  if it arrives at time $t$. It is also known as {\it remaining workload}, meaning the amount of work remains in the queue at time $t$ where the work is measured in the time it is needed to be served. We assume nothing about the inter-arrival times or the service process. The latter is considered an arbitrary sequence representing the workload that each packet brings with it to the system, namely, the time required to serve each packet. For simplicity, we assume that the system is empty at time $t=0$. Let $W_a(t)$ be a function of time representing the total work arrived during the interval $[0,t)$. Then, Reich's formula says that
\begin{equation}
\label{reichs}
W_q(t) = \sup_{0\leq s< t} \{W_a(t) - W_a(s) - t +s\}.
\end{equation}
If the queue is not empty at time $t$, the largest $s$ value among all the $s$ values that maximizes the
right-hand side of
(\ref{reichs}) corresponds to the point in time where the current (at time $t$) busy period started. The other optimal $s$ values correspond to the start times of earlier busy periods.
If the queue is empty at time $t$, then the largest $s$ value among all the optimal $s$ values is equal to $t$.

\subsubsection*{Homework \ref{gg1q}.\arabic{homework}}
\addtocounter{homework}{1} \addtocounter{tothomework}{1} Consider a queue that is empty at time 0 and its arrival process and the corresponding
service duration requirements are given in the following Table.
\begin{center}
\renewcommand{\arraystretch}{1.4}
\begin{tabular}{|c|c|c|c|}\hline
Arrival time & Service duration (work requirement) &
~~~$W_q(t^+)$~~~
&~~~optimal $s$~~~  \\ \hline 1 & 3 &  &  \\ \hline 3 & 4 &  & \\
\hline
                            4 & 3 &  & \\ \hline
                             9& 3 &  & \\ \hline
                             11 & 2 &  & \\ \hline
                             11.5 & 1 &  & \\ \hline
                             17 & 4 &  &  \\
                             \hline
                            \end{tabular}
                           \end{center}
 Plot the function $W_q(t)$ for every $t$, $0\leq t\leq 25$ and fill in the right values
 for $W_q(t^+)$ and the optimal $s$ for each time point in the Table. $~~~\Box$

 \subsection{Queue Size Versus Virtual Waiting Time}

Let us now consider  the queue size probability function at time $t$
$P(Q_t=n)$, for $n=0,1,2, \ldots ~$. Its complementary distribution
function is given by $P(Q_t>n)$.  Note that for a G/D/1 queue, we have
\cite{roberts96}
\begin{equation}\label{vwtvsqs}
   Q_t=\lceil W_q(t)\rceil,
\end{equation}
so if we consider $n$ integer, and consider the service time to be
equal to one unit of work, then for a G/D/1 queue we have the
following equality for the complementary distribution functions of
the virtual waiting time $P(W_q(t)>n)$ and the queue size
\cite{roberts96}
\begin{equation}\label{vwtvsqs1}
   P(Q_t>n)=P(W_q(t)>n), {\rm ~for~} n=0,1,2, \ldots ~.
\end{equation}

The steady-state probability $P(Q>x) = \lim_{t \rightarrow \infty} P(Q_t > x)$ is often called {the queue overflow probability}.

 \subsection{Is the Queue Overflow Probability a Good QoS Measure?}

In queueing theory, we often use the steady queue overflow probability as a measure of QoS which is clearly appropriate if the arrival process is Poisson. However, it is important to understand that this is not always true. Let us consider, for example, the SLB/D/1 queue and assume that a very large burst of say $10^{15}$ packets arrive at the SSQ at time 0 and no packet arrives for service afterward, and further assume that it takes 1 ms to serve each packet. Clearly, many of the packets will suffer a very long delay, but $P(Q>x)=0$ for $X>0$ because in steady state, the queue is empty. We admit that the SLB process is not practical, but it illustrates the point that if we have a queue with very bursty traffic, the queue overflow probability may not be appropriate to represent the QoS perceived by users.

It is appropriate to note here that recently the more modern concept of Quality of Experience (QoE) \cite{Fied10} became popular. While QoS is based on statistics measured in the network, QoE is based on directly asking customers about their experience. The reader is referred to \cite{Fied10} for a comparison between the two.

\subsection{G/GI/1 Queue and Its G/GI/1/$k$ Equivalent}
 Let us consider special cases of the G/G/1 and G/G/1/$k$ queues which we call G/GI/1 and  G/GI/1/$k$,
 respectively. The GI notation indicates that the service times are mutually independent and independent
 of the arrival process and the state of the queue.
 We consider two queueing systems: a G/GI/1 queue and a G/GI/1/$k$ queue which are statistically equal in every aspect except for the fact that the first has an infinite buffer and the second has a finite buffer.
 They both have the same arrival process the distribution of their service times and the relationship of
 service times to inter-arrival times are all statistically the same.

 In queueing theory, there are many cases where it is easier to obtain overflow probability estimations of
 the unlimited buffer queue G/GI/1, namely, the steady-state probability that the queue size $Q$ exceeds a
 threshold $k$, $P(Q>k)$, than to obtain the blocking probability, denoted $P_{loss}$, of its G/GI/1/$k$ equivalent.
 In practice, no buffer is of unlimited size, so the more important problem in applications is the blocking
 probability of a G/GI/1/$k$ queue.

 By applying Little's formula to the system defined by the single
 server we can observe that the mean delay of a G/GI/1/$k$ will be
 bounded above by the mean delay of its G/GI/1 equivalent. Notice
 that if we consider only the server system for both systems, we observe
 that they have the same mean delay
 (service time) and the one associated with the G/GI/1/$k$ has
 somewhat lower arrival rate due to the losses.
 The fact that only part of the customers offered to the G/GI/1
 enter the G/GI/1/$k$ equivalent also implies that the
 percentiles of the delay distribution of the G/GI/1/$k$ system will be lower than
 those of the G/GI/1 equivalent.

An interesting problem associated with these two equivalent queues is the following.
Given $P(Q>k)$ for a
G/GI/1 queue, what can we say about the blocking probability of the
G/GI/1/$k$ equivalent. Let us begin with two examples. First,
consider a discrete-time single-server queueing model where time is
divided into fixed-length intervals called slots. This example is a
discrete-time version of our earlier example where we demonstrated a
case of a finite buffer queue with arbitrarily low traffic and large
packet loss. Assume that the service time is deterministic and is
equal to a single slot. Let the arrival process be described as
follows: $10^9$ packets arrive at the first time slot and no packets
arrived later. Consider the case of $k=1$. In the finite buffer case
with a buffer size equal to $k$, almost all the $10^9$ packets that
arrived are lost because the buffer can store only one packet.
Therefore, $P_{loss}\approx 1$. However, for the case of an infinite
buffer where we are interested in $P(W_q>k)$, ($W_q = \lim_{t
\rightarrow \infty} W_q(t)$) the case is completely the opposite.
After the $10^9$ time slots that it takes to serve the initial burst
the queue is empty forever, so in steady-state $P(W_q >k)=0$.

In our second example, on the other hand, consider another
discrete-time queueing model with $k=10^9$ and a server that serves
$10^9$ customers -- all at once at the end of a time slot -- with
probability $1-10^{-9}$ and
 $10^{90}$ customers with probability $10^{-9}$.
The rare high service rate ensures stability. Assume that at a
beginning of every  time slot, $10^9+1$  customers arrive at the
buffer. This implies that one out of the arriving $10^9+1$ customers
is lost, thus $P_{loss}\approx 10^{-9}$, while $P(W_q>k)\approx 1$.
We conclude that $P_{loss}$ and $P(W_q>k)$ can be very different.

Wong \cite{Wong90} considered this problem in the context of an ATM
multiplexer fed by multiple deterministic flows (a queueing model
denoted N\raisebox{0.5mm}{$\ast$}D/D/1 and its finite buffer
equivalent) and obtained the following inequality.
\begin{equation}
\label{wong90} \rho P_{loss}\leq P(Q>k).
\end{equation}
Roberts et al. \cite{roberts96} argued that it can be generalized to G/D/1 and its G/D/1/$k$ equivalent.
This can be further generalized. The arguments are analogous to those made in \cite{Wong90}.
 Let $\lambda$ be the arrival rate and $\mu$ the service rate in both the G/GI/1 queue and its
  G/GI/1/$k$ equivalent, with $\rho=\lambda/\mu$.
Consider a continuous period of time, in our G/GI/1 queue, during which
$Q>k$ and that just before it begins and just after it ends $Q\leq k$,
and define such time period as {\em overflow period}.
Since the queue size at the beginning is the same as at the end of the
overflow period, the number of customers that joined the queue during an
overflow period must be equal to the number of customers served during the
overflow period, because the server is continuously busy during an overflow period.

Now consider a G/GI/1/$k$ queue that has the same realization of
arrivals and their work requirements as the G/GI/1 queue. Let us
argue that in the worst case, the number of lost customers in the
G/GI/1/$k$ queue is maximized if all customers that arrive during
overflow periods of the equivalent G/GI/1 queue are lost. If for a given G/GI/1
overflow period, not all arriving customers in the G/GI/1/$k$ queue
are lost, the losses are reduced from that maximum level without
increasing future losses because at the end of a G/GI/1 overflow
period, the number of customers in the equivalent G/GI/1/$k$ queue can never be
more than $k$.

Consider a long period of time of length $L$. The mean number of
lost customers in the G/GI/1/$k$ queue during this period of time of
length $L$ is $\lambda L P_{loss}$. This must be lower or equal to
the number of customers that arrived during the same period of time
during the G/GI/1 overflow periods. This must be equal to the number
of customers served during that period of time of length $L$ during
the G/GI/1 overflow periods, which is equal to $\mu L P(Q>k)$.

Therefore,
$$\lambda L P_{loss} \leq  \mu L P(Q>k)$$
and  (\ref{wong90}) follows. $~~~\Box$

\subsubsection*{Homework \ref{gg1q}.\arabic{homework}}
\addtocounter{homework}{1} \addtocounter{tothomework}{1} Show that
(\ref{wong90}) applies to an M/M/1 queue and its M/M/1/$N$
equivalent, and discuss how tight the bound is in this case for the
complete range of parameter values.
\subsubsection*{Guide}
Recall that for M/M/1/$N$, $$P_{loss}= \frac{\rho^N(1-\rho)}{1-\rho^{N+1}},$$
and for M/M/1, $$P(Q>N)=\rho^{N+1}(1-\rho) + \rho^{N+2}(1-\rho) + \rho^{N+3}(1-\rho)+ \ldots = \rho^{N+1}.~~~\Box $$

\subsubsection*{Homework \ref{gg1q}.\arabic{homework}}
\addtocounter{homework}{1} \addtocounter{tothomework}{1} Using the UNIX command {\it netstat}
collect a sequence of 100,000 numbers representing the number of
packets arriving recorded every second for consecutive 100,000
seconds. Assume that these numbers represent the amount of work,
measured in packets, which arrive at an SSQ during 100,000
consecutive seconds. Further, assume that packets arrive at the beginning of each second interval and that all packets are of the same length and therefore require the same service time. Write a simulation of an SSQ fed by this arrival process, and compute the Packet Loss Ratio (PLR) for a range of buffer sizes and the overflow probabilities for a range of thresholds.
PLRs are relevant in the case of a finite buffer queue and
overflow probabilities represent the probability of
exceeding a threshold in an infinite buffer queue.
Plot the results in two curves, one for the PLR and the other for the overflow probabilities times $\rho^{-1}$, and observe and discuss the relationship between the two. $~~~\Box$

\subsubsection*{Homework \ref{gg1q}.\arabic{homework}}
\addtocounter{homework}{1} \addtocounter{tothomework}{1} Consider the sequence of 100,000 numbers
you have collected. Let $E[A]$ be their average. Generate a sequence
of 100,000 independent random numbers governed by a Poisson
distribution with mean $\lambda = E[A]$. Use your SSQ simulation,
and compute the PLR for a range of buffer sizes and the overflow
probabilities for a range of thresholds. Compare your results to
those obtained in the previous Assignment, and try to explain the differences.
$~~~\Box$

\subsubsection*{Homework \ref{gg1q}.\arabic{homework}}
\addtocounter{homework}{1} \addtocounter{tothomework}{1}

In this exercise, the reader is asked to
repeat the previous homework assignment for the Bernoulli process.
Again, consider the sequence of 100,000 numbers you
 have collected. Let $E[A]$ be their average. Generate a
sequence of 100,000 independent random numbers governed by the
Bernoulli distribution with mean $p = E[A]$. Use your SSQ simulation
from Exercise 1, and compute the PLR for a range of buffer sizes
and the overflow probabilities for a range of thresholds. Compare
your results to those obtained previously, and discuss the
differences. $~~~\Box$

\newpage
\section{Multi-access Modeling and analyzes}
\label{multiaccess}
\setcounter{homework}{1} 

Multi-access systems are characterized by having multiple stations (or mobile terminals) where traffic is generated, queued, and then transmitted to the network through a common transmission medium. The stations are normally {\it distributed} in the sense that one station is not aware of how much data is queued for transmission in other stations. If we view a multi-access system as a queueing system, the transmission medium which can be viewed as the server that is not aware of the load in the stations \cite{Bertsekas02}. In Section \ref{multiaccessexample}, we have already discussed an example of a multiaccess system as an application of M/M/$\infty$ queueing model. In this chapter, we discuss several other examples of multi-access systems and their performance analyses using stochastic models. The coverage here is very limited given the vast literature on stochastic modeling and analysis of multi-access systems. The aim here is not to provide a comprehensive coverage of the topic, but instead to illustrate the application of concepts learned in this book to modeling and analysis of such systems through several examples.

\subsection{ALOHA}

The ALOHA multi-access system and the related ALOHA protocol were introduced in 1971 by Professor Norman Abramson of the University of Hawaii. Fundamentally, ALOHA's ideas have been used and are still being used in many systems and protocols including Ethernet and WiFi. Quoting \cite{Abramson09b}: ``Since the 1980’s ALOHA has been the primary random access mechanism utilized by mobile telephone networks, satellite data networks, DOCSIS based cable data networks, Ethernet, WiFi and WiMAX''. Given its importance, there are many publications on ALOHA. They include:
\cite{Abramson70,Abramson09,Abramson09b,Binder75,Kuo1981,Roberts75}.

We will consider two models that correspond to two versions of ALOHA: Pure ALOHA and Slotted ALOHA. In both models  packet sizes [bits] are assumed to be fixed and denoted $s$ and the channel bitrate [bits/sec] $r$ is also given. The packet transmission time (or {\em packet time} in short) is therefore given by $\tau = s/r$. For convenience and without loss of generality, we assume $\tau = 1$. There is no co-ordination between the users and if one user transmits a packet when the channel is busy transmitting another packet, both packets collide and are lost. The throughput of the channel $T_h$, $(0 \leq T_h \leq 1)$ is the long-run average of successful transmissions per time unit. In our models, we assume that packets (including retransmitted packets) are generated following a Poisson process with a rate of $\lambda$ per packet time.

Under {\bf Pure ALOHA}, packets are transmitted as soon as they are generated. If a generated packet does not collide with another packet that is being transmitted at the moment it is generated {\bf and} if no other packets are generated and transmitted during the packet time period since the packet is generated the transmission of the packet will be successful. Under {\bf Slotted ALOHA}, time is divided into consecutive fixed-length intervals called time slots or simply slots. Each slot is equal to one packet time. Then, when a packet is generated, it is not transmitted right away, but instead has to wait until the next slot boundary to start its transmission.
In this way, if multiple packets are generated within a slot they will all be transmitted (and collide) within the next time slot, but having their total transmission time limited to one slot has an efficiency benefit as the analyses below show.

\subsubsection{Throughput Analysis of Slotted ALOHA}

To find the throughput $T_h$ under Slotted ALOHA, notice that $T_h$ is equal to the mean number of successful transmissions during a time slot. Noticing that there can be no more than one successful transmission during a time slot, we obtain

\begin{eqnarray*}
T_h &
= &   1 \times P({\rm The ~ number ~ of ~ packets~ generated ~during~ a~ slot} = 1) \\ &
+ &   0 \times P({\rm The ~ number ~ of ~ packets~ generated ~during~ a~ slot} \neq 1 ).
\end{eqnarray*}

Since the number of arrivals during a time slot of length $\tau = 1$ follows a Poisson distribution with parameter $\lambda$, we obtain

\begin{equation}
\label{ThSLotALOHA} T_h = \lambda e^{-\lambda}.
\end{equation}

Another way to explain Equation (\ref{ThSLotALOHA}) is to observe that for an arriving test packet to be successfully transmitted, it should not collide with another packet. This requires that no other packet will arrive during the same time slot that the test packet arrives. The probability of such an event is equal to the probability that a Poisson random variable with parameter $\lambda$ takes the value of zero which in turn is equal to $e^{-\lambda}$. Since on average there are $\lambda$ arrivals per time slot, and the probability of each of them to be transmitted successfully is $e^{-\lambda}$, the throughput is the product of $\lambda$ and $e^{-\lambda}$ which is consistent with Equation (\ref{ThSLotALOHA}).

An interesting question here is what is the arrival rate $\lambda$ that maximizes the throughput? This is obtained by taking derivative of the throughput with relation to $\lambda$ and equating the result to zero. By Equation (\ref{ThSLotALOHA}), we obtain

$$ \frac{dTh}{d\lambda} = e^{-\lambda} + \lambda (-e^{-\lambda}) = e^{-\lambda} [1-\lambda]. $$

We obtain that $\lambda = 1$ is the value of $\lambda$ that maximized the throughput. Substituting $\lambda = 1$ back in Equation (\ref{ThSLotALOHA}), we obtain the maximal throughput to be given by

\begin{equation}
\label{ThSLotALOHAmax1} {\rm max} ~T_h  = \frac{1}{e} = 0.3679.
\end{equation}

The reader is reminded of the assumption that the arrivals follow a Poisson process. Under this assumption, the maximal throughput is 0.3679. According to Equation (\ref{ThSLotALOHA}), the throughput increases from zero (at $\lambda = 0$) until it reaches the maximum value of $1/e$ at (at $\lambda = 1$) and then descends and approaches zero as $\lambda$ approaches infinity. If the number of users and traffic increases, the number of collisions increases and we have the so-called {\em congestion collapse} where almost all the traffic comprises colliding packets.

Notice however that a throughput of 1 (100\% throughput) is achievable if there is only one station transmitting deterministically exactly one packet in every  slot. Such a deterministic arrival process (and very high throughput) is difficult to achieve if there are multiple stations transmitting data and they are not coordinated.

\subsubsection{Throughput Analysis of Pure ALOHA}

Under Pure ALOHA, for a test packet arrival not to collide, two events must
happen:
\begin{enumerate}
\item The inter-arrival with the previous arrival must be more
      than one time slot. This happens with probability $e^{-\lambda}$.
\item The inter-arrival  with the next arrival must be more than
      one time slot. This again happens with probability $e^{-\lambda}$.
      \end{enumerate}
These two events are independent so the probability of both
to happen is $(e^{-\lambda})^2 = e^{-2\lambda}.$

Since this probability of successful transmission applies to
one arrival and the arrival rate is $\lambda$ arrivals per time slot,
the throughput is obtained to be given by

\begin{equation}
\label{ThALOHA} T_h = \lambda e^{-2\lambda}.
\end{equation}

Taking derivative we obtain

$$ \frac{dTh}{d\lambda} = e^{-2\lambda} + \lambda (-2e^{-2\lambda}) = e^{-2\lambda} [1-2\lambda]. $$

We obtain that $\lambda = 1/2$ is the value of $\lambda$ that maximizes the throughput. Substituting $\lambda = 1/2$ back in Equation (\ref{ThALOHA}), we obtain the maximal throughput to be given by

\begin{equation}
\label{ThSLotALOHAmax2} {\rm max} ~T_h  = \frac{1}{2e} = 0.184.
\end{equation}

We observe that the maximal throughput obtained for Pure ALOHA is half that of Slotted ALOHA under the Poisson arrivals assumption. As discussed, Slotted ALOHA is more efficient than pure ALOHA as it reduces the probability of collision. According to Equation (\ref{ThALOHA}), the throughput increases from zero (at $\lambda = 0$) until it reaches a maximum value of $1/2e$ at (at $\lambda = 1/2$) and then descends and approaches zero as $\lambda$ approaches infinity. As in Slotted ALOHA, if the number of users and traffic increases, the number of collisions increases and we have congestion collapse where again almost all the traffic comprises colliding packets.

 Again, we notice that even for Pure ALOHA a throughput of 1 (100\% throughput) is achievable if there is only one station transmitting deterministically exactly one packet in every  slot, and again we comment that such a deterministic arrival process is difficult to achieve when we have multiple independent uncoordinated stations transmitting data and competing for the single channel.

\subsection{Random-access Channel}

A random-access channel (RACH) is a channel shared by multiple users of a mobile network for transmission of bursty data or for call set-up to reserve a channel that will be exclusively used by an individual user for a phone call. RACH has been used in a range of networks including TDMA/FDMA and CDMA.  In RACH, time is divided into fixed-length consecutive intervals called frames. Each frame is further subdivided into slots.
A user that tries to access the channel will attempt access by choosing a slot in a frame. If no other user attempts to access at the same slot in the same frame, the transmission is successful. Otherwise, if there is a collision with another slot, both collided packets fail and will re-attempt access in the next frame.

\subsubsection*{Homework \ref{multiaccess}.\arabic{homework}}
\addtocounter{homework}{1} \addtocounter{tothomework}{1}

Consider a random access channel (RACH) with 8 slots per frame. Assume that $n$ users are attempting to access the channel, by randomly choosing one of the 8 slots in a frame. If a collision occurs, any user that experiences the collision will attempt to access the channel again in the next frame by randomly choosing a slot in the next frame. Every user keeps retransmitting in this way until its  transmission is successful. Find the mean number of retransmissions for the following two cases:
\be
\item  Case 1: $n=2$
\item  Case 2: $n=3$.
\ee

Please notice that if a user successfully accesses the channel in the first attempt, the number of retransmissions for this user is zero.

\subsubsection*{Solution}

1. {\bf Case 1: $n=2$}

We can use the Law of Iterated Expectation by considering the following two mutually exclusive and exhaustive events:\\
	Event A: The two mobiles successfully access on the first attempt.\\
	Event B: The two mobiles collide in the first attempt.\\

We have

$$P(A)=\frac{7}{8}$$

and

$$P(B)=\frac{1}{8}.$$

Let $R_1$ be the number of retransmissions. The index 1 indicates Case 1.

By the law of iterated expectations, we have

                                                   $$E[R_1] = E[R_1\mid A]P(A)+E[R_1 \mid B]P(B) ~~~~~~~~~~~~~~~~~~~~~~~~~~~~~~~~(1) $$

where $E[R_1\mid A]$ is the mean number of retransmissions conditional on that Event $A$ occurs and $E[R_1\mid B]$ is the mean number of retransmissions conditional on that Event $B$ occurs.

We also know that if Event $A$ occurs, the number of retransmissions must be 0, so we have:

$$E[R_1 \mid A]=0.$$
We also know that
$$E[R_1 \mid B]=1+E[R_1]$$
because if there is a failure in the first attempt, the mean number of retransmissions is equal to 1 plus the unconditional mean of the number of retransmissions. Notice that after the first failure the system is at the same state it was before the failure and the remaining number of retransmissions is the same as it was before the failure.

Therefore by Eq. (1)

$$E[R_1]=(1+E[R_1)P(B)$$
or
$$E[R_1]=P(B)+E[R_1]P(B)$$
or
$$E[R_1]=(P(B))/(1-P(B))=(1/8)/(1-1/8)=(1/8)/(7/8)=1/7=0.142857.$$

\vspace{0.5 cm}

2. {\bf Case 2: $n=3$}

We can again use the Law of Iterated Expectation. Consider the following three mutually exclusive and exhaustive events:\\
	Event C: All three mobiles collide in the first attempt. \\
	Event D: The two mobiles collide in the first attempt, but a third one is successfully transmitted.\\
Event E: The three mobiles successfully access in the first attempt.

We have

$$P(C) = \frac{1}{8} \times \frac{1}{8} = \frac{1}{64}. $$

$$P(D) = \frac{1}{8} \times \frac{7}{8} + \frac{7}{8} \times \frac{2}{8} = \frac{21}{64}. $$

$$P(E) = \frac{7}{8} \times \frac{6}{8}  = \frac{42}{64}. $$

Let $R_2$ be the number of retransmissions. The index 2 indicates Case 2 with $n=3$.

By the law of iterated expectations, we have

                                                   $$E[R_2] = E[R_2\mid C]P(C)+ E[R_2 \mid D]P(D)+ E[R_2 \mid E]P(E) ~~~~~~~~~~~~~~~~(2) $$

We still use the notation $R_1$ to be the number of retransmissions in Case 1 with $n=2$.

We also know that if Event $E$ occurs, the number of retransmissions must be 0, so we have:

$$E[R_2 \mid E]=0.$$
We also know that
$$E[R_2 \mid D]=1+E[R_1]$$
because if two packets collide in the first attempt, the mean number of retransmissions is equal to 1 plus the unconditional mean of the number of retransmissions in Case 1.
We also have
$$E[R_2 \mid C]=1+E[R_2]$$
because if all three packets collide in the first attempt, the mean number of retransmissions is equal to 1 plus the unconditional mean of the number of retransmissions in Case 2.

Therefore by Eq. (2)

$$E[R_2]=0\times P(E) + (1+E[R_1])P(D) + (1+E[R_2)P(C) =  P(D)+ E[R_1]P(D) + P(C) + E[R_2]P(C)  $$
or
$$E[R_2] [1-P(C)]  = P(D)+ E[R_1]P(D) + P(C). $$
Thus,
$$E[R_2] = \frac{P(D)+ E[R_1]P(D) + P(C)}{1-P(C)} = \frac{21/64 + (1/7)*(21/64) + 1/64}{1-1/64)} = \frac{25}{63} = 0.396825397.$$

~~~~~~~~~~~~~~~~~~~~~~~~~~~~~~~~~~~~~~~~~~~~~~~~~~~~~~~~~~~~~~~~~~~ $~~~\Box$

\subsubsection*{Homework \ref{multiaccess}.\arabic{homework}}
\addtocounter{homework}{1} \addtocounter{tothomework}{1}

Consider a random access channel (RACH) with 8 slots per frame. Assume that in the first frame, there are two users that are attempting to access the channel by randomly and uniformly choosing one of the 8 slots in the frame. If no collision occurs in the first frame, they both successfully transmit without the need for retransmissions. If a collision occurs in the first frame, both users that experience the collision will re-attempt to access the channel again in the next frame by randomly choosing a slot. As long as a user keep colliding it will keep retransmitting until its  transmission is successful. Assume that in addition to the two users, there is a third user that attempts to transmit in the second frame. Assume that there are no other users in addition to the above-mentioned three users that are accessing the channel.

Consider the case that if the two users  collide in the first frame, they will have to contend with the third user attempting transmission in the second frame. In this case, any one of the three users follows the same access protocol, and it keeps retransmitting until its transmission is successful.

Please notice that if a user successfully accesses the channel in its first attempt, the number of retransmissions for this user is zero.

Find the mean number of retransmissions for any one of the two users that attempt transmitting in the first frame.

\subsubsection*{Guide and Answers}

Let $R_{2,1}$ be a random variable representing the number of retransmissions for any one of the two users that attempt transmitting in the first frame. Then, a third user attempts transmission in the second frame. In this question, it is required to derive the mean $E[R_{2,1}]$.

Let $R_2$ be the number of retransmissions if two users attempt simultaneously to access in the first frame and no other users are involved. The index 2 indicates that two users and only two users are attempting to access it.

Let $R_3$ be the number of retransmissions if three users attempt simultaneously to access in the first frame and no other users are involved. The index 3 indicates that three users and only three users are attempting to access it.

The solution of this problem goes through three steps:
\be
\item Obtaining   $E[R_2]$,
\item Obtaining  $E[R_3]$
\item Obtaining  $E[R_{2,1}]$
\ee

{\bf Step 1 Obtaining   $E[R_2]$}

To obtain $E[R_2]$, use the Law of Iterated Expectation by considering the following two mutually exclusive and exhaustive events:\\
	Event A: The two mobiles successfully access on the first attempt.\\
	Event B: The two mobiles collide on the first attempt.

$$E[R_2]=0.142857.$$

{\bf Step 2:  Obtaining  $E[R_3]$}

To obtain $E[R_3]$ again use the Law of Iterated Expectation. Consider the following three mutually exclusive and exhaustive events:\\
	Event C: All three mobiles collide on the first attempt. \\
	Event D: The two mobiles collide on the first attempt, but a third one is successfully transmitted.\\
Event E: The three mobiles successfully access on the first attempt.

$$E[R_3] = 0.396825397.$$

{\bf Step 3: Obtaining   $E[R_{2,1}]$}

To obtain $E[R_{2,1}]$, again use the Law of Iterated Expectation by considering the following two mutually exclusive and exhaustive events, namely, Event A and Event B that we defined in Step 1.

Following all these steps, the end result is:

$$E[R_{2,1}]  = 0.174603175. $$


~~~~~~~~~~~~~~~~~~~~~~~~~~~~~~~~~~~~~~~~~~~~~~~~~~~~~~~~~~~~~~~~~~~~~~~~ $~~~\Box$

\subsection{Binary Exponential Backoff (BEB)}

To avoid congestion collapse, access schemes use means to resolve congestion in cases where too many users attempt to access the medium. Protocols that do that are called Medium Access Control (MAC) Protocols. One technique used to resolve congestion is known as Binary Exponential Backoff (BEB).

In the MAC protocol of 802.11, truncated binary exponential backoff (BEB) is used to resolve contention. BEB aims to avoid a situation where many users repeatedly attempt to access a congested medium by determining how long these users should wait before accessing again if they fail due to collision. In this way, BEB limits access of users during congestion thus reducing the probability that users attempt access at the same time, and improving efficiency and resource utilization. In particular, BEB defines a sequence of {\it contention windows} $W_i, i=0, 1, 2, \ldots, M$, where $W_i$ represents the time measured in time slots that the user has to wait before it is allowed to transmit again after its $i$th consecutive failure. We set $W$ to be the initial window size, so $W_0=W$, and $M$ is the maximum allowed number of attempts per packet. After $M$ attempts, the packet is lost (its transmission is unsuccessful), and the user starts a new sequence of attempts with an initial window $W$ for a new packet. In particular, for BEB, the $\{W_i\}$ sequence is given by

\begin{equation}
\label{BEBCW}
W_i = \left\{\begin{array}{ll}
2^i W  & \mbox{~~~$0\leq i < r,$} \\
2^r W & \mbox{~~~$r\leq i \leq M.$}
\end{array} \right.
\end{equation}

\subsubsection{A Discrete-time Markov Chain Model of BEB \cite{chan17}}

The operation of BEB at a mobile terminal (MT) used by a customer/user for accessing the network can be modeled as a discrete-time Markov chain, where each state represents the backoff state at the MT. In backoff state $i$, $i=0,1,2 \ldots M$, the contention window size $W_i$  is given by (\ref{BEBCW}).



Please note that since we are considering the Markov chain to be in discrete time, we ignore here the time spent in each state, and we are interested in the frequency of visits in each state. The term {\it probability} of being in each state here represents the proportion of {\bf number of times} the system enters a state relative to the number of times the system enters all states.

Let $p$ be the probability that a packet transmission experiences a collision. Let $P_i$ be the probability that an MT is in state $i$. We are interested to find the $P_i$ values for $i=0, 1, 2, \ldots, M$.  The $P_i$s obey the following steady-state equations and the normalizing equation.

$$P_i = P_{i-1}p = p^i P_0, ~i=0, 1, 2, \ldots, M$$

$$\sum_{i=0}^M P_i = 1.$$

Solving these equations, we obtain

\begin{equation}
\label{BEBdis}
P_i = \frac{1-p}{1-p^{M+1}}p^i, ~~i=0, 1, 2, \ldots, M.
\end{equation}

\subsubsection*{Homework \ref{multiaccess}.\arabic{homework}}
\addtocounter{homework}{1} \addtocounter{tothomework}{1}

Derive the result of (\ref{BEBdis}).  $~~~\Box$

\subsubsection{A Continuous-time Markov Chain Model of BEB}

In the previous section, the aim was to study a discrete-time Markov chain where each state is a state of the BEB protocol and the steady state
probabilities represented the proportion of {\bf how many times} the system enters a given state. In this section,
 we are interested in {\bf the proportion of time} the system is in state $i$. This is equivalent to the probability denoted $\pi_i$ of picking a point in time at state $i$
 if one selects a point in time uniformly at random. To derive the $\pi_i$ values we first write:

 $$V_i=P_i W_i = \frac{1-p}{1-p^{M+1}}p^i W_i, ~~i=0, 1, 2, \ldots, M.$$

 Then, to obtain the $\pi_i$ values, we use the following normalization.

 $$\pi_i =\frac{V_i}{\sum_{i=0}^{M}V_i}.$$

 This will lead to the following result.
  \begingroup
  \Large
 \begin{equation}
\label{BEBpii}
\pi_i = \left\{\begin{array}{ll}
 \frac{(2p)^i}{\frac{[1-(2p)^{r+1}]}{1-2p}+\frac{2^rp^{r+1}[1-p^{M-r}]}{1-p}} & \mbox{~~~$0\leq i < r,$} \\
\frac{2^rp^i}{\frac{[1-(2p)^{r+1}]}{1-2p}+\frac{2^rp^{r+1}[1-p^{M-r}]}{1-p}} & \mbox{~~~$r\leq i \leq M.$}
\end{array} \right.
\end{equation}
 \endgroup

\subsubsection*{Homework \ref{multiaccess}.\arabic{homework}}
\addtocounter{homework}{1} \addtocounter{tothomework}{1}

Derive the result of (\ref{BEBpii}).  $~~~\Box$

\newpage
\section{Network Models and Applications}
\label{networks}

\setcounter{homework}{1} 

So far we have considered various queueing systems, but in each
case we have considered a single queueing system in isolation.
Very important and interesting models involve networks of queues.
One important application is the Internet itself. It may be viewed
as a network of queueing systems where all network elements, such
as routers and switches are connected and where the packets are
the customers served by the various network elements and are often
queued there waiting for service.

Queueing network models can be classified into two groups: (1)
open queueing networks, and (2) closed queueing networks. In
closed queueing networks the same customers stay in the network
all the time. No new customers join and no customer leaves the
network. Customers that complete their service in one queueing system go to another and then to another and so forth, and never leave the network. In open queueing systems, new customers from the outside of the network can join any queue, and after they
complete their service in the network, obtain service from a
required number of queueing systems, which may be different for different customers,
and may include service repetitions in some queues, they leave the network. In
this section, we will only consider open queueing networks.

\subsection{Jackson Networks}
\label{jackson} Consider an open network of queues. An important issue for such
a queueing network is the statistical characteristics of the output of such queues because, in queueing networks, the output of one queue may be the input of another.

Burke's Theorem states that, in steady-state, the output (departure)
process of M/M/1, M/M/$k$  or M/M/$\infty$ queue follows a Poisson process.
Because
no traffic is lost in such queues, the arrival rate must be equal to
the departure rate, then any M/M/1, M/M/$k$, or M/M/$\infty$ queue
with arrival rate of $\lambda$ will have a Poisson departure process
with rate $\lambda$ in steady-state.

Having information about the output processes, we will now consider
an example of a very simple queueing network made of two identical
single-server queues in series, in steady-state, where all the output of the first queue is the input of the second queue and all the customers that complete
service at the second queue leave the system. Let us assume that customers' arrivals into the first queue follow a Poisson
process with parameter $\lambda$. The service times required by each
of the arriving customers at the two queues are independent and
exponentially distributed with parameter $\mu$. This means that the
amount of time a customer requires in the first queue is independent
of the amount of time a customer requires in the second queue and
they are both independent of the arrival process into the first
queue. Since the output process of the first queue is Poisson with
parameter $\lambda$, and since the first queue is clearly an M/M/1
queue, we have here two identical M/M/1 queues in
series. This is an example of a network of queues where Burke's
theorem \cite{Burke56} leads immediately to a solution for queue
size and waiting time statistics. A class of networks that can be
easily analyzed this way is the class of the so-called acyclic
networks. These networks are characterized by the fact that a
customer never goes to the same queue twice for service.

If the network is not acyclic, the independence between inter-arrival times and between inter-arrival and service times do not hold any longer. This means that the queues are no longer Markovian. To illustrate this let us consider a single server queue with feedback
described as follows. Normally, a single node does not constitute a network, however,
this simple single queue example is sufficient to illustrate the feedback effect and related dependencies.
Customers arrive into the system from the outside
according to a Poisson process with parameter $\lambda$ and the service time is exponentially distributed
 with parameter $\mu$. Then, when the customer completes the service
 the customer returns to the end of the queue with probability $p$, and with
probability $(1-p)$, the customer leaves the system. Now assume that $\lambda$ is very small
and $\mu$ is very large. Assume also that $p>0.99$. This results in an arrival process
that is based on very infrequent original arrivals (from the outside) each of which
brings with it a burst of many feedback arrivals that are very close to each other.
Clearly, this is not a Poisson process. Furthermore, the
inter-arrivals of packets within a burst, most of which are feedback
from the queue output, are very much dependent on the service times,
so clearly we have a dependence between inter-arrival times and
service times.

Nevertheless, Jackson's theorem extends the simple result
applicable to an acyclic network of queues
 to networks that are not acyclic in terms of the queue size statistics. In other words, although
the queues are not M/M/1 (or M/M/$k$  or M/M/$\infty$), they behave in terms
of their queue-size statistics as if they are. In Homework 19.1 below we illustrate it for the case on an M/M/1 queue with feedback and despite the fact
that it does not obey the independence assumptions of M/M/1, its statistical behavior in terms of queue size probabilities is the same as that of M/M/1.

Consider a network of $N$ single-server queues with infinite buffers in steady-state.
Jackson's theorem also applies to multi-server queues, but let us consider
single-server queues for now.
For queue $i$, $i=1,2,3, ~\ldots, N$, the arrival process from the outside is Poisson with
rate $r_i$. We allow for $r_i=0$ for some queues, but there must be
at least one queue $j$, such that $r_j>0$. Once a customer completes
its service in queue $i$, it continues to queue $j$ with probability
$P_{ij}$, $i=1,2,3, ~\ldots, N$, or leaves the system with probability $1-\sum_{j=1}^N
P_{ij}$. Notice that we allow for $P_{ii}>0$ for some queues. That
is, we allow for a positive probability for customers to return to the
same queue they just exited.

Let $\lambda_j$ be the total arrival rate into queue $j$. These
arrival rates can be computed by solving the following set of
equations.
\begin{equation}
\label{jacks} \lambda_j=r_j+\sum_{i=1}^N \lambda_i P_{ij},
~~~~~~j=1,2,3, ~\ldots, N.
\end{equation}
The above set of equations can be solved uniquely, if every
customer eventually leaves the network. This means that the
routing probabilities $P_{ij}$ must be such that there is a sequence of positive routing probabilities and a final exit probability that create an exit path of positive probability from
each node.

The service times at the $j$th queue are assumed exponentially
distributed with parameter $\mu_j$. They are assumed to be mutually independent and
 also  independent of the arrival process at that queue.
Let $\rho_j$ be defined by
\begin{equation}
\rho_j=\frac{\lambda_j}{\mu_j} ~~~~~{\rm for}~j=1,2,3, ~\ldots, N.
\end{equation}
Assume that $$0 \leq \rho_j < 1  ~~~~~{\rm for}~j=1,2,3, ~\ldots, N.$$

Before we proceed with the derivation of the queueing performance results,  let us summarize the conditions of Jackson's Theorem which will henceforth be called {\it Jackson's Assumptions}.

Jackson's assumptions are:

\begin{enumerate}
\item The $N$-node network, where each node has an infinite buffer queue, is open (meaning that traffic can come from outside the network and goes out from the network), and any external arrivals (arrivals from outside the network) to node $j$  follow a Poisson process with rate $r_j$. We allow some of the external arrival rates to be zero, but we require that at least one of them to one queue must be positive.
\item All service times are exponentially distributed random variables that are independent of the arriving packets to queues and their service times in previous queues,
\item The service discipline at all queues is FIFO.
\item A customer completing service at queue $i$ will either move to queue $j$ with routing probability $P_{ij}$  (note that $i$   can be equal to $j$, so we have a self-loop), or leave the system with probability  $1-\sum^N_{j=1}P_{ij}$ which is positive for some queues.
\item All the queues must be stable - namely, $0 \leq \rho_j < 1  ~~~~~{\rm for}~j=1,2,3, ~\ldots, N$.
\end{enumerate}

Let $Q_j$ be the queue size of queue $j$. Then, according to
Jackson's Theorem, in steady-state, we have that
\begin{equation}
P(Q_1=k_1, Q_2=k_2, ~\ldots,~ Q_N=k_N)=P(k_1)P(k_2)P(k_3)\cdot
~\ldots~ \cdot P(k_N),
\end{equation}
where $P(k_i)=\rho_i^{k_i}(1-\rho_i),~~~~~{\rm for}~i=1,2,3,~
\ldots, N.$

{\bf  \underline{Comment:}} Although Jackson's theorem assumes that the arrival processes from \underline{the outside} follow Poisson processes, it does not assume that the input into every queue follows a Poisson process. Therefore, it does not assume that the queues are independent M/M/1 (or M/M/$k$ or M/M/$\infty$) queues. However, it turns out, according to Jackson's theorem that the steady-state joint probability distribution of the queue sizes of the $N$ queues is obtained \underline{as if} the
     queues are independent M/M/1 (or M/M/$k$ or M/M/$\infty$) queues. This result applies despite the fact that
 the network is cyclic (not acyclic) in which case we have demonstrated that the queues do not have to be M/M/1 queues.

Accordingly, the mean
queue-size of the $j$th queue is given by
\begin{equation}
E[Q_j]=\frac{\rho_j}{1-\rho_j}.
\end{equation}
The mean delay of a customer in the $j$th queue $E[D_j]$ defined as the time from the moment the customer joins the queue until it completes service, can
be obtained by Little's formula as follows.
\begin{equation}
E[D_j]=\frac{E[Q_j]}{\lambda_j}.
\end{equation}
Using Little's formula, by considering the entire queueing network as our system,
we can also derive the mean delay of an arbitrary customer $E[D]$:
\begin{equation}
E[D]=\frac{\sum_{j=1}^N E[Q_j]} {\sum_{j=1}^N r_j}.
\end{equation}
Let us now consider a network of two-queue in series  where all the traffic that completes service in queue 1 enters queue 2 and some of the traffic in queue 2 leaves the system while the rest enters queue 1. This example is similar to the above-mentioned example of a single queue with feedback.
Using our notation, let the arrivals from the outside follow Poisson processes with rates $r_1=10^{-8}$ and $r_2=0$ and let
$\mu_1=\mu_2=1$. Further, assume that the probability that a customer that completes service in queue 2 leaves the system is $10^{-3}$, so it enters queue 1 with probability $1-10^{-3}$.

Accordingly, $$\lambda_1 = r_1 + (1-10^{-3})\lambda_2$$
and $$\lambda_2=\lambda_1.$$ Thus, $$\lambda_1 = 10^{-8} +
(1-10^{-3})\lambda_1,$$ so $$\lambda_1=\lambda_2=10^{-5}$$ and
$$\rho_1=\rho_2=10^{-5},$$ so
$$E[Q_1]=E[Q_2]=\frac{10^{-5}}{1-10^{-5}}\approx10^{-5}$$ and
$$E[D_1]=E[D_2]\approx\frac{10^{-5}}{10^{-5}} = 1.$$
Recalling that the mean service time is equal to one, this means that negligible queueing delay is expected.
(The word
`negligible' is used instead of `zero' because of the approximation
$1-10^{-5} \approx 1 $ made above.) This result makes sense
intuitively. Although the feedbacked traffic is more bursty than
Poisson we are considering here the same packet that returns over
and over again and it is impossible for the same packet to wait in
the queue for itself to be served.

An open network of M/M/1, M/M/$k$ or M/M/$\infty$ queues described
above is called a Jackson Network. For such a network an exact solution is available. However, in most practical cases,
especially when we have to deal with the so-called loss networks
that comprise queues, such as M/M/$k$/$k$, where traffic is lost, we
have to make additional modeling assumptions and to rely on
approximations to evaluate performance measures, such as blocking probability,
or carried traffic. One approximation is the so-called
Reduced-Load Erlang Fixed-Point Approximation which is reasonably
accurate and useful for loss networks.

\subsubsection*{Homework \ref{networks}.\arabic{homework}}
\addtocounter{homework}{1} \addtocounter{tothomework}{1}

Illustrate Jackson's Theorem for the case of a single queue with feedback by deriving its queueing statistics in two different ways. In particular, consider the above-mentioned queue with feedback where arrivals from the outside follow a Poisson process with parameter $\lambda$, and the service time is exponentially distributed
 with parameter $\mu$. Upon service completion, the customer returns to the end of the queue with probability $p$, and with
probability $1-p$ the customer leaves the system. Solve this problem in two ways. The first one uses Jackson's Theorem, and the second by considering M/M/1 without feedback, but with longer service durations that account for the ``unsuccessful'' service completions due to feedback. Show that both methods lead to the same queue size distribution.

\subsubsection*{Solution}
In the first way, we use Jackson's Theorem, so we write
$$\lambda_1 = \lambda + p\lambda_1$$
and we obtain
$$\lambda_1 = \frac{\lambda}{1-p}.$$

Thus,
$$\rho=\frac{\lambda}{(1-p)\mu}.$$
In the second way, we consider an equivalent system with longer service durations, we let the feedback arrivals have preemptive resume priority over all other arrivals.
This priority regime will not change the queue size statistics.
Now we have that the service time comprises a geometric sum of exponential random variables which is also exponential. Recall our discussion in Section \ref{Lap_transforms} where we learnt that a geometric sum of IID exponentially distributed random variables is exponentially distributed. Accordingly, the longer service durations are exponentially distributed with mean service time  $$\frac{1}{\mu} \times \frac{1}{(1-p)}=\frac{1}{\mu (1-p)}.$$
As a result, we have an M/M/1 queue with arrival rate $\lambda$ and service rate $\mu (1-p)$.
Thus,
$$\rho=\frac{\lambda}{(1-p)\mu}.$$
This is exactly the same result for $\rho$ that we have obtained using Jackson's theorem.

Since in M/M/1, the queue size distribution, and therefore $E[Q]$ are functions of only $\rho$, then they are equal for both cases.

Then, the mean delay is obtained according to Little's formula by the ratio $E[Q]/\lambda$, so the mean delay is also the same for both cases.

This homework demonstrates the applicability of Jackson's theorem to the case of a single node with feedback. $~~~\Box$

\subsubsection*{Homework \ref{networks}.\arabic{homework}}
\addtocounter{homework}{1} \addtocounter{tothomework}{1} Consider a 6-node network of
single server queues, the service times in all the queues are independent and exponentially distributed with a service rate equal to one, i.e., $\mu_i=1$ for $i=1,2,3, \ldots, 6$.
The arrival processes from the outside into the different queues are independent Poisson processes with rates given by $r_1=0.6$, $r_2=0.5$,
and $r_i=0$ for $i=3,4,5,6$. The routing matrix is as follows
\begin{center}
\renewcommand{\arraystretch}{1.4}
\begin{tabular}{|c|c|c|c|c|c|c|}\hline
   & 1 & 2    & 3   &  4   & 5   & 6 \\ \hline
 1 & 0 & 0.4  & 0.6 &  0   & 0   & 0 \\ \hline
 2 & 0 & 0.1  & 0   &  0.7 & 0.2 & 0 \\ \hline
 3 & 0 & 0    & 0   &  0.3 & 0.7 & 0 \\ \hline
 4 & 0 & 0    & 0   &  0   & 0   & 0.6 \\ \hline
 5 & 0 & 0    & 0   &  0.3 & 0   & 0.2 \\ \hline
 6 & 0 & 0    & 0.3 &  0   & 0   & 0 \\ \hline
                            \end{tabular}
                           \end{center}
                           In this routing matrix every row gives the routing probabilities from a specific queue to other queues. The sum of the probabilities of each row is less or equal to 1. Assume that Jackson's conditions hold and the network is in steady state.

\begin{enumerate} \item Find the mean delay in each of the queues. \item Find the mean time a packet spends in the network from the moment it enters the network until it leaves the network. \item Find the probability that the entire  network is empty. $~~~\Box$
\end{enumerate}

\subsubsection*{Homework \ref{networks}.\arabic{homework}}
Consider a queueing network that processes arriving jobs and is composed of two infinite-buffer single server queues, which we call Queue 1, and Queue 2. The network is open, and any external arrivals to either Queue 1 or Queue 2 form Poisson processes with rate $r_1=20$  arrivals per second for Queue 1, and rate $r_1=24$ arrivals per second for Queue 2. The service times of the two queues are exponentially distributed random variables that are independent of the arriving packets to the queues and their service times in the previous queues.  The service rate of each of the queues is equal to  $\mu$ jobs per second (the service rate of the two queues are equal to each other). The service discipline at both queues are First In First Out (FIFO). A customer completing service at queue i (i = 1, 2) will either move to queue $j ( j = 1, 2)$ with routing probability $P_{ij}$ or leave the system with probability  $1-P_{i1}-P_{i2}.$   The routing probabilities (the $P_{ij}$  values) are given in the following table. Notice that self-loops are also included.
\begin{center}
\renewcommand{\arraystretch}{1.4}
\begin{tabular}{|c|c|c|}\hline
   & To Queue 1 & To Queue 2     \\ \hline
 From Queue 1 & 0.4 & 0.4   \\ \hline
 From Queue 2 & 0.5 & 0.2    \\ \hline
                             \end{tabular}
                           \end{center}
Both queues must be stable. Namely, their service rates are always higher than their respective arrival rates. Let $E[D_N]$  be the mean time that a randomly chosen job spends in this queueing network from the moment it arrives from outside the network until it leaves the network. First, derive a formula for  $E[D_N]$  as a function of $\mu$  for the numerical values provided in this question. Then, find the minimal value of $\mu$, called  $\mu^*$  such that the following quality of service constraint must be maintained: $E[D_N]\leq 100$ milliseconds.

\subsubsection*{Solution}
Because Jackson conditions hold, based on Jackson theorem, the relevant queueing performance results can be obtained as if the two queues are independent M/M/1 queues.
We start by solving the following equations for $\lambda_1$ and $\lambda_2$, the arrival rates of Queue 1 and Queue 2, respectively.\\
$$\lambda_1 = 20 + 0.4 \lambda_1  + 0.5 \lambda_2 $$
and
$$ \lambda_2 = 24 + 0.4 \lambda_1 + 0.2 \lambda_2.$$
These two equations yield:  $\lambda_1=100$ and $\lambda_2=80$.

For each of the two queues, we have
$$\rho_1 =\frac{\lambda_1}{\mu} ~~~~~~ {\rm and} ~~~~~~  \rho_2 =\frac{\lambda_2}{\mu}.$$
Note that the values of $\lambda_1=100$ and $\lambda_2=80$ have been obtained above.

Considering the network as a system, the total mean queue size as a function of $\mu$ in steady state is given by
$$E[Q](\mu)  = \frac{\rho_1}{1-\rho_1}+\frac{\rho_2}{1-\rho_2} = \frac{\lambda_1/\mu}{1-\lambda_1/\mu}+\frac{\lambda_2/\mu}{1-\lambda_2/\mu}
=  \frac{\lambda_1}{\mu-\lambda_1}+\frac{\lambda_2}{\mu-\lambda_2},  $$
and substituting $\lambda_1=100$ and $\lambda_2=80$, we obtain
$$E[Q](\mu) = \frac{100}{\mu-100}+\frac{80}{\mu-80}  = \frac{100\mu -8000 +80\mu - 8000}{(\mu-100)(\mu-80)} = \frac{180\mu - 16000}{(\mu-100)(\mu-80)}. $$
and invoking Little's formula, the mean end-to-end delay of a job is given by:
$$E[D_N] = \frac{E[Q](\mu)}{r_1+r_2}. $$
Substituting $r_1=20$, and $r_2=24$, we obtain
$$E[D_N] = \frac{E[Q](\mu)}{44} = \frac{45\mu - 4000}{11(\mu-100)(\mu-80)}. $$
Meeting QoS requirement, we have the condition
$$E[D_N] = \frac{45\mu - 4000}{11(\mu-100)(\mu-80)} \leq 0.1 s. $$
We are interested in the minimal value of $\mu$ that satisfies this inequality and also satisfies the additional condition that both queues are stable.
From the above inequality we obtain

$$  45\mu - 4000  \leq 1.1(\mu-100)(\mu-80), $$
or
$$  45\mu - 4000  \leq 1.1(\mu^2-180\mu + 8000), $$
or
$$ 0  \leq 1.1\mu^2- 243\mu + 12800. $$

The following quadratic equation
$$ 0  = 1.1\mu^2- 243\mu + 12800 $$
has two real solutions: $ \mu = 86.7$ and $ \mu = 134.2$.  Recall also that $\lambda_1=100$ and $\lambda_2=80$, so for stability we must have $\mu > 100$.
To satisfy both QoS and stability requirements, we must have $\mu \geq 134.2$. Therefore, the minimal value for $\mu$ that satisfies both requirements is: $ \mu^* = 134.2$ jobs per second.

{\bf Check: } \\
Substituting $\lambda_1=100$, $\lambda_2=80$, $r_1=20$, $r_2=24$ and $\mu = 134.2$ in
$$E[D_N] =
\frac{\frac{\lambda_1/\mu}{1-\lambda_1/\mu}+\frac{\lambda_2/\mu}{1-\lambda_2/\mu}}{r_1+r_2}, $$
gives
$$E[D_N] = \frac{\frac{100/134.2}{1-100/134.2}+\frac{80/134.8}{1-80/134.2}}{20+24}=0.1~ [{\rm second}].  ~~~\Box $$

\subsection{Computation of Blocking Probability in Circuit Switched Networks by the Erlang Fixed-Point Approximation}
Let us consider a circuit-switched network made of nodes (switching
centers) that are connected by links. Each link has a fixed number
of circuits. In order to make a call between two nodes: source and destination, a user
should reserve a free circuit in each consecutive link of a path
between the two nodes. Such reservation is successful if and only if
there exists a free circuit on each of the links of that path.

An important characteristic of a circuit-switched network is that once a user makes a reservation for a connection between
a source and a destination the capacity for this connection is exclusively available for the user of this connection and no other users can utilize this capacity for the entire duration of this connection holding time.

To evaluate the probability that a circuit reservation is blocked, we
first make the following simplifying assumptions:
\begin{enumerate}
\item all the links are independent, \item  the arrival process of
calls for each origin-destination pair is Poisson, and \item  the
arrival process seen by each link is Poisson. \end{enumerate}

Having made these assumptions, we now consider each link as an
independent M/M/$k$/$k$ system for which the blocking probability is
readily available by the Erlang B formula. We assume for simplicity that the number of circuits on each link is the same (equal to $k$) for simplicity of exposition, but an extension of the model to the case where the number of circuits is different for different links is not difficult.

Let $a_j$ be the total offered load to link $j$ from all the routes that
pass through link $j$. Recall that multiplexing of Poisson processes gives another
Poisson process, whose rate is the sum of the individual rates.
Then, the blocking probability
on link $j$ is obtained by

    \begin{equation} \label{Bj} B_{j}= E_{k}(a_{j}). \end{equation}

Now that we have the means to obtain the blocking
probability on each link, we can approximate the blocking probability of a call made on a given route.
Let $B(R)$ be the blocking probability of a call made on route $R$.
The route $R$ can be viewed  as an ordered set of links, so the route blocking probability is given by
    \begin{equation} \label{BR} B(R) =1-\prod_{i\in L_R}(1-B_{i}). \end{equation}
    Note that in the above equation, $L_R$ represents the set of links in route $R$.

 Let $A(R)$ be the offered traffic on route $R$ and let $a_j(R)$ be the total traffic offered to link $j$ from traffic that flows on route $R$. Then,  $a_j(R)$ can be computed by deducting from $A(R)$ the traffic lost due to congestion on links other than $j$. That is,

    \begin{equation} \label{ajR} a_j(R) = A(R) \prod_{i\in L_R;~i\neq j} ({1-B_{i}}) ~~{\rm for} ~j\in L_R, \end{equation}

    and $a_j(R)=0$ if $j$ is not in $L_R$.  This consideration to the reduced load due to blocking on other links gave rise to the name {\it reduced load approximation} to this procedure.

  Then, the total offered traffic on link $j$ is obtained by
    \begin{equation} \label{aj} a_{j} = \sum_{R\in \mathcal{R}}  a_j(R) \end{equation}
  where $\mathcal{R}$ is the set of all routes.

  These give a set of nonlinear equations that requires a fixed-point solution.
  Notice that Erlang B is non-linear.

To solve these equations, we start with an initial vector of
 $B_j$ values; for example, set  $B_j=0$ for all $j$. Since the $A(R)$ values are known, we use equation (\ref{ajR}) to obtain the  $a_j(R)$ values.
 Then, use equation (\ref{aj}) to obtain the  $a_j$ values, which can be substituted in equation (\ref{Bj}) to obtain a new set of values for the blocking probabilities. Then, the process repeats itself iteratively until the blocking probability values  obtained in one iteration is sufficiently close to those obtained in the previous iteration. Finally, when the algorithm converges, and we obtain the resulted $B_i$ values, we use (\ref{BR}) to obtain the blocking probability on each route $R$.

The above solution, based on the principles of the Reduced-Load and
Erlang Fixed-Point Approximations can be applied to many systems and networks. For example,
an application is an Optical Burst Switching (OBS) network is described in
\cite{ros2003} where bursts of data are moving between OBS nodes
each of which is modeled as an M/M/$k$/$k$ system.

We have discussed an approach to evaluate blocking probability for circuit-switched networks under the so-called fixed routing regime,
where a call is offered to a route, and if it is rejected, it is lost and cleared from the system.
There are, however, various other regimes involving alternate routing where
rejected calls offered to a given route can overflow to other routes. A similar Erlang
fixed-point approximation can be used for circuit switching with alternate routing. See \cite{girard90}.

\subsection{A Markov-chain Simulation of a Mobile Cellular Network}
\label{mobcellsim}
A mobile cellular network can be modeled as a network of M/M/$k$/$k$
systems by assuming that the number of channels in each cell is
fixed and equal to $k$, that new call generations in each cell
follows a Poisson process, that call holding times are exponentially
distributed, and the times, until a handover occurs in each cell, are
also exponentially distributed (e.g., \cite{Hiew00}). In the following, we describe how to simulate such a network.

 Variables and input parameters:\\
$N$ = total of M/M/$k$/$k$ Systems (cells) in the network; \\
$Q(i)$ = number of calls in progress  (queue size) in cell $i$ ;\\
 $B_p$ = estimation for the blocking probability; \\
$N_a(i)$ = number of call arrivals counted so far in cell $i$; \\
$N_b(i)$ = number of blocked new calls (arrivals) counted so far in cell $i$; \\
$N_h(i)$ = number of call handovers counted so far into cell $i$; \\
$N_d(i)$ = number of dropped calls that tried to handover so far into cell $i$; \\
$MAXN_a$ = maximal number of call arrivals - used as a stopping criterion; \\
$\mu$ =  1/(the mean call holding time) \\
$\lambda(i)$ = arrival rate of new calls in cell $i$;\\
$P(i,j)$ = the probability of a handoff from cell $i$ to
cell $j$ given that the call will handoff; \\
$\delta(i)$ = handover rate in cell $i$ per call = 1/(mean time a call stays in cell $i$ before it leaves the cell)\\
$P_B$ = Blocking probability estimation.\\
$P_D$ = Dropping probability estimation of dropped handover calls.\\
$Neib(i)$ = the set of neighboring cells of cell $i$. \\
$|Neib(i)|$= number of neighboring cells of cell $i$.

Again, we will repeatedly consider  $R(01)$ a uniform $U(0,1)$
random variate. A new value for $R(01)$ is generated every time it
is called.

To know if the next event is an arrival, we use the following
IF statement.

If $$R(01) \leq \frac{\sum_{i=1}^N \lambda(i)}{\sum_{i=1}^N
\lambda(i) + \sum_{i=1}^N Q(i)\mu+ \sum_{i=1}^N Q(i)\delta(i)},   $$
then the next event is an arrival. Else, to find out if it is a
departure (it could also be a handover) we use the following
IF statement without changing the $R(01)$ value. If $$R(01) \leq \frac{\sum_{i=1}^N
\lambda(i)+ \sum_{i=1}^N Q(i)\mu}{\sum_{i=1}^N \lambda(i) +
\sum_{i=1}^N Q(i)\mu+ \sum_{i=1}^N Q(i)\delta(i)},$$ then the next
event is a departure; else, it is a handover.

If the next event is an arrival, we need to know in which of the $N$
cells it occurs. To find out, We generate a new $R(0,1)$ value and use this same value in the following loop for each $i$, $i=1, 2, \ldots  N$.

For $i=1$ to $N$, do: If $$R(01)  \leq \frac{\sum_{j=1}^i
\lambda(j)}{\sum_{j=1}^N \lambda(j)},  $$  stop the loop. The
arrival occurs in cell $i$ and we increment the total number of arrivals in cell $i$ so far by
 $N_a(i)= N_a(i) +1$, and if $Q(i)<k$, then $Q(i)=Q(i)+1$,  else
the number of lost calls needs to be incremented, namely,
$N_b(i)=N_b(i)+1$.

Now we check if this arrival is the $MAXN_a$ arrival. If it is the $MAXN_a$ arrival, then we end the simulation. In particular, if
 $$\sum_{j=1}^N N_a(j) \geq MAXN_a,$$
the simulation ends, and we compute the blocking probability as
follows. $$P_B =\frac{\sum_{i=1}^N N_b(i)}{MAXN_a}.$$ In a similar way we can compute The dropping probability of handover calls as follows.
$$P_D =\frac{\sum_{i=1}^N N_d(i)}{\sum_{i=1}^N N_h(i)}.$$

If the next event is a  departure, we need to know in which of the
$N$ cells it occurs. To find out we use the following loop (and again, we use the same $R(0,1)$ value for the entire loop).

For $i=1$ to $N$, do: If $$R(01)  \leq \frac{\sum_{j=1}^i
Q(j)\mu}{\sum_{j=1}^N Q(j)\mu}
= \frac{\sum_{j=1}^i
Q(j)}{\sum_{j=1}^m Q(j)}.  $$ Then stop the loop. The
departure occurs in System (Cell) $i$, so $Q(i)=Q(i)-1$. Note that we do
not need to verify that $Q(i)>0$ (why?).

If the next event is a  handover, we need to know from which of the
$N$ cells it handovers out of. To find it out, we use the following
loop (and again, we use the same $R(0,1)$ value for the entire loop).

For $i=1$ to $N$, do: If $$R(01)  \leq \frac{\sum_{j=1}^i
Q(j)\delta(j)}{\sum_{j=1}^N Q(j)\delta(j)}.$$ Then, stop the loop.
The handover occurs out of cell $i$, so $Q(i)=Q(i)-1$. Note that
again we do not need to verify that $Q(i)>0$.

Then, to find out into which cell the call handover in, we use the
following look (and again, we use the same $R(0,1)$ value for the entire loop).

For $j=1$ to $|Neib(i)|$, $j \in Neib(i)$, do: If $$R(01) \leq \frac{\sum_{n=1}^j
P(i,n)}{\sum_{n=1}^{|Neib(i)|} P(i,n)}
= \sum_{n=1}^j
P(i,n), ~~~n\in Neib(i), $$
the call handovers into
cell $j$. Then, increment $N_h(j)$. If cell $j$ is full, namely $Q(j)=k$, the handover is dropped, so increment also $N_d(j)$. If $Q(j)<k$, increment $Q(j)$, i.e, $Q(j) = Q(j) + 1$.

To explain how we choose the next event in this procedure, notice that we use the following probabilities:
\begin{itemize}
\item the probability that the next event is an arrival, denoted $P_{arr}$;
\item the probability that the next event is a departure, denoted $P_{dep}$;
\item the probability that the next event is a handover, denoted $P_{ho}$.
\end{itemize}

They are given by
$$P_{arr} =  \frac{\sum_{i=1}^N \lambda(i)}{\sum_{i=1}^N
\lambda(i) + \sum_{i=1}^N Q(i)\mu+ \sum_{i=1}^N Q(i)\delta(i)}.$$

$$P_{dep} = \frac{\sum_{i=1}^N Q(i)\mu}{\sum_{i=1}^N
\lambda(i) + \sum_{i=1}^N Q(i)\mu+ \sum_{i=1}^N Q(i)\delta(i)}.$$

$$P_{ho} = \frac{\sum_{i=1}^N Q(i)\delta(i)}{\sum_{i=1}^N
\lambda(i) + \sum_{i=1}^N Q(i)\mu+ \sum_{i=1}^N Q(i)\delta(i)}.$$

Then, generating $R(01)$, the condition for an arrival to be the next event is $R(01)\leq P_{arr}$.

If this condition does not hold, i.e.,  $R(01) > P_{arr}$, then the condition for a departure to be the next event is $R(01)\leq P_{arr}  + P_{dep}$.

Then, if this condition does not hold,  i.e.,  $R(01) > P_{arr}  + P_{dep}$, then the next event must be a handover.

\subsubsection*{Homework \ref{trafmod}.\arabic{homework}}
\addtocounter{homework}{1} \addtocounter{tothomework}{1}

Consider the 49-cell cellular network model with a wrapped-around design depicted in Figure 2 of \cite{Wu17}. Choose your own input parameters for the number of channels per cell, arrival rates and mean holding times.  Use a Markov-chain simulation to approximate the overall blocking probability in this network.  $~~~\Box$


\newpage
\section{Stochastic Processes as Traffic Models}
\label{trafmod}

\setcounter{homework}{1} 

In general, the aim of
traffic modeling is to provide the network designer with a relatively simple means to characterize traffic load
on a network. Ideally, such means can be used to estimate
performance and to enable efficient provisioning of network
resources. modeling a traffic stream emitted from a source, or a
traffic stream that represents multiplexing of many Internet
traffic streams is part of traffic modeling. It is normally reduced to finding a stochastic process that behaves like the real traffic stream from the point of view of the way it affects network performance or provides QoS to customers.

\subsection{Parameter Fitting}
\label{parfit}
One way to choose such a stochastic process is by fitting its
statistical characteristics to those of the real traffic stream.
Consider time to be divided into fixed-length consecutive intervals,
and consider the number of packets arriving during each time
interval as the real traffic stream. Then, the model of this traffic
stream could be a stationary discrete-time stochastic process
$\{X_n, n \geq 0\}$, with similar statistical characteristics as
those of the real traffic stream. In this case, $X_n$ could be a
random variable representing the number of packets that arrive in
the nth interval. Let $S_n$ be a random variable representing the
number of packets arriving in $n$ consecutive intervals. We may
consider the following for fitting between the statistics of $\{X_n,
n \geq 0\}$ and those of the real traffic stream:
\begin{itemize}
\item The mean $E[X_n]$. \item The variance $Var[X_n]$. \item The
 AVR is discussed in Section \ref{gencons}. The AVR is related to the so-called Index
of Dispersion for Counts (IDC) \cite{HL86} as follows: the AVR is
equal to $E[X_n]$ times the IDC.
\end{itemize}

A stationary stochastic process $\{X_n, n\geq 0\}$, where
the autocorrelation function decays slower than exponential, is said to
be Long Range Dependent (LRD). Notice that if the autocovariance sum
$\sum_{k=1}^\infty Cov(X_1,X_k)$ is infinite, the autocorrelation function must decay slower than exponential, so the process is LRD.
In such processes, the use of AVR (or IDC) may not be appropriate
because it is not finite, so a time-dependent version of the IDC,
i.e., IDC(n) = $Var[S_n]/E[X_n]$ may be considered. Another
statistical characteristic that is suitable for LRD processes is the
so-called {\em Hurst parameter} denoted by $H$ for the range $0\leq
H < 1$ that satisfies \begin{equation} \lim_{n\rightarrow \infty}
\frac{Var[S_n]}{\alpha n^{2H}} =1. \end{equation}

Each of these statistical parameters has its respective
continuous-time counterpart. As the concepts are equivalent, we
do not present them here.  We will discuss a few examples of
 stochastic processes (out of many more available in the literature)
 that have been considered as traffic models.

\subsection{Poisson Process}

For many years the Poisson process has been used as a traffic
model for the arrival process of phone calls at a telephone
exchange. The Poisson process is characterized by one parameter
$\lambda$, and $\lambda t$ is the mean as well as the variance of
the number of occurrences during any time interval of length $t$.
Its memoryless nature makes it amenable to analysis as noticed
through the analyzes of the above-mentioned queueing systems. Its
ability to characterize telephone traffic well, being
characterized by a single parameter, and its memoryless nature,
which makes it so amenable to analysis, have made the Poisson
process useful in the design and dimensioning of telephone
networks.

By its nature, the Poisson process can accurately model events
generated by a large number of independent sources each of which
generate relatively sparsely spaced events. Such events could
include phone calls or the generation of Internet traffic flows. For
example, a download of a page could be considered such a traffic
flow. However, it cannot accurately model a packet traffic stream
generated by a single user or a small number of users. It is
important to note here that many textbooks and practitioners do
consider the Poisson process as a model of a packet traffic stream
(despite the inaccuracy it introduces) due to its nice analytical
properties.

Normally, the Poisson process is defined as a continuous-time
process. However, in many cases, it is used as a model for a
discrete sequence of a traffic stream by considering time to be
divided into fixed length intervals each of size one (i.e., $t=1$),
and simply to generate a sequence of independent random numbers
which are governed by a Poisson distribution with mean $\lambda$
where $\lambda$ is equal to the average of the sequence we try to
model. As we fit only one parameter here, namely the mean, such
a model will not have the same variance, and because of the
independence property of the Poisson process, it will not mimic the
autocorrelation function of the real process. In an assignment
below, you will be asked to demonstrate that such a process does not
lead to a similar queueing curve as the real traffic stream.

\subsection{Markov Modulated Poisson Process (MMPP)}

Traffic models based on MMPP have been used to model bursty traffic.
Due to its Markovian structure together with its versatility, the
MMPP can  capture bursty traffic statistics better than the Poisson
process and still be amenable to queueing analysis. The simplest
MMPP model is MMPP(2) with only four parameters: $\lambda_0$,
$\lambda_1$, $\delta_0$, and $\delta_1$.

Queueing models involving MMPP input have been analyzed in the 70s
and 80s using Z-transform \cite{YN71,ZR86a,ZR86b,ZR86c}. Neuts
developed matrix methods to analyze such queues \cite{neuts81}. For
applications of these matrix methods for Queueing models involving
MMPP and the use of MMPP in traffic modeling and its related
parameter fitting of MMPP the reader is referred to
\cite{FMH92,HL86,lr99,neuts89}.

\subsection{Autoregressive Gaussian Process}

A traffic model based on a Gaussian process can be described as a
traffic process where the amount of traffic generated within any
time interval has a Gaussian distribution. There are several ways
to represent a Gaussian process.
 The Gaussian auto-regressive is one of them. Also, in many engineering
applications, the Gaussian process is described as a continuous-time process.
 In this section, we shall define the process as a discrete time.

Let time be divided into fixed-length intervals. Let $X_n$ be a continuous random variable representing  the amount of work entering
the system during the $n$th interval.

According to the Gaussian Autoregressive model, we assume that
$X_n,~n=1,2,= 3 ~ \ldots$ is the so-called $k$th order autoregressive
process, defined by
\begin{equation}
X_n = a_1X_{n-1} +a_2X_{n-2} +~ \ldots~+ a_kX_{n-k}    b\tilde G_n,
\end{equation}
where $\tilde G_n$ is a sequence of IID
Gaussian random variables each with mean $\eta$ and variance 1, and $a_i$ $(i=1,2, ~ \ldots,~ k)$ and $b$ are
real numbers with $\vert a\vert < 1$.

In order to characterize real traffic, we will need to find the best
fit for the parameters $a_1, ~ \ldots,~a_k,~b,$ and $\eta$. On the other
hand, it has been shown in \cite{az94}, \cite{azatr94},
\cite{azitc94} that in any Gaussian process only three parameters
are sufficient to estimate queueing performance to a reasonable
degree of accuracy. It is therefore sufficient to reduce the complexity involved in fitting many parameters and use only the 1st
order autoregressive process, also called the AR(1) process. In this
case we assume that the $X_n$ process is given by
\begin{equation}
\label{ar1} X_n = aX_{n-1}  +  b\tilde G_n,
\end{equation}
where $\tilde G_n$ is again a sequence of IID Gaussian random
variables with mean $\eta$ and variance 1, and $a$ and $b$ are real
numbers with $\vert a\vert < 1$. Let $\lambda=E [X_n] $ and
$\sigma^2=Var[ X_n ]$. The AR(1) process was proposed in
\cite{mag88} as a model of a VBR traffic stream generated by a
single source of video telephony.

The $X_n$s can be negative with positive probability. This may seem
to hinder the application of this model to real traffic processes.
However, in modeling traffic, we are not necessarily interested in a
process that is similar in every detail to the real traffic. What
we are interested in is a process that has the property that when
it is fed into a queue, the queueing performance is sufficiently
close to that of the queue fed by the real traffic.

Fitting of the parameters $a$ $b$ and $\eta$ with measurable
(estimated) parameters of the process $\lambda$, $\sigma^2$ and $S$, are provided based on \cite{azatr94}:
\begin{equation}
\label{fita} a=\frac{S}{S+\sigma^2}
\end{equation}
\begin{equation}
\label{fitb} b=\sigma^2(1-a^2)
\end{equation}
\begin{equation}
\label{fiteta} \eta=\frac{(1-a)\lambda}{b}
\end{equation}
where $S$ is the autocovariance sum given by Eq.\@ (\ref{autosum}).

\subsection{Exponential Autoregressive (1) Process}

In the previous section, we considered an autoregressive process which is Gaussian.
What made it a Gaussian process was that the so-called
{\it innovation process}, which in the case of the previous section was the
sequence $b\tilde G_n$, was a sequence of Gaussian random variables.
Letting $D_n$ be a sequence of inter-arrival times, here we consider another AR(1) process called
{\it Exponential Autoregressive (1)} (EAR(1)) \cite{gl80}, defined as follows:
\begin{equation}
\label{ear1}
D_n = aD_{n-1}  +  I_n E_n,
\end{equation}
where $D_0=I_0$, $\{I_n\}$ is a sequence of IID Bernoulli random variables in which $P(I_n=1)=1-a$ and $P(I_n=0)=a$,
and $\{E_n\}$ is a sequence of IID exponential random variables with parameter $\lambda$.

The EAR(1) has many  nice and useful properties. The $\{D_n\}$
process is a sequence of exponential random variables with parameter
$\lambda$. These are IID only for the case $a=0$. That is, when
$a=0$, the $\{D_n\}$ is a sequence of inter-arrival times of a
Poisson process. The autocorrelation function of $\{D_n\}$ is given
by \begin{equation} \label{ear_autocovck} C_{EAR1}(k) = a^k.
\end{equation} It is very easy to simulate the $\{D_n\}$ process, so
it is useful to demonstrate by simulation the relationship between
correlation in the arrival process and queueing performance.

\subsubsection*{Homework \ref{trafmod}.\arabic{homework}}
\addtocounter{homework}{1} \addtocounter{tothomework}{1} Prove that $D_n$ is exponentially
distributed for all $n\geq 0$.

\subsubsection*{Guide}
Knowing that the statement is true for $D_0$, prove that
the statement is true for $D_1$. Let ${\mathcal{L}}_X(s)$ be the
Laplace transform of random variable $X$. By definition,
${\mathcal{L}}_X(s)=E[e^{-sX}]$, so
${\mathcal{L}}_{I_1E_1}(s)=E[e^{-sI_1E_1}]$. Thus, by
(\ref{meancondind}),
${\mathcal{L}}_{I_1E_1}(s)=P(I=1)E[e^{-SE_1}]+P(I=0)E[e^{-0}]=(1-a)\lambda/(\lambda+s)+a$.
By definition,
${\mathcal{L}}_{D_1}(s)=E[e^{-s(aD_0+I_1E_1)}]={\mathcal{L}}_{D_0}(as){\mathcal{L}}_{I_1E_1}(s)$.
Recall that $D_0$ is exponentially distributed with parameter
$\lambda$, so ${\mathcal{L}}_{D_0}(as)=\lambda/(\lambda+as)$. Use
the above to show that ${\mathcal{L}}_{D_1}(s)=\lambda/(\lambda+s)$.
This proves that $D_1$ is exponentially distributed. Use the
recursion to prove that $D_n$ is exponentially distributed for all
$n>1$. $~~~\Box$


\subsection{Poisson Pareto Burst Process}
\label{ppbpmodel}
Unlike the previous models, the Poisson Pareto Burst Process (PPBP) \cite{addie09,addie02,zukerman03}
is Long Range Dependent (LRD). The PPBP has been proposed as a
more realistic model for Internet traffic than its predecessors.
According to this model, bursts of data (e.g. files) are generated
in accordance with a Poisson process with parameter $\lambda$. The
size of any of these bursts has a Pareto distribution, and each of
them is transmitted at a fixed rate $r$. At any point in time, we
may have any number of sources transmitting at rate $r$
simultaneously because according to the model, new sources may
start transmission while others are active. If $m$ sources are
simultaneously active, the total rate equals $mr$. A further
generalization of this model is the case where the burst lengths
are generally distributed. In this case, the amount of work
introduced by this model as a function of time is equivalent to
the evolution of an M/G/$\infty$ queueing system. Having $m$
sources simultaneously active is equivalent to having $m$ servers
busy in an M/G/$\infty$ system. M/G/$\infty$ which is a name of a
queueing system has been used as an alternative name of PPBP \cite{kru98,par97}. The PPBP has also been called M/Pareto/$\infty$ or
simply M/Pareto \cite{addie98}.

Again, let time be divided into fixed-length intervals, and let
$X_n$ be a continuous random variable representing the amount of
work entering the system during the $n$th interval. For convenience,
we assume that the rate $r$ is the amount transmitted by a single
source within one time interval if the source was active during the
entire interval. We also assume that the Poisson rate $\lambda$ is
per time interval. That is, the total number of transmissions to
start in one time interval is $\lambda$.

To find the  mean of $X_n$ for the PPBP process, we consider the total amount of work generated in one time interval. The reader may
notice that the mean of  the total amount of work generated in one
time interval is equal to the mean of the amount of work transmitted
in one time interval. Recalling (\ref{meanpareto}),  we obtain 
\begin{equation}
\label{mAnpar} E [ X_n ] =\frac{ \lambda r\gamma\delta }{ \gamma - 1 }, ~~~~~~~~~\delta>0, \gamma>1.
\end{equation}
Also, another  important relationship for this model, which is
provided here
 without proof, is
\begin{equation}
\label{gamH} \gamma = 3- 2H,
\end{equation}
where $H$ is the Hurst parameter.

Having the last two equations, we are able to fit the overall mean
of the process ($ E [ X_n ]$) and the Hurst parameter of the process
with those measured in a real life process,  and generate traffic
based on the PPBP model. For fitting also the variance, the reader is referred to \cite{zukerman03}.

A similar model to PPBP, is the Poisson Lomax Burst Process (PLBP) considered in \cite{chen15,chen15b}. Lomax is a version of the Pareto distribution that allows burst lengths to start from zero. Accordingly, PLPB has the advantage of including very small flows that are prevalent on the Internet. It also retains the LRD property of PPBP. Therefore, PLBP is also a suitable model for Internet traffic.
Another LRD traffic model, which is also Gaussian, is the fractional Brownian motion   \cite{chen13,chen15a,chen15b}.

\subsubsection*{Homework \ref{trafmod}.\arabic{homework}}
\addtocounter{homework}{1} \addtocounter{tothomework}{1} Use the 100,000 numbers representing the number of
packets arriving recorded every second for consecutive 100,000
seconds you have collected in the assignments of Section \ref{gg1q} Using the UNIX command {\it netstat}. Again
assume that these numbers represent the amount of work,
measured in packets, which arrive at an SSQ during 100,000
consecutive time-intervals.
 Let $E[A]$ be their average.
Use your SSQ simulation of the
assignments of Section \ref{gg1q}, and compute
 the PLR, the correlation, and the
variance of the amount of work arrives in large intervals (each of
1000 packet-transmission times) for the various processes you have
considered and discussed the differences. $~~~\Box$
\subsubsection*{Homework \ref{trafmod}.\arabic{homework}}
\addtocounter{homework}{1} \addtocounter{tothomework}{1} Compare by simulations the effect of the
correlation parameter $a$ on the performance of the queues
EAR(1)/EAR(1)/1 versus their EAR(1)/M/1, M/EAR(1)/1 and M/M/1
equivalence. Demonstrate the effect of $a$ and $\rho$ on mean delay.
Use the ranges $0\leq a \leq 1$ and $0\leq \rho \leq 1$. $~~~\Box$

\newpage
\section*{The End of the Beginning}

It is appropriate now to recall Winston Churchill's famous quote: ``Now this
is not the end. It is not even the beginning of the end. But it is, perhaps, the
end of the beginning.'' In this book, the reader has been introduced to certain fundamental theories,
techniques, and numerical methods of queueing theory and related stochastic models as well as to certain
practical telecommunications applications. However, for someone who is interested in pursuing a research
career in this field, there is a
scope for far deeper and broader study of both the theory of queues as well as telecommunications and
other applications. For the last half a century, advances in telecommunications technologies have provided
queueing theorists with a wealth of interesting problems and research challenges and it is often said that
the telecommunications and information technologies actually revived queueing theory. However, this is only
part of the story. There are many other application areas of queueing theory. The fact is that many exceptional
queueing theorists also developed expertise in various real-life systems, operations, and technologies,  and have
made tremendous contributions to their design, operations and understanding. This dual relationship between
queueing theory and its applications will likely continue, so it is very much encouraged to develop
an understanding of real-life problems as well as queueing theory. And if the aim is to become an expert in both, it is
not the end of the beginning, but merely the beginning.


\end{document}